\documentclass[12pt]{article}

\usepackage{latexsym} 
\usepackage{amssymb}
\usepackage{amsmath}
\usepackage{url}
\usepackage{booktabs}

\usepackage{amsthm}

\theoremstyle{theorem}

\newtheorem{lemma}{Lemma}

\newcommand{\gl}{\operatorname{GL}}
\newcommand{\pg}{\operatorname{PG}}

\newcommand{\gf}{\operatorname{GF}}
\newcommand{\B}{{\mathcal B}}

\begin{document}

\title{New lower bounds on the size of $(n,r)$-arcs in $\pg(2,q)$}
\author{Michael Braun\\
Faculty of Computer Science\\
University of Applied Sciences, Darmstadt, Germany\\
\url{michael.braun@h-da.de}}
\maketitle

\begin{abstract}
An $(n,r)$-arc in $\pg(2,q)$ is a set of $n$ points such that each line contains at most $r$ of the selected points. It is well-known that $(n,r)$-arcs in $\pg(2,q)$ correspond to projective linear codes. Let $m_r(2,q)$ denote the maximal number $n$ of points for which an $(n,r)$-arc in $\pg(2,q)$ exists. In this paper we obtain improved lower bounds on $m_r(2,q)$ by explicitly constructing $(n,r)$-arcs. Some of the constructed $(n,r)$-arcs correspond to linear codes meeting the Griesmer bound. All results are obtained by integer linear programming.
\end{abstract}

\section{Introduction}

A linear $[n,k,d]_q$ code is a $k$-dimensional subspace of the $n$-dimensional row vector space $\gf(q)^n$ over the finite field with $q$ elements such that the minimal Hamming distance of two distinct codewords is exactly $d$. A linear code is called projective if and only if any two columns of the generator matrix are pairwise linearly independent.

\begin{lemma}[\cite{Gri60}]
The parameters of a linear $[n,k,d]_q$ code satisfy the Griesmer bound $n\ge\sum_{i=0}^{k-1}\lceil d/q^i\rceil$.
\end{lemma}

An $(n,r)$-arc in $\pg(2,q)$ is a set $\B$ of points in $\pg(2,q)$ such that each line in $\pg(2,q)$ contains at most $r$ elements of $\B$ and such that there is at least one line containing exactly $r$ elements of $\B$.

If we write the points of an $(n,r)$-arc in $\pg(2,q)$ as columns of a generator matrix we can define a projective linear $[n,3,n-r]_q$ code. In addition, by taking the complement any $(n,r)$-arc in $\pg(2,q)$ is equivalent to a so-called $(q^2+q+1-n,q+1-r)$-blocking set in $\pg(2,q)$ (see \cite{HS01}).

A main goal in studying $(n,r)$-arcs in $\pg(2,q)$ is the determination of $m_r(2,q)$ which  is the maximum number $n$ such that an $(n,r)$-arc in $\pg(2,q)$ exists.

The values $m_r(2,q)$ for $q\le 9$ are exactly determined (see \cite{BH05}) and listed in Table~\ref{tab:qle9}. 

\begin{table}[!htbp]
\centering
\caption{Exact values $m_r(2,q)$ for $q\le 9$}\label{tab:qle9}
\begin{tabular}{rrrrrrr}
\toprule
$r\backslash q$ & 3 & 4 & 5 & 7 & 8 & 9 \\
\midrule
2&4&6&6&8&10&10\\
3&&9&11&15&17&17\\
4&&&16&22&28&28\\
5&&&&29&33&37\\
6&&&&36&42&48\\
7&&&&&49&55\\
8&&&&&&65\\
\bottomrule
\end{tabular}
\end{table}

In addition, Tables~\ref{tab:qge11} and \ref{tab:qge23} show lower and upper bound on $m_r(2,q)$. For $q\le 19$ we refer to the online table \lq\lq Three-dimensional linear codes\rq\rq{} hosted by Ball (see \cite{Bal}) whereas Table~\ref{tab:qge23} for $23\le q\le 31$ is compiled from several sources:

\begin{itemize}
\item $q=23$. \cite{DM11}, $r=3$ \cite{DM17a}.
\item $q=25$. \cite{DM13}, $r=24$ \cite[Table 5.1]{HS01}.
\item $q=27$. \cite{DM13}, $r=26$ \cite[Theorem 3.5]{DGM+11}.
\item $q=29$. \cite{Das}.
\item $q=31$. \cite{Das08,DM17b}, $r=21,22,23,24$ \cite{DM18},  $r=28,29$ \cite{DM17c}, $r=30$ \cite{CH}.
\end{itemize}

In this article many improvements of lower bounds on $m_r(2,q)$ are given indicated by the boldface numbers in Tables~\ref{tab:qge11} and \ref{tab:qge23}. Furthermore, all underlined entries define $(n,r)$-arcs in $\pg(2,q)$ corresponding to new projective linear $[n,3,n-r]_q$ codes whose parameters meet the Griesmer bound.

\begin{table}[!htbp]
\centering
\caption{Lower and upper bounds on $m_r(2,q)$ for $11 \le q\le 19$}\label{tab:qge11}
\begin{tabular}{lllllllllll}
\toprule
$r\backslash q$ & 11 & 13 & 16 & 17 & 19\\
\midrule
2&12&14&18&18&20\\
3&21&23&28&28\dots33&31\dots39\\
4&32&38\dots40&52&48\dots52&52\dots58\\
5&43\dots45&49\dots53&65&61\dots69&68\dots77\\
6&56&64\dots66&78\dots82&79\dots86&\underline{\textbf{87}}\dots96\\
7&67&79&93\dots97&95\dots103&105\dots115\\
8&78&92&120&114\dots120&126\dots134\\
9&89\dots90&105&129\dots130&137&147\dots153\\
10&100\dots102&118\dots119&142\dots148&154&172\\
11&&132\dots133&159\dots164&\underline{\textbf{167}}\dots171&191\\
12&&145\dots147&180\dots181&184\dots189&204\dots210\\
13&&&195\dots199&205\dots207&225\dots230\\
14&&&210\dots214&221\dots225&244\dots250\\
15&&&231&239\dots243&267\dots270\\
16&&&&256\dots261&286\dots290\\
17&&&&&305\dots310\\
18&&&&&325\dots330\\
\bottomrule
\end{tabular}
\end{table}

\begin{table}[!htbp]
\centering
\caption{Lower and upper bounds on $m_r(2,q)$ for $23 \le q\le 31$}\label{tab:qge23}
\begin{tabular}{llllll}
\toprule
$r\backslash q$ & 23 & 25 & 27 & 29 & 31\\
\midrule
2&24&26&28&30&32\\
3&37\dots47&38\dots51&42\dots55&\underline{\textbf{44}}\dots59&46\dots63\\
4&58\dots70&64\dots76&66\dots82&70\dots88&75\dots94\\
5&79\dots93&85\dots103&88\dots109&94\dots117&100\dots125\\
6&102\dots116&126&116\dots138&126\dots146&132\dots156\\
7&124\dots139&135\dots151&142\dots163&\textbf{148}\dots175&158\dots187\\
8&\textbf{147}\dots162&168\dots177&169\dots191&181\dots204&193\dots218\\
9&169\dots185&189\dots203&198\dots223&\textbf{208}\dots233&217\dots249\\
10&192\dots208&210\dots230&223\dots246&\textbf{234}\dots262&252\dots280\\
11&223\dots231&231\dots254&253\dots274&\textbf{262}\dots291&\textbf{282}\dots311\\
12&254&266\dots280&281\dots302&\textbf{300}\dots320&312\dots342\\
13&277&301\dots306&315\dots330&325\dots349&\textbf{348}\dots373\\
14&291\dots300&326\dots332&352\dots358&361\dots378&\textbf{378}\dots404\\
15&313\dots324&341\dots358&379\dots386&407&423\dots435\\
16&\textbf{336}\dots348&\textbf{366}\dots385&394\dots414&436&466\\
17&361\dots372&\underline{\textbf{393}}\dots411&\textbf{421}\dots442&452\dots465&497\\
18&\underline{\textbf{387}}\dots396&\textbf{416}\dots437&468\dots471&476\dots495&514\dots528\\
19&415\dots420&447\dots464&483\dots499&\textbf{507}\dots525&\textbf{539}\dots560\\
20&437\dots444&480&\underline{\textbf{510}}\dots527&\textbf{534}\dots555&\textbf{567}\dots592\\
21&461\dots468&501\dots517&\underline{\textbf{540}}\dots556&\textbf{565}\dots585&\textbf{597}\dots624\\
22&484\dots492&527\dots544&\textbf{561}\dots584&\textbf{595}\dots615&\textbf{631}\dots656\\
23&&558\dots571&\underline{\textbf{595}}\dots613&\textbf{628}\dots645&\textbf{663}\dots688\\
24&&589&624\dots642&662\dots675&\textbf{698}\dots720\\
25&&&652\dots671&\underline{\textbf{695}}\dots705&\textbf{733}\dots752\\
26&&&677\dots699&725\dots735&\textbf{768}\dots784\\
27&&&&755\dots765&\underline{\textbf{805}}\dots816\\
28&&&&784\dots795&837\dots848\\
29&&&&&869\dots880\\
30&&&&&901\dots913\\
\bottomrule
\end{tabular}
\end{table}

\section{Construction approach}

We use the construction approach of $(n,r)$-arcs in $\pg(2,q)$ with prescribed groups of linear automorphisms using integer linear programming described in \cite{BKW05}. We briefly summarize the construction.

An $(n,r)$-arc $\B$ in $\pg(2,q)$ admits a subgroup $G\le\gl(3,q)$ as a group of linear automorphisms if and only if $\B$ consists of $G$-orbits on the set of points in $\pg(2,q)$ (which correspond to the $1$-dimensional subspaces of $\gf(q)^3$). Let $m^G_r(2,q)$ denote the maximal size $n$ of an $(n,r)$-arc in $\pg(2,q)$ admitting $G\le\gl(3,q)$ as a group of linear automorphisms. Obviously, for any $G\le\gl(3,q)$ we get a lower bound $m^G_r(2,q)\le  m_r(2,q)$.

Let $M^G$ denote the $G$-incidence matrix whose rows belong to the $G$-orbits on the set of $2$-dimensional subspaces of $\gf(q)^3$ and whose columns correspond to the orbits on the set of $1$-dimensional subspaces. The entry $m_{ij}$ of $M^G$ counts the number of subspaces in the $j$th column orbit which are contained in a fixed subspace of the $i$th row orbit. Let $w=(w_1,\ldots,w_t)^T$ denote the vector whose entry $w_j$ is the cardinality of the $j$th column orbit. For the trivial subgroup $G=\{1\}$ the matrix $M^{\{1\}}$ is the incidence matrix between lines and points of $\pg(2,q)$.

If $u=(1,\ldots,1)^T$ any binary vector $x=(x_1,\ldots,x_t)^T$ with $M^Gx\le r u$ defines a $(w^Tx,r)$-arc in $\pg(2,q)$ admitting $G\le \gl(3,q)$ as a group of linear automorphisms. Consequently, we obtain the following integer linear programming
\[
m^G_r(2,q)=\max_{x\in\{0,1\}^t}\{w^Tx\mid M^Gx\le r u\}.
\]

\section{New lower bounds}

In the following we list the cyclic subgroups of linear automorphisms defining the $(n,r)$-arcs in $\pg(2,q)$ with improved size given in Tables~\ref{tab:qge11} and \ref{tab:qge23}. All $(n,r)$-arcs in $\pg(2,q)$ were computed with Gurobi (see \cite{Gur}) as ILP solver. We list all constructed $(n,r)$-arcs in the appendix. 

In order to find the cyclic groups to be prescribed as groups of linear automorphisms we randomly chose invertible $3\times 3$ matrices. Furthermore, for $q=11$, $13$, and $17$ we considered representatives of all conjugacy classes of elements of $\gl(3,q)$. 

In the following each of the constructed new $(n,r)$-arcs in $\pg(2,q)$ will be presented by a vector $(q,r,n;a_{11},a_{12},a_{13},a_{21},a_{22},a_{23},a_{31},a_{32},a_{33};time)$. The entries $a_{11},\ldots,a_{33}$ define the generator of the prescribed cyclic group:
\[
\begin{pmatrix}
a_{11} & a_{12} & a_{13}\\
a_{21} & a_{22} & a_{23}\\
a_{31} & a_{32} & a_{33}
\end{pmatrix}.
\]
The rightmost entry $time$ indicates the runtime in seconds needed to compute the corresponding $(n,r)$-arc on a personal computer (1.2 GHz Intel Core m3 processor).

Field elements of the prime field $\gf(p)$ are represented by integers $0\le a < p$ where elements are added and multiplied modulo $p$. In extension fields $\gf(p^e)$ the elements $\sum_{i=0}^{e-1}a_ix^i$ are given by the numbers $\sum_{i=0}^{e-1}a_ip^i$ where elements are added and multiplied modulo a given irreducible polynomial $f(x)\in\gf(p)[x]$ of degree $e$. For the fields $\gf(16)$, $\gf(25)$, and $\gf(27)$ we use the irreducible polynomials $x^4+x^3+1\in\gf(2)[x]$, $x^2+x+2\in\gf(5)[x]$, and $x^3+2x+1\in\gf(3)[x]$.

$(17,11,167;0,12,10,5,1,12,1,9,6;11)$ %

$(19,6,87;0,18,13,12,16,13,13,14,9;12)$ %

$(23,8,147;21,17,17,3,5,22,17,5,16;44)$ %

$(23,16,336;2,14,0,7,17,5,20,19,19;22)$ %

$(23,18,387;10,21,8,0,12,6,15,7,9;3)$ %

$(25,16,366;10,16,19,1,19,14,12,20,13;2)$ %

$(25,17,393;15,19,8,6,8,10,16,16,18;10)$ %

$(25,18,416;15,12,24,9,24,19,24,1,2;44)$ %

$(27,17,421;4,24,3,13,6,5,26,3,3;5)$ %

$(27,20,510;5,7,0,15,0,19,2,3,22;0)$ %

$(27,21,540;1,17,15,23,22,17,7,17,25;21)$ %

$(27,22,561;26,0,25,9,6,23,20,21,13;10)$ %

$(27,23,595;14,10,10,26,6,24,20,11,19;276)$ %

$(29,3,44;9,16,22,22,7,27,4,21,16;0)$ %

$(29,7,148;1,13,3,28,1,20,6,1,9;9)$ %

$(29,9,208;18,25,4,11,11,12,26,17,23;2)$ %

$(29,10,234;26,4,3,1,21,0,2,20,19;13)$ %

$(29,11,262;4,24,18,27,11,12,26,7,17;50)$ %

$(29,12,300;22,14,11,23,11,4,15,26,16;19)$ %

$(29,19,507;10,19,28,22,22,18,12,8,5;1)$ %

$(29,20,534;16,2,8,10,11,25,16,18,26;11)$ %

$(29,21,565;25,25,21,23,17,21,1,27,6;38)$ %

$(29,22,595;8,12,21,14,25,24,22,5,6;39)$ %

$(29,23,628;24,16,26,10,0,20,20,18,11;39)$ %

$(29,25,695;25,25,21,23,17,21,1,27,6;29)$ %

$(31,11,282;2,19,1,30,10,22,8,17,23;34)$ %

$(31,13,348;27,20,30,26,26,8,22,20,11;89)$ %

$(31,14,378;19,16,26,5,13,19,3,20,19;2)$ %

$(31,19,539;21,20,1,12,15,6,27,25,22;240)$

$(31,20,567;0,7,29,6,8,21,13,12,27;292)$ %

$(31,21,597;29,12,3,30,28,2,30,6,3;10)$ %

$(31,22,631;26,8,1,25,14,23,17,15,4;14)$ %

$(31,23,663;21,10,23,3,5,23,30,28,20;7)$ %

$(31,24,698;23,10,11,12,18,9,6,18,26;35)$ %

$(31,25,733;30,30,25,3,0,2,24,14,18;22)$ %

$(31,26,768;28,6,30,19,28,8,21,25,4;1)$

$(31,27,805;28,18,2,10,12,0,15,30,7;17)$

\section{Acknowledgement}

I would like to thank Rumen Daskalov and the anonymous referees for their valuable comments.

\section*{Appendix}

\subsection*{$m_{11}(2,17)\ge 167$}
$(0,1,1)$, $(1,6,3)$, $(0,1,3)$, $(1,11,14)$, $(0,1,4)$, $(1,15,16)$, $(0,1,12)$, 
$(1,2,1)$, $(0,1,7)$, $(1,0,0)$, $(0,1,11)$, $(1,16,8)$, $(0,1,14)$, $(1,1,9)$, 
$(0,1,13)$, $(1,12,6)$, $(0,1,15)$, $(1,5,11)$, $(1,0,1)$, $(1,0,16)$, $(1,16,9)$, 
$(1,1,8)$, $(1,0,2)$, $(1,4,10)$, $(1,15,1)$, $(1,5,2)$, $(1,0,3)$, $(1,11,8)$, 
$(1,14,10)$, $(1,12,0)$, $(1,0,6)$, $(1,1,6)$, $(1,11,3)$, $(1,2,15)$, $(1,0,8)$, 
$(1,7,14)$, $(1,9,4)$, $(1,8,6)$, $(1,0,9)$, $(1,9,11)$, $(1,8,13)$, $(1,10,3)$, 
$(1,0,10)$, $(1,14,12)$, $(1,7,5)$, $(1,15,4)$, $(1,0,11)$, $(1,15,2)$, $(1,6,14)$, 
$(1,16,11)$, $(1,1,0)$, $(1,9,15)$, $(1,8,12)$, $(1,6,5)$, $(1,1,1)$, $(1,7,10)$, 
$(1,7,4)$, $(1,12,4)$, $(1,1,2)$, $(1,2,6)$, $(1,6,13)$, $(1,10,10)$, $(1,1,4)$, 
$(1,3,0)$, $(1,4,14)$, $(1,7,2)$, $(1,1,11)$, $(1,12,14)$, $(1,14,9)$, $(1,14,15)$, 
$(1,1,13)$, $(1,10,9)$, $(1,12,10)$, $(1,3,14)$, $(1,1,15)$, $(1,15,13)$, $(1,10,11)$, 
$(1,5,8)$, $(1,1,16)$, $(1,14,2)$, $(1,9,3)$, $(1,8,16)$, $(1,2,2)$, $(1,15,15)$, 
$(1,14,8)$, $(1,3,9)$, $(1,2,4)$, $(1,16,2)$, $(1,12,9)$, $(1,7,6)$, $(1,2,9)$, 
$(1,11,16)$, $(1,7,3)$, $(1,4,4)$, $(1,2,10)$, $(1,5,9)$, $(1,6,12)$, $(1,14,5)$, 
$(1,2,11)$, $(1,7,0)$, $(1,5,4)$, $(1,5,16)$, $(1,2,12)$, $(1,4,5)$, $(1,4,13)$, 
$(1,10,8)$, $(1,2,14)$, $(1,3,15)$, $(1,16,1)$, $(1,9,1)$, $(1,8,14)$, $(1,3,16)$, 
$(1,9,10)$, $(1,8,10)$, $(1,11,11)$, $(1,4,0)$, $(1,14,16)$, $(1,15,14)$, $(1,11,6)$, 
$(1,4,6)$, $(1,5,15)$, $(1,9,0)$, $(1,8,8)$, $(1,4,8)$, $(1,6,0)$, $(1,7,1)$, 
$(1,14,4)$, $(1,4,9)$, $(1,16,3)$, $(1,6,10)$, $(1,6,15)$, $(1,4,16)$, $(1,12,12)$, 
$(1,16,5)$, $(1,16,14)$, $(1,5,1)$, $(1,11,9)$, $(1,5,12)$, $(1,12,5)$, $(1,11,15)$, 
$(1,6,2)$, $(1,5,14)$, $(1,15,10)$, $(1,9,16)$, $(1,8,15)$, $(1,6,4)$, $(1,11,13)$, 
$(1,10,6)$, $(1,7,11)$, $(1,6,8)$, $(1,12,16)$, $(1,7,7)$, $(1,11,4)$, $(1,15,11)$, 
$(1,16,15)$, $(1,8,2)$, $(1,16,0)$, $(1,11,12)$, $(1,9,5)$, $(1,8,5)$, $(1,9,12)$, 
$(1,10,0)$, $(1,15,6)$, $(1,12,1)$, $(1,12,13)$, $(1,10,2)$, $(1,14,13)$.

\subsection*{$m_{6}(2,19)\ge 87$}
$(0,0,1)$, $(1,1,8)$, $(1,7,10)$, $(1,5,16)$, $(1,7,16)$, $(1,6,18)$, $(0,1,1)$, 
$(1,4,13)$, $(1,13,7)$, $(1,13,15)$, $(1,17,11)$, $(1,15,7)$, $(0,1,4)$, $(1,14,17)$, 
$(1,9,9)$, $(1,2,14)$, $(1,4,8)$, $(1,17,13)$, $(0,1,13)$, $(1,8,7)$, $(1,14,16)$, 
$(1,16,17)$, $(1,6,7)$, $(1,1,3)$, $(1,0,2)$, $(1,0,18)$, $(1,3,7)$, $(1,11,13)$, 
$(1,12,3)$, $(1,9,7)$, $(1,0,9)$, $(1,5,6)$, $(1,13,5)$, $(1,0,15)$, $(1,11,3)$, 
$(1,2,11)$, $(1,2,8)$, $(1,13,8)$, $(1,14,14)$, $(1,0,16)$, $(1,8,14)$, $(1,6,14)$, 
$(1,1,18)$, $(1,3,15)$, $(1,4,0)$, $(1,1,9)$, $(1,6,16)$, $(1,2,5)$, $(1,15,8)$, 
$(1,4,11)$, $(1,8,9)$, $(1,17,12)$, $(1,10,6)$, $(1,2,15)$, $(1,10,8)$, $(1,9,3)$, 
$(1,16,3)$, $(1,15,9)$, $(1,12,0)$, $(1,3,13)$, $(1,15,15)$, $(1,18,6)$, $(1,12,5)$, 
$(1,4,10)$, $(1,14,18)$, $(1,3,14)$, $(1,16,6)$, $(1,16,5)$, $(1,13,17)$, $(1,15,13)$, 
$(1,11,18)$, $(1,5,0)$, $(1,12,10)$, $(1,17,6)$, $(1,5,5)$, $(1,8,11)$, $(1,18,17)$, 
$(1,7,6)$, $(1,9,5)$, $(1,8,13)$, $(1,9,10)$, $(1,11,11)$, $(1,11,0)$, $(1,14,9)$, 
$(1,18,3)$, $(1,16,11)$, $(1,17,1)$.

\subsection*{$m_{8}(2,23)\ge 147$}
$(0,1,2)$, $(1,19,12)$, $(1,13,16)$, $(1,18,2)$, $(1,10,13)$, $(1,3,12)$, $(0,1,22)$, 
$(0,1,9)$, $(1,20,14)$, $(1,20,19)$, $(1,9,18)$, $(1,13,6)$, $(1,7,6)$, $(1,18,20)$, 
$(0,1,12)$, $(1,11,19)$, $(1,8,4)$, $(1,6,8)$, $(1,8,10)$, $(1,2,2)$, $(1,4,11)$, 
$(0,1,15)$, $(1,14,2)$, $(1,15,7)$, $(1,12,5)$, $(1,22,8)$, $(1,18,1)$, $(1,0,15)$, 
$(0,1,16)$, $(1,8,13)$, $(1,3,15)$, $(1,19,13)$, $(1,3,14)$, $(1,15,17)$, $(1,10,5)$, 
$(0,1,17)$, $(1,18,10)$, $(1,1,1)$, $(1,11,17)$, $(1,21,18)$, $(1,8,16)$, $(1,11,4)$, 
$(0,1,18)$, $(1,10,17)$, $(1,11,2)$, $(1,20,1)$, $(1,1,11)$, $(1,19,11)$, $(1,9,6)$, 
$(1,0,19)$, $(1,16,1)$, $(1,20,4)$, $(1,2,9)$, $(1,4,10)$, $(1,6,6)$, $(1,13,22)$, 
$(1,1,17)$, $(1,12,22)$, $(1,18,15)$, $(1,21,20)$, $(1,13,1)$, $(1,15,1)$, $(1,14,5)$, 
$(1,2,4)$, $(1,4,20)$, $(1,14,10)$, $(1,18,18)$, $(1,12,8)$, $(1,2,20)$, $(1,4,12)$, 
$(1,2,5)$, $(1,4,19)$, $(1,21,2)$, $(1,16,15)$, $(1,6,2)$, $(1,21,9)$, $(1,15,18)$, 
$(1,2,14)$, $(1,4,8)$, $(1,20,13)$, $(1,3,7)$, $(1,18,14)$, $(1,22,6)$, $(1,6,11)$, 
$(1,2,16)$, $(1,4,6)$, $(1,9,19)$, $(1,7,13)$, $(1,3,22)$, $(1,20,12)$, $(1,7,22)$, 
$(1,2,19)$, $(1,4,17)$, $(1,19,1)$, $(1,19,8)$, $(1,16,12)$, $(1,10,19)$, $(1,10,9)$, 
$(1,3,5)$, $(1,14,18)$, $(1,6,20)$, $(1,12,15)$, $(1,7,17)$, $(1,6,17)$, $(1,15,13)$, 
$(1,3,9)$, $(1,12,20)$, $(1,21,4)$, $(1,9,13)$, $(1,3,10)$, $(1,16,16)$, $(1,14,13)$, 
$(1,6,14)$, $(1,17,4)$, $(1,13,12)$, $(1,19,10)$, $(1,11,11)$, $(1,6,16)$, $(1,10,1)$, 
$(1,7,7)$, $(1,9,12)$, $(1,22,7)$, $(1,19,4)$, $(1,22,4)$, $(1,21,10)$, $(1,10,10)$, 
$(1,7,16)$, $(1,12,7)$, $(1,13,18)$, $(1,22,11)$, $(1,8,17)$, $(1,8,11)$, $(1,14,20)$, 
$(1,8,14)$, $(1,13,7)$, $(1,11,15)$, $(1,22,12)$, $(1,9,20)$, $(1,15,19)$, $(1,13,5)$, 
$(1,8,22)$, $(1,22,9)$, $(1,20,7)$, $(1,14,8)$, $(1,22,2)$, $(1,12,6)$, $(1,16,7)$.

\subsection*{$m_{16}(2,23)\ge 336$}
$(0,0,1)$, $(0,1,13)$, $(1,19,19)$, $(1,13,5)$, $(0,1,2)$, $(1,20,9)$, $(1,4,7)$, 
$(1,13,21)$, $(0,1,3)$, $(1,22,12)$, $(1,15,2)$, $(1,13,18)$, $(0,1,5)$, $(1,3,18)$, 
$(1,18,9)$, $(1,13,4)$, $(0,1,6)$, $(1,5,21)$, $(1,12,18)$, $(1,13,17)$, $(0,1,7)$, 
$(1,7,1)$, $(1,2,10)$, $(1,13,19)$, $(0,1,8)$, $(1,9,4)$, $(1,5,17)$, $(1,13,8)$, 
$(0,1,9)$, $(1,11,7)$, $(1,14,15)$, $(1,13,3)$, $(0,1,12)$, $(1,17,16)$, $(1,10,21)$, 
$(1,13,20)$, $(0,1,14)$, $(1,21,22)$, $(1,3,20)$, $(1,13,13)$, $(0,1,16)$, $(1,2,5)$, 
$(1,16,12)$, $(1,13,0)$, $(0,1,17)$, $(1,4,8)$, $(1,0,13)$, $(1,13,7)$, $(0,1,18)$, 
$(1,6,11)$, $(1,11,8)$, $(1,13,2)$, $(0,1,19)$, $(1,8,14)$, $(1,17,22)$, $(1,13,14)$, 
$(0,1,21)$, $(1,12,20)$, $(1,9,11)$, $(1,13,6)$, $(1,0,4)$, $(1,2,2)$, $(1,4,17)$, 
$(1,21,10)$, $(1,0,7)$, $(1,21,19)$, $(1,8,16)$, $(1,7,7)$, $(1,0,8)$, $(1,12,17)$, 
$(1,15,20)$, $(1,8,22)$, $(1,0,10)$, $(1,17,13)$, $(1,20,13)$, $(1,15,12)$, $(1,0,11)$, 
$(1,8,11)$, $(1,9,10)$, $(1,2,1)$, $(1,0,14)$, $(1,4,5)$, $(1,16,14)$, $(1,22,2)$, 
$(1,0,15)$, $(1,18,3)$, $(1,6,5)$, $(1,16,4)$, $(1,0,16)$, $(1,9,1)$, $(1,18,2)$, 
$(1,1,9)$, $(1,0,17)$, $(1,0,22)$, $(1,1,12)$, $(1,11,21)$, $(1,0,18)$, $(1,14,20)$, 
$(1,0,19)$, $(1,5,18)$, $(1,7,22)$, $(1,18,11)$, $(1,0,20)$, $(1,19,16)$, $(1,12,15)$, 
$(1,19,3)$, $(1,0,21)$, $(1,10,14)$, $(1,10,4)$, $(1,9,14)$, $(1,1,2)$, $(1,5,12)$, 
$(1,20,7)$, $(1,10,16)$, $(1,1,4)$, $(1,20,0)$, $(1,8,13)$, $(1,22,18)$, $(1,1,5)$, 
$(1,16,17)$, $(1,1,8)$, $(1,4,22)$, $(1,2,16)$, $(1,14,9)$, $(1,1,13)$, $(1,7,15)$, 
$(1,5,3)$, $(1,5,19)$, $(1,1,14)$, $(1,3,9)$, $(1,6,14)$, $(1,20,10)$, $(1,1,17)$, 
$(1,14,14)$, $(1,11,0)$, $(1,21,14)$, $(1,1,18)$, $(1,10,8)$, $(1,14,10)$, $(1,3,11)$, 
$(1,1,20)$, $(1,2,19)$, $(1,3,4)$, $(1,7,4)$, $(1,1,21)$, $(1,21,13)$, $(1,18,8)$, 
$(1,8,8)$, $(1,1,22)$, $(1,17,7)$, $(1,17,20)$, $(1,12,1)$, $(1,2,3)$, $(1,8,0)$, 
$(1,18,12)$, $(1,5,15)$, $(1,2,9)$, $(1,9,13)$, $(1,12,16)$, $(1,17,0)$, $(1,2,12)$, 
$(1,21,8)$, $(1,11,9)$, $(1,12,12)$, $(1,2,13)$, $(1,2,14)$, $(1,6,20)$, $(1,15,14)$, 
$(1,2,15)$, $(1,10,3)$, $(1,2,17)$, $(1,18,15)$, $(1,20,3)$, $(1,22,11)$, $(1,2,18)$, 
$(1,22,21)$, $(1,17,5)$, $(1,16,7)$, $(1,2,22)$, $(1,15,22)$, $(1,10,2)$, $(1,18,16)$, 
$(1,3,1)$, $(1,3,21)$, $(1,22,15)$, $(1,8,3)$, $(1,3,2)$, $(1,12,0)$, $(1,3,3)$, 
$(1,21,2)$, $(1,21,1)$, $(1,15,15)$, $(1,3,7)$, $(1,11,10)$, $(1,11,22)$, $(1,22,4)$, 
$(1,3,10)$, $(1,15,16)$, $(1,4,16)$, $(1,11,18)$, $(1,3,12)$, $(1,10,20)$, $(1,6,21)$, 
$(1,18,7)$, $(1,3,13)$, $(1,19,22)$, $(1,14,18)$, $(1,19,12)$, $(1,3,14)$, $(1,5,1)$, 
$(1,17,14)$, $(1,9,8)$, $(1,3,19)$, $(1,4,11)$, $(1,7,12)$, $(1,6,16)$, $(1,4,0)$, 
$(1,12,8)$, $(1,10,1)$, $(1,11,17)$, $(1,4,1)$, $(1,22,0)$, $(1,20,21)$, $(1,14,22)$, 
$(1,4,2)$, $(1,9,15)$, $(1,11,3)$, $(1,18,21)$, $(1,4,9)$, $(1,10,5)$, $(1,16,13)$, 
$(1,6,1)$, $(1,4,10)$, $(1,20,20)$, $(1,17,15)$, $(1,21,3)$, $(1,4,12)$, $(1,17,4)$, 
$(1,4,20)$, $(1,5,9)$, $(1,15,11)$, $(1,22,20)$, $(1,5,4)$, $(1,22,10)$, $(1,12,9)$, 
$(1,8,21)$, $(1,5,8)$, $(1,21,20)$, $(1,14,2)$, $(1,10,10)$, $(1,5,10)$, $(1,9,2)$, 
$(1,6,7)$, $(1,22,13)$, $(1,5,16)$, $(1,19,17)$, $(1,20,4)$, $(1,19,18)$, $(1,5,22)$, 
$(1,6,9)$, $(1,6,3)$, $(1,10,18)$, $(1,15,0)$, $(1,11,15)$, $(1,6,10)$, $(1,8,15)$, 
$(1,12,3)$, $(1,20,16)$, $(1,6,12)$, $(1,14,1)$, $(1,8,7)$, $(1,6,17)$, $(1,6,13)$, 
$(1,17,17)$, $(1,22,16)$, $(1,21,11)$, $(1,6,18)$, $(1,9,5)$, $(1,16,22)$, $(1,12,10)$, 
$(1,7,0)$, $(1,10,22)$, $(1,20,19)$, $(1,20,8)$, $(1,7,3)$, $(1,9,9)$, $(1,14,19)$, 
$(1,21,15)$, $(1,7,11)$, $(1,14,5)$, $(1,16,19)$, $(1,10,7)$, $(1,7,17)$, $(1,12,2)$, 
$(1,22,19)$, $(1,14,12)$, $(1,7,20)$, $(1,11,12)$, $(1,9,19)$, $(1,9,0)$, $(1,8,4)$, 
$(1,21,5)$, $(1,16,5)$, $(1,16,16)$, $(1,8,9)$, $(1,19,2)$, $(1,15,8)$, $(1,19,4)$, 
$(1,8,18)$, $(1,20,15)$, $(1,9,3)$, $(1,17,12)$, $(1,10,11)$, $(1,12,7)$, $(1,12,4)$, 
$(1,18,13)$, $(1,11,2)$, $(1,19,11)$, $(1,18,1)$, $(1,19,13)$, $(1,14,8)$, $(1,22,5)$, 
$(1,16,15)$, $(1,15,3)$, $(1,14,17)$, $(1,17,9)$, $(1,18,10)$, $(1,18,0)$, $(1,15,7)$, 
$(1,18,14)$, $(1,17,11)$, $(1,19,0)$, $(1,22,7)$, $(1,19,8)$, $(1,20,11)$, $(1,21,21)$.

\subsection*{$m_{18}(2,23)\ge 387$}
$(0,1,0)$, $(1,17,8)$, $(1,4,4)$, $(1,17,3)$, $(0,1,1)$, $(1,3,18)$, $(1,19,6)$, 
$(1,4,17)$, $(0,1,3)$, $(1,16,12)$, $(1,16,1)$, $(1,22,10)$, $(0,1,4)$, $(1,15,16)$, 
$(1,0,5)$, $(1,19,15)$, $(0,1,6)$, $(0,1,19)$, $(1,22,11)$, $(1,1,22)$, $(0,1,7)$, 
$(1,1,3)$, $(1,11,8)$, $(1,7,12)$, $(0,1,8)$, $(1,21,15)$, $(1,18,12)$, $(1,12,19)$, 
$(0,1,9)$, $(1,20,19)$, $(1,21,17)$, $(1,18,9)$, $(0,1,10)$, $(1,8,21)$, $(1,15,7)$, 
$(1,10,7)$, $(0,1,11)$, $(1,10,13)$, $(1,7,9)$, $(1,0,16)$, $(0,1,14)$, $(1,13,1)$, 
$(1,20,0)$, $(1,15,14)$, $(0,1,15)$, $(1,11,9)$, $(1,10,14)$, $(1,2,5)$, $(0,1,17)$, 
$(1,6,6)$, $(0,1,18)$, $(1,7,2)$, $(1,8,3)$, $(1,14,8)$, $(0,1,20)$, $(1,2,22)$, 
$(1,14,13)$, $(1,5,0)$, $(0,1,22)$, $(1,4,14)$, $(1,6,15)$, $(1,8,18)$, $(1,0,0)$, 
$(1,0,13)$, $(1,14,6)$, $(1,16,14)$, $(1,0,1)$, $(1,8,9)$, $(1,19,3)$, $(1,19,8)$, 
$(1,0,2)$, $(1,4,11)$, $(1,12,21)$, $(1,14,18)$, $(1,0,3)$, $(1,10,8)$, $(1,21,11)$, 
$(1,18,10)$, $(1,0,8)$, $(1,5,22)$, $(1,22,15)$, $(1,10,3)$, $(1,0,9)$, $(1,13,18)$, 
$(1,15,10)$, $(1,4,15)$, $(1,0,10)$, $(1,16,5)$, $(1,20,7)$, $(1,2,19)$, $(1,0,11)$, 
$(1,11,19)$, $(1,8,5)$, $(1,1,21)$, $(1,0,12)$, $(1,15,17)$, $(1,6,20)$, $(1,5,13)$, 
$(1,0,14)$, $(1,12,7)$, $(1,11,17)$, $(1,3,17)$, $(1,0,15)$, $(1,6,10)$, $(1,7,1)$, 
$(1,11,1)$, $(1,0,17)$, $(1,7,21)$, $(1,1,0)$, $(1,13,20)$, $(1,0,20)$, $(1,21,14)$, 
$(1,18,22)$, $(1,6,11)$, $(1,0,21)$, $(1,2,12)$, $(1,5,16)$, $(1,12,22)$, $(1,0,22)$, 
$(1,20,3)$, $(1,3,8)$, $(1,8,7)$, $(1,1,2)$, $(1,1,17)$, $(1,19,10)$, $(1,2,9)$, 
$(1,1,5)$, $(1,21,22)$, $(1,18,14)$, $(1,19,13)$, $(1,1,7)$, $(1,3,6)$, $(1,12,15)$, 
$(1,20,20)$, $(1,1,9)$, $(1,6,1)$, $(1,13,11)$, $(1,5,7)$, $(1,1,10)$, $(1,5,18)$, 
$(1,5,20)$, $(1,4,0)$, $(1,1,12)$, $(1,7,7)$, $(1,21,2)$, $(1,18,6)$, $(1,1,13)$, 
$(1,16,15)$, $(1,1,14)$, $(1,10,2)$, $(1,22,21)$, $(1,6,14)$, $(1,1,15)$, $(1,22,5)$, 
$(1,3,5)$, $(1,13,17)$, $(1,1,16)$, $(1,15,9)$, $(1,16,22)$, $(1,3,16)$, $(1,1,18)$, 
$(1,2,0)$, $(1,4,1)$, $(1,10,19)$, $(1,1,19)$, $(1,12,14)$, $(1,15,3)$, $(1,15,8)$, 
$(1,2,1)$, $(1,12,6)$, $(1,8,16)$, $(1,11,11)$, $(1,2,3)$, $(1,6,8)$, $(1,19,21)$, 
$(1,16,19)$, $(1,2,7)$, $(1,7,0)$, $(1,2,11)$, $(1,22,18)$, $(1,13,12)$, $(1,4,9)$, 
$(1,2,13)$, $(1,21,3)$, $(1,18,8)$, $(1,13,5)$, $(1,2,14)$, $(1,13,21)$, $(1,7,3)$, 
$(1,12,8)$, $(1,2,15)$, $(1,2,17)$, $(1,20,11)$, $(1,14,2)$, $(1,2,18)$, $(1,11,14)$, 
$(1,4,10)$, $(1,3,12)$, $(1,2,21)$, $(1,8,15)$, $(1,21,1)$, $(1,18,13)$, $(1,3,1)$, 
$(1,15,21)$, $(1,3,2)$, $(1,7,5)$, $(1,7,19)$, $(1,6,12)$, $(1,3,3)$, $(1,20,8)$, 
$(1,22,6)$, $(1,5,11)$, $(1,3,4)$, $(1,17,2)$, $(1,19,4)$, $(1,17,0)$, $(1,3,7)$, 
$(1,22,12)$, $(1,21,13)$, $(1,18,1)$, $(1,3,13)$, $(1,10,11)$, $(1,13,0)$, $(1,19,2)$, 
$(1,3,14)$, $(1,5,1)$, $(1,14,16)$, $(1,7,13)$, $(1,3,15)$, $(1,14,19)$, $(1,10,21)$, 
$(1,10,16)$, $(1,3,19)$, $(1,16,0)$, $(1,8,12)$, $(1,8,14)$, $(1,3,22)$, $(1,4,22)$, 
$(1,16,2)$, $(1,12,18)$, $(1,4,2)$, $(1,11,20)$, $(1,16,18)$, $(1,10,15)$, $(1,4,5)$, 
$(1,15,1)$, $(1,19,18)$, $(1,14,11)$, $(1,4,6)$, $(1,21,7)$, $(1,18,18)$, $(1,16,9)$, 
$(1,4,7)$, $(1,19,5)$, $(1,10,18)$, $(1,12,13)$, $(1,4,8)$, $(1,14,0)$, $(1,6,18)$, 
$(1,22,3)$, $(1,4,18)$, $(1,8,17)$, $(1,4,19)$, $(1,13,22)$, $(1,11,18)$, $(1,6,19)$, 
$(1,4,20)$, $(1,7,16)$, $(1,20,18)$, $(1,20,5)$, $(1,5,2)$, $(1,16,10)$, $(1,21,20)$, 
$(1,18,21)$, $(1,5,4)$, $(1,17,7)$, $(1,12,4)$, $(1,17,19)$, $(1,5,5)$, $(1,8,11)$, 
$(1,16,6)$, $(1,7,22)$, $(1,5,9)$, $(1,15,13)$, $(1,8,2)$, $(1,20,2)$, $(1,5,15)$, 
$(1,7,14)$, $(1,22,9)$, $(1,12,9)$, $(1,5,17)$, $(1,11,2)$, $(1,13,16)$, $(1,8,1)$, 
$(1,6,0)$, $(1,10,6)$, $(1,15,20)$, $(1,12,12)$, $(1,6,2)$, $(1,6,7)$, $(1,20,15)$, 
$(1,19,16)$, $(1,6,4)$, $(1,17,10)$, $(1,8,4)$, $(1,17,5)$, $(1,6,5)$, $(1,16,16)$, 
$(1,14,21)$, $(1,13,6)$, $(1,6,9)$, $(1,11,0)$, $(1,12,0)$, $(1,16,11)$, $(1,6,21)$, 
$(1,12,17)$, $(1,22,13)$, $(1,20,10)$, $(1,7,4)$, $(1,17,20)$, $(1,22,4)$, $(1,17,1)$, 
$(1,7,10)$, $(1,14,22)$, $(1,15,19)$, $(1,20,9)$, $(1,7,17)$, $(1,15,6)$, $(1,11,21)$, 
$(1,22,22)$, $(1,7,20)$, $(1,19,11)$, $(1,21,16)$, $(1,18,19)$, $(1,8,6)$, $(1,13,9)$, 
$(1,21,21)$, $(1,18,20)$, $(1,8,10)$, $(1,22,0)$, $(1,22,16)$, $(1,19,7)$, $(1,8,22)$, 
$(1,10,12)$, $(1,14,10)$, $(1,20,17)$, $(1,10,4)$, $(1,17,17)$, $(1,13,4)$, $(1,17,6)$, 
$(1,10,20)$, $(1,20,1)$, $(1,16,3)$, $(1,13,8)$, $(1,11,4)$, $(1,17,11)$, $(1,15,4)$, 
$(1,17,13)$, $(1,11,6)$, $(1,20,13)$, $(1,13,13)$, $(1,22,19)$, $(1,11,16)$, $(1,21,6)$, 
$(1,18,2)$, $(1,14,14)$, $(1,11,22)$, $(1,19,20)$, $(1,12,1)$, $(1,21,5)$, $(1,18,0)$, 
$(1,20,22)$, $(1,12,11)$, $(1,19,12)$, $(1,22,20)$, $(1,14,15)$, $(1,12,16)$, $(1,13,10)$, 
$(1,14,4)$, $(1,17,9)$, $(1,14,17)$, $(1,21,12)$, $(1,18,15)$, $(1,15,0)$, $(1,17,4)$, 
$(1,18,3)$, $(1,21,8)$.

\subsection*{$m_{16}(2,25)\ge 366$}
$(0,0,1)$, $(1,18,10)$, $(1,21,11)$, $(1,16,6)$, $(1,21,8)$, $(1,14,8)$, $(1,21,22)$, 
$(1,24,21)$, $(1,21,19)$, $(1,3,3)$, $(0,1,2)$, $(1,3,22)$, $(1,18,19)$, $(1,13,6)$, 
$(1,7,22)$, $(1,4,17)$, $(1,10,8)$, $(1,11,16)$, $(1,4,11)$, $(1,20,5)$, $(0,1,3)$, 
$(1,2,14)$, $(1,10,22)$, $(1,3,6)$, $(1,18,11)$, $(1,23,19)$, $(1,4,19)$, $(1,7,1)$, 
$(1,7,8)$, $(1,9,15)$, $(0,1,4)$, $(1,15,6)$, $(0,1,5)$, $(1,0,18)$, $(1,13,24)$, 
$(1,9,6)$, $(1,17,12)$, $(0,1,6)$, $(1,1,1)$, $(1,6,1)$, $(1,24,6)$, $(1,1,3)$, 
$(1,5,15)$, $(1,16,2)$, $(1,0,19)$, $(1,11,4)$, $(1,7,14)$, $(0,1,8)$, $(1,9,11)$, 
$(1,14,23)$, $(1,20,6)$, $(1,15,14)$, $(1,1,2)$, $(1,2,16)$, $(1,14,13)$, $(1,8,7)$, 
$(1,10,7)$, $(0,1,9)$, $(1,8,3)$, $(1,5,2)$, $(1,17,6)$, $(1,3,1)$, $(1,18,21)$, 
$(1,19,4)$, $(1,20,10)$, $(1,12,3)$, $(1,5,8)$, $(0,1,13)$, $(1,20,12)$, $(1,1,6)$, 
$(1,23,6)$, $(1,11,18)$, $(1,24,24)$, $(1,6,12)$, $(1,22,18)$, $(1,16,24)$, $(0,1,18)$, 
$(0,1,14)$, $(1,12,21)$, $(1,3,9)$, $(1,18,6)$, $(1,12,17)$, $(1,13,3)$, $(1,5,13)$, 
$(1,1,8)$, $(1,19,21)$, $(1,2,10)$, $(0,1,15)$, $(1,4,5)$, $(1,16,16)$, $(1,6,6)$, 
$(1,6,23)$, $(1,17,16)$, $(1,11,7)$, $(1,2,22)$, $(1,1,14)$, $(1,24,17)$, $(0,1,16)$, 
$(1,7,15)$, $(1,2,5)$, $(1,11,6)$, $(1,14,15)$, $(1,19,1)$, $(1,8,10)$, $(1,17,14)$, 
$(1,15,20)$, $(1,14,19)$, $(0,1,20)$, $(1,13,9)$, $(1,8,4)$, $(1,14,6)$, $(1,2,2)$, 
$(1,6,20)$, $(1,15,3)$, $(1,5,23)$, $(1,14,1)$, $(1,23,4)$, $(0,1,23)$, $(1,11,13)$, 
$(1,9,3)$, $(1,5,6)$, $(1,0,4)$, $(1,2,7)$, $(1,17,1)$, $(1,10,2)$, $(1,13,2)$, 
$(1,22,11)$, $(0,1,24)$, $(1,22,8)$, $(1,12,20)$, $(1,19,6)$, $(1,19,10)$, $(1,9,10)$, 
$(1,3,15)$, $(1,18,3)$, $(1,5,5)$, $(1,6,21)$, $(1,0,0)$, $(1,6,13)$, $(1,20,0)$, 
$(1,24,22)$, $(1,19,0)$, $(1,17,10)$, $(1,16,0)$, $(1,4,1)$, $(1,2,0)$, $(1,2,23)$, 
$(1,0,8)$, $(1,4,2)$, $(1,15,9)$, $(1,13,15)$, $(1,23,18)$, $(1,3,24)$, $(1,18,16)$, 
$(1,20,22)$, $(1,11,1)$, $(1,11,9)$, $(1,0,11)$, $(1,8,8)$, $(1,6,3)$, $(1,5,4)$, 
$(1,8,9)$, $(1,20,17)$, $(1,3,5)$, $(1,18,15)$, $(1,10,14)$, $(1,17,7)$, $(1,0,16)$, 
$(1,22,16)$, $(1,3,4)$, $(1,18,18)$, $(1,6,24)$, $(1,19,9)$, $(1,10,10)$, $(1,6,8)$, 
$(1,8,11)$, $(1,9,1)$, $(1,0,17)$, $(1,24,11)$, $(1,2,18)$, $(1,19,24)$, $(1,16,10)$, 
$(1,15,16)$, $(1,19,14)$, $(1,13,10)$, $(1,20,2)$, $(1,10,4)$, $(1,0,21)$, $(1,12,14)$, 
$(1,0,24)$, $(1,20,21)$, $(1,14,2)$, $(1,9,23)$, $(1,12,22)$, $(1,10,15)$, $(1,22,15)$, 
$(1,14,11)$, $(1,4,9)$, $(1,1,18)$, $(1,1,5)$, $(1,12,8)$, $(1,22,1)$, $(1,13,18)$, 
$(1,15,24)$, $(1,2,1)$, $(1,3,7)$, $(1,18,9)$, $(1,17,21)$, $(1,1,16)$, $(1,1,9)$, 
$(1,14,0)$, $(1,10,3)$, $(1,5,0)$, $(1,20,11)$, $(1,13,0)$, $(1,14,14)$, $(1,6,0)$, 
$(1,9,5)$, $(1,23,0)$, $(1,1,11)$, $(1,17,4)$, $(1,9,4)$, $(1,19,13)$, $(1,14,10)$, 
$(1,11,19)$, $(1,20,16)$, $(1,7,18)$, $(1,10,24)$, $(1,14,5)$, $(1,1,12)$, $(1,11,22)$, 
$(1,13,23)$, $(1,3,21)$, $(1,18,5)$, $(1,4,12)$, $(1,12,2)$, $(1,20,14)$, $(1,2,8)$, 
$(1,13,4)$, $(1,1,15)$, $(1,8,23)$, $(1,15,12)$, $(1,20,8)$, $(1,17,13)$, $(1,17,5)$, 
$(1,22,3)$, $(1,5,12)$, $(1,3,16)$, $(1,18,2)$, $(1,1,21)$, $(1,7,12)$, $(1,8,14)$, 
$(1,12,9)$, $(1,23,22)$, $(1,20,7)$, $(1,19,12)$, $(1,15,10)$, $(1,4,4)$, $(1,6,19)$, 
$(1,1,24)$, $(1,21,14)$, $(1,4,15)$, $(1,21,17)$, $(1,19,2)$, $(1,21,10)$, $(1,23,9)$, 
$(1,21,2)$, $(1,8,22)$, $(1,21,18)$, $(1,2,3)$, $(1,5,11)$, $(1,11,14)$, $(1,22,5)$, 
$(1,10,9)$, $(1,2,20)$, $(1,9,2)$, $(1,4,18)$, $(1,4,24)$, $(1,17,3)$, $(1,5,9)$, 
$(1,24,7)$, $(1,15,15)$, $(1,6,22)$, $(1,16,22)$, $(1,2,24)$, $(1,8,18)$, $(1,3,13)$, 
$(1,18,1)$, $(1,14,20)$, $(1,3,18)$, $(1,18,24)$, $(1,24,8)$, $(1,24,14)$, $(1,14,16)$, 
$(1,16,8)$, $(1,16,9)$, $(1,4,21)$, $(1,11,15)$, $(1,15,4)$, $(1,12,16)$, $(1,9,12)$, 
$(1,24,3)$, $(1,5,22)$, $(1,9,13)$, $(1,16,14)$, $(1,8,12)$, $(1,4,22)$, $(1,22,4)$, 
$(1,4,23)$, $(1,13,14)$, $(1,16,1)$, $(1,19,3)$, $(1,5,10)$, $(1,10,5)$, $(1,9,22)$, 
$(1,15,11)$, $(1,15,13)$, $(1,24,9)$, $(1,5,16)$, $(1,23,3)$, $(1,6,15)$, $(1,13,13)$, 
$(1,6,18)$, $(1,21,24)$, $(1,12,4)$, $(1,21,5)$, $(1,15,19)$, $(1,21,7)$, $(1,7,5)$, 
$(1,21,1)$, $(1,20,13)$, $(1,21,21)$, $(1,7,2)$, $(1,11,0)$, $(1,16,21)$, $(1,9,0)$, 
$(1,13,20)$, $(1,12,0)$, $(1,8,16)$, $(1,8,0)$, $(1,22,19)$, $(1,24,0)$, $(1,8,24)$, 
$(1,23,15)$, $(1,9,21)$, $(1,10,21)$, $(1,24,1)$, $(1,12,11)$, $(1,12,7)$, $(1,16,23)$, 
$(1,11,5)$, $(1,17,18)$, $(1,10,16)$, $(1,15,21)$, $(1,17,11)$, $(1,19,15)$, $(1,12,10)$, 
$(1,24,18)$, $(1,23,24)$, $(1,10,23)$, $(1,16,4)$, $(1,13,11)$, $(1,11,10)$, $(1,19,7)$, 
$(1,16,5)$, $(1,24,2)$.

\subsection*{$m_{17}(2,25)\ge 393$}
$(0,1,0)$, $(1,2,12)$, $(0,1,11)$, $(1,23,22)$, $(1,8,20)$, $(1,20,11)$, $(1,1,17)$, 
$(1,5,23)$, $(0,1,3)$, $(1,5,13)$, $(1,7,14)$, $(1,16,16)$, $(1,10,3)$, $(1,9,10)$, 
$(1,18,17)$, $(1,15,18)$, $(0,1,4)$, $(1,24,11)$, $(1,13,16)$, $(1,8,1)$, $(1,2,22)$, 
$(1,1,0)$, $(1,9,17)$, $(1,12,7)$, $(0,1,5)$, $(1,12,20)$, $(1,5,17)$, $(1,3,14)$, 
$(0,1,6)$, $(1,19,7)$, $(1,4,10)$, $(1,21,8)$, $(1,11,8)$, $(1,8,2)$, $(1,4,17)$, 
$(1,1,10)$, $(0,1,7)$, $(1,6,2)$, $(1,12,5)$, $(1,22,0)$, $(1,5,5)$, $(1,18,14)$, 
$(1,3,17)$, $(1,9,21)$, $(0,1,8)$, $(1,17,4)$, $(1,10,13)$, $(1,17,13)$, $(1,23,9)$, 
$(1,15,5)$, $(1,24,17)$, $(1,13,9)$, $(0,1,9)$, $(1,8,5)$, $(1,6,3)$, $(1,1,20)$, 
$(1,17,6)$, $(1,22,7)$, $(1,22,17)$, $(1,6,20)$, $(0,1,10)$, $(1,1,23)$, $(1,18,12)$, 
$(1,10,7)$, $(1,20,19)$, $(1,3,21)$, $(1,15,17)$, $(1,8,24)$, $(0,1,12)$, $(1,22,8)$, 
$(1,0,21)$, $(1,0,3)$, $(1,18,11)$, $(1,7,19)$, $(1,21,17)$, $(1,7,22)$, $(0,1,13)$, 
$(1,10,17)$, $(1,16,15)$, $(0,1,22)$, $(1,15,21)$, $(1,17,1)$, $(1,2,17)$, $(1,19,16)$, 
$(0,1,16)$, $(1,11,6)$, $(1,23,8)$, $(1,7,9)$, $(1,22,4)$, $(1,16,18)$, $(1,12,17)$, 
$(1,23,4)$, $(0,1,17)$, $(1,4,15)$, $(1,11,24)$, $(1,9,23)$, $(1,0,12)$, $(1,4,9)$, 
$(1,6,17)$, $(1,20,3)$, $(0,1,20)$, $(1,18,18)$, $(1,22,22)$, $(1,4,6)$, $(1,16,1)$, 
$(1,2,13)$, $(1,11,17)$, $(1,4,11)$, $(0,1,21)$, $(1,9,24)$, $(1,21,11)$, $(1,12,21)$, 
$(1,7,15)$, $(1,19,22)$, $(1,20,17)$, $(1,11,5)$, $(1,0,0)$, $(1,21,20)$, $(1,22,24)$, 
$(1,11,9)$, $(1,3,15)$, $(1,22,20)$, $(1,3,16)$, $(1,20,10)$, $(1,0,8)$, $(1,23,13)$, 
$(1,3,6)$, $(1,17,2)$, $(1,9,1)$, $(1,2,21)$, $(1,9,4)$, $(1,3,3)$, $(1,0,14)$, 
$(1,22,19)$, $(1,19,13)$, $(1,8,8)$, $(1,5,18)$, $(1,1,14)$, $(1,7,5)$, $(1,5,7)$, 
$(1,0,16)$, $(1,7,18)$, $(1,24,23)$, $(1,13,12)$, $(1,24,13)$, $(1,13,5)$, $(1,12,6)$, 
$(1,6,16)$, $(1,0,18)$, $(1,3,10)$, $(1,8,16)$, $(1,22,6)$, $(1,1,11)$, $(1,19,18)$, 
$(1,17,7)$, $(1,17,22)$, $(1,0,19)$, $(1,5,0)$, $(1,17,14)$, $(1,6,0)$, $(1,18,2)$, 
$(1,21,13)$, $(1,19,1)$, $(1,2,19)$, $(1,1,1)$, $(1,2,7)$, $(1,7,7)$, $(1,12,2)$, 
$(1,18,9)$, $(1,8,3)$, $(1,8,23)$, $(1,1,12)$, $(1,1,4)$, $(1,4,19)$, $(1,8,7)$, 
$(1,23,16)$, $(1,1,6)$, $(1,3,13)$, $(1,18,7)$, $(1,24,15)$, $(1,13,13)$, $(1,6,15)$, 
$(1,8,21)$, $(1,19,10)$, $(1,1,8)$, $(1,16,0)$, $(1,20,7)$, $(1,19,20)$, $(1,10,5)$, 
$(1,7,24)$, $(1,8,13)$, $(1,15,20)$, $(1,1,15)$, $(1,21,3)$, $(1,16,7)$, $(1,6,8)$, 
$(1,17,23)$, $(1,12,12)$, $(1,8,19)$, $(1,17,15)$, $(1,1,18)$, $(1,6,4)$, $(1,21,7)$, 
$(1,3,11)$, $(1,23,0)$, $(1,16,21)$, $(1,8,22)$, $(1,7,13)$, $(1,1,19)$, $(1,23,10)$, 
$(1,11,7)$, $(1,18,21)$, $(1,11,16)$, $(1,17,0)$, $(1,8,4)$, $(1,9,8)$, $(1,1,24)$, 
$(1,24,16)$, $(1,13,7)$, $(1,11,3)$, $(1,7,6)$, $(1,20,5)$, $(1,8,15)$, $(1,18,0)$, 
$(1,2,0)$, $(1,23,1)$, $(1,2,15)$, $(1,5,10)$, $(1,22,14)$, $(1,9,15)$, $(1,7,4)$, 
$(1,7,1)$, $(1,2,2)$, $(1,20,23)$, $(1,15,12)$, $(1,22,16)$, $(1,12,24)$, $(1,17,5)$, 
$(1,11,10)$, $(1,10,1)$, $(1,2,4)$, $(1,15,2)$, $(1,7,21)$, $(1,23,11)$, $(1,21,10)$, 
$(1,6,9)$, $(1,10,22)$, $(1,6,1)$, $(1,2,5)$, $(1,21,9)$, $(1,20,18)$, $(1,4,13)$, 
$(1,9,2)$, $(1,19,24)$, $(1,18,15)$, $(1,4,1)$, $(1,2,8)$, $(1,22,15)$, $(1,17,8)$, 
$(1,24,6)$, $(1,13,23)$, $(1,5,12)$, $(1,16,14)$, $(1,18,1)$, $(1,2,14)$, $(1,12,3)$, 
$(1,21,1)$, $(1,2,23)$, $(1,24,12)$, $(1,13,10)$, $(1,3,18)$, $(1,5,1)$, $(1,2,20)$, 
$(1,6,21)$, $(1,6,13)$, $(1,16,19)$, $(1,11,20)$, $(1,12,18)$, $(1,19,8)$, $(1,12,1)$, 
$(1,2,24)$, $(1,3,22)$, $(1,24,5)$, $(1,13,2)$, $(1,17,19)$, $(1,20,8)$, $(1,6,11)$, 
$(1,11,1)$, $(1,3,5)$, $(1,16,22)$, $(1,21,4)$, $(1,24,0)$, $(1,13,3)$, $(1,23,5)$, 
$(1,9,11)$, $(1,16,8)$, $(1,3,12)$, $(1,3,23)$, $(1,19,12)$, $(1,11,18)$, $(1,16,9)$, 
$(1,24,10)$, $(1,13,19)$, $(1,7,23)$, $(1,3,20)$, $(1,21,15)$, $(1,15,19)$, $(1,18,4)$, 
$(1,18,20)$, $(1,16,2)$, $(1,23,18)$, $(1,10,16)$, $(1,3,24)$, $(1,10,24)$, $(1,5,8)$, 
$(1,19,14)$, $(1,15,11)$, $(1,20,15)$, $(1,20,6)$, $(1,12,22)$, $(1,4,2)$, $(1,22,3)$, 
$(1,24,22)$, $(1,13,21)$, $(1,4,21)$, $(1,24,18)$, $(1,13,15)$, $(1,16,11)$, $(1,4,4)$, 
$(1,19,11)$, $(1,5,4)$, $(1,23,14)$, $(1,4,22)$, $(1,23,23)$, $(1,17,12)$, $(1,21,5)$, 
$(1,4,5)$, $(1,20,13)$, $(1,22,10)$, $(1,7,10)$, $(1,4,24)$, $(1,12,19)$, $(1,16,6)$, 
$(1,9,19)$, $(1,5,2)$, $(1,15,14)$, $(1,5,20)$, $(1,11,14)$, $(1,24,4)$, $(1,13,6)$, 
$(1,10,20)$, $(1,23,3)$, $(1,5,16)$, $(1,9,14)$, $(1,20,9)$, $(1,17,16)$, $(1,19,21)$, 
$(1,16,24)$, $(1,22,21)$, $(1,21,2)$, $(1,5,19)$, $(1,6,14)$, $(1,19,6)$, $(1,24,9)$, 
$(1,13,11)$, $(1,9,9)$, $(1,11,0)$, $(1,11,22)$, $(1,6,19)$, $(1,24,24)$, $(1,13,4)$, 
$(1,20,16)$, $(1,21,23)$, $(1,23,12)$, $(1,20,24)$, $(1,7,3)$, $(1,6,22)$, $(1,22,11)$, 
$(1,18,23)$, $(1,9,12)$, $(1,10,10)$, $(1,24,21)$, $(1,13,0)$, $(1,15,8)$, $(1,12,0)$, 
$(1,12,10)$, $(1,20,2)$, $(1,16,4)$, $(1,17,9)$, $(1,18,24)$, $(1,23,6)$, $(1,18,3)$, 
$(1,21,6)$.

\subsection*{$m_{18}(2,25)\ge 416$}
$(0,0,1)$, $(1,17,5)$, $(1,9,23)$, $(1,0,1)$, $(1,2,23)$, $(1,23,3)$, $(1,10,11)$, 
$(1,3,6)$, $(1,13,2)$, $(1,24,14)$, $(1,8,21)$, $(1,19,1)$, $(0,1,1)$, $(1,12,5)$, 
$(1,0,13)$, $(1,14,15)$, $(1,10,12)$, $(1,15,16)$, $(1,24,13)$, $(1,1,19)$, $(1,7,16)$, 
$(1,0,7)$, $(1,10,7)$, $(1,24,11)$, $(0,1,5)$, $(1,21,5)$, $(1,10,0)$, $(1,18,2)$, 
$(1,7,21)$, $(1,14,11)$, $(1,22,0)$, $(1,11,13)$, $(1,10,6)$, $(1,16,13)$, $(1,20,14)$, 
$(1,15,3)$, $(0,1,6)$, $(1,20,5)$, $(1,20,17)$, $(1,8,18)$, $(1,24,7)$, $(1,9,22)$, 
$(1,16,8)$, $(1,0,23)$, $(0,1,17)$, $(1,19,5)$, $(1,17,18)$, $(1,5,8)$, $(0,1,8)$, 
$(1,8,5)$, $(1,21,3)$, $(1,9,13)$, $(1,19,9)$, $(1,17,23)$, $(1,3,22)$, $(1,14,21)$, 
$(1,12,10)$, $(1,22,6)$, $(1,5,16)$, $(1,23,14)$, $(0,1,10)$, $(1,6,5)$, $(0,1,11)$, 
$(1,11,5)$, $(1,6,15)$, $(1,13,20)$, $(1,12,24)$, $(1,12,9)$, $(1,4,1)$, $(1,2,10)$, 
$(1,19,13)$, $(1,17,2)$, $(0,1,14)$, $(1,10,5)$, $(1,7,1)$, $(1,2,16)$, $(1,16,16)$, 
$(1,4,8)$, $(1,6,2)$, $(1,10,17)$, $(1,5,12)$, $(1,4,18)$, $(1,11,17)$, $(1,13,19)$, 
$(0,1,15)$, $(1,24,5)$, $(1,3,16)$, $(1,15,17)$, $(1,0,11)$, $(1,5,13)$, $(1,8,10)$, 
$(1,23,24)$, $(1,9,20)$, $(1,5,11)$, $(1,15,23)$, $(1,16,0)$, $(0,1,16)$, $(1,0,5)$, 
$(1,4,2)$, $(1,20,24)$, $(1,3,4)$, $(1,1,10)$, $(1,11,15)$, $(0,1,21)$, $(1,1,5)$, 
$(1,23,20)$, $(1,24,4)$, $(1,20,13)$, $(0,1,19)$, $(1,3,5)$, $(1,18,4)$, $(0,1,20)$, 
$(1,4,5)$, $(1,16,7)$, $(1,23,9)$, $(1,13,0)$, $(1,6,4)$, $(1,18,16)$, $(1,6,1)$, 
$(1,2,22)$, $(0,1,23)$, $(1,14,5)$, $(1,11,11)$, $(1,19,22)$, $(1,17,22)$, $(1,19,0)$, 
$(1,17,12)$, $(1,16,20)$, $(1,15,0)$, $(1,7,19)$, $(1,12,2)$, $(1,1,20)$, $(0,1,24)$, 
$(1,13,5)$, $(1,13,8)$, $(1,3,11)$, $(1,23,1)$, $(1,2,1)$, $(1,2,18)$, $(1,7,22)$, 
$(1,23,15)$, $(1,15,24)$, $(1,14,22)$, $(1,3,24)$, $(1,0,2)$, $(1,16,18)$, $(1,1,13)$, 
$(1,12,13)$, $(1,11,22)$, $(1,24,9)$, $(1,18,6)$, $(1,4,9)$, $(1,0,14)$, $(1,24,8)$, 
$(1,4,17)$, $(1,21,22)$, $(1,0,3)$, $(1,3,10)$, $(1,12,8)$, $(1,7,8)$, $(1,16,3)$, 
$(1,12,21)$, $(1,23,18)$, $(1,22,18)$, $(1,4,24)$, $(1,24,21)$, $(1,7,20)$, $(1,22,13)$, 
$(1,0,4)$, $(1,11,24)$, $(1,15,18)$, $(1,19,15)$, $(1,17,11)$, $(1,4,7)$, $(1,5,4)$, 
$(1,5,22)$, $(1,18,11)$, $(1,24,2)$, $(1,15,8)$, $(1,18,10)$, $(1,0,6)$, $(1,19,14)$, 
$(1,17,14)$, $(1,1,2)$, $(1,21,9)$, $(1,21,16)$, $(1,13,24)$, $(1,11,10)$, $(1,20,0)$, 
$(1,24,1)$, $(1,2,15)$, $(1,20,6)$, $(1,0,9)$, $(1,22,4)$, $(1,16,1)$, $(1,2,3)$, 
$(1,18,24)$, $(1,18,0)$, $(1,12,3)$, $(1,19,4)$, $(1,17,21)$, $(1,24,16)$, $(1,22,10)$, 
$(1,13,21)$, $(1,0,10)$, $(1,15,1)$, $(1,2,21)$, $(1,11,12)$, $(1,22,17)$, $(1,20,10)$, 
$(1,16,19)$, $(1,13,11)$, $(1,14,17)$, $(1,24,3)$, $(1,13,1)$, $(1,2,2)$, $(1,0,15)$, 
$(1,7,17)$, $(1,11,20)$, $(1,14,10)$, $(1,6,16)$, $(1,8,4)$, $(1,3,0)$, $(1,9,24)$, 
$(1,6,13)$, $(1,24,23)$, $(1,14,2)$, $(1,23,4)$, $(1,0,17)$, $(1,23,16)$, $(1,3,9)$, 
$(1,24,20)$, $(1,19,7)$, $(1,17,24)$, $(1,0,21)$, $(1,9,21)$, $(1,22,15)$, $(1,22,23)$, 
$(1,7,4)$, $(1,6,17)$, $(1,19,2)$, $(1,17,3)$, $(1,7,3)$, $(1,24,22)$, $(1,11,4)$, 
$(1,8,7)$, $(1,1,0)$, $(1,4,10)$, $(1,4,11)$, $(1,12,15)$, $(1,6,21)$, $(1,3,8)$, 
$(1,22,7)$, $(1,14,0)$, $(1,15,7)$, $(1,20,18)$, $(1,23,10)$, $(1,5,3)$, $(1,1,1)$, 
$(1,2,24)$, $(1,5,19)$, $(1,16,4)$, $(1,7,12)$, $(1,5,10)$, $(1,7,9)$, $(1,19,19)$, 
$(1,17,13)$, $(1,13,4)$, $(1,4,20)$, $(1,1,12)$, $(1,1,3)$, $(1,6,23)$, $(1,15,13)$, 
$(1,15,20)$, $(1,13,17)$, $(1,16,21)$, $(1,6,22)$, $(1,22,1)$, $(1,2,11)$, $(1,8,11)$, 
$(1,1,9)$, $(1,9,7)$, $(1,1,7)$, $(1,22,2)$, $(1,3,15)$, $(1,4,22)$, $(1,8,3)$, 
$(1,11,16)$, $(1,9,3)$, $(1,23,17)$, $(1,1,23)$, $(1,19,21)$, $(1,17,0)$, $(1,22,21)$, 
$(1,1,8)$, $(1,8,14)$, $(1,6,10)$, $(1,9,10)$, $(1,21,18)$, $(1,18,23)$, $(1,8,16)$, 
$(1,5,2)$, $(1,9,16)$, $(1,12,20)$, $(1,16,23)$, $(1,11,6)$, $(1,1,15)$, $(1,5,0)$, 
$(1,19,17)$, $(1,17,8)$, $(1,10,24)$, $(1,19,24)$, $(1,17,6)$, $(1,6,18)$, $(1,12,14)$, 
$(1,18,17)$, $(1,10,13)$, $(1,21,0)$, $(1,1,18)$, $(1,13,6)$, $(1,10,1)$, $(1,2,14)$, 
$(1,12,1)$, $(1,2,7)$, $(1,21,20)$, $(1,8,20)$, $(1,10,8)$, $(1,9,15)$, $(1,18,7)$, 
$(1,13,23)$, $(1,1,22)$, $(1,15,19)$, $(1,23,8)$, $(1,19,16)$, $(1,17,15)$, $(1,8,13)$, 
$(1,18,18)$, $(1,9,11)$, $(1,7,10)$, $(1,3,23)$, $(1,6,12)$, $(1,20,9)$, $(1,2,4)$, 
$(1,14,9)$, $(1,14,13)$, $(1,5,17)$, $(1,12,7)$, $(1,8,15)$, $(1,12,4)$, $(1,9,2)$, 
$(1,13,22)$, $(1,20,23)$, $(1,5,20)$, $(1,11,1)$, $(1,2,6)$, $(1,3,7)$, $(1,11,14)$, 
$(1,4,16)$, $(1,7,0)$, $(1,21,6)$, $(1,23,22)$, $(1,7,14)$, $(1,13,9)$, $(1,5,15)$, 
$(1,21,2)$, $(1,5,1)$, $(1,2,17)$, $(1,9,0)$, $(1,20,16)$, $(1,14,23)$, $(1,18,9)$, 
$(1,16,10)$, $(1,18,20)$, $(1,20,22)$, $(1,13,7)$, $(1,7,13)$, $(1,7,24)$, $(1,8,1)$, 
$(1,3,3)$, $(1,5,18)$, $(1,10,20)$, $(1,6,9)$, $(1,23,12)$, $(1,23,21)$, $(1,16,2)$, 
$(1,22,19)$, $(1,10,9)$, $(1,20,8)$, $(1,11,0)$, $(1,3,20)$, $(1,3,19)$, $(1,8,0)$, 
$(1,16,22)$, $(1,9,17)$, $(1,4,4)$, $(1,23,23)$, $(1,20,7)$, $(1,12,16)$, $(1,20,12)$, 
$(1,19,12)$, $(1,17,4)$, $(1,12,19)$, $(1,4,3)$, $(1,15,9)$, $(1,8,24)$, $(1,20,15)$, 
$(1,9,9)$, $(1,10,21)$, $(1,8,23)$, $(1,22,11)$, $(1,9,12)$, $(1,14,24)$, $(1,4,23)$, 
$(1,10,19)$, $(1,11,2)$, $(1,16,15)$.

\subsection*{$m_{17}(2,27)\ge 421$}
$(0,0,1)$, $(1,9,1)$, $(1,2,3)$, $(1,8,23)$, $(1,10,6)$, $(1,1,26)$, $(1,26,19)$, 
$(0,1,0)$, $(1,21,15)$, $(1,22,9)$, $(1,15,25)$, $(1,11,22)$, $(1,4,16)$, $(1,5,18)$, 
$(0,1,1)$, $(1,25,19)$, $(1,21,20)$, $(1,4,2)$, $(1,14,14)$, $(1,12,3)$, $(1,16,13)$, 
$(0,1,3)$, $(1,8,8)$, $(1,5,4)$, $(1,24,1)$, $(1,1,8)$, $(1,10,18)$, $(1,2,15)$, 
$(0,1,6)$, $(1,14,2)$, $(1,0,22)$, $(1,18,5)$, $(1,12,0)$, $(1,14,25)$, $(1,7,2)$, 
$(0,1,7)$, $(1,5,10)$, $(1,23,1)$, $(1,14,4)$, $(1,16,17)$, $(1,5,21)$, $(1,13,1)$, 
$(0,1,9)$, $(1,1,6)$, $(1,6,21)$, $(1,2,24)$, $(1,15,19)$, $(1,21,7)$, $(1,14,5)$, 
$(0,1,10)$, $(1,18,20)$, $(1,7,13)$, $(1,26,11)$, $(1,9,11)$, $(1,26,10)$, $(1,9,21)$, 
$(0,1,13)$, $(1,24,4)$, $(1,4,12)$, $(1,0,17)$, $(1,17,3)$, $(1,6,19)$, $(1,21,11)$, 
$(0,1,14)$, $(1,23,18)$, $(1,13,15)$, $(1,25,18)$, $(1,25,15)$, $(1,24,24)$, $(1,6,7)$, 
$(0,1,15)$, $(1,19,17)$, $(1,18,19)$, $(1,21,6)$, $(1,7,16)$, $(1,23,20)$, $(1,15,9)$, 
$(0,1,16)$, $(1,17,24)$, $(1,10,17)$, $(1,20,12)$, $(1,19,7)$, $(1,8,14)$, $(1,17,17)$, 
$(0,1,20)$, $(1,2,21)$, $(1,20,0)$, $(1,23,16)$, $(0,1,25)$, $(1,0,9)$, $(1,1,14)$, 
$(0,1,22)$, $(1,3,7)$, $(1,15,24)$, $(1,12,21)$, $(1,26,4)$, $(1,9,13)$, $(1,12,6)$, 
$(0,1,23)$, $(1,26,16)$, $(1,9,25)$, $(1,3,9)$, $(1,5,13)$, $(1,16,1)$, $(1,3,22)$, 
$(1,0,0)$, $(1,24,12)$, $(1,26,26)$, $(1,9,16)$, $(1,22,4)$, $(1,14,9)$, $(1,23,0)$, 
$(1,0,1)$, $(1,6,16)$, $(1,1,13)$, $(1,24,9)$, $(1,2,20)$, $(1,5,14)$, $(1,5,1)$, 
$(1,0,5)$, $(1,25,20)$, $(1,16,21)$, $(1,18,10)$, $(1,24,8)$, $(1,17,0)$, $(1,12,24)$, 
$(1,0,6)$, $(1,8,2)$, $(1,20,18)$, $(1,15,15)$, $(1,3,11)$, $(1,10,21)$, $(1,0,7)$, 
$(1,0,14)$, $(1,14,13)$, $(1,3,24)$, $(1,19,2)$, $(1,17,5)$, $(1,7,8)$, $(1,24,17)$, 
$(1,0,18)$, $(1,10,15)$, $(1,14,10)$, $(1,17,26)$, $(1,8,12)$, $(1,22,10)$, $(1,20,11)$, 
$(1,0,19)$, $(1,21,22)$, $(1,23,22)$, $(1,16,7)$, $(1,6,25)$, $(1,20,25)$, $(1,16,5)$, 
$(1,0,20)$, $(1,3,26)$, $(1,0,23)$, $(1,22,0)$, $(1,5,7)$, $(1,25,1)$, $(1,11,10)$, 
$(1,0,21)$, $(1,12,18)$, $(1,11,15)$, $(1,1,5)$, $(1,18,13)$, $(1,2,23)$, $(1,13,13)$, 
$(1,1,3)$, $(1,25,2)$, $(1,2,25)$, $(1,19,6)$, $(1,19,18)$, $(1,4,15)$, $(1,15,22)$, 
$(1,1,9)$, $(1,12,5)$, $(1,2,10)$, $(1,26,18)$, $(1,9,15)$, $(1,13,16)$, $(1,19,12)$, 
$(1,1,10)$, $(1,19,3)$, $(1,2,11)$, $(1,11,7)$, $(1,16,2)$, $(1,16,8)$, $(1,4,23)$, 
$(1,1,17)$, $(1,16,26)$, $(1,2,18)$, $(1,22,15)$, $(1,12,13)$, $(1,20,14)$, $(1,20,21)$, 
$(1,1,20)$, $(1,11,11)$, $(1,2,17)$, $(1,14,16)$, $(1,25,23)$, $(1,6,9)$, $(1,16,4)$, 
$(1,1,21)$, $(1,17,16)$, $(1,2,2)$, $(1,23,9)$, $(1,3,19)$, $(1,21,19)$, $(1,21,8)$, 
$(1,1,22)$, $(1,5,0)$, $(1,2,1)$, $(1,17,25)$, $(1,13,4)$, $(1,7,7)$, $(1,17,13)$, 
$(1,1,24)$, $(1,23,24)$, $(1,2,14)$, $(1,4,11)$, $(1,17,20)$, $(1,23,3)$, $(1,12,20)$, 
$(1,1,25)$, $(1,22,7)$, $(1,2,22)$, $(1,7,20)$, $(1,14,22)$, $(1,19,25)$, $(1,23,26)$, 
$(1,3,0)$, $(1,3,6)$, $(1,11,17)$, $(1,10,13)$, $(1,15,0)$, $(1,16,23)$, $(1,8,24)$, 
$(1,3,1)$, $(1,4,20)$, $(1,20,19)$, $(1,21,3)$, $(1,18,18)$, $(1,19,15)$, $(1,5,17)$, 
$(1,3,4)$, $(1,17,4)$, $(1,18,23)$, $(1,25,16)$, $(1,4,5)$, $(1,13,9)$, $(1,11,3)$, 
$(1,3,10)$, $(1,12,14)$, $(1,6,22)$, $(1,24,19)$, $(1,21,4)$, $(1,3,25)$, $(1,6,10)$, 
$(1,3,12)$, $(1,16,11)$, $(1,7,26)$, $(1,19,5)$, $(1,22,24)$, $(1,25,4)$, $(1,12,8)$, 
$(1,3,13)$, $(1,13,7)$, $(1,23,4)$, $(1,6,11)$, $(1,19,14)$, $(1,7,12)$, $(1,20,2)$, 
$(1,3,15)$, $(1,8,26)$, $(1,14,20)$, $(1,22,18)$, $(1,7,15)$, $(1,10,7)$, $(1,13,18)$, 
$(1,4,0)$, $(1,26,13)$, $(1,9,3)$, $(1,7,24)$, $(1,4,18)$, $(1,24,15)$, $(1,20,3)$, 
$(1,4,6)$, $(1,25,11)$, $(1,15,3)$, $(1,15,17)$, $(1,13,0)$, $(1,20,10)$, $(1,14,6)$, 
$(1,4,9)$, $(1,8,6)$, $(1,26,3)$, $(1,9,5)$, $(1,15,5)$, $(1,10,16)$, $(1,17,11)$, 
$(1,4,13)$, $(1,4,25)$, $(1,10,3)$, $(1,24,3)$, $(1,23,10)$, $(1,8,22)$, $(1,13,21)$, 
$(1,4,19)$, $(1,21,0)$, $(1,6,3)$, $(1,5,8)$, $(1,14,1)$, $(1,26,5)$, $(1,9,7)$, 
$(1,5,11)$, $(1,18,1)$, $(1,20,26)$, $(1,11,14)$, $(1,22,8)$, $(1,6,14)$, $(1,13,23)$, 
$(1,5,16)$, $(1,20,1)$, $(1,8,9)$, $(1,24,2)$, $(1,18,24)$, $(1,16,10)$, $(1,10,2)$, 
$(1,5,20)$, $(1,6,1)$, $(1,10,12)$, $(1,25,22)$, $(1,14,17)$, $(1,8,11)$, $(1,12,25)$, 
$(1,5,23)$, $(1,26,1)$, $(1,9,17)$, $(1,19,8)$, $(1,11,2)$, $(1,19,21)$, $(1,25,21)$, 
$(1,6,4)$, $(1,11,26)$, $(1,13,24)$, $(1,22,6)$, $(1,23,23)$, $(1,18,2)$, $(1,15,2)$, 
$(1,6,20)$, $(1,26,17)$, $(1,9,6)$, $(1,16,12)$, $(1,12,4)$, $(1,15,18)$, $(1,14,15)$, 
$(1,6,24)$, $(1,8,16)$, $(1,18,26)$, $(1,10,9)$, $(1,22,21)$, $(1,10,5)$, $(1,8,17)$, 
$(1,7,10)$, $(1,22,23)$, $(1,16,20)$, $(1,17,6)$, $(1,17,7)$, $(1,11,21)$, $(1,15,26)$, 
$(1,7,18)$, $(1,8,15)$, $(1,7,23)$, $(1,12,19)$, $(1,21,16)$, $(1,10,14)$, $(1,11,9)$, 
$(1,7,19)$, $(1,21,26)$, $(1,25,26)$, $(1,22,11)$, $(1,22,19)$, $(1,21,24)$, $(1,24,13)$, 
$(1,7,22)$, $(1,11,4)$, $(1,20,6)$, $(1,24,25)$, $(1,14,23)$, $(1,15,8)$, $(1,20,8)$, 
$(1,8,13)$, $(1,11,5)$, $(1,11,25)$, $(1,12,12)$, $(1,10,23)$, $(1,19,23)$, $(1,24,22)$, 
$(1,9,0)$, $(1,17,22)$, $(1,15,14)$, $(1,18,14)$, $(1,23,25)$, $(1,25,8)$, $(1,26,6)$, 
$(1,9,12)$, $(1,23,14)$, $(1,24,10)$, $(1,11,16)$, $(1,24,7)$, $(1,19,9)$, $(1,26,23)$, 
$(1,9,14)$, $(1,25,0)$, $(1,13,26)$, $(1,12,22)$, $(1,22,2)$, $(1,13,17)$, $(1,26,25)$, 
$(1,18,7)$.

\subsection*{$m_{20}(2,27)\ge 510$}
$(0,1,0)$, $(1,0,14)$, $(1,11,21)$, $(1,10,19)$, $(1,21,19)$, $(1,6,20)$, $(1,26,3)$, 
$(1,15,10)$, $(1,5,4)$, $(1,24,8)$, $(1,18,14)$, $(1,21,12)$, $(1,1,8)$, $(0,1,2)$, 
$(1,8,8)$, $(1,16,14)$, $(1,10,20)$, $(1,2,1)$, $(1,1,19)$, $(0,1,7)$, $(1,19,6)$, 
$(1,9,5)$, $(1,10,24)$, $(1,26,17)$, $(1,14,4)$, $(1,1,6)$, $(0,1,3)$, $(1,12,7)$, 
$(1,2,24)$, $(1,10,16)$, $(1,0,11)$, $(1,23,17)$, $(1,18,22)$, $(1,23,5)$, $(1,2,0)$, 
$(1,6,13)$, $(1,3,5)$, $(1,15,5)$, $(1,1,25)$, $(0,1,5)$, $(1,11,19)$, $(1,25,11)$, 
$(1,10,8)$, $(1,13,14)$, $(1,18,0)$, $(1,13,7)$, $(1,6,19)$, $(1,24,13)$, $(1,15,18)$, 
$(1,2,8)$, $(1,20,14)$, $(1,1,15)$, $(0,1,8)$, $(1,23,12)$, $(1,14,23)$, $(1,10,21)$, 
$(1,6,12)$, $(1,4,12)$, $(1,15,21)$, $(1,4,21)$, $(1,18,17)$, $(1,3,0)$, $(1,14,26)$, 
$(1,25,13)$, $(1,1,11)$, $(0,1,9)$, $(1,14,24)$, $(1,18,2)$, $(1,10,2)$, $(1,3,16)$, 
$(1,0,13)$, $(1,12,0)$, $(1,5,23)$, $(1,16,18)$, $(1,7,6)$, $(1,5,9)$, $(1,5,17)$, 
$(1,1,2)$, $(0,1,11)$, $(1,10,9)$, $(1,12,4)$, $(1,10,5)$, $(1,23,9)$, $(1,16,9)$, 
$(1,20,24)$, $(1,8,20)$, $(1,19,7)$, $(1,4,23)$, $(1,1,12)$, $(0,1,21)$, $(1,1,21)$, 
$(0,1,12)$, $(1,26,11)$, $(1,22,9)$, $(1,10,1)$, $(1,25,7)$, $(1,3,6)$, $(1,6,17)$, 
$(1,2,26)$, $(1,15,1)$, $(1,26,12)$, $(1,17,0)$, $(1,8,21)$, $(1,1,22)$, $(0,1,14)$, 
$(1,22,5)$, $(1,7,17)$, $(1,10,18)$, $(1,16,10)$, $(1,13,25)$, $(1,3,23)$, $(1,25,0)$, 
$(1,21,9)$, $(1,17,7)$, $(1,13,3)$, $(1,16,26)$, $(1,1,5)$, $(0,1,15)$, $(1,2,4)$, 
$(1,13,12)$, $(1,10,22)$, $(1,14,21)$, $(1,12,10)$, $(1,19,20)$, $(1,1,24)$, $(0,1,20)$, 
$(1,21,22)$, $(1,20,18)$, $(1,10,3)$, $(1,1,18)$, $(0,1,18)$, $(1,25,1)$, $(1,5,26)$, 
$(1,10,26)$, $(1,9,8)$, $(1,8,14)$, $(1,4,18)$, $(1,26,2)$, $(1,23,19)$, $(1,19,5)$, 
$(1,16,16)$, $(1,0,22)$, $(1,1,26)$, $(0,1,22)$, $(1,5,0)$, $(1,23,20)$, $(1,10,17)$, 
$(1,17,20)$, $(1,24,16)$, $(1,0,2)$, $(1,7,18)$, $(1,13,26)$, $(1,22,15)$, $(1,23,4)$, 
$(1,11,24)$, $(1,1,4)$, $(0,1,24)$, $(1,13,17)$, $(1,19,10)$, $(1,10,12)$, $(1,5,6)$, 
$(1,20,21)$, $(1,21,16)$, $(1,0,25)$, $(1,25,3)$, $(1,23,2)$, $(1,4,25)$, $(1,12,25)$, 
$(1,1,17)$, $(0,1,25)$, $(1,17,23)$, $(1,4,15)$, $(1,10,0)$, $(1,11,25)$, $(1,7,8)$, 
$(1,22,14)$, $(1,20,8)$, $(1,4,14)$, $(1,25,19)$, $(1,8,22)$, $(1,24,19)$, $(1,1,9)$, 
$(0,1,26)$, $(1,9,2)$, $(1,21,1)$, $(1,10,23)$, $(1,19,3)$, $(1,2,7)$, $(1,25,8)$, 
$(1,14,14)$, $(1,26,23)$, $(1,5,16)$, $(1,0,19)$, $(1,13,10)$, $(1,1,7)$, $(1,0,1)$, 
$(1,2,6)$, $(1,21,20)$, $(1,12,9)$, $(1,3,3)$, $(1,19,8)$, $(1,26,14)$, $(1,3,26)$, 
$(1,3,9)$, $(1,4,19)$, $(1,16,4)$, $(1,9,10)$, $(1,25,16)$, $(1,0,3)$, $(1,18,19)$, 
$(1,14,17)$, $(1,5,12)$, $(1,8,13)$, $(1,17,21)$, $(1,25,26)$, $(1,17,10)$, $(1,16,2)$, 
$(1,19,17)$, $(1,22,8)$, $(1,21,14)$, $(1,22,16)$, $(1,0,4)$, $(1,26,13)$, $(1,19,13)$, 
$(1,20,10)$, $(1,14,9)$, $(1,13,20)$, $(1,23,25)$, $(1,21,18)$, $(1,18,1)$, $(1,16,3)$, 
$(1,25,9)$, $(1,18,13)$, $(1,18,16)$, $(1,0,5)$, $(1,22,7)$, $(1,17,12)$, $(1,2,20)$, 
$(1,18,21)$, $(1,20,6)$, $(1,16,7)$, $(1,22,13)$, $(1,13,6)$, $(1,17,6)$, $(1,19,22)$, 
$(1,25,4)$, $(1,4,16)$, $(1,0,6)$, $(1,15,8)$, $(1,3,14)$, $(1,8,7)$, $(1,20,17)$, 
$(1,12,19)$, $(1,20,0)$, $(1,18,11)$, $(1,12,24)$, $(1,12,18)$, $(1,4,13)$, $(1,8,26)$, 
$(1,6,16)$, $(1,0,10)$, $(1,24,23)$, $(1,20,4)$, $(1,17,17)$, $(1,21,11)$, $(1,2,15)$, 
$(1,8,3)$, $(1,11,20)$, $(1,22,12)$, $(1,6,7)$, $(1,9,21)$, $(1,15,9)$, $(1,23,16)$, 
$(1,0,12)$, $(1,16,15)$, $(1,16,21)$, $(1,11,0)$, $(1,26,18)$, $(1,8,12)$, $(1,18,24)$, 
$(1,2,21)$, $(1,7,23)$, $(1,2,13)$, $(1,24,25)$, $(1,16,1)$, $(1,8,16)$, $(1,0,15)$, 
$(1,4,22)$, $(1,11,10)$, $(1,18,6)$, $(1,4,10)$, $(1,4,11)$, $(1,21,13)$, $(1,16,24)$, 
$(1,21,3)$, $(1,24,5)$, $(1,6,3)$, $(1,13,0)$, $(1,15,16)$, $(1,0,16)$, $(1,0,17)$, 
$(1,8,1)$, $(1,12,8)$, $(1,24,14)$, $(1,13,5)$, $(1,24,4)$, $(1,2,19)$, $(1,26,15)$, 
$(1,25,17)$, $(1,11,15)$, $(1,3,20)$, $(1,11,18)$, $(1,19,16)$, $(1,0,18)$, $(1,17,25)$, 
$(1,26,6)$, $(1,26,24)$, $(1,24,7)$, $(1,23,0)$, $(1,22,1)$, $(1,4,9)$, $(1,14,15)$, 
$(1,7,1)$, $(1,11,1)$, $(1,23,22)$, $(1,26,16)$, $(1,0,21)$, $(1,5,5)$, $(1,25,15)$, 
$(1,19,23)$, $(1,11,22)$, $(1,16,23)$, $(1,24,2)$, $(1,12,22)$, $(1,9,22)$, $(1,3,25)$, 
$(1,20,15)$, $(1,24,12)$, $(1,11,16)$, $(1,0,24)$, $(1,20,9)$, $(1,2,25)$, $(1,14,22)$, 
$(1,17,8)$, $(1,7,14)$, $(1,14,6)$, $(1,14,0)$, $(1,20,19)$, $(1,22,26)$, $(1,2,23)$, 
$(1,3,17)$, $(1,13,16)$, $(1,2,5)$, $(1,17,18)$, $(1,9,24)$, $(1,4,8)$, $(1,17,14)$, 
$(1,17,9)$, $(1,22,18)$, $(1,5,19)$, $(1,12,12)$, $(1,22,25)$, $(1,22,21)$, $(1,19,25)$, 
$(1,8,2)$, $(1,2,10)$, $(1,12,5)$, $(1,8,19)$, $(1,5,1)$, $(1,21,0)$, $(1,4,4)$, 
$(1,23,3)$, $(1,20,12)$, $(1,13,2)$, $(1,26,21)$, $(1,16,20)$, $(1,14,8)$, $(1,9,14)$, 
$(1,2,14)$, $(1,19,4)$, $(1,21,25)$, $(1,13,19)$, $(1,4,3)$, $(1,3,1)$, $(1,24,20)$, 
$(1,8,23)$, $(1,14,25)$, $(1,23,24)$, $(1,3,15)$, $(1,26,10)$, $(1,11,8)$, $(1,2,17)$, 
$(1,4,24)$, $(1,22,23)$, $(1,12,26)$, $(1,26,4)$, $(1,20,23)$, $(1,9,23)$, $(1,18,5)$, 
$(1,19,26)$, $(1,18,20)$, $(1,17,13)$, $(1,23,1)$, $(1,13,13)$, $(1,2,18)$, $(1,22,0)$, 
$(1,24,10)$, $(1,21,10)$, $(1,23,26)$, $(1,12,3)$, $(1,21,15)$, $(1,21,6)$, $(1,11,3)$, 
$(1,17,3)$, $(1,12,23)$, $(1,7,9)$, $(1,25,10)$, $(1,3,4)$, $(1,5,2)$, $(1,22,10)$, 
$(1,26,0)$, $(1,7,5)$, $(1,26,20)$, $(1,21,2)$, $(1,25,12)$, $(1,20,3)$, $(1,5,21)$, 
$(1,17,1)$, $(1,20,5)$, $(1,20,1)$, $(1,3,11)$, $(1,18,25)$, $(1,5,20)$, $(1,13,22)$, 
$(1,5,13)$, $(1,6,23)$, $(1,23,7)$, $(1,15,2)$, $(1,8,0)$, $(1,25,2)$, $(1,13,4)$, 
$(1,8,9)$, $(1,8,6)$, $(1,3,18)$, $(1,21,21)$, $(1,14,2)$, $(1,24,3)$, $(1,4,17)$, 
$(1,20,7)$, $(1,19,18)$, $(1,23,15)$, $(1,9,7)$, $(1,11,7)$, $(1,5,10)$, $(1,3,24)$, 
$(1,16,12)$, $(1,4,7)$, $(1,24,9)$, $(1,26,8)$, $(1,12,14)$, $(1,6,2)$, $(1,17,15)$, 
$(1,15,19)$, $(1,11,23)$, $(1,13,24)$, $(1,14,1)$, $(1,22,17)$, $(1,16,5)$, $(1,5,11)$, 
$(1,5,8)$, $(1,19,14)$, $(1,12,1)$, $(1,18,10)$, $(1,9,19)$, $(1,17,2)$, $(1,14,13)$, 
$(1,21,7)$, $(1,26,22)$, $(1,18,3)$, $(1,22,22)$, $(1,8,24)$, $(1,23,13)$, $(1,5,24)$, 
$(1,11,6)$, $(1,8,15)$, $(1,13,15)$, $(1,11,5)$, $(1,11,2)$, $(1,20,26)$, $(1,23,6)$, 
$(1,24,26)$, $(1,14,12)$, $(1,19,9)$, $(1,24,21)$, $(1,24,15)$, $(1,6,4)$, $(1,12,20)$, 
$(1,15,26)$, $(1,24,22)$, $(1,22,20)$, $(1,25,20)$, $(1,16,0)$, $(1,12,15)$, $(1,14,11)$, 
$(1,11,4)$, $(1,14,7)$, $(1,16,25)$, $(1,17,4)$, $(1,19,0)$, $(1,25,24)$.

\subsection*{$m_{21}(2,27)\ge 540$}
$(0,1,0)$, $(1,2,1)$, $(1,21,21)$, $(1,21,24)$, $(1,25,26)$, $(1,25,20)$, $(1,11,13)$, 
$(0,1,2)$, $(1,22,11)$, $(1,23,26)$, $(1,16,0)$, $(1,3,14)$, $(1,26,24)$, $(1,24,3)$, 
$(0,1,3)$, $(1,12,6)$, $(1,18,0)$, $(1,4,2)$, $(1,15,12)$, $(1,19,14)$, $(1,7,17)$, 
$(0,1,5)$, $(1,8,16)$, $(1,20,5)$, $(1,26,14)$, $(1,23,0)$, $(1,20,9)$, $(1,14,7)$, 
$(0,1,6)$, $(1,11,17)$, $(1,6,22)$, $(1,15,5)$, $(1,4,9)$, $(1,13,19)$, $(1,2,6)$, 
$(0,1,7)$, $(0,1,8)$, $(1,26,0)$, $(1,17,4)$, $(1,5,6)$, $(1,17,17)$, $(1,5,10)$, 
$(1,23,14)$, $(0,1,9)$, $(1,15,18)$, $(1,22,19)$, $(1,2,12)$, $(1,11,2)$, $(1,24,22)$, 
$(1,4,23)$, $(0,1,11)$, $(1,16,25)$, $(1,12,17)$, $(1,20,8)$, $(1,20,15)$, $(1,7,6)$, 
$(1,3,25)$, $(0,1,12)$, $(1,17,23)$, $(1,5,3)$, $(1,11,19)$, $(1,2,22)$, $(1,17,11)$, 
$(1,5,18)$, $(0,1,13)$, $(1,4,24)$, $(1,8,24)$, $(1,1,17)$, $(1,9,6)$, $(1,14,26)$, 
$(1,15,26)$, $(0,1,14)$, $(1,21,13)$, $(1,3,1)$, $(1,25,16)$, $(1,21,7)$, $(1,16,13)$, 
$(1,25,1)$, $(0,1,15)$, $(1,25,5)$, $(1,24,15)$, $(1,14,13)$, $(1,8,1)$, $(1,22,8)$, 
$(1,21,9)$, $(0,1,16)$, $(1,23,15)$, $(1,13,12)$, $(1,7,23)$, $(1,12,18)$, $(1,6,2)$, 
$(1,26,8)$, $(0,1,17)$, $(1,1,3)$, $(1,4,8)$, $(1,3,4)$, $(1,16,10)$, $(1,15,15)$, 
$(1,9,11)$, $(0,1,18)$, $(1,13,4)$, $(1,11,25)$, $(1,22,22)$, $(1,24,19)$, $(1,2,25)$, 
$(1,6,10)$, $(0,1,21)$, $(1,18,19)$, $(1,2,9)$, $(1,18,3)$, $(1,19,11)$, $(1,11,5)$, 
$(1,19,22)$, $(0,1,22)$, $(1,24,7)$, $(1,7,20)$, $(1,23,20)$, $(1,26,21)$, $(1,12,21)$, 
$(1,22,16)$, $(0,1,23)$, $(1,19,26)$, $(1,19,7)$, $(1,9,26)$, $(1,1,24)$, $(1,18,16)$, 
$(1,18,24)$, $(0,1,25)$, $(1,14,2)$, $(1,1,14)$, $(1,13,15)$, $(1,6,8)$, $(1,9,0)$, 
$(1,8,12)$, $(0,1,26)$, $(1,9,12)$, $(1,26,11)$, $(1,6,25)$, $(1,13,25)$, $(1,23,3)$, 
$(1,1,2)$, $(1,1,0)$, $(1,1,13)$, $(1,16,1)$, $(1,13,6)$, $(1,21,18)$, $(1,10,6)$, 
$(1,25,15)$, $(1,1,5)$, $(1,12,19)$, $(1,2,10)$, $(1,14,20)$, $(1,16,19)$, $(1,2,21)$, 
$(1,16,14)$, $(1,1,6)$, $(1,15,14)$, $(1,24,23)$, $(1,11,16)$, $(1,26,22)$, $(1,3,2)$, 
$(1,6,5)$, $(1,1,7)$, $(1,3,9)$, $(1,17,18)$, $(1,5,12)$, $(1,20,0)$, $(1,15,19)$, 
$(1,2,17)$, $(1,1,8)$, $(1,22,25)$, $(1,1,20)$, $(1,10,23)$, $(1,4,21)$, $(1,17,3)$, 
$(1,5,7)$, $(1,1,9)$, $(1,8,17)$, $(1,7,25)$, $(1,4,19)$, $(1,2,20)$, $(1,7,16)$, 
$(1,13,21)$, $(1,1,10)$, $(1,17,24)$, $(1,5,14)$, $(1,6,9)$, $(1,7,2)$, $(1,22,26)$, 
$(1,7,10)$, $(1,1,11)$, $(1,24,21)$, $(1,20,7)$, $(1,24,4)$, $(1,13,14)$, $(1,14,14)$, 
$(1,4,20)$, $(1,1,19)$, $(1,2,0)$, $(1,10,5)$, $(1,8,25)$, $(1,3,11)$, $(1,12,0)$, 
$(1,26,23)$, $(1,1,21)$, $(1,14,5)$, $(1,23,2)$, $(1,18,17)$, $(1,23,17)$, $(1,24,20)$, 
$(1,19,8)$, $(1,1,22)$, $(1,26,7)$, $(1,15,11)$, $(1,3,8)$, $(1,24,8)$, $(1,6,18)$, 
$(1,11,9)$, $(1,1,23)$, $(1,25,11)$, $(1,22,9)$, $(1,12,13)$, $(1,12,1)$, $(1,16,17)$, 
$(1,21,11)$, $(1,1,25)$, $(1,7,18)$, $(1,19,17)$, $(1,22,14)$, $(1,22,4)$, $(1,8,5)$, 
$(1,18,18)$, $(1,1,26)$, $(1,21,2)$, $(1,8,15)$, $(1,25,24)$, $(1,17,5)$, $(1,5,13)$, 
$(1,24,1)$, $(1,2,3)$, $(1,24,12)$, $(1,17,10)$, $(1,5,8)$, $(1,10,16)$, $(1,24,17)$, 
$(1,8,19)$, $(1,2,4)$, $(1,9,8)$, $(1,7,7)$, $(1,26,2)$, $(1,9,7)$, $(1,10,19)$, 
$(1,2,19)$, $(1,2,5)$, $(1,13,23)$, $(1,22,15)$, $(1,16,15)$, $(1,4,3)$, $(1,23,16)$, 
$(1,9,19)$, $(1,2,7)$, $(1,22,18)$, $(1,4,0)$, $(1,9,17)$, $(1,13,9)$, $(1,16,12)$, 
$(1,23,19)$, $(1,2,11)$, $(1,20,6)$, $(1,6,16)$, $(1,14,16)$, $(1,7,8)$, $(1,17,6)$, 
$(1,5,19)$, $(1,2,14)$, $(1,8,13)$, $(1,15,1)$, $(1,15,3)$, $(1,21,15)$, $(1,12,14)$, 
$(1,25,19)$, $(1,2,16)$, $(1,4,25)$, $(1,10,26)$, $(1,22,12)$, $(1,14,18)$, $(1,3,24)$, 
$(1,14,19)$, $(1,2,18)$, $(1,26,15)$, $(1,18,20)$, $(1,24,14)$, $(1,3,21)$, $(1,19,15)$, 
$(1,20,19)$, $(1,2,24)$, $(1,19,0)$, $(1,12,3)$, $(1,10,20)$, $(1,6,26)$, $(1,18,21)$, 
$(1,6,19)$, $(1,3,0)$, $(1,11,20)$, $(1,3,22)$, $(1,8,21)$, $(1,11,4)$, $(1,7,22)$, 
$(1,9,9)$, $(1,3,3)$, $(1,22,17)$, $(1,18,2)$, $(1,13,16)$, $(1,12,2)$, $(1,19,4)$, 
$(1,3,16)$, $(1,3,6)$, $(1,7,13)$, $(1,7,1)$, $(1,24,5)$, $(1,21,3)$, $(1,9,18)$, 
$(1,25,3)$, $(1,3,7)$, $(1,4,14)$, $(1,12,10)$, $(1,3,26)$, $(1,23,25)$, $(1,17,8)$, 
$(1,5,24)$, $(1,3,10)$, $(1,18,7)$, $(1,14,0)$, $(1,14,15)$, $(1,19,18)$, $(1,15,10)$, 
$(1,13,2)$, $(1,3,15)$, $(1,10,9)$, $(1,19,9)$, $(1,6,23)$, $(1,17,14)$, $(1,5,25)$, 
$(1,18,12)$, $(1,3,18)$, $(1,13,11)$, $(1,8,11)$, $(1,17,12)$, $(1,5,20)$, $(1,8,6)$, 
$(1,23,21)$, $(1,3,23)$, $(1,21,6)$, $(1,13,20)$, $(1,25,4)$, $(1,14,21)$, $(1,13,13)$, 
$(1,10,1)$, $(1,4,1)$, $(1,7,15)$, $(1,21,23)$, $(1,18,4)$, $(1,25,7)$, $(1,23,23)$, 
$(1,19,13)$, $(1,4,4)$, $(1,18,23)$, $(1,8,8)$, $(1,8,14)$, $(1,19,16)$, $(1,16,4)$, 
$(1,6,12)$, $(1,4,6)$, $(1,11,8)$, $(1,18,22)$, $(1,10,3)$, $(1,20,10)$, $(1,19,25)$, 
$(1,8,2)$, $(1,4,7)$, $(1,21,25)$, $(1,26,5)$, $(1,25,10)$, $(1,8,3)$, $(1,26,13)$, 
$(1,13,1)$, $(1,4,11)$, $(1,4,13)$, $(1,14,1)$, $(1,11,0)$, $(1,21,22)$, $(1,7,0)$, 
$(1,25,8)$, $(1,4,17)$, $(1,19,5)$, $(1,10,11)$, $(1,26,16)$, $(1,10,14)$, $(1,18,11)$, 
$(1,9,10)$, $(1,4,18)$, $(1,24,18)$, $(1,23,4)$, $(1,4,26)$, $(1,6,6)$, $(1,20,3)$, 
$(1,12,24)$, $(1,5,1)$, $(1,19,24)$, $(1,21,20)$, $(1,20,25)$, $(1,25,21)$, $(1,18,26)$, 
$(1,17,13)$, $(1,5,5)$, $(1,7,21)$, $(1,10,12)$, $(1,15,8)$, $(1,23,5)$, $(1,15,4)$, 
$(1,17,20)$, $(1,5,9)$, $(1,21,26)$, $(1,14,8)$, $(1,25,12)$, $(1,9,24)$, $(1,22,13)$, 
$(1,17,1)$, $(1,5,15)$, $(1,26,25)$, $(1,16,24)$, $(1,7,4)$, $(1,15,0)$, $(1,16,16)$, 
$(1,17,26)$, $(1,5,17)$, $(1,12,15)$, $(1,8,7)$, $(1,13,7)$, $(1,20,17)$, $(1,11,3)$, 
$(1,17,22)$, $(1,6,0)$, $(1,22,10)$, $(1,20,16)$, $(1,22,20)$, $(1,9,3)$, $(1,15,21)$, 
$(1,8,0)$, $(1,6,7)$, $(1,19,12)$, $(1,24,6)$, $(1,16,11)$, $(1,16,7)$, $(1,18,10)$, 
$(1,7,12)$, $(1,6,13)$, $(1,23,1)$, $(1,14,11)$, $(1,21,4)$, $(1,23,9)$, $(1,25,25)$, 
$(1,15,16)$, $(1,6,15)$, $(1,24,26)$, $(1,10,17)$, $(1,10,0)$, $(1,7,24)$, $(1,11,15)$, 
$(1,23,22)$, $(1,7,26)$, $(1,20,11)$, $(1,13,17)$, $(1,15,24)$, $(1,26,10)$, $(1,22,23)$, 
$(1,15,23)$, $(1,8,9)$, $(1,9,20)$, $(1,18,15)$, $(1,22,21)$, $(1,26,6)$, $(1,19,6)$, 
$(1,26,12)$, $(1,8,26)$, $(1,12,9)$, $(1,23,11)$, $(1,10,24)$, $(1,9,2)$, $(1,16,2)$, 
$(1,20,18)$, $(1,9,13)$, $(1,20,1)$, $(1,23,10)$, $(1,21,5)$, $(1,11,12)$, $(1,25,22)$, 
$(1,13,10)$, $(1,9,25)$, $(1,19,23)$, $(1,14,9)$, $(1,24,10)$, $(1,24,0)$, $(1,18,6)$, 
$(1,12,23)$, $(1,10,2)$, $(1,23,13)$, $(1,22,1)$, $(1,10,25)$, $(1,21,14)$, $(1,15,7)$, 
$(1,25,6)$, $(1,10,8)$, $(1,13,3)$, $(1,13,26)$, $(1,24,16)$, $(1,15,9)$, $(1,20,24)$, 
$(1,16,9)$, $(1,10,13)$, $(1,11,1)$, $(1,12,22)$, $(1,21,12)$, $(1,22,3)$, $(1,25,14)$, 
$(1,20,12)$, $(1,10,22)$, $(1,14,4)$, $(1,26,17)$, $(1,26,26)$, $(1,11,18)$, $(1,16,22)$, 
$(1,11,24)$, $(1,12,8)$, $(1,15,13)$, $(1,18,1)$, $(1,26,18)$, $(1,21,16)$, $(1,19,10)$, 
$(1,25,18)$.

\subsection*{$m_{22}(2,27)\ge 561$}
$(0,0,1)$, $(1,17,20)$, $(0,1,10)$, $(1,24,15)$, $(1,3,21)$, $(1,18,12)$, $(1,14,20)$, 
$(0,1,0)$, $(0,1,5)$, $(1,16,24)$, $(1,10,23)$, $(1,19,14)$, $(1,18,2)$, $(1,18,20)$, 
$(0,1,1)$, $(1,12,16)$, $(1,5,14)$, $(1,14,19)$, $(1,17,4)$, $(1,18,23)$, $(1,22,20)$, 
$(0,1,2)$, $(1,10,3)$, $(1,18,1)$, $(1,13,20)$, $(0,1,12)$, $(1,18,21)$, $(1,6,20)$, 
$(0,1,3)$, $(1,8,6)$, $(1,8,10)$, $(1,6,6)$, $(1,21,13)$, $(1,18,11)$, $(1,4,20)$, 
$(0,1,4)$, $(1,22,5)$, $(1,21,6)$, $(1,4,2)$, $(1,13,17)$, $(1,18,24)$, $(1,21,20)$, 
$(0,1,6)$, $(1,26,13)$, $(1,24,5)$, $(1,25,17)$, $(1,23,1)$, $(1,18,5)$, $(1,3,20)$, 
$(0,1,7)$, $(1,15,22)$, $(1,13,19)$, $(1,12,18)$, $(1,12,5)$, $(1,18,18)$, $(1,7,20)$, 
$(0,1,8)$, $(1,0,17)$, $(1,2,15)$, $(1,8,7)$, $(1,15,16)$, $(1,18,17)$, $(1,19,20)$, 
$(0,1,11)$, $(1,23,7)$, $(1,19,11)$, $(1,7,8)$, $(1,2,7)$, $(1,18,15)$, $(1,20,20)$, 
$(0,1,13)$, $(1,20,19)$, $(1,22,16)$, $(1,9,21)$, $(1,5,9)$, $(1,18,13)$, $(1,10,20)$, 
$(0,1,14)$, $(1,1,10)$, $(1,11,3)$, $(1,3,0)$, $(1,0,10)$, $(1,18,19)$, $(1,24,20)$, 
$(0,1,15)$, $(1,25,11)$, $(1,25,12)$, $(1,22,11)$, $(1,16,19)$, $(1,18,8)$, $(1,23,20)$, 
$(0,1,16)$, $(1,3,23)$, $(1,14,2)$, $(1,17,25)$, $(1,24,24)$, $(1,18,14)$, $(1,2,20)$, 
$(0,1,17)$, $(1,9,1)$, $(1,0,25)$, $(1,21,9)$, $(1,25,0)$, $(1,18,26)$, $(1,0,20)$, 
$(0,1,22)$, $(1,13,9)$, $(1,23,26)$, $(1,23,10)$, $(1,20,26)$, $(1,18,6)$, $(1,12,20)$, 
$(0,1,23)$, $(1,19,26)$, $(1,9,13)$, $(1,11,22)$, $(1,11,15)$, $(1,18,22)$, $(1,1,20)$, 
$(0,1,24)$, $(1,4,25)$, $(1,26,22)$, $(1,5,1)$, $(1,4,6)$, $(1,18,7)$, $(1,11,20)$, 
$(0,1,25)$, $(1,21,0)$, $(1,12,9)$, $(1,2,4)$, $(1,26,12)$, $(1,18,0)$, $(1,16,20)$, 
$(1,0,0)$, $(1,4,18)$, $(1,3,9)$, $(1,12,8)$, $(1,25,3)$, $(1,24,22)$, $(1,4,8)$, 
$(1,0,1)$, $(1,26,21)$, $(1,2,3)$, $(1,22,9)$, $(1,10,22)$, $(1,9,11)$, $(1,2,16)$, 
$(1,0,2)$, $(1,25,16)$, $(1,7,24)$, $(1,1,0)$, $(1,14,16)$, $(1,17,13)$, $(1,25,7)$, 
$(1,0,5)$, $(1,23,4)$, $(1,22,4)$, $(1,2,14)$, $(1,13,26)$, $(1,8,9)$, $(1,7,23)$, 
$(1,0,6)$, $(1,13,22)$, $(1,15,5)$, $(1,8,3)$, $(1,11,12)$, $(1,20,7)$, $(1,16,9)$, 
$(1,0,7)$, $(1,22,26)$, $(1,14,26)$, $(1,17,5)$, $(1,22,2)$, $(1,15,18)$, $(1,6,24)$, 
$(1,0,9)$, $(1,8,25)$, $(1,19,8)$, $(1,4,11)$, $(1,9,5)$, $(1,13,15)$, $(1,1,11)$, 
$(1,0,11)$, $(1,6,3)$, $(1,23,14)$, $(1,11,13)$, $(1,6,25)$, $(1,14,23)$, $(1,17,17)$, 
$(1,0,12)$, $(1,0,23)$, $(1,16,12)$, $(1,6,17)$, $(1,16,18)$, $(1,4,14)$, $(1,22,10)$, 
$(1,0,14)$, $(1,5,8)$, $(1,11,18)$, $(1,23,23)$, $(1,2,0)$, $(1,21,5)$, $(1,11,5)$, 
$(1,0,18)$, $(1,14,0)$, $(1,17,22)$, $(1,14,21)$, $(1,17,11)$, $(1,26,17)$, $(1,10,6)$, 
$(1,0,19)$, $(1,11,10)$, $(1,13,16)$, $(1,20,12)$, $(1,15,1)$, $(1,16,8)$, $(1,8,19)$, 
$(1,0,21)$, $(1,1,1)$, $(1,5,2)$, $(1,10,2)$, $(1,20,24)$, $(1,23,24)$, $(1,26,3)$, 
$(1,0,22)$, $(1,19,6)$, $(1,1,23)$, $(1,7,19)$, $(1,7,15)$, $(1,5,19)$, $(1,21,14)$, 
$(1,0,24)$, $(1,3,13)$, $(1,20,15)$, $(1,13,10)$, $(1,24,13)$, $(1,2,2)$, $(1,19,22)$, 
$(1,0,26)$, $(1,15,7)$, $(1,21,21)$, $(1,26,4)$, $(1,8,8)$, $(1,10,25)$, $(1,12,25)$, 
$(1,1,2)$, $(1,22,1)$, $(1,12,24)$, $(1,15,12)$, $(1,1,25)$, $(1,16,6)$, $(1,11,9)$, 
$(1,1,3)$, $(1,10,10)$, $(1,23,5)$, $(1,7,1)$, $(1,9,9)$, $(1,22,21)$, $(1,4,16)$, 
$(1,1,4)$, $(1,3,22)$, $(1,25,18)$, $(1,11,21)$, $(1,24,6)$, $(1,23,11)$, $(1,10,7)$, 
$(1,1,7)$, $(1,26,7)$, $(1,4,10)$, $(1,12,14)$, $(1,6,16)$, $(1,2,8)$, $(1,9,23)$, 
$(1,1,9)$, $(1,20,11)$, $(1,1,14)$, $(1,21,17)$, $(1,14,17)$, $(1,17,23)$, $(1,12,4)$, 
$(1,1,12)$, $(1,8,12)$, $(1,16,3)$, $(1,16,4)$, $(1,7,3)$, $(1,8,13)$, $(1,26,0)$, 
$(1,1,15)$, $(1,19,24)$, $(1,22,22)$, $(1,6,9)$, $(1,19,15)$, $(1,6,0)$, $(1,8,2)$, 
$(1,1,16)$, $(1,12,0)$, $(1,24,11)$, $(1,22,6)$, $(1,20,5)$, $(1,20,25)$, $(1,3,5)$, 
$(1,1,17)$, $(1,25,14)$, $(1,20,6)$, $(1,26,26)$, $(1,25,23)$, $(1,15,10)$, $(1,2,10)$, 
$(1,1,18)$, $(1,9,26)$, $(1,21,12)$, $(1,2,18)$, $(1,16,13)$, $(1,14,1)$, $(1,17,1)$, 
$(1,1,19)$, $(1,13,23)$, $(1,26,1)$, $(1,14,22)$, $(1,17,0)$, $(1,13,14)$, $(1,25,25)$, 
$(1,1,26)$, $(1,24,18)$, $(1,13,7)$, $(1,13,3)$, $(1,5,7)$, $(1,3,17)$, $(1,21,3)$, 
$(1,2,1)$, $(1,11,16)$, $(1,19,10)$, $(1,8,4)$, $(1,13,11)$, $(1,14,15)$, $(1,17,18)$, 
$(1,2,5)$, $(1,10,0)$, $(1,2,17)$, $(1,12,10)$, $(1,7,2)$, $(1,4,21)$, $(1,8,24)$, 
$(1,2,6)$, $(1,24,19)$, $(1,9,12)$, $(1,25,26)$, $(1,16,0)$, $(1,3,25)$, $(1,4,7)$, 
$(1,2,11)$, $(1,3,24)$, $(1,6,19)$, $(1,5,5)$, $(1,24,12)$, $(1,26,19)$, $(1,16,2)$, 
$(1,2,13)$, $(1,7,18)$, $(1,24,4)$, $(1,19,25)$, $(1,11,17)$, $(1,22,8)$, $(1,3,15)$, 
$(1,2,19)$, $(1,15,6)$, $(1,22,15)$, $(1,23,8)$, $(1,14,6)$, $(1,17,3)$, $(1,19,3)$, 
$(1,2,23)$, $(1,25,8)$, $(1,5,13)$, $(1,15,9)$, $(1,15,14)$, $(1,7,9)$, $(1,13,1)$, 
$(1,2,24)$, $(1,5,11)$, $(1,12,11)$, $(1,11,0)$, $(1,9,3)$, $(1,15,2)$, $(1,5,26)$, 
$(1,3,2)$, $(1,16,23)$, $(1,21,24)$, $(1,24,8)$, $(1,21,22)$, $(1,7,14)$, $(1,23,12)$, 
$(1,3,3)$, $(1,9,17)$, $(1,9,4)$, $(1,24,17)$, $(1,6,5)$, $(1,15,19)$, $(1,25,22)$, 
$(1,3,4)$, $(1,23,6)$, $(1,15,8)$, $(1,24,1)$, $(1,20,0)$, $(1,6,4)$, $(1,6,1)$, 
$(1,3,7)$, $(1,7,25)$, $(1,6,11)$, $(1,24,10)$, $(1,14,9)$, $(1,17,9)$, $(1,21,19)$, 
$(1,3,8)$, $(1,19,9)$, $(1,3,12)$, $(1,24,21)$, $(1,16,15)$, $(1,24,16)$, $(1,10,26)$, 
$(1,3,16)$, $(1,13,5)$, $(1,4,22)$, $(1,24,25)$, $(1,15,4)$, $(1,14,13)$, $(1,17,16)$, 
$(1,3,19)$, $(1,26,24)$, $(1,20,9)$, $(1,24,7)$, $(1,8,18)$, $(1,22,18)$, $(1,26,2)$, 
$(1,4,5)$, $(1,14,11)$, $(1,17,21)$, $(1,7,22)$, $(1,21,10)$, $(1,10,1)$, $(1,15,3)$, 
$(1,4,9)$, $(1,6,18)$, $(1,20,17)$, $(1,5,21)$, $(1,13,0)$, $(1,7,6)$, $(1,10,16)$, 
$(1,4,12)$, $(1,23,13)$, $(1,16,14)$, $(1,26,18)$, $(1,15,11)$, $(1,20,23)$, $(1,4,13)$, 
$(1,4,24)$, $(1,7,21)$, $(1,25,2)$, $(1,6,7)$, $(1,10,24)$, $(1,13,2)$, $(1,11,4)$, 
$(1,4,26)$, $(1,20,4)$, $(1,9,0)$, $(1,19,19)$, $(1,20,3)$, $(1,14,18)$, $(1,17,15)$, 
$(1,5,0)$, $(1,20,21)$, $(1,19,17)$, $(1,15,25)$, $(1,10,12)$, $(1,21,15)$, $(1,16,25)$, 
$(1,5,3)$, $(1,6,15)$, $(1,25,9)$, $(1,26,5)$, $(1,12,13)$, $(1,19,13)$, $(1,21,26)$, 
$(1,5,6)$, $(1,16,1)$, $(1,22,13)$, $(1,5,17)$, $(1,7,7)$, $(1,14,8)$, $(1,17,12)$, 
$(1,5,12)$, $(1,19,5)$, $(1,16,21)$, $(1,20,10)$, $(1,25,15)$, $(1,10,4)$, $(1,11,7)$, 
$(1,5,16)$, $(1,25,4)$, $(1,12,15)$, $(1,21,1)$, $(1,6,10)$, $(1,26,11)$, $(1,19,18)$, 
$(1,6,21)$, $(1,11,26)$, $(1,12,19)$, $(1,10,17)$, $(1,8,0)$, $(1,25,13)$, $(1,13,8)$, 
$(1,6,22)$, $(1,22,12)$, $(1,7,10)$, $(1,16,22)$, $(1,12,26)$, $(1,26,25)$, $(1,23,16)$, 
$(1,6,23)$, $(1,9,25)$, $(1,20,1)$, $(1,7,12)$, $(1,11,14)$, $(1,16,5)$, $(1,9,7)$, 
$(1,7,13)$, $(1,12,3)$, $(1,20,8)$, $(1,12,2)$, $(1,14,10)$, $(1,17,24)$, $(1,11,25)$, 
$(1,8,11)$, $(1,21,4)$, $(1,15,21)$, $(1,12,7)$, $(1,9,22)$, $(1,10,13)$, $(1,22,19)$, 
$(1,8,14)$, $(1,12,21)$, $(1,22,23)$, $(1,11,23)$, $(1,10,18)$, $(1,25,6)$, $(1,9,14)$, 
$(1,8,15)$, $(1,9,2)$, $(1,23,3)$, $(1,25,24)$, $(1,19,16)$, $(1,23,18)$, $(1,21,25)$, 
$(1,9,6)$, $(1,14,24)$, $(1,17,10)$, $(1,9,8)$, $(1,26,8)$, $(1,13,4)$, $(1,21,2)$, 
$(1,10,21)$.

\subsection*{$m_{23}(2,27)\ge 595$}
$(0,0,1)$, $(1,9,15)$, $(1,20,12)$, $(0,1,0)$, $(1,26,21)$, $(1,3,2)$, $(0,1,1)$, 
$(1,4,0)$, $(1,5,0)$, $(0,1,2)$, $(0,1,24)$, $(1,13,19)$, $(0,1,3)$, $(1,7,25)$, 
$(1,1,4)$, $(0,1,5)$, $(1,22,11)$, $(1,21,11)$, $(0,1,6)$, $(1,3,12)$, $(1,26,15)$, 
$(0,1,7)$, $(1,23,8)$, $(1,0,5)$, $(0,1,8)$, $(1,8,13)$, $(1,8,6)$, $(0,1,9)$, 
$(1,20,10)$, $(1,9,23)$, $(0,1,10)$, $(1,19,22)$, $(1,25,16)$, $(0,1,13)$, $(1,17,5)$, 
$(1,6,8)$, $(0,1,14)$, $(1,6,1)$, $(1,17,24)$, $(0,1,15)$, $(1,0,26)$, $(1,23,9)$, 
$(0,1,16)$, $(1,10,3)$, $(1,14,18)$, $(0,1,18)$, $(1,2,2)$, $(1,11,21)$, $(0,1,20)$, 
$(1,25,6)$, $(1,19,13)$, $(0,1,21)$, $(1,11,18)$, $(1,2,3)$, $(0,1,22)$, $(1,18,7)$, 
$(1,18,14)$, $(0,1,23)$, $(1,16,17)$, $(1,12,20)$, $(0,1,25)$, $(1,15,20)$, $(1,24,17)$, 
$(0,1,26)$, $(1,21,23)$, $(1,22,10)$, $(1,0,1)$, $(1,22,14)$, $(1,12,1)$, $(1,0,2)$, 
$(1,15,21)$, $(1,0,4)$, $(1,0,7)$, $(1,1,6)$, $(1,16,0)$, $(1,0,8)$, $(1,12,0)$, 
$(1,22,20)$, $(1,0,9)$, $(1,21,16)$, $(1,4,3)$, $(1,0,10)$, $(1,19,20)$, $(1,18,18)$, 
$(1,0,11)$, $(1,7,12)$, $(1,26,19)$, $(1,0,12)$, $(1,26,3)$, $(1,7,26)$, $(1,0,13)$, 
$(1,16,19)$, $(1,1,10)$, $(1,0,14)$, $(1,10,13)$, $(1,10,23)$, $(1,0,15)$, $(1,17,26)$, 
$(1,3,24)$, $(1,0,17)$, $(1,25,8)$, $(1,24,7)$, $(1,0,18)$, $(1,24,1)$, $(1,25,13)$, 
$(1,0,20)$, $(1,5,25)$, $(1,8,5)$, $(1,0,21)$, $(1,14,5)$, $(1,11,2)$, $(1,0,22)$, 
$(1,4,18)$, $(1,21,14)$, $(1,0,23)$, $(1,11,11)$, $(1,14,22)$, $(1,0,24)$, $(1,2,4)$, 
$(1,13,16)$, $(1,0,25)$, $(1,6,17)$, $(1,6,11)$, $(1,1,1)$, $(1,5,19)$, $(1,25,22)$, 
$(1,1,2)$, $(1,14,21)$, $(1,6,6)$, $(1,1,3)$, $(1,22,23)$, $(1,24,3)$, $(1,1,5)$, 
$(1,8,1)$, $(1,18,4)$, $(1,1,7)$, $(1,25,5)$, $(1,5,16)$, $(1,1,8)$, $(1,21,20)$, 
$(1,15,0)$, $(1,1,12)$, $(1,26,8)$, $(1,11,9)$, $(1,1,13)$, $(1,3,22)$, $(1,2,15)$, 
$(1,1,15)$, $(1,2,10)$, $(1,3,8)$, $(1,1,16)$, $(1,13,18)$, $(1,1,26)$, $(1,1,17)$, 
$(1,11,12)$, $(1,26,14)$, $(1,1,18)$, $(1,10,9)$, $(1,17,11)$, $(1,1,20)$, $(1,19,14)$, 
$(1,4,24)$, $(1,1,21)$, $(1,23,26)$, $(1,12,2)$, $(1,1,22)$, $(1,24,2)$, $(1,22,21)$, 
$(1,1,23)$, $(1,4,25)$, $(1,19,23)$, $(1,1,24)$, $(1,15,6)$, $(1,21,5)$, $(1,1,25)$, 
$(1,18,11)$, $(1,8,17)$, $(1,2,0)$, $(1,19,26)$, $(1,15,14)$, $(1,2,1)$, $(1,8,18)$, 
$(1,4,13)$, $(1,2,5)$, $(1,17,17)$, $(1,17,23)$, $(1,2,6)$, $(1,2,7)$, $(1,10,0)$, 
$(1,2,8)$, $(1,22,3)$, $(1,5,4)$, $(1,2,9)$, $(1,5,14)$, $(1,22,19)$, $(1,2,13)$, 
$(1,15,2)$, $(1,19,21)$, $(1,2,14)$, $(1,6,12)$, $(1,26,17)$, $(1,2,16)$, $(1,14,1)$, 
$(1,14,25)$, $(1,2,17)$, $(1,7,6)$, $(1,12,16)$, $(1,2,18)$, $(1,20,11)$, $(1,9,9)$, 
$(1,2,20)$, $(1,16,23)$, $(1,21,1)$, $(1,2,22)$, $(1,21,9)$, $(1,16,5)$, $(1,2,23)$, 
$(1,12,22)$, $(1,7,11)$, $(1,2,24)$, $(1,25,10)$, $(1,18,3)$, $(1,2,26)$, $(1,24,25)$, 
$(1,23,10)$, $(1,3,0)$, $(1,11,15)$, $(1,5,1)$, $(1,3,4)$, $(1,21,15)$, $(1,10,11)$, 
$(1,3,5)$, $(1,20,15)$, $(1,9,14)$, $(1,3,6)$, $(1,6,15)$, $(1,24,18)$, $(1,3,9)$, 
$(1,25,15)$, $(1,16,20)$, $(1,3,10)$, $(1,12,15)$, $(1,14,8)$, $(1,3,13)$, $(1,18,15)$, 
$(1,7,22)$, $(1,3,14)$, $(1,14,15)$, $(1,12,5)$, $(1,3,16)$, $(1,22,15)$, $(1,15,23)$, 
$(1,3,17)$, $(1,13,15)$, $(1,8,19)$, $(1,3,20)$, $(1,7,15)$, $(1,18,9)$, $(1,3,21)$, 
$(1,8,15)$, $(1,13,2)$, $(1,3,25)$, $(1,4,15)$, $(1,23,3)$, $(1,4,5)$, $(1,13,25)$, 
$(1,7,20)$, $(1,4,6)$, $(1,17,8)$, $(1,25,11)$, $(1,4,7)$, $(1,24,13)$, $(1,12,7)$, 
$(1,4,8)$, $(1,18,12)$, $(1,26,9)$, $(1,4,9)$, $(1,22,22)$, $(1,10,5)$, $(1,4,10)$, 
$(1,16,24)$, $(1,11,3)$, $(1,4,11)$, $(1,6,10)$, $(1,15,1)$, $(1,4,12)$, $(1,26,5)$, 
$(1,18,13)$, $(1,4,14)$, $(1,20,4)$, $(1,9,4)$, $(1,4,16)$, $(1,15,16)$, $(1,6,19)$, 
$(1,4,17)$, $(1,10,26)$, $(1,22,17)$, $(1,4,19)$, $(1,14,6)$, $(1,4,26)$, $(1,4,21)$, 
$(1,11,7)$, $(1,16,2)$, $(1,4,22)$, $(1,12,17)$, $(1,24,10)$, $(1,5,2)$, $(1,24,21)$, 
$(1,5,11)$, $(1,5,3)$, $(1,16,13)$, $(1,14,17)$, $(1,5,5)$, $(1,15,12)$, $(1,26,13)$, 
$(1,5,6)$, $(1,14,11)$, $(1,16,8)$, $(1,5,7)$, $(1,12,9)$, $(1,18,23)$, $(1,5,12)$, 
$(1,26,23)$, $(1,15,24)$, $(1,5,13)$, $(1,23,20)$, $(1,17,16)$, $(1,5,17)$, $(1,21,18)$, 
$(1,19,4)$, $(1,5,18)$, $(1,18,24)$, $(1,12,25)$, $(1,5,20)$, $(1,13,10)$, $(1,6,18)$, 
$(1,5,21)$, $(1,10,16)$, $(1,7,2)$, $(1,5,22)$, $(1,6,3)$, $(1,13,6)$, $(1,5,23)$, 
$(1,20,26)$, $(1,9,26)$, $(1,5,26)$, $(1,7,4)$, $(1,10,7)$, $(1,6,0)$, $(1,16,14)$, 
$(1,18,26)$, $(1,6,5)$, $(1,23,0)$, $(1,8,3)$, $(1,6,7)$, $(1,20,22)$, $(1,9,24)$, 
$(1,6,9)$, $(1,8,7)$, $(1,23,23)$, $(1,6,14)$, $(1,19,9)$, $(1,7,19)$, $(1,6,21)$, 
$(1,12,23)$, $(1,10,2)$, $(1,6,23)$, $(1,18,8)$, $(1,16,18)$, $(1,6,24)$, $(1,10,19)$, 
$(1,12,8)$, $(1,6,26)$, $(1,21,13)$, $(1,11,13)$, $(1,7,0)$, $(1,23,18)$, $(1,22,18)$, 
$(1,7,3)$, $(1,25,0)$, $(1,7,5)$, $(1,7,9)$, $(1,14,24)$, $(1,21,24)$, $(1,7,10)$, 
$(1,24,11)$, $(1,17,1)$, $(1,7,13)$, $(1,21,10)$, $(1,14,19)$, $(1,7,14)$, $(1,17,25)$, 
$(1,24,6)$, $(1,7,16)$, $(1,20,20)$, $(1,9,13)$, $(1,7,17)$, $(1,8,22)$, $(1,11,10)$, 
$(1,7,23)$, $(1,16,6)$, $(1,15,4)$, $(1,7,24)$, $(1,22,2)$, $(1,23,21)$, $(1,8,4)$, 
$(1,19,0)$, $(1,12,10)$, $(1,8,11)$, $(1,10,14)$, $(1,16,3)$, $(1,8,12)$, $(1,26,11)$, 
$(1,25,19)$, $(1,8,14)$, $(1,24,16)$, $(1,21,8)$, $(1,8,16)$, $(1,16,26)$, $(1,10,18)$, 
$(1,8,20)$, $(1,22,24)$, $(1,14,14)$, $(1,8,21)$, $(1,15,18)$, $(1,17,2)$, $(1,8,23)$, 
$(1,25,12)$, $(1,26,24)$, $(1,8,24)$, $(1,17,22)$, $(1,15,7)$, $(1,8,25)$, $(1,20,8)$, 
$(1,9,22)$, $(1,8,26)$, $(1,14,7)$, $(1,22,4)$, $(1,9,0)$, $(1,15,13)$, $(1,20,25)$, 
$(1,9,2)$, $(1,18,21)$, $(1,20,19)$, $(1,9,5)$, $(1,24,19)$, $(1,20,18)$, $(1,9,6)$, 
$(1,10,24)$, $(1,20,5)$, $(1,9,8)$, $(1,17,4)$, $(1,20,3)$, $(1,9,10)$, $(1,11,6)$, 
$(1,20,9)$, $(1,9,11)$, $(1,13,20)$, $(1,20,7)$, $(1,9,12)$, $(1,26,10)$, $(1,20,16)$, 
$(1,9,16)$, $(1,12,11)$, $(1,20,0)$, $(1,9,19)$, $(1,21,26)$, $(1,20,6)$, $(1,9,20)$, 
$(1,25,1)$, $(1,20,23)$, $(1,9,21)$, $(1,22,8)$, $(1,20,2)$, $(1,10,1)$, $(1,23,22)$, 
$(1,18,19)$, $(1,10,4)$, $(1,11,1)$, $(1,24,5)$, $(1,10,8)$, $(1,13,13)$, $(1,19,8)$, 
$(1,10,10)$, $(1,25,7)$, $(1,15,10)$, $(1,10,17)$, $(1,19,24)$, $(1,13,0)$, $(1,10,21)$, 
$(1,18,10)$, $(1,23,2)$, $(1,10,25)$, $(1,15,25)$, $(1,25,9)$, $(1,11,0)$, $(1,25,23)$, 
$(1,23,19)$, $(1,11,4)$, $(1,12,13)$, $(1,17,9)$, $(1,11,8)$, $(1,23,16)$, $(1,25,14)$, 
$(1,11,16)$, $(1,18,0)$, $(1,13,11)$, $(1,11,19)$, $(1,22,25)$, $(1,22,7)$, $(1,11,22)$, 
$(1,17,20)$, $(1,12,26)$, $(1,11,24)$, $(1,19,18)$, $(1,11,25)$, $(1,12,6)$, $(1,13,12)$, 
$(1,26,1)$, $(1,12,12)$, $(1,26,26)$, $(1,13,3)$, $(1,12,14)$, $(1,21,6)$, $(1,25,25)$, 
$(1,12,21)$, $(1,25,3)$, $(1,21,2)$, $(1,13,1)$, $(1,21,3)$, $(1,17,7)$, $(1,13,9)$, 
$(1,24,0)$, $(1,14,9)$, $(1,13,14)$, $(1,25,2)$, $(1,13,21)$, $(1,13,17)$, $(1,15,9)$, 
$(1,23,1)$, $(1,13,22)$, $(1,23,4)$, $(1,15,22)$, $(1,13,26)$, $(1,17,10)$, $(1,21,25)$, 
$(1,14,0)$, $(1,17,12)$, $(1,26,22)$, $(1,14,3)$, $(1,15,19)$, $(1,18,20)$, $(1,14,4)$, 
$(1,25,18)$, $(1,25,17)$, $(1,14,20)$, $(1,23,11)$, $(1,19,7)$, $(1,14,23)$, $(1,19,17)$, 
$(1,23,5)$, $(1,14,26)$, $(1,18,1)$, $(1,15,8)$, $(1,16,1)$, $(1,16,16)$, $(1,23,25)$, 
$(1,16,4)$, $(1,17,19)$, $(1,19,19)$, $(1,16,9)$, $(1,19,2)$, $(1,17,21)$, $(1,16,12)$, 
$(1,26,0)$, $(1,24,22)$, $(1,17,6)$, $(1,18,17)$, $(1,22,16)$, $(1,18,16)$, $(1,24,4)$, 
$(1,24,24)$, $(1,19,5)$, $(1,19,11)$, $(1,24,26)$, $(1,19,12)$, $(1,26,6)$, $(1,22,1)$, 
$(1,21,4)$, $(1,23,12)$, $(1,26,20)$, $(1,21,12)$, $(1,26,18)$, $(1,23,7)$, $(1,23,6)$.

\subsection*{$m_{3}(2,29)\ge 44$}
$(0,1,11)$, $(1,5,2)$, $(1,27,10)$, $(1,2,28)$, $(1,2,26)$, $(1,25,14)$, $(1,26,11)$, 
$(1,0,11)$, $(1,0,16)$, $(1,26,20)$, $(1,2,0)$, $(1,3,28)$, $(1,22,23)$, $(1,5,12)$, 
$(1,24,8)$, $(1,0,27)$, $(1,15,24)$, $(1,11,23)$, $(1,1,7)$, $(1,3,10)$, $(1,25,16)$, 
$(1,7,21)$, $(1,1,3)$, $(1,13,11)$, $(1,5,8)$, $(1,3,23)$, $(1,7,0)$, $(1,20,7)$, 
$(1,20,16)$, $(1,1,10)$, $(1,23,12)$, $(1,24,4)$, $(1,9,20)$, $(1,28,6)$, $(1,10,12)$, 
$(1,24,0)$, $(1,7,4)$, $(1,25,2)$, $(1,22,7)$, $(1,10,21)$, $(1,26,25)$, $(1,15,27)$, 
$(1,22,13)$, $(1,21,20)$.

\subsection*{$m_{7}(2,29)\ge 148$}
$(0,0,1)$, $(1,26,3)$, $(0,1,14)$, $(1,3,6)$, $(0,1,25)$, $(1,8,23)$, $(0,1,6)$, 
$(1,17,13)$, $(0,1,19)$, $(1,10,24)$, $(0,1,0)$, $(1,9,9)$, $(0,1,11)$, $(1,13,11)$, 
$(1,0,1)$, $(1,12,11)$, $(1,9,20)$, $(1,15,27)$, $(1,2,2)$, $(1,3,21)$, $(1,1,16)$, 
$(1,22,16)$, $(1,5,18)$, $(1,4,7)$, $(1,18,10)$, $(1,18,14)$, $(1,24,13)$, $(1,20,15)$, 
$(1,0,11)$, $(1,9,21)$, $(1,28,25)$, $(1,1,17)$, $(1,3,27)$, $(1,4,4)$, $(1,16,19)$, 
$(1,21,27)$, $(1,22,22)$, $(1,11,22)$, $(1,27,10)$, $(1,22,13)$, $(1,7,0)$, $(1,7,20)$, 
$(1,0,18)$, $(1,6,2)$, $(1,21,14)$, $(1,16,7)$, $(1,24,1)$, $(1,5,16)$, $(1,12,24)$, 
$(1,20,9)$, $(1,2,19)$, $(1,18,8)$, $(1,13,10)$, $(1,26,12)$, $(1,27,17)$, $(1,23,25)$, 
$(1,0,21)$, $(1,7,18)$, $(1,18,1)$, $(1,11,20)$, $(1,4,23)$, $(1,24,12)$, $(1,2,22)$, 
$(1,1,15)$, $(1,10,26)$, $(1,6,3)$, $(1,7,10)$, $(1,15,22)$, $(1,19,16)$, $(1,8,4)$, 
$(1,0,26)$, $(1,4,28)$, $(1,13,18)$, $(1,26,10)$, $(1,19,21)$, $(1,25,24)$, $(1,24,9)$, 
$(1,0,27)$, $(1,14,14)$, $(1,12,4)$, $(1,5,24)$, $(1,22,9)$, $(1,12,13)$, $(1,11,18)$, 
$(1,13,28)$, $(1,26,11)$, $(1,9,26)$, $(1,24,10)$, $(1,25,5)$, $(1,3,14)$, $(1,19,2)$, 
$(1,1,3)$, $(1,19,4)$, $(1,18,26)$, $(1,12,15)$, $(1,8,21)$, $(1,11,0)$, $(1,19,12)$, 
$(1,15,2)$, $(1,4,19)$, $(1,28,23)$, $(1,6,20)$, $(1,5,26)$, $(1,27,16)$, $(1,10,14)$, 
$(1,1,28)$, $(1,14,13)$, $(1,22,24)$, $(1,6,9)$, $(1,25,11)$, $(1,9,25)$, $(1,13,5)$, 
$(1,26,19)$, $(1,3,0)$, $(1,16,14)$, $(1,14,20)$, $(1,27,5)$, $(1,28,27)$, $(1,12,12)$, 
$(1,2,15)$, $(1,7,4)$, $(1,17,8)$, $(1,25,13)$, $(1,15,7)$, $(1,11,6)$, $(1,23,11)$, 
$(1,9,5)$, $(1,20,24)$, $(1,17,9)$, $(1,3,1)$, $(1,14,22)$, $(1,13,6)$, $(1,26,28)$, 
$(1,4,6)$, $(1,5,25)$, $(1,19,28)$, $(1,11,28)$, $(1,17,27)$, $(1,16,16)$, $(1,25,2)$, 
$(1,25,6)$, $(1,22,15)$, $(1,18,17)$, $(1,5,21)$, $(1,17,2)$, $(1,15,26)$, $(1,21,4)$, 
$(1,17,6)$.

\subsection*{$m_{9}(2,29)\ge 208$}
$(0,1,0)$, $(1,19,3)$, $(1,2,1)$, $(1,26,8)$, $(1,4,18)$, $(1,26,10)$, $(1,13,5)$, 
$(1,22,28)$, $(1,23,0)$, $(1,27,12)$, $(1,21,24)$, $(1,8,7)$, $(1,11,21)$, $(0,1,2)$, 
$(1,16,23)$, $(1,17,6)$, $(1,3,6)$, $(1,0,12)$, $(1,23,16)$, $(1,24,23)$, $(1,0,17)$, 
$(1,17,18)$, $(1,3,16)$, $(1,10,25)$, $(1,8,11)$, $(1,28,9)$, $(0,1,6)$, $(1,23,15)$, 
$(1,18,16)$, $(1,14,12)$, $(1,28,25)$, $(1,21,20)$, $(1,22,25)$, $(1,2,18)$, $(1,14,27)$, 
$(1,19,23)$, $(1,18,19)$, $(1,8,2)$, $(1,23,4)$, $(0,1,8)$, $(1,9,2)$, $(1,4,21)$, 
$(1,12,3)$, $(1,11,14)$, $(1,10,13)$, $(1,16,2)$, $(1,12,23)$, $(1,25,23)$, $(1,4,11)$, 
$(1,11,17)$, $(1,8,13)$, $(1,9,19)$, $(0,1,15)$, $(1,6,22)$, $(1,13,24)$, $(1,9,4)$, 
$(1,8,24)$, $(1,27,8)$, $(0,1,28)$, $(1,11,8)$, $(1,9,13)$, $(1,1,26)$, $(1,22,16)$, 
$(1,8,17)$, $(1,4,14)$, $(0,1,22)$, $(1,5,19)$, $(1,22,27)$, $(1,19,20)$, $(1,13,17)$, 
$(1,2,0)$, $(1,12,6)$, $(1,27,16)$, $(1,12,4)$, $(1,5,6)$, $(1,1,10)$, $(1,8,26)$, 
$(1,21,2)$, $(0,1,23)$, $(1,26,24)$, $(1,15,15)$, $(1,23,9)$, $(1,7,8)$, $(1,11,11)$, 
$(1,5,13)$, $(1,23,14)$, $(1,5,25)$, $(1,28,7)$, $(1,4,15)$, $(1,8,1)$, $(1,15,25)$, 
$(1,0,0)$, $(1,28,24)$, $(1,28,5)$, $(1,18,14)$, $(1,0,8)$, $(1,12,10)$, $(1,6,17)$, 
$(1,5,21)$, $(1,6,3)$, $(1,14,28)$, $(1,15,6)$, $(1,12,17)$, $(1,16,11)$, $(1,0,21)$, 
$(1,4,3)$, $(1,21,22)$, $(1,23,18)$, $(1,21,5)$, $(1,1,27)$, $(1,19,14)$, $(1,9,27)$, 
$(1,19,25)$, $(1,11,22)$, $(1,12,13)$, $(1,7,28)$, $(1,7,12)$, $(1,1,1)$, $(1,18,23)$, 
$(1,12,28)$, $(1,26,5)$, $(1,17,25)$, $(1,3,12)$, $(1,11,6)$, $(1,14,13)$, $(1,25,5)$, 
$(1,15,18)$, $(1,15,28)$, $(1,4,7)$, $(1,23,23)$, $(1,1,5)$, $(1,28,20)$, $(1,16,10)$, 
$(1,2,28)$, $(1,18,11)$, $(1,21,14)$, $(1,17,2)$, $(1,3,28)$, $(1,16,27)$, $(1,19,5)$, 
$(1,14,0)$, $(1,1,19)$, $(1,6,15)$, $(1,1,14)$, $(1,11,28)$, $(1,25,13)$, $(1,22,4)$, 
$(1,10,7)$, $(1,10,16)$, $(1,4,1)$, $(1,16,5)$, $(1,5,12)$, $(1,26,4)$, $(1,22,21)$, 
$(1,5,3)$, $(1,16,18)$, $(1,1,16)$, $(1,2,22)$, $(1,27,4)$, $(1,14,2)$, $(1,14,9)$, 
$(1,23,11)$, $(1,10,26)$, $(1,12,21)$, $(1,10,3)$, $(1,22,17)$, $(1,18,25)$, $(1,18,9)$, 
$(1,13,20)$, $(1,3,0)$, $(1,17,8)$, $(1,3,19)$, $(1,21,19)$, $(1,18,24)$, $(1,7,10)$, 
$(1,5,27)$, $(1,19,8)$, $(1,28,23)$, $(1,7,21)$, $(1,14,3)$, $(1,9,24)$, $(1,17,9)$, 
$(1,3,9)$, $(1,5,4)$, $(1,4,0)$, $(1,13,18)$, $(1,25,9)$, $(1,28,21)$, $(1,27,3)$, 
$(1,6,10)$, $(1,22,2)$, $(1,25,19)$, $(1,11,26)$, $(1,28,13)$, $(1,17,12)$, $(1,3,10)$, 
$(1,25,1)$, $(1,17,14)$, $(1,3,24)$, $(1,26,11)$, $(1,25,7)$, $(1,20,8)$, $(1,21,1)$, 
$(1,9,0)$, $(1,10,11)$, $(1,26,27)$, $(1,19,6)$, $(1,17,7)$.

\subsection*{$m_{10}(2,29)\ge 234$}
$(0,1,10)$, $(1,10,13)$, $(1,23,4)$, $(1,18,26)$, $(1,1,22)$, $(1,25,7)$, $(1,2,13)$, 
$(1,28,27)$, $(1,6,11)$, $(1,19,26)$, $(1,28,1)$, $(1,5,7)$, $(1,15,10)$, $(1,0,2)$, 
$(1,10,23)$, $(1,5,26)$, $(1,6,2)$, $(1,9,7)$, $(1,25,1)$, $(1,7,20)$, $(1,13,0)$, 
$(1,5,16)$, $(1,11,0)$, $(1,0,4)$, $(1,13,28)$, $(1,11,16)$, $(1,0,5)$, $(1,17,25)$, 
$(1,27,28)$, $(1,5,2)$, $(1,21,25)$, $(1,27,13)$, $(1,12,23)$, $(1,4,23)$, $(1,12,18)$, 
$(1,9,10)$, $(1,9,28)$, $(1,16,18)$, $(1,11,3)$, $(1,0,8)$, $(1,18,17)$, $(1,15,19)$, 
$(1,16,2)$, $(1,2,11)$, $(1,8,15)$, $(1,16,8)$, $(1,20,24)$, $(1,11,26)$, $(1,0,26)$, 
$(1,12,7)$, $(1,28,7)$, $(1,11,2)$, $(1,0,15)$, $(1,9,2)$, $(1,19,22)$, $(1,1,2)$, 
$(1,28,21)$, $(1,10,27)$, $(1,4,24)$, $(1,1,21)$, $(1,23,17)$, $(1,25,17)$, $(1,11,14)$, 
$(1,0,23)$, $(1,11,15)$, $(1,0,19)$, $(1,7,18)$, $(1,25,12)$, $(1,19,2)$, $(1,8,20)$, 
$(1,12,10)$, $(1,20,22)$, $(1,7,25)$, $(1,27,14)$, $(1,23,27)$, $(1,24,10)$, $(1,10,9)$, 
$(1,11,9)$, $(1,0,21)$, $(1,15,12)$, $(1,14,11)$, $(1,24,2)$, $(1,1,24)$, $(1,15,28)$, 
$(1,8,9)$, $(1,17,22)$, $(1,20,12)$, $(1,24,22)$, $(1,5,27)$, $(1,7,19)$, $(1,11,27)$, 
$(1,1,4)$, $(1,24,12)$, $(1,21,10)$, $(1,16,11)$, $(1,15,22)$, $(1,12,20)$, $(1,6,12)$, 
$(1,18,27)$, $(1,12,1)$, $(1,24,0)$, $(1,2,3)$, $(1,1,5)$, $(1,5,20)$, $(1,1,9)$, 
$(1,7,10)$, $(1,28,15)$, $(1,1,20)$, $(1,17,18)$, $(1,13,25)$, $(1,27,16)$, $(1,13,26)$, 
$(1,17,16)$, $(1,16,7)$, $(1,8,25)$, $(1,27,17)$, $(1,20,18)$, $(1,1,27)$, $(1,13,9)$, 
$(1,12,16)$, $(1,18,4)$, $(1,19,14)$, $(1,6,19)$, $(1,2,14)$, $(1,4,1)$, $(1,18,19)$, 
$(1,6,23)$, $(1,23,22)$, $(1,2,1)$, $(1,9,4)$, $(1,2,4)$, $(1,23,24)$, $(1,23,3)$, 
$(1,15,5)$, $(1,6,18)$, $(1,16,0)$, $(1,6,1)$, $(1,21,24)$, $(1,19,1)$, $(1,19,11)$, 
$(1,18,25)$, $(1,27,4)$, $(1,17,9)$, $(1,2,7)$, $(1,5,19)$, $(1,3,9)$, $(1,5,25)$, 
$(1,27,11)$, $(1,10,17)$, $(1,8,3)$, $(1,22,5)$, $(1,20,19)$, $(1,25,15)$, $(1,18,18)$, 
$(1,21,26)$, $(1,15,24)$, $(1,2,15)$, $(1,20,28)$, $(1,10,4)$, $(1,12,11)$, $(1,17,24)$, 
$(1,6,9)$, $(1,25,20)$, $(1,10,1)$, $(1,6,28)$, $(1,28,17)$, $(1,18,1)$, $(1,25,21)$, 
$(1,8,4)$, $(1,2,22)$, $(1,10,22)$, $(1,15,17)$, $(1,4,27)$, $(1,8,27)$, $(1,25,18)$, 
$(1,20,15)$, $(1,17,13)$, $(1,4,21)$, $(1,16,9)$, $(1,18,0)$, $(1,16,25)$, $(1,27,21)$, 
$(1,2,24)$, $(1,16,14)$, $(1,12,15)$, $(1,7,21)$, $(1,3,19)$, $(1,8,13)$, $(1,14,9)$, 
$(1,20,14)$, $(1,5,10)$, $(1,22,13)$, $(1,18,10)$, $(1,5,17)$, $(1,28,28)$, $(1,3,20)$, 
$(1,19,27)$, $(1,9,6)$, $(1,6,24)$, $(1,26,19)$, $(1,22,1)$, $(1,28,26)$, $(1,23,21)$, 
$(1,7,7)$, $(1,7,23)$, $(1,17,4)$, $(1,28,22)$, $(1,9,23)$, $(1,4,7)$, $(1,17,14)$, 
$(1,15,13)$, $(1,19,3)$, $(1,7,5)$, $(1,24,11)$, $(1,7,14)$, $(1,10,26)$, $(1,21,0)$, 
$(1,23,7)$, $(1,16,12)$, $(1,25,26)$, $(1,4,10)$, $(1,8,28)$, $(1,21,20)$, $(1,20,17)$, 
$(1,23,0)$, $(1,10,28)$, $(1,19,25)$, $(1,27,15)$, $(1,15,3)$, $(1,25,5)$, $(1,28,3)$, 
$(1,19,5)$, $(1,23,13)$, $(1,21,18)$.

\subsection*{$m_{11}(2,29)\ge 262$}
$(0,1,2)$, $(1,3,6)$, $(1,19,12)$, $(1,10,21)$, $(1,21,17)$, $(1,28,28)$, $(1,6,3)$, 
$(1,16,26)$, $(1,15,0)$, $(1,12,10)$, $(1,24,6)$, $(1,9,21)$, $(1,21,28)$, $(1,5,6)$, 
$(0,1,16)$, $(1,0,14)$, $(1,19,11)$, $(1,1,26)$, $(1,20,23)$, $(1,28,23)$, $(1,24,19)$, 
$(1,15,25)$, $(1,14,17)$, $(1,1,19)$, $(1,20,3)$, $(1,5,19)$, $(1,10,18)$, $(1,2,10)$, 
$(0,1,19)$, $(1,2,28)$, $(1,19,22)$, $(1,14,22)$, $(1,22,11)$, $(1,28,27)$, $(1,25,7)$, 
$(1,25,6)$, $(1,27,28)$, $(1,9,23)$, $(1,1,25)$, $(1,20,12)$, $(1,18,20)$, $(1,9,20)$, 
$(0,1,20)$, $(1,11,4)$, $(1,19,4)$, $(1,11,14)$, $(1,16,18)$, $(1,28,6)$, $(1,1,5)$, 
$(1,20,1)$, $(1,3,1)$, $(1,2,5)$, $(1,12,26)$, $(1,6,5)$, $(1,11,11)$, $(1,24,0)$, 
$(1,0,0)$, $(1,14,21)$, $(1,21,10)$, $(1,18,26)$, $(1,3,26)$, $(1,9,13)$, $(1,24,27)$, 
$(1,24,2)$, $(1,16,16)$, $(1,2,23)$, $(1,3,27)$, $(1,27,17)$, $(1,6,11)$, $(1,18,3)$, 
$(1,0,5)$, $(1,0,20)$, $(1,4,12)$, $(1,25,1)$, $(1,13,19)$, $(1,17,22)$, $(1,11,8)$, 
$(1,2,6)$, $(1,10,23)$, $(1,25,20)$, $(1,7,11)$, $(1,16,4)$, $(1,17,21)$, $(1,21,22)$, 
$(1,0,7)$, $(1,10,0)$, $(1,9,8)$, $(1,3,5)$, $(1,14,27)$, $(1,5,23)$, $(1,26,21)$, 
$(1,21,21)$, $(1,21,15)$, $(1,11,13)$, $(1,8,7)$, $(1,18,9)$, $(1,4,25)$, $(1,6,14)$, 
$(1,0,13)$, $(1,16,17)$, $(1,14,4)$, $(1,14,3)$, $(1,23,12)$, $(1,1,4)$, $(1,20,10)$, 
$(1,13,4)$, $(1,4,1)$, $(1,27,21)$, $(1,21,13)$, $(1,9,1)$, $(1,23,8)$, $(1,7,1)$, 
$(1,0,16)$, $(1,8,4)$, $(1,15,9)$, $(1,17,13)$, $(1,10,24)$, $(1,15,27)$, $(1,9,14)$, 
$(1,15,1)$, $(1,6,18)$, $(1,6,25)$, $(1,17,0)$, $(1,26,0)$, $(1,18,14)$, $(1,27,2)$, 
$(1,0,19)$, $(1,3,14)$, $(1,6,22)$, $(1,7,28)$, $(1,7,0)$, $(1,6,6)$, $(1,25,24)$, 
$(1,6,0)$, $(1,2,13)$, $(1,14,0)$, $(1,10,28)$, $(1,13,11)$, $(1,13,20)$, $(1,23,25)$, 
$(1,0,26)$, $(1,17,15)$, $(1,5,17)$, $(1,24,17)$, $(1,11,3)$, $(1,13,3)$, $(1,12,5)$, 
$(1,9,10)$, $(1,17,10)$, $(1,23,19)$, $(1,14,12)$, $(1,7,25)$, $(1,11,5)$, $(1,15,13)$, 
$(1,1,1)$, $(1,20,20)$, $(1,27,14)$, $(1,23,21)$, $(1,21,23)$, $(1,27,22)$, $(1,2,24)$, 
$(1,27,5)$, $(1,26,4)$, $(1,25,11)$, $(1,3,25)$, $(1,3,2)$, $(1,8,16)$, $(1,11,9)$, 
$(1,1,3)$, $(1,20,18)$, $(1,22,27)$, $(1,14,10)$, $(1,9,27)$, $(1,16,20)$, $(1,18,17)$, 
$(1,15,14)$, $(1,10,3)$, $(1,18,22)$, $(1,23,6)$, $(1,18,10)$, $(1,12,17)$, $(1,25,27)$, 
$(1,1,13)$, $(1,20,22)$, $(1,25,25)$, $(1,8,22)$, $(1,5,9)$, $(1,25,19)$, $(1,23,13)$, 
$(1,25,21)$, $(1,21,20)$, $(1,12,19)$, $(1,12,15)$, $(1,15,20)$, $(1,3,22)$, $(1,26,20)$, 
$(1,1,14)$, $(1,20,28)$, $(1,16,2)$, $(1,11,16)$, $(1,10,17)$, $(1,7,21)$, $(1,21,3)$, 
$(1,3,23)$, $(1,4,28)$, $(1,23,10)$, $(1,10,14)$, $(1,5,5)$, $(1,5,8)$, $(1,13,24)$, 
$(1,2,9)$, $(1,9,7)$, $(1,11,12)$, $(1,3,11)$, $(1,17,18)$, $(1,18,12)$, $(1,16,13)$, 
$(1,12,12)$, $(1,26,19)$, $(1,26,23)$, $(1,10,5)$, $(1,17,14)$, $(1,8,27)$, $(1,8,9)$, 
$(1,2,11)$, $(1,2,26)$, $(1,8,20)$, $(1,13,9)$, $(1,24,11)$, $(1,14,5)$, $(1,4,22)$, 
$(1,15,17)$, $(1,16,24)$, $(1,7,3)$, $(1,25,8)$, $(1,11,26)$, $(1,26,10)$, $(1,15,7)$, 
$(1,24,23)$, $(1,4,9)$, $(1,5,4)$, $(1,27,23)$, $(1,22,17)$, $(1,26,8)$, $(1,5,3)$, 
$(1,15,28)$, $(1,26,1)$, $(1,8,10)$, $(1,26,9)$, $(1,27,6)$, $(1,8,19)$, $(1,16,9)$, 
$(1,8,28)$, $(1,10,19)$, $(1,27,7)$, $(1,17,26)$, $(1,12,21)$, $(1,21,1)$, $(1,24,4)$, 
$(1,18,27)$, $(1,16,25)$, $(1,17,28)$.

\subsection*{$m_{12}(2,29)\ge 300$}
$(0,1,2)$, $(1,11,0)$, $(1,14,20)$, $(1,18,1)$, $(1,13,22)$, $(1,2,14)$, $(1,14,1)$, 
$(1,17,23)$, $(1,18,5)$, $(1,27,7)$, $(1,0,8)$, $(1,15,10)$, $(1,20,17)$, $(1,2,23)$, 
$(1,15,0)$, $(0,1,5)$, $(1,16,7)$, $(1,6,27)$, $(1,2,2)$, $(1,10,5)$, $(1,13,15)$, 
$(1,8,2)$, $(1,24,10)$, $(1,2,4)$, $(1,17,27)$, $(1,24,8)$, $(1,0,17)$, $(1,20,14)$, 
$(1,17,3)$, $(1,16,4)$, $(0,1,7)$, $(1,17,20)$, $(1,0,25)$, $(1,1,22)$, $(1,14,18)$, 
$(1,22,0)$, $(1,3,27)$, $(1,26,27)$, $(1,17,14)$, $(1,26,9)$, $(1,4,8)$, $(1,13,9)$, 
$(1,20,7)$, $(1,5,19)$, $(1,14,25)$, $(0,1,10)$, $(1,10,16)$, $(1,18,2)$, $(1,14,23)$, 
$(1,26,28)$, $(1,7,25)$, $(1,24,9)$, $(1,16,0)$, $(1,8,8)$, $(1,6,20)$, $(1,20,8)$, 
$(1,10,22)$, $(1,20,20)$, $(0,1,18)$, $(1,6,22)$, $(0,1,11)$, $(1,9,3)$, $(1,15,1)$, 
$(1,28,4)$, $(1,5,25)$, $(1,17,18)$, $(1,11,16)$, $(1,23,16)$, $(1,4,15)$, $(1,1,1)$, 
$(1,15,8)$, $(1,17,11)$, $(1,20,6)$, $(1,9,4)$, $(1,27,19)$, $(0,1,14)$, $(1,13,26)$, 
$(1,26,24)$, $(1,12,5)$, $(1,9,9)$, $(1,26,3)$, $(1,18,10)$, $(1,0,9)$, $(1,6,26)$, 
$(1,9,14)$, $(1,10,8)$, $(1,14,24)$, $(1,20,15)$, $(1,14,7)$, $(1,12,17)$, $(0,1,16)$, 
$(1,21,14)$, $(1,13,10)$, $(1,11,25)$, $(1,8,13)$, $(1,12,7)$, $(1,16,20)$, $(1,9,13)$, 
$(1,0,22)$, $(1,8,16)$, $(1,16,8)$, $(1,27,16)$, $(1,20,25)$, $(1,21,17)$, $(1,9,5)$, 
$(0,1,19)$, $(1,0,2)$, $(1,4,7)$, $(1,24,26)$, $(1,12,26)$, $(1,24,16)$, $(1,22,19)$, 
$(1,3,20)$, $(1,16,23)$, $(1,5,22)$, $(1,13,8)$, $(1,12,23)$, $(1,20,18)$, $(1,16,14)$, 
$(1,1,2)$, $(0,1,20)$, $(1,2,28)$, $(1,7,8)$, $(1,9,7)$, $(1,20,23)$, $(1,4,1)$, 
$(1,21,24)$, $(0,1,23)$, $(1,22,27)$, $(1,3,26)$, $(1,22,8)$, $(1,3,4)$, $(1,20,13)$, 
$(1,15,25)$, $(1,4,14)$, $(0,1,28)$, $(1,12,13)$, $(1,21,3)$, $(1,5,0)$, $(1,4,0)$, 
$(1,28,19)$, $(1,27,23)$, $(1,12,24)$, $(1,1,13)$, $(1,23,15)$, $(1,21,8)$, $(1,2,18)$, 
$(1,20,24)$, $(1,11,11)$, $(1,2,6)$, $(1,0,5)$, $(1,16,5)$, $(1,12,28)$, $(1,7,10)$, 
$(1,17,10)$, $(1,14,3)$, $(1,13,24)$, $(1,22,3)$, $(1,3,23)$, $(1,13,16)$, $(1,9,0)$, 
$(1,2,25)$, $(1,0,15)$, $(1,22,4)$, $(1,3,15)$, $(1,0,6)$, $(1,18,24)$, $(1,18,27)$, 
$(1,15,24)$, $(1,27,2)$, $(1,6,16)$, $(1,28,1)$, $(1,10,14)$, $(1,24,27)$, $(1,13,5)$, 
$(1,23,1)$, $(1,11,7)$, $(1,22,10)$, $(1,3,16)$, $(1,15,11)$, $(1,0,7)$, $(1,26,13)$, 
$(1,6,0)$, $(1,23,9)$, $(1,10,4)$, $(1,11,26)$, $(1,21,4)$, $(1,8,11)$, $(1,8,6)$, 
$(1,13,0)$, $(1,21,5)$, $(1,22,14)$, $(1,3,13)$, $(1,27,10)$, $(1,28,26)$, $(1,0,23)$, 
$(1,2,17)$, $(1,11,4)$, $(1,6,1)$, $(1,6,13)$, $(1,10,24)$, $(1,10,17)$, $(1,27,25)$, 
$(1,9,20)$, $(1,13,27)$, $(1,27,22)$, $(1,15,28)$, $(1,7,20)$, $(1,21,26)$, $(1,11,22)$, 
$(1,1,10)$, $(1,16,27)$, $(1,16,19)$, $(1,21,15)$, $(1,27,13)$, $(1,26,4)$, $(1,22,23)$, 
$(1,3,9)$, $(1,18,18)$, $(1,4,13)$, $(1,5,23)$, $(1,11,3)$, $(1,26,5)$, $(1,21,27)$, 
$(1,18,9)$, $(1,1,11)$, $(1,21,6)$, $(1,15,27)$, $(1,5,16)$, $(1,12,20)$, $(1,8,24)$, 
$(1,14,22)$, $(1,28,13)$, $(1,1,28)$, $(1,7,1)$, $(1,5,20)$, $(1,23,4)$, $(1,16,22)$, 
$(1,18,0)$, $(1,17,26)$, $(1,1,19)$, $(1,4,2)$, $(1,8,25)$, $(1,22,24)$, $(1,3,1)$, 
$(1,24,3)$, $(1,11,18)$, $(1,15,26)$, $(1,14,5)$, $(1,14,2)$, $(1,5,26)$, $(1,8,10)$, 
$(1,15,15)$, $(1,24,25)$, $(1,27,1)$, $(1,2,7)$, $(1,4,11)$, $(1,23,6)$, $(1,17,25)$, 
$(1,6,28)$, $(1,7,23)$, $(1,10,11)$, $(1,12,6)$, $(1,28,14)$, $(1,4,4)$, $(1,9,16)$, 
$(1,14,14)$, $(1,23,7)$, $(1,28,0)$, $(1,16,24)$, $(1,2,19)$, $(1,12,3)$, $(1,4,26)$, 
$(1,27,28)$, $(1,7,27)$, $(1,11,15)$, $(1,4,18)$, $(1,17,2)$, $(1,9,15)$, $(1,18,20)$, 
$(1,11,20)$, $(1,24,22)$, $(1,26,22)$, $(1,27,3)$, $(1,8,9)$, $(1,3,5)$, $(1,6,19)$, 
$(1,24,19)$, $(1,10,10)$, $(1,21,0)$, $(1,12,16)$, $(1,17,5)$, $(1,15,18)$, $(1,6,18)$, 
$(1,26,10)$, $(1,10,13)$, $(1,22,9)$, $(1,3,28)$, $(1,7,2)$, $(1,22,1)$.

\subsection*{$m_{19}(2,29)\ge 507$}
$(0,0,1)$, $(1,11,24)$, $(1,0,16)$, $(1,16,4)$, $(1,2,8)$, $(1,27,22)$, $(1,25,6)$, 
$(1,26,20)$, $(1,10,16)$, $(1,24,23)$, $(1,19,0)$, $(1,4,21)$, $(1,20,13)$, $(0,1,1)$, 
$(1,28,12)$, $(1,27,8)$, $(1,24,28)$, $(1,13,18)$, $(1,24,17)$, $(1,28,2)$, $(1,2,4)$, 
$(1,15,9)$, $(1,2,21)$, $(1,10,6)$, $(1,3,9)$, $(0,1,19)$, $(0,1,2)$, $(1,0,13)$, 
$(1,21,13)$, $(1,21,19)$, $(1,10,10)$, $(1,1,27)$, $(1,4,5)$, $(1,15,3)$, $(1,0,1)$, 
$(1,27,18)$, $(1,1,12)$, $(1,17,3)$, $(1,4,28)$, $(0,1,5)$, $(1,8,21)$, $(1,1,20)$, 
$(1,3,23)$, $(1,18,12)$, $(1,15,2)$, $(1,23,28)$, $(1,8,8)$, $(1,9,0)$, $(1,19,12)$, 
$(1,28,23)$, $(1,7,28)$, $(1,16,24)$, $(0,1,7)$, $(1,22,6)$, $(1,13,10)$, $(1,19,13)$, 
$(1,21,20)$, $(1,22,4)$, $(1,16,18)$, $(1,10,19)$, $(1,13,6)$, $(1,28,26)$, $(1,9,26)$, 
$(1,5,4)$, $(1,27,1)$, $(0,1,8)$, $(1,23,7)$, $(1,9,23)$, $(1,22,22)$, $(0,1,22)$, 
$(1,25,9)$, $(1,8,19)$, $(1,20,16)$, $(1,11,3)$, $(1,20,20)$, $(1,11,15)$, $(1,12,1)$, 
$(1,26,11)$, $(0,1,10)$, $(1,16,0)$, $(1,18,1)$, $(1,10,15)$, $(1,7,2)$, $(1,5,24)$, 
$(1,22,10)$, $(1,4,15)$, $(1,27,27)$, $(1,0,5)$, $(1,5,19)$, $(1,19,27)$, $(1,15,5)$, 
$(0,1,11)$, $(1,13,26)$, $(1,7,15)$, $(1,13,24)$, $(1,15,4)$, $(1,18,7)$, $(1,12,4)$, 
$(1,5,6)$, $(1,20,2)$, $(1,11,6)$, $(1,14,13)$, $(1,21,22)$, $(1,1,0)$, $(0,1,16)$, 
$(1,26,10)$, $(1,20,9)$, $(1,11,18)$, $(1,25,21)$, $(0,1,21)$, $(1,3,16)$, $(1,25,0)$, 
$(1,1,17)$, $(1,5,16)$, $(1,7,8)$, $(1,26,24)$, $(1,17,14)$, $(0,1,17)$, $(1,19,3)$, 
$(1,16,22)$, $(1,8,9)$, $(1,16,26)$, $(1,14,10)$, $(1,15,0)$, $(1,24,9)$, $(1,6,10)$, 
$(1,7,3)$, $(1,12,24)$, $(1,2,26)$, $(1,19,23)$, $(0,1,18)$, $(1,27,11)$, $(1,23,21)$, 
$(1,25,2)$, $(1,3,1)$, $(1,6,16)$, $(1,10,26)$, $(0,1,20)$, $(1,24,8)$, $(1,10,27)$, 
$(1,26,5)$, $(1,27,7)$, $(1,2,19)$, $(0,1,23)$, $(1,7,20)$, $(1,2,24)$, $(1,28,11)$, 
$(1,22,13)$, $(1,21,12)$, $(1,6,12)$, $(1,9,28)$, $(1,18,28)$, $(1,22,7)$, $(1,17,11)$, 
$(1,24,0)$, $(1,14,15)$, $(0,1,27)$, $(1,9,22)$, $(1,15,18)$, $(1,6,3)$, $(1,1,15)$, 
$(1,10,13)$, $(1,21,21)$, $(1,23,18)$, $(1,3,20)$, $(1,12,14)$, $(1,27,14)$, $(1,18,15)$, 
$(1,6,8)$, $(0,1,28)$, $(1,6,19)$, $(1,5,7)$, $(1,5,0)$, $(1,17,19)$, $(1,9,21)$, 
$(1,7,1)$, $(1,3,24)$, $(1,4,7)$, $(1,14,1)$, $(1,0,3)$, $(1,15,8)$, $(1,28,20)$, 
$(1,0,0)$, $(1,8,7)$, $(1,15,14)$, $(1,7,13)$, $(1,21,24)$, $(1,13,23)$, $(1,13,14)$, 
$(1,14,9)$, $(1,24,4)$, $(1,8,4)$, $(1,23,23)$, $(1,4,6)$, $(1,2,6)$, $(1,0,4)$, 
$(1,6,15)$, $(1,14,18)$, $(1,14,6)$, $(1,10,1)$, $(1,22,12)$, $(1,20,19)$, $(1,11,28)$, 
$(1,23,9)$, $(1,27,3)$, $(1,13,9)$, $(1,4,1)$, $(1,23,19)$, $(1,0,7)$, $(1,1,6)$, 
$(1,23,11)$, $(1,1,19)$, $(1,27,26)$, $(1,10,17)$, $(1,19,10)$, $(1,12,12)$, $(1,13,1)$, 
$(1,16,28)$, $(1,27,17)$, $(1,4,4)$, $(1,10,22)$, $(1,0,8)$, $(1,25,26)$, $(1,24,7)$, 
$(1,4,16)$, $(1,9,20)$, $(1,0,26)$, $(1,2,2)$, $(1,6,21)$, $(1,26,23)$, $(1,3,18)$, 
$(1,9,15)$, $(1,4,20)$, $(1,8,18)$, $(1,0,9)$, $(1,10,28)$, $(1,4,0)$, $(1,6,14)$, 
$(1,0,17)$, $(1,7,11)$, $(1,27,24)$, $(1,20,0)$, $(1,11,11)$, $(1,17,2)$, $(1,7,18)$, 
$(1,4,24)$, $(1,5,12)$, $(1,0,11)$, $(1,12,20)$, $(1,16,10)$, $(1,23,26)$, $(1,20,14)$, 
$(1,11,19)$, $(1,1,22)$, $(1,9,2)$, $(1,28,13)$, $(1,21,14)$, $(1,28,1)$, $(1,4,2)$, 
$(1,19,11)$, $(1,0,12)$, $(1,26,22)$, $(1,26,28)$, $(1,19,1)$, $(1,5,9)$, $(1,23,14)$, 
$(1,16,12)$, $(1,25,7)$, $(1,25,28)$, $(1,12,16)$, $(1,18,16)$, $(1,4,19)$, $(1,18,9)$, 
$(1,0,19)$, $(1,24,1)$, $(1,17,6)$, $(1,5,15)$, $(1,18,23)$, $(1,2,1)$, $(1,24,26)$, 
$(1,18,3)$, $(1,8,26)$, $(1,13,19)$, $(1,2,11)$, $(1,4,14)$, $(1,1,4)$, $(1,0,21)$, 
$(1,19,21)$, $(1,3,4)$, $(1,3,17)$, $(1,16,3)$, $(1,19,6)$, $(1,22,8)$, $(1,5,8)$, 
$(1,16,15)$, $(1,20,11)$, $(1,11,12)$, $(1,4,11)$, $(1,13,28)$, $(1,0,24)$, $(1,9,3)$, 
$(1,25,3)$, $(1,9,11)$, $(1,26,16)$, $(1,15,27)$, $(1,12,5)$, $(1,28,17)$, $(1,15,20)$, 
$(1,26,0)$, $(1,25,20)$, $(1,4,22)$, $(1,7,16)$, $(1,0,27)$, $(1,23,5)$, $(1,22,15)$, 
$(1,28,21)$, $(1,28,7)$, $(1,3,3)$, $(1,14,23)$, $(1,23,10)$, $(1,6,7)$, $(1,19,8)$, 
$(1,8,2)$, $(1,4,27)$, $(1,16,5)$, $(1,1,1)$, $(1,25,4)$, $(1,14,0)$, $(1,22,16)$, 
$(1,26,6)$, $(1,7,21)$, $(1,15,24)$, $(1,3,2)$, $(1,26,19)$, $(1,28,15)$, $(1,25,10)$, 
$(1,13,8)$, $(1,1,16)$, $(1,1,7)$, $(1,13,17)$, $(1,26,7)$, $(1,6,13)$, $(1,21,4)$, 
$(1,24,3)$, $(1,5,13)$, $(1,21,6)$, $(1,25,17)$, $(1,2,9)$, $(1,17,16)$, $(1,6,2)$, 
$(1,22,2)$, $(1,1,8)$, $(1,20,7)$, $(1,11,20)$, $(1,18,8)$, $(1,19,9)$, $(1,14,17)$, 
$(1,1,26)$, $(1,16,21)$, $(1,9,14)$, $(1,3,7)$, $(1,20,21)$, $(1,11,27)$, $(1,3,5)$, 
$(1,1,10)$, $(1,24,22)$, $(1,2,22)$, $(1,20,12)$, $(1,11,0)$, $(1,2,28)$, $(1,14,20)$, 
$(1,7,19)$, $(1,14,24)$, $(1,12,18)$, $(1,12,27)$, $(1,24,5)$, $(1,9,1)$, $(1,1,24)$, 
$(1,14,28)$, $(1,5,2)$, $(1,25,22)$, $(1,18,26)$, $(1,23,16)$, $(1,28,18)$, $(1,17,18)$, 
$(1,23,13)$, $(1,21,0)$, $(1,26,2)$, $(1,12,3)$, $(1,7,12)$, $(1,2,3)$, $(1,22,19)$, 
$(1,12,17)$, $(1,20,22)$, $(1,11,17)$, $(1,10,18)$, $(1,20,1)$, $(1,11,22)$, $(1,6,23)$, 
$(1,24,27)$, $(1,27,5)$, $(1,12,2)$, $(1,8,15)$, $(1,2,10)$, $(1,8,28)$, $(1,10,12)$, 
$(1,12,15)$, $(1,22,28)$, $(1,25,13)$, $(1,21,2)$, $(1,27,16)$, $(1,19,22)$, $(1,24,13)$, 
$(1,21,23)$, $(1,16,7)$, $(1,26,12)$, $(1,2,14)$, $(1,23,8)$, $(1,23,1)$, $(1,18,13)$, 
$(1,21,27)$, $(1,20,5)$, $(1,11,21)$, $(1,5,17)$, $(1,3,21)$, $(1,24,10)$, $(1,3,19)$, 
$(1,15,13)$, $(1,21,8)$, $(1,2,23)$, $(1,25,15)$, $(1,5,14)$, $(1,9,16)$, $(1,8,14)$, 
$(1,5,10)$, $(1,10,20)$, $(1,23,3)$, $(1,26,17)$, $(1,24,12)$, $(1,22,20)$, $(1,27,28)$, 
$(1,9,10)$, $(1,3,6)$, $(1,17,22)$, $(1,17,4)$, $(1,12,11)$, $(1,25,8)$, $(1,6,6)$, 
$(1,15,16)$, $(1,13,15)$, $(1,23,22)$, $(1,12,10)$, $(1,14,4)$, $(1,26,21)$, $(1,13,11)$, 
$(1,3,10)$, $(1,22,1)$, $(1,7,23)$, $(1,17,24)$, $(1,10,8)$, $(1,15,17)$, $(1,8,23)$, 
$(1,14,2)$, $(1,17,17)$, $(1,27,19)$, $(1,8,3)$, $(1,18,21)$, $(1,18,0)$, $(1,3,14)$, 
$(1,19,2)$, $(1,13,0)$, $(1,10,0)$, $(1,16,8)$, $(1,7,4)$, $(1,7,24)$, $(1,19,24)$, 
$(1,27,6)$, $(1,6,18)$, $(1,16,14)$, $(1,22,21)$, $(1,14,3)$, $(1,5,28)$, $(1,6,4)$, 
$(1,13,20)$, $(1,6,20)$, $(1,20,6)$, $(1,11,4)$, $(1,28,10)$, $(1,18,17)$, $(1,9,6)$, 
$(1,9,27)$, $(1,14,5)$, $(1,13,12)$, $(1,18,18)$, $(1,12,9)$, $(1,15,7)$, $(1,22,18)$, 
$(1,19,26)$, $(1,26,18)$, $(1,18,22)$, $(1,16,11)$, $(1,15,15)$, $(1,19,17)$, $(1,18,27)$, 
$(1,25,5)$, $(1,16,13)$, $(1,21,16)$.

\subsection*{$m_{20}(2,29)\ge 534$}
$(0,1,0)$, $(1,20,9)$, $(1,21,17)$, $(1,13,7)$, $(1,14,5)$, $(1,10,22)$, $(1,10,8)$, 
$(1,9,11)$, $(1,6,0)$, $(1,11,21)$, $(1,6,2)$, $(1,20,4)$, $(1,11,16)$, $(1,17,1)$, 
$(0,1,1)$, $(1,21,16)$, $(1,22,24)$, $(1,2,17)$, $(1,2,8)$, $(1,0,10)$, $(1,16,21)$, 
$(1,19,23)$, $(1,26,24)$, $(1,3,23)$, $(1,3,10)$, $(1,6,22)$, $(1,17,0)$, $(1,8,25)$, 
$(0,1,2)$, $(1,5,20)$, $(1,6,28)$, $(1,4,2)$, $(1,20,18)$, $(1,15,28)$, $(1,7,16)$, 
$(1,4,5)$, $(1,25,17)$, $(1,15,20)$, $(1,22,27)$, $(1,27,24)$, $(1,8,24)$, $(1,7,18)$, 
$(0,1,5)$, $(1,6,27)$, $(1,7,6)$, $(1,22,12)$, $(1,8,21)$, $(1,5,16)$, $(1,13,0)$, 
$(1,14,17)$, $(1,16,12)$, $(1,7,22)$, $(1,19,6)$, $(1,13,13)$, $(1,14,8)$, $(1,27,13)$, 
$(0,1,6)$, $(1,27,0)$, $(1,28,8)$, $(1,23,19)$, $(1,17,26)$, $(1,27,25)$, $(1,23,12)$, 
$(1,21,8)$, $(1,1,23)$, $(1,13,6)$, $(1,14,0)$, $(1,9,14)$, $(1,24,20)$, $(1,12,24)$, 
$(0,1,7)$, $(0,1,9)$, $(1,22,23)$, $(1,23,2)$, $(1,20,27)$, $(1,19,11)$, $(1,19,27)$, 
$(1,22,5)$, $(1,0,6)$, $(1,17,19)$, $(1,24,25)$, $(1,0,18)$, $(1,21,11)$, $(1,23,13)$, 
$(1,28,20)$, $(0,1,12)$, $(1,23,1)$, $(1,24,9)$, $(1,9,8)$, $(1,7,14)$, $(1,9,15)$, 
$(1,28,18)$, $(1,10,18)$, $(1,8,14)$, $(1,16,27)$, $(1,26,26)$, $(1,7,0)$, $(1,0,26)$, 
$(1,19,15)$, $(0,1,13)$, $(1,7,5)$, $(1,8,13)$, $(1,11,22)$, $(1,25,24)$, $(1,24,4)$, 
$(1,19,13)$, $(1,24,0)$, $(1,7,7)$, $(1,28,24)$, $(1,16,14)$, $(1,28,2)$, $(1,20,21)$, 
$(1,18,8)$, $(0,1,14)$, $(1,8,12)$, $(1,9,20)$, $(1,0,3)$, $(1,13,27)$, $(1,14,21)$, 
$(1,25,26)$, $(1,5,12)$, $(1,27,2)$, $(1,20,26)$, $(1,13,22)$, $(1,14,20)$, $(1,26,5)$, 
$(1,9,3)$, $(0,1,15)$, $(1,16,10)$, $(1,17,18)$, $(1,28,25)$, $(1,4,22)$, $(1,21,12)$, 
$(1,15,14)$, $(1,27,21)$, $(1,13,20)$, $(1,14,13)$, $(1,18,28)$, $(1,18,19)$, $(1,16,22)$, 
$(1,24,21)$, $(0,1,16)$, $(1,19,2)$, $(1,20,10)$, $(1,24,26)$, $(1,26,2)$, $(1,20,5)$, 
$(1,4,24)$, $(1,28,28)$, $(1,15,5)$, $(1,19,19)$, $(1,9,23)$, $(1,5,15)$, $(1,5,3)$, 
$(1,26,6)$, $(0,1,18)$, $(1,26,22)$, $(1,27,1)$, $(1,5,9)$, $(1,0,23)$, $(1,8,8)$, 
$(1,17,28)$, $(1,11,25)$, $(1,10,28)$, $(1,21,4)$, $(1,17,21)$, $(1,23,25)$, $(1,18,7)$, 
$(1,21,0)$, $(0,1,20)$, $(1,13,18)$, $(1,14,26)$, $(1,3,24)$, $(1,11,13)$, $(1,22,19)$, 
$(1,26,4)$, $(1,26,14)$, $(1,11,6)$, $(1,9,7)$, $(1,27,4)$, $(1,2,23)$, $(1,27,12)$, 
$(1,22,7)$, $(0,1,22)$, $(1,17,17)$, $(1,18,25)$, $(1,17,6)$, $(1,21,25)$, $(1,11,0)$, 
$(1,21,27)$, $(1,8,4)$, $(1,4,15)$, $(1,6,15)$, $(1,15,7)$, $(1,4,8)$, $(1,22,6)$, 
$(1,15,16)$, $(0,1,23)$, $(1,15,3)$, $(1,16,11)$, $(1,10,15)$, $(1,16,19)$, $(1,2,24)$, 
$(1,9,1)$, $(1,17,9)$, $(1,22,25)$, $(1,22,11)$, $(1,21,20)$, $(1,3,1)$, $(1,10,9)$, 
$(1,4,26)$, $(0,1,24)$, $(1,3,6)$, $(1,4,14)$, $(1,26,11)$, $(1,15,12)$, $(1,6,23)$, 
$(1,24,19)$, $(1,13,10)$, $(1,14,27)$, $(1,2,16)$, $(1,28,11)$, $(1,26,17)$, $(1,25,27)$, 
$(1,25,28)$, $(0,1,25)$, $(1,2,28)$, $(1,3,7)$, $(1,8,1)$, $(1,27,9)$, $(1,16,6)$, 
$(1,18,6)$, $(1,3,27)$, $(1,23,3)$, $(1,10,14)$, $(1,2,3)$, $(1,11,28)$, $(1,19,14)$, 
$(1,5,4)$, $(0,1,26)$, $(1,24,8)$, $(1,25,16)$, $(1,27,18)$, $(1,24,17)$, $(1,28,3)$, 
$(1,5,2)$, $(1,20,1)$, $(1,28,9)$, $(1,8,0)$, $(1,23,5)$, $(1,22,18)$, $(1,6,10)$, 
$(1,10,10)$, $(1,0,0)$, $(1,26,1)$, $(1,13,9)$, $(1,14,19)$, $(1,4,3)$, $(1,19,24)$, 
$(1,15,8)$, $(1,10,24)$, $(1,23,27)$, $(1,4,12)$, $(1,16,13)$, $(1,15,13)$, $(1,13,3)$, 
$(1,14,12)$, $(1,0,1)$, $(1,22,9)$, $(1,28,15)$, $(1,24,28)$, $(1,6,21)$, $(1,26,12)$, 
$(1,10,19)$, $(1,28,12)$, $(1,13,16)$, $(1,14,22)$, $(1,13,17)$, $(1,14,28)$, $(1,24,18)$, 
$(1,23,15)$, $(1,0,4)$, $(1,18,17)$, $(1,10,2)$, $(1,20,7)$, $(1,22,20)$, $(1,5,19)$, 
$(1,19,5)$, $(1,16,20)$, $(1,27,14)$, $(1,6,14)$, $(1,25,1)$, $(1,23,9)$, $(1,11,24)$, 
$(1,26,16)$, $(1,0,5)$, $(1,5,14)$, $(1,7,24)$, $(1,27,22)$, $(1,26,27)$, $(1,17,15)$, 
$(1,2,25)$, $(1,25,14)$, $(1,21,19)$, $(1,7,15)$, $(1,17,2)$, $(1,20,25)$, $(1,3,21)$, 
$(1,21,24)$, $(1,0,7)$, $(1,7,10)$, $(1,26,20)$, $(1,4,10)$, $(1,23,0)$, $(1,15,6)$, 
$(1,5,1)$, $(1,18,9)$, $(1,7,21)$, $(1,11,19)$, $(1,27,8)$, $(1,4,4)$, $(1,10,20)$, 
$(1,15,22)$, $(1,0,8)$, $(1,28,26)$, $(1,25,8)$, $(1,16,15)$, $(1,8,10)$, $(1,0,11)$, 
$(1,27,28)$, $(1,27,3)$, $(1,9,0)$, $(1,16,24)$, $(1,4,0)$, $(1,24,23)$, $(1,15,11)$, 
$(1,25,6)$, $(1,0,12)$, $(1,24,5)$, $(1,3,5)$, $(1,13,21)$, $(1,14,6)$, $(1,28,21)$, 
$(1,22,10)$, $(1,4,28)$, $(1,0,22)$, $(1,15,23)$, $(1,2,22)$, $(1,8,2)$, $(1,20,2)$, 
$(1,20,14)$, $(1,0,13)$, $(1,4,16)$, $(1,11,14)$, $(1,15,17)$, $(1,17,4)$, $(1,16,25)$, 
$(1,9,27)$, $(1,21,7)$, $(1,24,2)$, $(1,20,28)$, $(1,23,23)$, $(1,6,3)$, $(1,28,5)$, 
$(1,27,26)$, $(1,0,17)$, $(1,0,24)$, $(1,6,12)$, $(1,18,11)$, $(1,2,14)$, $(1,13,26)$, 
$(1,14,16)$, $(1,8,6)$, $(1,6,17)$, $(1,5,13)$, $(1,3,11)$, $(1,5,18)$, $(1,9,16)$, 
$(1,18,23)$, $(1,0,19)$, $(1,11,2)$, $(1,20,6)$, $(1,7,4)$, $(1,25,18)$, $(1,18,5)$, 
$(1,26,7)$, $(1,19,18)$, $(1,17,3)$, $(1,17,25)$, $(1,7,25)$, $(1,21,10)$, $(1,21,6)$, 
$(1,16,3)$, $(1,0,21)$, $(1,10,4)$, $(1,22,1)$, $(1,19,9)$, $(1,27,7)$, $(1,25,22)$, 
$(1,11,11)$, $(1,2,2)$, $(1,20,19)$, $(1,25,0)$, $(1,3,22)$, $(1,9,5)$, $(1,11,12)$, 
$(1,3,26)$, $(1,16,4)$, $(1,13,8)$, $(1,14,18)$, $(1,19,7)$, $(1,2,26)$, $(1,18,21)$, 
$(1,9,18)$, $(1,2,4)$, $(1,7,3)$, $(1,3,3)$, $(1,24,7)$, $(1,18,27)$, $(1,6,8)$, 
$(1,3,20)$, $(1,28,23)$, $(1,7,20)$, $(1,24,16)$, $(1,2,19)$, $(1,22,17)$, $(1,11,23)$, 
$(1,28,0)$, $(1,2,9)$, $(1,5,5)$, $(1,2,15)$, $(1,23,16)$, $(1,7,13)$, $(1,10,17)$, 
$(1,3,0)$, $(1,27,19)$, $(1,3,28)$, $(1,22,28)$, $(1,25,13)$, $(1,7,2)$, $(1,20,8)$, 
$(1,6,1)$, $(1,2,21)$, $(1,17,22)$, $(1,6,25)$, $(1,6,24)$, $(1,21,15)$, $(1,9,22)$, 
$(1,3,16)$, $(1,19,16)$, $(1,23,17)$, $(1,19,17)$, $(1,4,21)$, $(1,8,3)$, $(1,15,26)$, 
$(1,23,20)$, $(1,4,7)$, $(1,17,7)$, $(1,15,25)$, $(1,24,13)$, $(1,21,14)$, $(1,19,28)$, 
$(1,8,11)$, $(1,28,10)$, $(1,19,12)$, $(1,11,18)$, $(1,11,20)$, $(1,13,4)$, $(1,14,4)$, 
$(1,5,11)$, $(1,4,19)$, $(1,18,1)$, $(1,15,9)$, $(1,26,28)$, $(1,5,26)$, $(1,24,14)$, 
$(1,22,0)$, $(1,10,3)$, $(1,25,21)$, $(1,28,13)$, $(1,17,23)$, $(1,9,25)$, $(1,19,8)$, 
$(1,21,22)$, $(1,5,7)$, $(1,10,12)$, $(1,19,10)$, $(1,18,15)$, $(1,13,12)$, $(1,14,15)$, 
$(1,7,27)$, $(1,13,23)$, $(1,14,23)$, $(1,11,1)$, $(1,16,9)$, $(1,23,4)$, $(1,27,10)$, 
$(1,7,19)$, $(1,5,27)$, $(1,16,17)$, $(1,6,5)$, $(1,23,8)$, $(1,26,0)$, $(1,18,2)$, 
$(1,20,16)$, $(1,6,13)$, $(1,16,18)$, $(1,16,1)$, $(1,25,9)$, $(1,10,16)$, $(1,6,20)$, 
$(1,8,17)$, $(1,24,6)$, $(1,10,23)$, $(1,17,16)$, $(1,26,25)$, $(1,28,17)$, $(1,21,1)$, 
$(1,26,9)$, $(1,8,19)$, $(1,8,26)$, $(1,21,3)$, $(1,18,10)$, $(1,25,23)$, $(1,22,14)$, 
$(1,18,18)$, $(1,22,13)$, $(1,9,12)$, $(1,23,6)$, $(1,23,28)$, $(1,13,15)$, $(1,14,24)$, 
$(1,17,20)$, $(1,25,25)$.

\subsection*{$m_{21}(2,29)\ge 565$}
$(0,0,1)$, $(1,1,21)$, $(1,6,10)$, $(1,2,17)$, $(1,7,25)$, $(0,1,0)$, $(1,3,15)$, 
$(1,11,3)$, $(1,24,23)$, $(1,17,20)$, $(0,1,1)$, $(1,21,19)$, $(1,15,9)$, $(1,10,6)$, 
$(1,0,14)$, $(0,1,2)$, $(1,13,14)$, $(1,22,5)$, $(1,23,28)$, $(1,13,22)$, $(0,1,3)$, 
$(1,22,16)$, $(1,20,2)$, $(1,9,11)$, $(1,12,8)$, $(0,1,5)$, $(1,17,2)$, $(1,23,21)$, 
$(1,5,2)$, $(1,18,5)$, $(0,1,6)$, $(1,19,25)$, $(1,3,20)$, $(1,14,15)$, $(1,20,4)$, 
$(0,1,8)$, $(1,25,7)$, $(1,28,14)$, $(1,17,0)$, $(1,15,21)$, $(0,1,10)$, $(1,8,0)$, 
$(1,27,27)$, $(1,22,4)$, $(1,16,6)$, $(0,1,11)$, $(1,20,22)$, $(1,13,6)$, $(1,26,13)$, 
$(1,5,26)$, $(0,1,12)$, $(1,15,8)$, $(1,16,25)$, $(1,3,12)$, $(1,9,24)$, $(0,1,14)$, 
$(0,1,26)$, $(1,18,28)$, $(1,4,7)$, $(1,19,19)$, $(0,1,15)$, $(1,2,18)$, $(1,0,1)$, 
$(1,6,26)$, $(1,28,0)$, $(0,1,17)$, $(1,11,20)$, $(1,19,15)$, $(1,27,8)$, $(1,23,17)$, 
$(0,1,18)$, $(1,23,13)$, $(1,26,11)$, $(1,0,27)$, $(1,25,16)$, $(0,1,20)$, $(1,6,6)$, 
$(1,7,26)$, $(1,1,22)$, $(1,2,13)$, $(0,1,21)$, $(1,26,4)$, $(1,17,12)$, $(1,7,21)$, 
$(1,1,28)$, $(0,1,22)$, $(1,12,17)$, $(1,8,13)$, $(1,11,1)$, $(1,14,7)$, $(0,1,27)$, 
$(1,10,23)$, $(1,24,8)$, $(1,28,3)$, $(1,8,10)$, $(0,1,28)$, $(1,28,27)$, $(1,5,23)$, 
$(1,21,9)$, $(1,22,3)$, $(1,0,3)$, $(1,28,19)$, $(1,25,4)$, $(1,14,23)$, $(1,25,27)$, 
$(1,0,4)$, $(1,22,13)$, $(1,17,5)$, $(1,10,21)$, $(1,23,1)$, $(1,0,7)$, $(1,2,22)$, 
$(1,23,26)$, $(1,6,19)$, $(1,19,7)$, $(1,0,8)$, $(1,7,27)$, $(1,1,7)$, $(1,22,27)$, 
$(1,27,24)$, $(1,0,10)$, $(1,10,1)$, $(1,26,22)$, $(1,24,28)$, $(1,5,28)$, $(1,0,11)$, 
$(1,13,4)$, $(1,20,1)$, $(1,16,24)$, $(1,16,26)$, $(1,0,12)$, $(1,19,10)$, $(1,19,12)$, 
$(1,20,26)$, $(1,24,14)$, $(1,0,15)$, $(1,23,14)$, $(1,15,27)$, $(1,12,22)$, $(1,20,20)$, 
$(1,0,17)$, $(1,18,9)$, $(1,11,13)$, $(1,18,25)$, $(1,3,2)$, $(1,0,19)$, $(1,17,8)$, 
$(1,5,21)$, $(1,13,8)$, $(1,13,16)$, $(1,0,20)$, $(1,24,15)$, $(1,9,6)$, $(1,17,10)$, 
$(1,0,21)$, $(1,0,22)$, $(1,11,2)$, $(1,16,16)$, $(1,28,1)$, $(1,22,17)$, $(1,0,23)$, 
$(1,26,17)$, $(1,4,3)$, $(1,11,7)$, $(1,8,9)$, $(1,0,25)$, $(1,3,23)$, $(1,2,25)$, 
$(1,3,3)$, $(1,6,12)$, $(1,0,26)$, $(1,27,18)$, $(1,27,11)$, $(1,25,14)$, $(1,21,4)$, 
$(1,0,28)$, $(1,15,6)$, $(1,14,9)$, $(1,19,11)$, $(1,14,0)$, $(1,1,1)$, $(1,27,16)$, 
$(1,21,20)$, $(1,8,15)$, $(1,7,9)$, $(1,1,3)$, $(1,14,4)$, $(1,18,18)$, $(1,13,23)$, 
$(1,7,4)$, $(1,1,6)$, $(1,25,3)$, $(1,13,5)$, $(1,16,22)$, $(1,7,23)$, $(1,1,8)$, 
$(1,10,7)$, $(1,14,25)$, $(1,3,7)$, $(1,7,16)$, $(1,1,13)$, $(1,13,12)$, $(1,12,14)$, 
$(1,26,9)$, $(1,7,14)$, $(1,1,16)$, $(1,16,17)$, $(1,15,16)$, $(1,11,14)$, $(1,7,19)$, 
$(1,1,18)$, $(1,19,22)$, $(1,5,19)$, $(1,22,20)$, $(1,7,15)$, $(1,1,20)$, $(1,9,15)$, 
$(1,4,28)$, $(1,28,18)$, $(1,7,24)$, $(1,1,24)$, $(1,17,9)$, $(1,8,21)$, $(1,9,5)$, 
$(1,7,6)$, $(1,1,25)$, $(1,3,5)$, $(1,17,27)$, $(1,19,21)$, $(1,7,18)$, $(1,1,26)$, 
$(1,4,26)$, $(1,7,1)$, $(1,1,27)$, $(1,7,2)$, $(1,2,0)$, $(1,17,22)$, $(1,6,21)$, 
$(1,27,4)$, $(1,15,0)$, $(1,2,2)$, $(1,12,9)$, $(1,6,7)$, $(1,23,20)$, $(1,16,15)$, 
$(1,2,3)$, $(1,16,2)$, $(1,6,13)$, $(1,15,23)$, $(1,5,24)$, $(1,2,5)$, $(1,27,19)$, 
$(1,6,28)$, $(1,18,11)$, $(1,27,6)$, $(1,2,7)$, $(1,26,28)$, $(1,6,20)$, $(1,17,15)$, 
$(1,28,21)$, $(1,2,8)$, $(1,2,12)$, $(1,6,5)$, $(1,6,1)$, $(1,4,9)$, $(1,2,10)$, 
$(1,22,6)$, $(1,6,16)$, $(1,24,16)$, $(1,25,5)$, $(1,2,11)$, $(1,18,13)$, $(1,6,17)$, 
$(1,25,12)$, $(1,11,27)$, $(1,2,14)$, $(1,25,8)$, $(1,6,24)$, $(1,21,28)$, $(1,14,14)$, 
$(1,2,19)$, $(1,5,14)$, $(1,6,18)$, $(1,5,5)$, $(1,8,11)$, $(1,2,20)$, $(1,21,15)$, 
$(1,6,25)$, $(1,3,13)$, $(1,22,18)$, $(1,2,23)$, $(1,15,11)$, $(1,6,3)$, $(1,12,6)$, 
$(1,13,28)$, $(1,2,26)$, $(1,14,20)$, $(1,6,15)$, $(1,8,22)$, $(1,18,16)$, $(1,2,28)$, 
$(1,20,24)$, $(1,6,22)$, $(1,22,24)$, $(1,10,12)$, $(1,3,1)$, $(1,19,6)$, $(1,9,17)$, 
$(1,9,1)$, $(1,11,25)$, $(1,3,6)$, $(1,15,19)$, $(1,26,18)$, $(1,26,16)$, $(1,10,25)$, 
$(1,3,9)$, $(1,27,9)$, $(1,17,26)$, $(1,22,21)$, $(1,26,25)$, $(1,3,10)$, $(1,20,10)$, 
$(1,8,5)$, $(1,24,4)$, $(1,24,25)$, $(1,3,11)$, $(1,23,22)$, $(1,10,0)$, $(1,20,9)$, 
$(1,21,25)$, $(1,3,16)$, $(1,4,4)$, $(1,5,27)$, $(1,28,28)$, $(1,9,25)$, $(1,3,17)$, 
$(1,13,11)$, $(1,11,12)$, $(1,27,22)$, $(1,8,25)$, $(1,3,18)$, $(1,24,26)$, $(1,15,2)$, 
$(1,25,10)$, $(1,27,25)$, $(1,3,21)$, $(1,21,14)$, $(1,12,24)$, $(1,14,2)$, $(1,17,25)$, 
$(1,3,22)$, $(1,25,1)$, $(1,23,11)$, $(1,8,24)$, $(1,28,25)$, $(1,3,24)$, $(1,12,7)$, 
$(1,25,6)$, $(1,5,6)$, $(1,4,25)$, $(1,4,2)$, $(1,10,20)$, $(1,25,15)$, $(1,12,16)$, 
$(1,27,12)$, $(1,4,8)$, $(1,9,4)$, $(1,27,14)$, $(1,10,22)$, $(1,28,7)$, $(1,4,10)$, 
$(1,17,16)$, $(1,13,21)$, $(1,24,9)$, $(1,20,18)$, $(1,4,16)$, $(1,22,9)$, $(1,18,4)$, 
$(1,19,24)$, $(1,15,14)$, $(1,4,17)$, $(1,12,23)$, $(1,11,22)$, $(1,26,3)$, $(1,25,22)$, 
$(1,4,19)$, $(1,13,10)$, $(1,15,20)$, $(1,22,15)$, $(1,24,27)$, $(1,4,20)$, $(1,26,15)$, 
$(1,21,17)$, $(1,16,4)$, $(1,11,5)$, $(1,4,27)$, $(1,24,12)$, $(1,23,16)$, $(1,14,10)$, 
$(1,13,24)$, $(1,5,0)$, $(1,10,4)$, $(1,12,11)$, $(1,22,26)$, $(1,13,9)$, $(1,5,1)$, 
$(1,15,1)$, $(1,18,12)$, $(1,15,4)$, $(1,21,7)$, $(1,5,3)$, $(1,20,27)$, $(1,23,8)$, 
$(1,25,23)$, $(1,16,1)$, $(1,5,7)$, $(1,24,13)$, $(1,21,27)$, $(1,9,10)$, $(1,5,11)$, 
$(1,5,10)$, $(1,12,26)$, $(1,19,17)$, $(1,24,24)$, $(1,11,24)$, $(1,5,12)$, $(1,26,6)$, 
$(1,28,4)$, $(1,28,20)$, $(1,18,15)$, $(1,5,20)$, $(1,11,15)$, $(1,15,26)$, $(1,20,28)$, 
$(1,10,17)$, $(1,8,1)$, $(1,17,18)$, $(1,25,21)$, $(1,22,22)$, $(1,14,8)$, $(1,8,2)$, 
$(1,19,14)$, $(1,16,28)$, $(1,8,4)$, $(1,25,2)$, $(1,8,12)$, $(1,16,20)$, $(1,12,15)$, 
$(1,10,19)$, $(1,10,26)$, $(1,8,18)$, $(1,12,28)$, $(1,23,0)$, $(1,12,5)$, $(1,9,16)$, 
$(1,8,19)$, $(1,14,24)$, $(1,27,13)$, $(1,20,7)$, $(1,15,18)$, $(1,8,27)$, $(1,13,26)$, 
$(1,17,24)$, $(1,18,21)$, $(1,18,19)$, $(1,9,3)$, $(1,23,24)$, $(1,24,6)$, $(1,27,23)$, 
$(1,18,10)$, $(1,9,8)$, $(1,19,5)$, $(1,25,9)$, $(1,15,12)$, $(1,13,1)$, $(1,9,9)$, 
$(1,26,2)$, $(1,20,23)$, $(1,20,19)$, $(1,27,3)$, $(1,9,13)$, $(1,14,3)$, $(1,21,26)$, 
$(1,10,5)$, $(1,22,23)$, $(1,9,14)$, $(1,24,7)$, $(1,16,11)$, $(1,26,10)$, $(1,16,18)$, 
$(1,9,28)$, $(1,16,27)$, $(1,18,17)$, $(1,14,28)$, $(1,11,9)$, $(1,10,2)$, $(1,13,18)$, 
$(1,14,18)$, $(1,28,6)$, $(1,19,1)$, $(1,10,11)$, $(1,17,17)$, $(1,17,21)$, $(1,11,21)$, 
$(1,28,12)$, $(1,10,15)$, $(1,16,10)$, $(1,25,0)$, $(1,18,8)$, $(1,24,20)$, $(1,10,18)$, 
$(1,11,4)$, $(1,13,17)$, $(1,23,7)$, $(1,11,17)$, $(1,11,0)$, $(1,21,24)$, $(1,13,19)$, 
$(1,28,17)$, $(1,21,2)$, $(1,11,11)$, $(1,20,5)$, $(1,14,12)$, $(1,14,17)$, $(1,22,19)$, 
$(1,12,10)$, $(1,18,14)$, $(1,28,11)$, $(1,24,18)$, $(1,22,7)$, $(1,12,12)$, $(1,24,3)$, 
$(1,18,0)$, $(1,22,10)$, $(1,28,23)$, $(1,12,18)$, $(1,16,8)$, $(1,20,8)$, $(1,26,26)$, 
$(1,18,6)$, $(1,12,21)$, $(1,17,11)$, $(1,16,21)$, $(1,19,27)$, $(1,17,13)$, $(1,14,1)$, 
$(1,24,5)$, $(1,26,14)$, $(1,27,10)$, $(1,26,27)$, $(1,15,28)$, $(1,27,26)$, $(1,16,13)$, 
$(1,23,9)$, $(1,28,16)$, $(1,18,7)$, $(1,21,5)$, $(1,20,13)$, $(1,25,20)$, $(1,20,15)$, 
$(1,19,8)$, $(1,21,11)$, $(1,28,2)$, $(1,26,19)$, $(1,23,5)$.

\subsection*{$m_{22}(2,29)\ge 595$}
$(0,1,1)$, $(1,5,10)$, $(1,24,8)$, $(0,1,7)$, $(1,20,22)$, $(1,4,14)$, $(1,22,5)$, 
$(0,1,2)$, $(1,18,3)$, $(1,5,5)$, $(1,2,10)$, $(1,13,15)$, $(1,6,6)$, $(1,6,10)$, 
$(0,1,3)$, $(1,4,15)$, $(1,8,7)$, $(1,14,7)$, $(1,17,19)$, $(1,27,9)$, $(1,23,21)$, 
$(0,1,4)$, $(1,7,0)$, $(1,3,23)$, $(1,18,6)$, $(1,14,16)$, $(1,18,16)$, $(1,4,7)$, 
$(0,1,5)$, $(1,0,6)$, $(1,12,0)$, $(1,20,20)$, $(1,8,10)$, $(1,2,22)$, $(1,12,19)$, 
$(0,1,6)$, $(1,9,19)$, $(1,25,28)$, $(1,27,11)$, $(1,19,21)$, $(1,21,4)$, $(1,8,13)$, 
$(0,1,8)$, $(1,12,4)$, $(1,1,12)$, $(1,6,9)$, $(1,4,6)$, $(1,25,17)$, $(1,13,6)$, 
$(0,1,9)$, $(1,10,14)$, $(1,23,17)$, $(1,23,12)$, $(1,15,17)$, $(1,17,20)$, $(1,16,25)$, 
$(0,1,10)$, $(1,17,8)$, $(1,17,13)$, $(1,4,24)$, $(1,9,11)$, $(1,19,12)$, $(1,7,26)$, 
$(0,1,11)$, $(1,21,17)$, $(1,6,25)$, $(1,12,22)$, $(1,6,8)$, $(1,15,28)$, $(1,10,16)$, 
$(0,1,12)$, $(1,27,16)$, $(1,0,21)$, $(1,8,23)$, $(1,10,12)$, $(1,10,19)$, $(1,28,14)$, 
$(0,1,13)$, $(1,8,24)$, $(1,22,26)$, $(1,24,19)$, $(1,3,5)$, $(1,16,24)$, $(1,15,9)$, 
$(0,1,14)$, $(1,28,11)$, $(1,19,24)$, $(1,13,0)$, $(1,27,0)$, $(1,24,21)$, $(1,11,3)$, 
$(0,1,17)$, $(1,22,12)$, $(1,20,15)$, $(1,10,8)$, $(1,18,20)$, $(1,22,0)$, $(1,17,12)$, 
$(0,1,18)$, $(1,24,2)$, $(1,15,2)$, $(1,28,18)$, $(1,24,26)$, $(1,20,8)$, $(1,1,17)$, 
$(0,1,19)$, $(1,15,18)$, $(1,18,4)$, $(1,15,14)$, $(1,0,2)$, $(1,14,3)$, $(1,25,24)$, 
$(0,1,21)$, $(1,2,25)$, $(1,14,11)$, $(1,19,13)$, $(1,12,14)$, $(1,9,23)$, $(1,5,23)$, 
$(0,1,22)$, $(1,6,5)$, $(1,10,18)$, $(1,0,25)$, $(1,11,13)$, $(1,3,18)$, $(1,9,0)$, 
$(0,1,23)$, $(1,13,28)$, $(1,2,3)$, $(1,17,28)$, $(1,21,23)$, $(1,1,26)$, $(1,21,18)$, 
$(0,1,26)$, $(1,14,23)$, $(1,16,22)$, $(1,3,17)$, $(1,5,7)$, $(1,28,5)$, $(1,14,22)$, 
$(0,1,27)$, $(1,23,7)$, $(1,21,6)$, $(1,5,2)$, $(1,22,24)$, $(1,7,2)$, $(1,18,28)$, 
$(1,0,0)$, $(1,9,10)$, $(1,27,14)$, $(1,4,21)$, $(1,24,13)$, $(1,16,9)$, $(1,9,13)$, 
$(1,0,3)$, $(1,20,12)$, $(1,11,7)$, $(1,12,28)$, $(1,0,7)$, $(1,24,18)$, $(1,21,15)$, 
$(1,0,5)$, $(1,23,2)$, $(1,4,13)$, $(1,2,12)$, $(1,9,2)$, $(1,2,15)$, $(1,11,23)$, 
$(1,0,8)$, $(1,16,6)$, $(1,1,28)$, $(1,14,8)$, $(1,6,23)$, $(1,21,11)$, $(1,19,5)$, 
$(1,0,10)$, $(1,15,19)$, $(1,13,26)$, $(1,17,7)$, $(1,13,3)$, $(1,28,8)$, $(1,8,8)$, 
$(1,0,11)$, $(1,19,25)$, $(1,23,5)$, $(1,18,26)$, $(1,10,24)$, $(1,5,22)$, $(1,0,26)$, 
$(1,0,15)$, $(1,21,28)$, $(1,28,9)$, $(1,16,17)$, $(1,21,5)$, $(1,20,28)$, $(1,18,0)$, 
$(1,0,16)$, $(1,2,14)$, $(1,24,0)$, $(1,25,14)$, $(1,14,25)$, $(1,13,2)$, $(1,13,4)$, 
$(1,0,17)$, $(1,18,9)$, $(1,6,3)$, $(1,27,23)$, $(1,27,21)$, $(1,17,21)$, $(1,1,2)$, 
$(1,0,18)$, $(1,27,8)$, $(1,14,21)$, $(1,7,20)$, $(1,23,20)$, $(1,14,14)$, $(1,20,10)$, 
$(1,0,19)$, $(1,6,20)$, $(1,5,8)$, $(1,3,2)$, $(1,8,9)$, $(1,8,0)$, $(1,16,19)$, 
$(1,0,22)$, $(1,13,16)$, $(1,8,22)$, $(1,5,11)$, $(1,19,19)$, $(1,4,10)$, $(1,7,3)$, 
$(1,1,4)$, $(1,25,3)$, $(1,12,16)$, $(1,17,5)$, $(1,25,4)$, $(1,20,23)$, $(1,23,24)$, 
$(1,1,5)$, $(1,8,16)$, $(1,21,20)$, $(1,7,19)$, $(1,11,16)$, $(1,15,12)$, $(1,2,28)$, 
$(1,1,6)$, $(1,9,5)$, $(1,28,7)$, $(1,23,14)$, $(1,16,20)$, $(1,18,7)$, $(1,1,13)$, 
$(1,1,7)$, $(1,16,15)$, $(1,22,14)$, $(1,2,26)$, $(1,4,22)$, $(1,28,0)$, $(1,10,3)$, 
$(1,1,8)$, $(1,2,24)$, $(1,25,25)$, $(1,1,10)$, $(1,20,0)$, $(1,14,4)$, $(1,4,0)$, 
$(1,1,11)$, $(1,19,11)$, $(1,19,3)$, $(1,14,15)$, $(1,15,25)$, $(1,7,6)$, $(1,28,12)$, 
$(1,1,15)$, $(1,5,20)$, $(1,27,13)$, $(1,18,21)$, $(1,9,26)$, $(1,13,25)$, $(1,15,20)$, 
$(1,1,16)$, $(1,13,19)$, $(1,1,24)$, $(1,3,13)$, $(1,17,15)$, $(1,9,22)$, $(1,8,2)$, 
$(1,1,18)$, $(1,10,23)$, $(1,2,18)$, $(1,8,6)$, $(1,24,9)$, $(1,5,19)$, $(1,5,15)$, 
$(1,1,19)$, $(1,18,22)$, $(1,11,22)$, $(1,20,24)$, $(1,6,12)$, $(1,8,14)$, $(1,11,18)$, 
$(1,1,21)$, $(1,22,7)$, $(1,5,0)$, $(1,11,25)$, $(1,5,17)$, $(1,12,17)$, $(1,27,26)$, 
$(1,1,23)$, $(1,4,2)$, $(1,20,26)$, $(1,27,20)$, $(1,10,21)$, $(1,4,11)$, $(1,19,22)$, 
$(1,1,25)$, $(1,24,14)$, $(1,8,11)$, $(1,19,8)$, $(1,22,19)$, $(1,23,18)$, $(1,22,9)$, 
$(1,2,2)$, $(1,7,10)$, $(1,14,17)$, $(1,28,15)$, $(1,18,12)$, $(1,24,25)$, $(1,6,26)$, 
$(1,2,5)$, $(1,6,14)$, $(1,17,3)$, $(1,22,11)$, $(1,19,14)$, $(1,27,7)$, $(1,3,21)$, 
$(1,2,6)$, $(1,16,3)$, $(1,11,2)$, $(1,24,22)$, $(1,10,25)$, $(1,17,9)$, $(1,14,20)$, 
$(1,2,8)$, $(1,21,12)$, $(1,13,12)$, $(1,25,13)$, $(1,23,22)$, $(1,9,28)$, $(1,16,4)$, 
$(1,2,9)$, $(1,22,8)$, $(1,12,7)$, $(1,2,17)$, $(1,20,16)$, $(1,7,11)$, $(1,19,9)$, 
$(1,2,11)$, $(1,19,20)$, $(1,20,18)$, $(1,13,5)$, $(1,3,11)$, $(1,19,26)$, $(1,18,17)$, 
$(1,2,19)$, $(1,15,7)$, $(1,4,25)$, $(1,18,18)$, $(1,17,10)$, $(1,8,5)$, $(1,17,25)$, 
$(1,2,20)$, $(1,3,26)$, $(1,22,28)$, $(1,4,28)$, $(1,15,6)$, $(1,21,14)$, $(1,21,22)$, 
$(1,3,0)$, $(1,4,16)$, $(1,10,15)$, $(1,20,13)$, $(1,15,3)$, $(1,9,3)$, $(1,10,17)$, 
$(1,3,4)$, $(1,13,22)$, $(1,22,18)$, $(1,25,20)$, $(1,25,2)$, $(1,10,7)$, $(1,6,4)$, 
$(1,3,9)$, $(1,15,4)$, $(1,23,11)$, $(1,19,0)$, $(1,28,22)$, $(1,25,9)$, $(1,18,14)$, 
$(1,3,10)$, $(1,11,11)$, $(1,19,10)$, $(1,3,24)$, $(1,20,17)$, $(1,16,2)$, $(1,23,23)$, 
$(1,3,16)$, $(1,3,25)$, $(1,27,12)$, $(1,5,21)$, $(1,13,9)$, $(1,24,5)$, $(1,4,12)$, 
$(1,3,19)$, $(1,21,8)$, $(1,24,4)$, $(1,24,7)$, $(1,18,23)$, $(1,8,28)$, $(1,25,15)$, 
$(1,3,22)$, $(1,14,13)$, $(1,11,8)$, $(1,9,15)$, $(1,17,26)$, $(1,12,15)$, $(1,12,9)$, 
$(1,4,3)$, $(1,4,20)$, $(1,6,28)$, $(1,24,17)$, $(1,14,24)$, $(1,21,7)$, $(1,11,24)$, 
$(1,4,4)$, $(1,16,10)$, $(1,4,26)$, $(1,23,9)$, $(1,25,18)$, $(1,6,18)$, $(1,7,12)$, 
$(1,4,9)$, $(1,12,23)$, $(1,11,4)$, $(1,28,20)$, $(1,15,5)$, $(1,15,23)$, $(1,25,8)$, 
$(1,4,18)$, $(1,5,24)$, $(1,23,16)$, $(1,12,8)$, $(1,13,14)$, $(1,18,15)$, $(1,28,17)$, 
$(1,5,3)$, $(1,16,14)$, $(1,25,12)$, $(1,16,5)$, $(1,24,3)$, $(1,7,4)$, $(1,12,26)$, 
$(1,5,9)$, $(1,20,11)$, $(1,19,28)$, $(1,17,2)$, $(1,5,12)$, $(1,21,3)$, $(1,15,0)$, 
$(1,5,13)$, $(1,6,7)$, $(1,9,16)$, $(1,23,13)$, $(1,10,2)$, $(1,9,8)$, $(1,20,5)$, 
$(1,5,14)$, $(1,7,28)$, $(1,6,24)$, $(1,11,20)$, $(1,21,9)$, $(1,27,15)$, $(1,7,21)$, 
$(1,5,26)$, $(1,25,0)$, $(1,21,13)$, $(1,22,16)$, $(1,9,4)$, $(1,14,18)$, $(1,14,28)$, 
$(1,5,28)$, $(1,8,20)$, $(1,28,4)$, $(1,21,19)$, $(1,17,17)$, $(1,10,10)$, $(1,25,10)$, 
$(1,6,15)$, $(1,13,24)$, $(1,16,28)$, $(1,20,25)$, $(1,22,25)$, $(1,21,16)$, $(1,16,23)$, 
$(1,6,19)$, $(1,22,13)$, $(1,28,23)$, $(1,9,20)$, $(1,12,3)$, $(1,8,18)$, $(1,23,3)$, 
$(1,6,21)$, $(1,28,25)$, $(1,20,7)$, $(1,25,22)$, $(1,24,12)$, $(1,28,6)$, $(1,27,4)$, 
$(1,7,13)$, $(1,9,25)$, $(1,9,24)$, $(1,24,6)$, $(1,22,22)$, $(1,15,10)$, $(1,12,13)$, 
$(1,7,14)$, $(1,10,6)$, $(1,8,26)$, $(1,7,23)$, $(1,20,19)$, $(1,24,23)$, $(1,13,17)$, 
$(1,7,15)$, $(1,27,2)$, $(1,25,21)$, $(1,15,15)$, $(1,14,10)$, $(1,22,4)$, $(1,10,5)$, 
$(1,8,4)$, $(1,18,2)$, $(1,27,19)$, $(1,12,21)$, $(1,25,16)$, $(1,11,26)$, $(1,16,26)$, 
$(1,8,17)$, $(1,15,13)$, $(1,8,19)$, $(1,9,6)$, $(1,18,19)$, $(1,16,12)$, $(1,14,0)$, 
$(1,9,7)$, $(1,15,11)$, $(1,19,17)$, $(1,11,28)$, $(1,11,14)$, $(1,13,13)$, $(1,21,10)$, 
$(1,9,9)$, $(1,11,10)$, $(1,18,25)$, $(1,25,23)$, $(1,22,3)$, $(1,18,13)$, $(1,20,4)$, 
$(1,10,0)$, $(1,22,6)$, $(1,13,8)$, $(1,23,26)$, $(1,11,5)$, $(1,27,6)$, $(1,11,21)$, 
$(1,11,6)$, $(1,17,11)$, $(1,19,18)$, $(1,16,16)$, $(1,28,10)$, $(1,28,16)$, $(1,22,2)$, 
$(1,12,6)$, $(1,15,24)$, $(1,27,5)$, $(1,13,21)$, $(1,14,12)$, $(1,27,18)$, $(1,20,21)$, 
$(1,17,23)$, $(1,24,16)$, $(1,27,28)$, $(1,22,21)$, $(1,27,22)$, $(1,21,21)$, $(1,23,10)$.

\subsection*{$m_{23}(2,29)\ge 628$}
$(0,1,0)$, $(1,0,12)$, $(1,13,26)$, $(1,17,2)$, $(0,1,3)$, $(1,21,28)$, $(1,9,20)$, 
$(1,14,15)$, $(0,1,5)$, $(1,13,15)$, $(1,6,1)$, $(1,1,23)$, $(0,1,6)$, $(1,27,16)$, 
$(1,20,22)$, $(1,6,11)$, $(0,1,7)$, $(1,1,10)$, $(1,21,9)$, $(1,26,21)$, $(0,1,8)$, 
$(1,9,23)$, $(1,15,0)$, $(1,13,0)$, $(0,1,9)$, $(1,10,21)$, $(1,28,5)$, $(1,5,25)$, 
$(0,1,10)$, $(1,23,24)$, $(1,26,2)$, $(1,8,12)$, $(0,1,11)$, $(1,28,14)$, $(1,0,21)$, 
$(1,15,1)$, $(0,1,13)$, $(1,24,22)$, $(1,12,10)$, $(1,3,24)$, $(0,1,14)$, $(1,16,9)$, 
$(1,27,18)$, $(1,9,27)$, $(0,1,15)$, $(0,1,27)$, $(1,14,13)$, $(1,16,16)$, $(0,1,17)$, 
$(1,11,19)$, $(1,8,4)$, $(1,11,28)$, $(0,1,18)$, $(1,18,5)$, $(1,24,28)$, $(1,27,7)$, 
$(0,1,19)$, $(1,7,27)$, $(1,10,7)$, $(1,18,17)$, $(0,1,20)$, $(1,12,17)$, $(1,3,11)$, 
$(1,23,5)$, $(0,1,21)$, $(1,25,20)$, $(1,22,25)$, $(1,20,18)$, $(0,1,23)$, $(1,5,2)$, 
$(1,23,12)$, $(1,24,20)$, $(0,1,24)$, $(1,8,25)$, $(1,7,17)$, $(1,28,22)$, $(0,1,26)$, 
$(1,15,11)$, $(1,2,24)$, $(1,2,9)$, $(0,1,28)$, $(1,2,8)$, $(1,18,19)$, $(1,22,19)$, 
$(1,0,0)$, $(1,27,25)$, $(1,26,14)$, $(1,0,14)$, $(1,0,1)$, $(1,18,7)$, $(1,5,15)$, 
$(1,20,14)$, $(1,0,2)$, $(1,6,12)$, $(1,23,10)$, $(1,15,14)$, $(1,0,3)$, $(1,24,19)$, 
$(1,14,27)$, $(1,25,14)$, $(1,0,4)$, $(1,22,15)$, $(1,21,17)$, $(1,23,14)$, $(1,0,5)$, 
$(1,9,18)$, $(1,18,13)$, $(1,8,14)$, $(1,0,7)$, $(1,21,13)$, $(1,19,24)$, $(1,18,14)$, 
$(1,0,10)$, $(1,23,17)$, $(1,8,19)$, $(1,27,14)$, $(1,0,11)$, $(1,26,23)$, $(1,28,7)$, 
$(1,13,14)$, $(1,0,13)$, $(1,11,22)$, $(1,25,3)$, $(1,11,14)$, $(1,0,15)$, $(1,17,5)$, 
$(1,20,6)$, $(1,2,14)$, $(1,0,16)$, $(1,8,16)$, $(1,16,20)$, $(1,19,14)$, $(1,0,17)$, 
$(1,1,2)$, $(1,10,12)$, $(1,6,14)$, $(1,0,18)$, $(1,7,14)$, $(1,0,20)$, $(1,16,3)$, 
$(1,7,8)$, $(1,1,14)$, $(1,0,22)$, $(1,10,20)$, $(1,15,9)$, $(1,9,14)$, $(1,0,26)$, 
$(1,2,4)$, $(1,6,26)$, $(1,14,14)$, $(1,0,27)$, $(1,28,27)$, $(1,2,11)$, $(1,10,14)$, 
$(1,0,28)$, $(1,5,10)$, $(1,24,21)$, $(1,21,14)$, $(1,1,1)$, $(1,11,17)$, $(1,15,7)$, 
$(1,11,7)$, $(1,1,3)$, $(1,6,21)$, $(1,5,18)$, $(1,19,5)$, $(1,1,4)$, $(1,26,5)$, 
$(1,18,24)$, $(1,25,18)$, $(1,1,5)$, $(1,16,13)$, $(1,1,6)$, $(1,27,10)$, $(1,25,25)$, 
$(1,17,20)$, $(1,1,7)$, $(1,14,3)$, $(1,10,27)$, $(1,9,22)$, $(1,1,11)$, $(1,8,2)$, 
$(1,21,12)$, $(1,21,19)$, $(1,1,13)$, $(1,9,7)$, $(1,17,28)$, $(1,15,6)$, $(1,1,15)$, 
$(1,25,0)$, $(1,7,10)$, $(1,8,15)$, $(1,1,16)$, $(1,24,24)$, $(1,22,8)$, $(1,28,10)$, 
$(1,1,17)$, $(1,13,27)$, $(1,24,0)$, $(1,5,23)$, $(1,1,18)$, $(1,15,8)$, $(1,14,11)$, 
$(1,24,11)$, $(1,1,19)$, $(1,18,23)$, $(1,19,20)$, $(1,13,21)$, $(1,1,22)$, $(1,5,16)$, 
$(1,9,2)$, $(1,20,12)$, $(1,1,24)$, $(1,1,25)$, $(1,2,1)$, $(1,23,4)$, $(1,1,26)$, 
$(1,12,22)$, $(1,3,26)$, $(1,18,27)$, $(1,2,0)$, $(1,24,1)$, $(1,28,1)$, $(1,6,20)$, 
$(1,2,3)$, $(1,20,13)$, $(1,18,2)$, $(1,14,12)$, $(1,2,5)$, $(1,14,2)$, $(1,5,12)$, 
$(1,25,1)$, $(1,2,7)$, $(1,25,27)$, $(1,6,9)$, $(1,22,4)$, $(1,2,10)$, $(1,17,22)$, 
$(1,7,6)$, $(1,7,19)$, $(1,2,12)$, $(1,27,21)$, $(1,16,8)$, $(1,24,2)$, $(1,2,13)$, 
$(1,21,10)$, $(1,27,4)$, $(1,10,16)$, $(1,2,15)$, $(1,15,28)$, $(1,8,3)$, $(1,23,3)$, 
$(1,2,17)$, $(1,12,8)$, $(1,3,18)$, $(1,27,28)$, $(1,2,18)$, $(1,11,11)$, $(1,21,22)$, 
$(1,11,15)$, $(1,2,19)$, $(1,16,25)$, $(1,23,16)$, $(1,19,7)$, $(1,2,22)$, $(1,13,5)$, 
$(1,12,20)$, $(1,3,23)$, $(1,2,25)$, $(1,22,7)$, $(1,9,0)$, $(1,8,18)$, $(1,2,26)$, 
$(1,26,24)$, $(1,10,26)$, $(1,5,21)$, $(1,2,27)$, $(1,7,23)$, $(1,20,25)$, $(1,13,13)$, 
$(1,2,28)$, $(1,19,16)$, $(1,25,10)$, $(1,26,0)$, $(1,3,0)$, $(1,9,26)$, $(1,22,11)$, 
$(1,12,26)$, $(1,3,2)$, $(1,28,12)$, $(1,18,18)$, $(1,12,3)$, $(1,3,4)$, $(1,16,1)$, 
$(1,20,0)$, $(1,12,21)$, $(1,3,5)$, $(1,6,16)$, $(1,7,1)$, $(1,12,28)$, $(1,3,7)$, 
$(1,20,24)$, $(1,5,19)$, $(1,12,0)$, $(1,3,8)$, $(1,12,7)$, $(1,3,10)$, $(1,5,3)$, 
$(1,16,7)$, $(1,12,25)$, $(1,3,13)$, $(1,24,18)$, $(1,23,2)$, $(1,12,12)$, $(1,3,15)$, 
$(1,19,11)$, $(1,17,27)$, $(1,12,6)$, $(1,3,16)$, $(1,21,8)$, $(1,10,3)$, $(1,12,1)$, 
$(1,3,17)$, $(1,7,0)$, $(1,15,16)$, $(1,12,18)$, $(1,3,21)$, $(1,22,21)$, $(1,25,13)$, 
$(1,12,15)$, $(1,3,27)$, $(1,13,20)$, $(1,27,24)$, $(1,12,23)$, $(1,3,28)$, $(1,25,2)$, 
$(1,9,12)$, $(1,12,16)$, $(1,5,1)$, $(1,27,19)$, $(1,23,18)$, $(1,21,18)$, $(1,5,5)$, 
$(1,26,10)$, $(1,9,19)$, $(1,9,3)$, $(1,5,6)$, $(1,15,27)$, $(1,16,4)$, $(1,17,13)$, 
$(1,5,7)$, $(1,6,4)$, $(1,25,22)$, $(1,16,19)$, $(1,5,8)$, $(1,13,9)$, $(1,8,17)$, 
$(1,25,23)$, $(1,5,9)$, $(1,10,11)$, $(1,19,10)$, $(1,19,1)$, $(1,5,11)$, $(1,11,20)$, 
$(1,5,13)$, $(1,22,3)$, $(1,22,16)$, $(1,27,11)$, $(1,5,17)$, $(1,17,16)$, $(1,17,6)$, 
$(1,23,6)$, $(1,5,20)$, $(1,8,22)$, $(1,14,0)$, $(1,26,17)$, $(1,5,22)$, $(1,21,23)$, 
$(1,6,13)$, $(1,13,8)$, $(1,5,26)$, $(1,7,13)$, $(1,27,26)$, $(1,8,9)$, $(1,6,5)$, 
$(1,19,15)$, $(1,22,20)$, $(1,18,20)$, $(1,6,6)$, $(1,28,11)$, $(1,14,9)$, $(1,17,19)$, 
$(1,6,7)$, $(1,27,5)$, $(1,28,21)$, $(1,8,10)$, $(1,6,15)$, $(1,8,7)$, $(1,10,18)$, 
$(1,14,16)$, $(1,6,17)$, $(1,16,26)$, $(1,21,15)$, $(1,24,26)$, $(1,6,19)$, $(1,20,21)$, 
$(1,7,3)$, $(1,9,11)$, $(1,6,23)$, $(1,24,16)$, $(1,13,4)$, $(1,7,9)$, $(1,6,28)$, 
$(1,11,25)$, $(1,18,0)$, $(1,11,13)$, $(1,7,2)$, $(1,16,12)$, $(1,16,27)$, $(1,22,1)$, 
$(1,7,5)$, $(1,22,17)$, $(1,10,28)$, $(1,26,19)$, $(1,7,11)$, $(1,18,4)$, $(1,9,4)$, 
$(1,24,10)$, $(1,7,12)$, $(1,17,8)$, $(1,26,6)$, $(1,19,2)$, $(1,7,16)$, $(1,11,3)$, 
$(1,28,25)$, $(1,11,24)$, $(1,7,18)$, $(1,13,24)$, $(1,23,21)$, $(1,14,23)$, $(1,7,20)$, 
$(1,26,1)$, $(1,15,3)$, $(1,27,9)$, $(1,7,21)$, $(1,19,0)$, $(1,14,8)$, $(1,15,13)$, 
$(1,7,24)$, $(1,14,20)$, $(1,21,2)$, $(1,18,12)$, $(1,7,28)$, $(1,21,21)$, $(1,20,7)$, 
$(1,23,20)$, $(1,8,0)$, $(1,18,11)$, $(1,16,21)$, $(1,10,24)$, $(1,8,1)$, $(1,22,22)$, 
$(1,8,23)$, $(1,13,19)$, $(1,13,28)$, $(1,16,23)$, $(1,8,24)$, $(1,17,1)$, $(1,10,6)$, 
$(1,20,3)$, $(1,8,28)$, $(1,28,24)$, $(1,15,4)$, $(1,19,8)$, $(1,9,5)$, $(1,21,26)$, 
$(1,28,16)$, $(1,28,26)$, $(1,9,6)$, $(1,26,9)$, $(1,23,1)$, $(1,17,25)$, $(1,9,9)$, 
$(1,25,24)$, $(1,21,24)$, $(1,16,17)$, $(1,9,10)$, $(1,28,8)$, $(1,22,27)$, $(1,23,15)$, 
$(1,9,16)$, $(1,10,17)$, $(1,24,4)$, $(1,20,20)$, $(1,9,21)$, $(1,11,2)$, $(1,17,12)$, 
$(1,11,6)$, $(1,9,24)$, $(1,19,27)$, $(1,27,13)$, $(1,14,1)$, $(1,9,25)$, $(1,27,23)$, 
$(1,18,15)$, $(1,13,22)$, $(1,10,1)$, $(1,25,26)$, $(1,25,28)$, $(1,23,26)$, $(1,10,2)$, 
$(1,27,12)$, $(1,22,0)$, $(1,19,19)$, $(1,10,8)$, $(1,21,25)$, $(1,14,22)$, $(1,20,28)$, 
$(1,10,9)$, $(1,11,8)$, $(1,27,8)$, $(1,11,5)$, $(1,10,15)$, $(1,16,2)$, $(1,15,12)$, 
$(1,10,25)$, $(1,10,22)$, $(1,22,18)$, $(1,17,21)$, $(1,24,6)$, $(1,11,0)$, $(1,16,24)$, 
$(1,11,10)$, $(1,20,5)$, $(1,11,4)$, $(1,23,27)$, $(1,11,21)$, $(1,26,20)$, $(1,11,18)$, 
$(1,22,10)$, $(1,11,27)$, $(1,19,17)$, $(1,13,6)$, $(1,25,8)$, $(1,20,23)$, $(1,17,10)$, 
$(1,13,10)$, $(1,22,13)$, $(1,15,23)$, $(1,21,1)$, $(1,13,11)$, $(1,15,15)$, $(1,23,23)$, 
$(1,16,5)$, $(1,13,12)$, $(1,14,7)$, $(1,28,23)$, $(1,27,2)$, $(1,13,18)$, $(1,20,26)$, 
$(1,19,23)$, $(1,26,26)$, $(1,13,25)$, $(1,28,3)$, $(1,17,23)$, $(1,22,6)$, $(1,14,5)$, 
$(1,23,8)$, $(1,23,0)$, $(1,20,27)$, $(1,14,10)$, $(1,14,17)$, $(1,19,25)$, $(1,16,11)$, 
$(1,25,17)$, $(1,15,18)$, $(1,28,20)$, $(1,20,8)$, $(1,25,4)$, $(1,15,19)$, $(1,25,16)$, 
$(1,18,21)$, $(1,27,27)$, $(1,15,21)$, $(1,17,15)$, $(1,27,6)$, $(1,24,7)$, $(1,15,24)$, 
$(1,24,5)$, $(1,25,19)$, $(1,19,22)$, $(1,15,26)$, $(1,20,19)$, $(1,16,15)$, $(1,26,3)$, 
$(1,26,22)$, $(1,24,8)$, $(1,16,28)$, $(1,18,28)$, $(1,17,11)$, $(1,28,6)$, $(1,17,0)$, 
$(1,21,6)$, $(1,17,17)$, $(1,18,6)$, $(1,18,1)$, $(1,26,8)$, $(1,18,3)$, $(1,21,0)$, 
$(1,25,9)$, $(1,24,27)$, $(1,18,16)$, $(1,23,9)$, $(1,19,4)$, $(1,28,0)$, $(1,23,22)$, 
$(1,27,22)$, $(1,19,13)$, $(1,20,17)$, $(1,20,15)$, $(1,28,9)$, $(1,20,9)$, $(1,21,7)$, 
$(1,22,2)$, $(1,26,12)$, $(1,24,25)$, $(1,25,12)$, $(1,28,2)$.

\subsection*{$m_{25}(2,29)\ge 695$}
$(0,0,1)$, $(1,1,21)$, $(1,6,10)$, $(1,2,17)$, $(1,7,25)$, $(0,1,0)$, $(1,3,15)$, 
$(1,11,3)$, $(1,24,23)$, $(1,17,20)$, $(0,1,1)$, $(1,21,19)$, $(1,15,9)$, $(1,10,6)$, 
$(1,0,14)$, $(0,1,3)$, $(1,22,16)$, $(1,20,2)$, $(1,9,11)$, $(1,12,8)$, $(0,1,4)$, 
$(1,27,1)$, $(1,12,19)$, $(1,12,25)$, $(1,3,27)$, $(0,1,5)$, $(1,17,2)$, $(1,23,21)$, 
$(1,5,2)$, $(1,18,5)$, $(0,1,7)$, $(1,5,9)$, $(1,25,24)$, $(1,25,18)$, $(1,10,9)$, 
$(0,1,8)$, $(1,25,7)$, $(1,28,14)$, $(1,17,0)$, $(1,15,21)$, $(0,1,10)$, $(1,8,0)$, 
$(1,27,27)$, $(1,22,4)$, $(1,16,6)$, $(0,1,12)$, $(1,15,8)$, $(1,16,25)$, $(1,3,12)$, 
$(1,9,24)$, $(0,1,13)$, $(1,0,24)$, $(1,9,0)$, $(1,13,20)$, $(1,27,15)$, $(0,1,14)$, 
$(0,1,26)$, $(1,18,28)$, $(1,4,7)$, $(1,19,19)$, $(0,1,15)$, $(1,2,18)$, $(1,0,1)$, 
$(1,6,26)$, $(1,28,0)$, $(0,1,17)$, $(1,11,20)$, $(1,19,15)$, $(1,27,8)$, $(1,23,17)$, 
$(0,1,18)$, $(1,23,13)$, $(1,26,11)$, $(1,0,27)$, $(1,25,16)$, $(0,1,19)$, $(1,7,3)$, 
$(1,1,17)$, $(1,20,14)$, $(1,11,23)$, $(0,1,21)$, $(1,26,4)$, $(1,17,12)$, $(1,7,21)$, 
$(1,1,28)$, $(0,1,22)$, $(1,12,17)$, $(1,8,13)$, $(1,11,1)$, $(1,14,7)$, $(0,1,23)$, 
$(1,14,11)$, $(1,10,16)$, $(1,18,24)$, $(1,26,1)$, $(0,1,25)$, $(1,9,26)$, $(1,2,4)$, 
$(1,8,16)$, $(1,6,11)$, $(0,1,27)$, $(1,10,23)$, $(1,24,8)$, $(1,28,3)$, $(1,8,10)$, 
$(0,1,28)$, $(1,28,27)$, $(1,5,23)$, $(1,21,9)$, $(1,22,3)$, $(1,0,0)$, $(1,16,7)$, 
$(1,10,24)$, $(1,7,5)$, $(1,1,5)$, $(1,0,4)$, $(1,22,13)$, $(1,17,5)$, $(1,10,21)$, 
$(1,23,1)$, $(1,0,7)$, $(1,2,22)$, $(1,23,26)$, $(1,6,19)$, $(1,19,7)$, $(1,0,8)$, 
$(1,7,27)$, $(1,1,7)$, $(1,22,27)$, $(1,27,24)$, $(1,0,9)$, $(1,14,5)$, $(1,12,2)$, 
$(1,5,4)$, $(1,9,22)$, $(1,0,10)$, $(1,10,1)$, $(1,26,22)$, $(1,24,28)$, $(1,5,28)$, 
$(1,0,11)$, $(1,13,4)$, $(1,20,1)$, $(1,16,24)$, $(1,16,26)$, $(1,0,12)$, $(1,19,10)$, 
$(1,19,12)$, $(1,20,26)$, $(1,24,14)$, $(1,0,13)$, $(1,8,28)$, $(1,22,8)$, $(1,4,18)$, 
$(1,15,13)$, $(1,0,15)$, $(1,23,14)$, $(1,15,27)$, $(1,12,22)$, $(1,20,20)$, $(1,0,16)$, 
$(1,12,3)$, $(1,0,18)$, $(1,21,12)$, $(1,18,23)$, $(1,0,17)$, $(1,18,9)$, $(1,11,13)$, 
$(1,18,25)$, $(1,3,2)$, $(1,0,20)$, $(1,24,15)$, $(1,9,6)$, $(1,17,10)$, $(1,0,21)$, 
$(1,0,22)$, $(1,11,2)$, $(1,16,16)$, $(1,28,1)$, $(1,22,17)$, $(1,0,23)$, $(1,26,17)$, 
$(1,4,3)$, $(1,11,7)$, $(1,8,9)$, $(1,0,26)$, $(1,27,18)$, $(1,27,11)$, $(1,25,14)$, 
$(1,21,4)$, $(1,0,28)$, $(1,15,6)$, $(1,14,9)$, $(1,19,11)$, $(1,14,0)$, $(1,1,0)$, 
$(1,24,11)$, $(1,19,9)$, $(1,23,10)$, $(1,7,7)$, $(1,1,3)$, $(1,14,4)$, $(1,18,18)$, 
$(1,13,23)$, $(1,7,4)$, $(1,1,4)$, $(1,23,19)$, $(1,9,12)$, $(1,25,19)$, $(1,7,13)$, 
$(1,1,6)$, $(1,25,3)$, $(1,13,5)$, $(1,16,22)$, $(1,7,23)$, $(1,1,8)$, $(1,10,7)$, 
$(1,14,25)$, $(1,3,7)$, $(1,7,16)$, $(1,1,9)$, $(1,12,20)$, $(1,28,15)$, $(1,14,13)$, 
$(1,7,8)$, $(1,1,10)$, $(1,11,28)$, $(1,25,13)$, $(1,24,0)$, $(1,7,17)$, $(1,1,11)$, 
$(1,15,25)$, $(1,3,8)$, $(1,18,2)$, $(1,7,22)$, $(1,1,12)$, $(1,5,18)$, $(1,23,2)$, 
$(1,27,28)$, $(1,7,12)$, $(1,1,13)$, $(1,13,12)$, $(1,12,14)$, $(1,26,9)$, $(1,7,14)$, 
$(1,1,15)$, $(1,26,24)$, $(1,22,11)$, $(1,21,1)$, $(1,7,20)$, $(1,1,16)$, $(1,16,17)$, 
$(1,15,16)$, $(1,11,14)$, $(1,7,19)$, $(1,1,18)$, $(1,19,22)$, $(1,5,19)$, $(1,22,20)$, 
$(1,7,15)$, $(1,1,20)$, $(1,9,15)$, $(1,4,28)$, $(1,28,18)$, $(1,7,24)$, $(1,1,23)$, 
$(1,8,23)$, $(1,10,3)$, $(1,15,3)$, $(1,7,11)$, $(1,1,24)$, $(1,17,9)$, $(1,8,21)$, 
$(1,9,5)$, $(1,7,6)$, $(1,1,25)$, $(1,3,5)$, $(1,17,27)$, $(1,19,21)$, $(1,7,18)$, 
$(1,1,26)$, $(1,4,26)$, $(1,7,1)$, $(1,1,27)$, $(1,7,2)$, $(1,2,0)$, $(1,17,22)$, 
$(1,6,21)$, $(1,27,4)$, $(1,15,0)$, $(1,2,1)$, $(1,28,10)$, $(1,6,9)$, $(1,16,19)$, 
$(1,24,19)$, $(1,2,2)$, $(1,12,9)$, $(1,6,7)$, $(1,23,20)$, $(1,16,15)$, $(1,2,3)$, 
$(1,16,2)$, $(1,6,13)$, $(1,15,23)$, $(1,5,24)$, $(1,2,5)$, $(1,27,19)$, $(1,6,28)$, 
$(1,18,11)$, $(1,27,6)$, $(1,2,7)$, $(1,26,28)$, $(1,6,20)$, $(1,17,15)$, $(1,28,21)$, 
$(1,2,9)$, $(1,24,17)$, $(1,6,8)$, $(1,11,10)$, $(1,23,4)$, $(1,2,10)$, $(1,22,6)$, 
$(1,6,16)$, $(1,24,16)$, $(1,25,5)$, $(1,2,11)$, $(1,18,13)$, $(1,6,17)$, $(1,25,12)$, 
$(1,11,27)$, $(1,2,14)$, $(1,25,8)$, $(1,6,24)$, $(1,21,28)$, $(1,14,14)$, $(1,2,15)$, 
$(1,13,0)$, $(1,6,2)$, $(1,14,27)$, $(1,26,20)$, $(1,2,16)$, $(1,9,7)$, $(1,6,23)$, 
$(1,13,2)$, $(1,21,3)$, $(1,2,19)$, $(1,5,14)$, $(1,6,18)$, $(1,5,5)$, $(1,8,11)$, 
$(1,2,20)$, $(1,21,15)$, $(1,6,25)$, $(1,3,13)$, $(1,22,18)$, $(1,2,21)$, $(1,4,23)$, 
$(1,6,14)$, $(1,20,3)$, $(1,17,1)$, $(1,2,24)$, $(1,19,4)$, $(1,6,4)$, $(1,26,8)$, 
$(1,12,13)$, $(1,2,26)$, $(1,14,20)$, $(1,6,15)$, $(1,8,22)$, $(1,18,16)$, $(1,2,27)$, 
$(1,11,18)$, $(1,6,0)$, $(1,9,18)$, $(1,20,17)$, $(1,2,28)$, $(1,20,24)$, $(1,6,22)$, 
$(1,22,24)$, $(1,10,12)$, $(1,3,0)$, $(1,28,13)$, $(1,16,14)$, $(1,13,25)$, $(1,3,25)$, 
$(1,3,1)$, $(1,19,6)$, $(1,9,17)$, $(1,9,1)$, $(1,11,25)$, $(1,3,6)$, $(1,15,19)$, 
$(1,26,18)$, $(1,26,16)$, $(1,10,25)$, $(1,3,9)$, $(1,27,9)$, $(1,17,26)$, $(1,22,21)$, 
$(1,26,25)$, $(1,3,10)$, $(1,20,10)$, $(1,8,5)$, $(1,24,4)$, $(1,24,25)$, $(1,3,11)$, 
$(1,23,22)$, $(1,10,0)$, $(1,20,9)$, $(1,21,25)$, $(1,3,16)$, $(1,4,4)$, $(1,5,27)$, 
$(1,28,28)$, $(1,9,25)$, $(1,3,17)$, $(1,13,11)$, $(1,11,12)$, $(1,27,22)$, $(1,8,25)$, 
$(1,3,18)$, $(1,24,26)$, $(1,15,2)$, $(1,25,10)$, $(1,27,25)$, $(1,3,19)$, $(1,16,23)$, 
$(1,22,28)$, $(1,19,3)$, $(1,22,25)$, $(1,3,21)$, $(1,21,14)$, $(1,12,24)$, $(1,14,2)$, 
$(1,17,25)$, $(1,3,24)$, $(1,12,7)$, $(1,25,6)$, $(1,5,6)$, $(1,4,25)$, $(1,3,26)$, 
$(1,5,8)$, $(1,14,19)$, $(1,4,0)$, $(1,23,25)$, $(1,3,28)$, $(1,26,5)$, $(1,13,7)$, 
$(1,18,26)$, $(1,25,25)$, $(1,4,1)$, $(1,25,28)$, $(1,17,19)$, $(1,20,21)$, $(1,12,0)$, 
$(1,4,2)$, $(1,10,20)$, $(1,25,15)$, $(1,12,16)$, $(1,27,12)$, $(1,4,6)$, $(1,16,0)$, 
$(1,14,6)$, $(1,23,12)$, $(1,21,13)$, $(1,4,8)$, $(1,9,4)$, $(1,27,14)$, $(1,10,22)$, 
$(1,28,7)$, $(1,4,10)$, $(1,17,16)$, $(1,13,21)$, $(1,24,9)$, $(1,20,18)$, $(1,4,13)$, 
$(1,27,2)$, $(1,24,1)$, $(1,13,13)$, $(1,10,10)$, $(1,4,14)$, $(1,14,26)$, $(1,9,23)$, 
$(1,28,26)$, $(1,23,3)$, $(1,4,16)$, $(1,22,9)$, $(1,18,4)$, $(1,19,24)$, $(1,15,14)$, 
$(1,4,17)$, $(1,12,23)$, $(1,11,22)$, $(1,26,3)$, $(1,25,22)$, $(1,4,19)$, $(1,13,10)$, 
$(1,15,20)$, $(1,22,15)$, $(1,24,27)$, $(1,4,22)$, $(1,21,22)$, $(1,19,18)$, $(1,18,27)$, 
$(1,16,9)$, $(1,4,27)$, $(1,24,12)$, $(1,23,16)$, $(1,14,10)$, $(1,13,24)$, $(1,5,0)$, 
$(1,10,4)$, $(1,12,11)$, $(1,22,26)$, $(1,13,9)$, $(1,5,1)$, $(1,15,1)$, $(1,18,12)$, 
$(1,15,4)$, $(1,21,7)$, $(1,5,3)$, $(1,20,27)$, $(1,23,8)$, $(1,25,23)$, $(1,16,1)$, 
$(1,5,7)$, $(1,24,13)$, $(1,21,27)$, $(1,9,10)$, $(1,5,11)$, $(1,5,12)$, $(1,26,6)$, 
$(1,28,4)$, $(1,28,20)$, $(1,18,15)$, $(1,5,13)$, $(1,19,16)$, $(1,26,23)$, $(1,19,0)$, 
$(1,8,3)$, $(1,5,17)$, $(1,27,17)$, $(1,11,6)$, $(1,11,8)$, $(1,22,14)$, $(1,5,20)$, 
$(1,11,15)$, $(1,15,26)$, $(1,20,28)$, $(1,10,17)$, $(1,5,22)$, $(1,16,12)$, $(1,17,7)$, 
$(1,27,21)$, $(1,14,16)$, $(1,8,1)$, $(1,17,18)$, $(1,25,21)$, $(1,22,22)$, $(1,14,8)$, 
$(1,8,2)$, $(1,19,14)$, $(1,16,28)$, $(1,8,4)$, $(1,25,2)$, $(1,8,6)$, $(1,8,7)$, 
$(1,20,12)$, $(1,28,9)$, $(1,13,27)$, $(1,8,8)$, $(1,15,22)$, $(1,9,27)$, $(1,21,0)$, 
$(1,26,12)$, $(1,8,12)$, $(1,16,20)$, $(1,12,15)$, $(1,10,19)$, $(1,10,26)$, $(1,8,14)$, 
$(1,23,6)$, $(1,18,20)$, $(1,23,15)$, $(1,22,1)$, $(1,8,18)$, $(1,12,28)$, $(1,23,0)$, 
$(1,12,5)$, $(1,9,16)$, $(1,8,26)$, $(1,21,10)$, $(1,21,8)$, $(1,17,28)$, $(1,24,21)$, 
$(1,8,27)$, $(1,13,26)$, $(1,17,24)$, $(1,18,21)$, $(1,18,19)$, $(1,9,3)$, $(1,23,24)$, 
$(1,24,6)$, $(1,27,23)$, $(1,18,10)$, $(1,9,9)$, $(1,26,2)$, $(1,20,23)$, $(1,20,19)$, 
$(1,27,3)$, $(1,9,13)$, $(1,14,3)$, $(1,21,26)$, $(1,10,5)$, $(1,22,23)$, $(1,9,14)$, 
$(1,24,7)$, $(1,16,11)$, $(1,26,10)$, $(1,16,18)$, $(1,9,20)$, $(1,15,15)$, $(1,19,20)$, 
$(1,13,15)$, $(1,25,11)$, $(1,10,2)$, $(1,13,18)$, $(1,14,18)$, $(1,28,6)$, $(1,19,1)$, 
$(1,10,11)$, $(1,17,17)$, $(1,17,21)$, $(1,11,21)$, $(1,28,12)$, $(1,10,15)$, $(1,16,10)$, 
$(1,25,0)$, $(1,18,8)$, $(1,24,20)$, $(1,11,0)$, $(1,21,24)$, $(1,13,19)$, $(1,28,17)$, 
$(1,21,2)$, $(1,11,11)$, $(1,20,5)$, $(1,14,12)$, $(1,14,17)$, $(1,22,19)$, $(1,11,16)$, 
$(1,17,6)$, $(1,21,21)$, $(1,25,17)$, $(1,19,26)$, $(1,11,26)$, $(1,12,27)$, $(1,15,5)$, 
$(1,15,17)$, $(1,28,5)$, $(1,12,1)$, $(1,21,23)$, $(1,27,7)$, $(1,13,3)$, $(1,19,28)$, 
$(1,12,12)$, $(1,24,3)$, $(1,18,0)$, $(1,22,10)$, $(1,28,23)$, $(1,12,18)$, $(1,16,8)$, 
$(1,20,8)$, $(1,26,26)$, $(1,18,6)$, $(1,14,1)$, $(1,24,5)$, $(1,26,14)$, $(1,27,10)$, 
$(1,26,27)$, $(1,18,7)$, $(1,21,5)$, $(1,20,13)$, $(1,25,20)$, $(1,20,15)$, $(1,19,8)$, 
$(1,21,11)$, $(1,28,2)$, $(1,26,19)$, $(1,23,5)$, $(1,19,13)$, $(1,28,22)$, $(1,27,5)$, 
$(1,28,24)$, $(1,20,16)$.

\subsection*{$m_{11}(2,31)\ge 282$}
$(0,1,8)$, $(1,0,4)$, $(1,30,27)$, $(1,18,24)$, $(1,17,16)$, $(0,1,9)$, $(1,3,8)$, 
$(1,10,30)$, $(1,3,0)$, $(1,11,1)$, $(1,11,6)$, $(0,1,10)$, $(1,9,16)$, $(1,23,11)$, 
$(1,12,2)$, $(1,15,11)$, $(1,1,19)$, $(1,13,10)$, $(1,12,22)$, $(1,19,9)$, $(0,1,11)$, 
$(1,27,9)$, $(1,20,13)$, $(1,13,16)$, $(1,21,26)$, $(1,8,13)$, $(1,2,24)$, $(1,10,18)$, 
$(1,9,4)$, $(1,1,27)$, $(1,30,7)$, $(1,26,22)$, $(1,7,4)$, $(1,27,27)$, $(1,10,25)$, 
$(0,1,12)$, $(0,1,22)$, $(1,6,12)$, $(1,11,19)$, $(1,0,20)$, $(1,27,10)$, $(1,19,8)$, 
$(1,7,12)$, $(1,1,0)$, $(1,27,13)$, $(1,28,22)$, $(1,6,1)$, $(1,15,12)$, $(1,16,13)$, 
$(1,15,23)$, $(0,1,13)$, $(1,17,6)$, $(1,4,3)$, $(1,30,6)$, $(1,20,8)$, $(1,26,2)$, 
$(1,12,0)$, $(1,2,2)$, $(1,17,8)$, $(1,15,29)$, $(1,2,0)$, $(1,9,15)$, $(1,8,5)$, 
$(1,24,26)$, $(1,28,24)$, $(0,1,20)$, $(1,2,17)$, $(1,2,25)$, $(1,14,30)$, $(1,29,15)$, 
$(1,3,4)$, $(1,24,27)$, $(1,6,10)$, $(1,0,15)$, $(1,23,8)$, $(1,24,21)$, $(1,8,0)$, 
$(1,14,11)$, $(0,1,21)$, $(1,18,28)$, $(0,1,26)$, $(1,15,24)$, $(1,26,9)$, $(1,6,11)$, 
$(1,28,28)$, $(1,30,3)$, $(1,16,9)$, $(1,28,23)$, $(1,21,10)$, $(1,11,24)$, $(1,13,26)$, 
$(1,3,18)$, $(1,18,15)$, $(1,30,28)$, $(1,24,7)$, $(0,1,29)$, $(1,29,22)$, $(1,9,10)$, 
$(1,21,4)$, $(1,9,27)$, $(1,20,16)$, $(1,20,18)$, $(1,29,25)$, $(1,28,29)$, $(1,16,7)$, 
$(1,1,23)$, $(1,11,14)$, $(1,25,22)$, $(1,17,3)$, $(1,1,10)$, $(1,0,9)$, $(1,1,28)$, 
$(1,2,1)$, $(1,1,22)$, $(1,23,21)$, $(1,23,22)$, $(1,0,14)$, $(1,25,9)$, $(1,8,29)$, 
$(1,9,19)$, $(1,7,17)$, $(1,28,3)$, $(1,19,15)$, $(1,9,30)$, $(1,24,30)$, $(1,0,17)$, 
$(1,18,21)$, $(1,28,19)$, $(1,28,8)$, $(1,20,28)$, $(1,19,0)$, $(1,10,8)$, $(1,13,17)$, 
$(1,9,0)$, $(1,17,21)$, $(1,27,1)$, $(1,13,9)$, $(1,11,25)$, $(1,20,3)$, $(1,18,13)$, 
$(1,0,25)$, $(1,10,17)$, $(1,30,18)$, $(1,13,2)$, $(1,6,9)$, $(1,9,7)$, $(1,6,29)$, 
$(1,14,6)$, $(1,14,10)$, $(1,13,20)$, $(1,2,21)$, $(1,15,2)$, $(1,10,3)$, $(1,8,24)$, 
$(1,21,6)$, $(1,0,26)$, $(1,6,15)$, $(1,29,3)$, $(1,24,25)$, $(1,18,12)$, $(1,13,29)$, 
$(1,8,3)$, $(1,26,29)$, $(1,1,15)$, $(1,12,12)$, $(1,24,22)$, $(1,11,16)$, $(1,12,16)$, 
$(1,28,20)$, $(1,14,12)$, $(1,1,2)$, $(1,5,26)$, $(1,30,22)$, $(1,14,25)$, $(1,16,21)$, 
$(1,11,23)$, $(1,7,28)$, $(1,12,15)$, $(1,26,21)$, $(1,24,9)$, $(1,17,17)$, $(1,16,6)$, 
$(1,6,23)$, $(1,17,24)$, $(1,29,17)$, $(1,1,6)$, $(1,19,29)$, $(1,15,22)$, $(1,18,6)$, 
$(1,7,2)$, $(1,23,30)$, $(1,13,1)$, $(1,29,2)$, $(1,13,14)$, $(1,25,14)$, $(1,25,4)$, 
$(1,23,27)$, $(1,17,28)$, $(1,6,19)$, $(1,18,30)$, $(1,1,11)$, $(1,3,30)$, $(1,6,22)$, 
$(1,10,13)$, $(1,17,30)$, $(1,26,24)$, $(1,16,3)$, $(1,23,23)$, $(1,3,11)$, $(1,30,8)$, 
$(1,23,15)$, $(1,14,0)$, $(1,24,3)$, $(1,9,26)$, $(1,27,25)$, $(1,2,12)$, $(1,12,24)$, 
$(1,20,24)$, $(1,15,21)$, $(1,20,30)$, $(1,30,12)$, $(1,30,21)$, $(1,12,10)$, $(1,8,30)$, 
$(1,18,17)$, $(1,23,12)$, $(1,29,16)$, $(1,16,17)$, $(1,10,23)$, $(1,23,3)$, $(1,3,21)$, 
$(1,29,6)$, $(1,18,19)$, $(1,6,13)$, $(1,7,26)$, $(1,19,2)$, $(1,11,17)$, $(1,27,11)$, 
$(1,16,24)$, $(1,8,23)$, $(1,28,1)$, $(1,26,13)$, $(1,30,29)$, $(1,29,29)$, $(1,27,28)$, 
$(1,3,26)$, $(1,14,7)$, $(1,12,30)$, $(1,21,8)$, $(1,28,14)$, $(1,25,0)$, $(1,13,18)$, 
$(1,26,11)$, $(1,23,16)$, $(1,7,7)$, $(1,21,16)$, $(1,29,28)$, $(1,7,8)$, $(1,8,18)$, 
$(1,17,0)$, $(1,7,13)$, $(1,29,10)$, $(1,15,16)$, $(1,17,25)$, $(1,11,28)$, $(1,9,17)$, 
$(1,15,14)$, $(1,25,28)$, $(1,10,6)$, $(1,27,16)$, $(1,8,15)$, $(1,27,15)$, $(1,24,2)$, 
$(1,24,13)$, $(1,10,21)$, $(1,7,29)$, $(1,20,9)$, $(1,16,8)$, $(1,14,23)$, $(1,15,0)$, 
$(1,22,12)$, $(1,22,25)$.

\subsection*{$m_{13}(2,31)\ge 348$}
$(0,1,1)$, $(1,23,0)$, $(1,3,5)$, $(1,1,27)$, $(1,1,3)$, $(1,13,1)$, $(1,0,2)$, 
$(1,24,3)$, $(1,30,1)$, $(1,22,28)$, $(1,28,30)$, $(1,29,18)$, $(0,1,8)$, $(1,23,9)$, 
$(1,13,28)$, $(1,14,17)$, $(1,14,5)$, $(1,16,29)$, $(1,30,9)$, $(1,26,27)$, $(1,15,13)$, 
$(1,6,0)$, $(1,16,21)$, $(1,11,30)$, $(0,1,9)$, $(1,23,8)$, $(1,27,23)$, $(1,19,6)$, 
$(1,26,14)$, $(1,22,23)$, $(1,10,25)$, $(1,11,2)$, $(1,25,5)$, $(1,18,21)$, $(1,24,27)$, 
$(1,4,14)$, $(0,1,11)$, $(1,23,26)$, $(1,9,25)$, $(1,9,28)$, $(1,15,29)$, $(1,8,6)$, 
$(1,18,0)$, $(1,4,11)$, $(1,21,2)$, $(1,14,14)$, $(1,27,6)$, $(1,21,13)$, $(0,1,14)$, 
$(1,23,29)$, $(1,15,14)$, $(1,21,14)$, $(1,29,24)$, $(1,24,21)$, $(1,5,29)$, $(1,28,20)$, 
$(1,12,3)$, $(1,11,1)$, $(1,26,13)$, $(1,0,27)$, $(0,1,15)$, $(1,23,6)$, $(1,11,11)$, 
$(1,8,24)$, $(1,30,17)$, $(1,5,9)$, $(1,24,20)$, $(1,3,30)$, $(1,18,23)$, $(1,16,2)$, 
$(1,6,29)$, $(1,24,11)$, $(0,1,17)$, $(1,23,7)$, $(1,4,29)$, $(1,11,5)$, $(1,24,28)$, 
$(1,27,18)$, $(1,3,12)$, $(1,13,26)$, $(1,13,27)$, $(1,7,25)$, $(1,2,26)$, $(1,10,10)$, 
$(0,1,19)$, $(1,23,12)$, $(1,10,18)$, $(1,16,25)$, $(1,7,23)$, $(1,1,13)$, $(1,21,10)$, 
$(1,21,29)$, $(1,4,28)$, $(1,19,15)$, $(1,1,2)$, $(1,27,9)$, $(0,1,22)$, $(1,23,24)$, 
$(1,26,30)$, $(1,5,12)$, $(1,3,20)$, $(1,28,17)$, $(1,7,15)$, $(1,25,15)$, $(1,11,10)$, 
$(1,1,30)$, $(1,19,0)$, $(1,6,23)$, $(0,1,28)$, $(1,23,21)$, $(1,17,0)$, $(1,18,2)$, 
$(1,13,12)$, $(1,9,5)$, $(1,11,18)$, $(1,7,16)$, $(1,9,24)$, $(1,25,10)$, $(1,5,5)$, 
$(1,26,20)$, $(1,0,1)$, $(1,18,12)$, $(1,21,6)$, $(1,0,15)$, $(1,7,13)$, $(1,29,5)$, 
$(1,13,10)$, $(1,26,21)$, $(1,29,2)$, $(1,11,8)$, $(1,3,15)$, $(1,10,3)$, $(1,0,3)$, 
$(1,15,1)$, $(1,11,15)$, $(1,4,23)$, $(1,3,26)$, $(1,29,4)$, $(1,27,24)$, $(1,11,9)$, 
$(1,5,0)$, $(1,21,20)$, $(1,15,25)$, $(1,16,30)$, $(1,0,7)$, $(1,1,22)$, $(1,24,25)$, 
$(1,19,22)$, $(1,13,9)$, $(1,29,22)$, $(1,9,6)$, $(1,4,22)$, $(1,16,19)$, $(1,2,22)$, 
$(1,27,4)$, $(1,28,22)$, $(1,0,13)$, $(1,27,14)$, $(1,18,18)$, $(1,25,3)$, $(1,25,1)$, 
$(1,29,25)$, $(1,15,12)$, $(1,28,4)$, $(1,24,30)$, $(1,7,28)$, $(1,4,21)$, $(1,15,10)$, 
$(1,0,17)$, $(1,10,24)$, $(1,19,14)$, $(1,5,25)$, $(1,12,20)$, $(1,29,14)$, $(1,6,3)$, 
$(1,1,1)$, $(1,4,18)$, $(1,9,18)$, $(1,21,30)$, $(1,27,2)$, $(1,0,22)$, $(1,28,28)$, 
$(1,17,22)$, $(1,6,27)$, $(1,9,22)$, $(1,29,21)$, $(1,25,22)$, $(1,12,16)$, $(1,21,22)$, 
$(1,3,17)$, $(1,30,22)$, $(1,11,23)$, $(1,0,28)$, $(1,29,11)$, $(1,0,29)$, $(1,26,0)$, 
$(1,1,24)$, $(1,16,16)$, $(1,5,4)$, $(1,29,20)$, $(1,26,23)$, $(1,25,14)$, $(1,19,27)$, 
$(1,30,6)$, $(1,17,6)$, $(1,1,9)$, $(1,0,30)$, $(1,25,17)$, $(1,27,13)$, $(1,9,2)$, 
$(1,8,2)$, $(1,29,13)$, $(1,28,25)$, $(1,6,5)$, $(1,2,23)$, $(1,13,29)$, $(1,10,26)$, 
$(1,4,7)$, $(1,1,21)$, $(1,18,26)$, $(1,2,3)$, $(1,20,1)$, $(1,16,4)$, $(1,20,7)$, 
$(1,19,2)$, $(1,20,26)$, $(1,26,18)$, $(1,20,29)$, $(1,18,13)$, $(1,20,30)$, $(1,24,9)$, 
$(1,20,6)$, $(1,2,17)$, $(1,8,18)$, $(1,27,5)$, $(1,6,7)$, $(1,11,4)$, $(1,18,11)$, 
$(1,6,25)$, $(1,26,24)$, $(1,9,0)$, $(1,5,17)$, $(1,19,20)$, $(1,4,27)$, $(1,2,28)$, 
$(1,30,23)$, $(1,4,17)$, $(1,25,7)$, $(1,15,3)$, $(1,27,1)$, $(1,30,29)$, $(1,12,15)$, 
$(1,6,6)$, $(1,19,25)$, $(1,10,15)$, $(1,18,28)$, $(1,3,2)$, $(1,16,18)$, $(1,8,17)$, 
$(1,4,6)$, $(1,30,21)$, $(1,19,21)$, $(1,27,12)$, $(1,25,27)$, $(1,6,11)$, $(1,7,9)$, 
$(1,21,25)$, $(1,8,10)$, $(1,3,8)$, $(1,19,10)$, $(1,24,18)$, $(1,15,24)$, $(1,15,27)$, 
$(1,5,28)$, $(1,22,20)$, $(1,13,0)$, $(1,30,12)$, $(1,18,29)$, $(1,19,23)$, $(1,7,21)$, 
$(1,3,24)$, $(1,17,5)$, $(1,5,11)$, $(1,16,20)$, $(1,19,13)$, $(1,3,29)$, $(1,9,16)$, 
$(1,28,26)$, $(1,21,0)$, $(1,13,3)$, $(1,28,1)$, $(1,25,9)$, $(1,3,27)$, $(1,26,12)$, 
$(1,30,30)$, $(1,30,26)$, $(1,7,24)$, $(1,8,11)$, $(1,5,10)$, $(1,15,20)$, $(1,25,26)$, 
$(1,28,19)$, $(1,10,14)$, $(1,9,30)$, $(1,4,26)$, $(1,8,9)$, $(1,16,23)$, $(1,27,27)$, 
$(1,13,16)$, $(1,16,28)$, $(1,12,24)$, $(1,18,25)$, $(1,5,18)$, $(1,11,20)$, $(1,10,9)$, 
$(1,17,11)$, $(1,5,1)$, $(1,20,20)$, $(1,6,21)$, $(1,20,16)$, $(1,30,5)$, $(1,20,13)$, 
$(1,14,26)$, $(1,20,15)$, $(1,12,17)$, $(1,20,23)$, $(1,17,24)$, $(1,20,3)$, $(1,6,18)$, 
$(1,28,6)$, $(1,16,12)$, $(1,24,14)$, $(1,25,19)$, $(1,7,17)$, $(1,18,17)$, $(1,26,10)$, 
$(1,10,28)$, $(1,13,5)$, $(1,21,9)$, $(1,6,19)$, $(1,8,3)$, $(1,26,1)$, $(1,9,27)$, 
$(1,24,15)$, $(1,12,6)$, $(1,13,30)$, $(1,21,5)$, $(1,30,4)$.

\subsection*{$m_{14}(2,31)\ge 378$}
$(0,0,1)$, $(1,21,21)$, $(1,13,8)$, $(1,16,12)$, $(1,12,19)$, $(1,20,2)$, $(1,29,2)$, 
$(1,6,4)$, $(1,2,22)$, $(1,5,9)$, $(1,28,20)$, $(1,30,16)$, $(1,3,16)$, $(1,9,29)$, 
$(1,15,27)$, $(0,1,0)$, $(1,26,9)$, $(1,1,11)$, $(1,15,7)$, $(1,21,18)$, $(1,0,14)$, 
$(1,19,16)$, $(1,19,30)$, $(1,25,3)$, $(1,15,2)$, $(1,21,0)$, $(1,11,28)$, $(1,18,25)$, 
$(1,23,23)$, $(1,27,22)$, $(0,1,1)$, $(1,17,12)$, $(1,28,12)$, $(1,2,4)$, $(0,1,24)$, 
$(1,25,30)$, $(1,3,26)$, $(1,20,1)$, $(1,17,15)$, $(1,16,23)$, $(1,6,28)$, $(1,25,11)$, 
$(1,14,4)$, $(1,18,3)$, $(1,18,18)$, $(0,1,3)$, $(1,8,15)$, $(1,27,20)$, $(1,24,21)$, 
$(1,23,4)$, $(1,17,10)$, $(1,7,8)$, $(1,0,23)$, $(1,18,29)$, $(1,11,11)$, $(1,1,27)$, 
$(1,26,12)$, $(1,25,23)$, $(1,22,19)$, $(1,10,11)$, $(0,1,10)$, $(1,15,23)$, $(1,21,6)$, 
$(1,22,11)$, $(1,11,26)$, $(1,28,22)$, $(1,24,9)$, $(1,5,2)$, $(1,4,19)$, $(1,22,25)$, 
$(1,20,6)$, $(1,24,10)$, $(1,19,7)$, $(1,29,16)$, $(0,1,28)$, $(0,1,12)$, $(1,22,0)$, 
$(1,19,22)$, $(1,3,9)$, $(1,14,5)$, $(1,3,6)$, $(1,16,14)$, $(1,4,0)$, $(1,16,1)$, 
$(1,19,24)$, $(1,7,22)$, $(1,23,9)$, $(1,10,14)$, $(1,27,8)$, $(1,17,21)$, $(0,1,16)$, 
$(1,12,24)$, $(1,17,7)$, $(1,10,13)$, $(1,30,17)$, $(1,23,25)$, $(1,18,5)$, $(1,8,8)$, 
$(1,12,7)$, $(1,14,12)$, $(1,27,26)$, $(1,4,21)$, $(1,30,26)$, $(1,24,27)$, $(1,26,25)$, 
$(0,1,17)$, $(1,14,13)$, $(1,0,19)$, $(1,27,5)$, $(1,25,21)$, $(1,4,24)$, $(1,23,29)$, 
$(1,28,17)$, $(1,1,8)$, $(1,28,27)$, $(1,26,1)$, $(1,18,4)$, $(1,20,20)$, $(1,11,6)$, 
$(1,20,12)$, $(0,1,26)$, $(1,29,8)$, $(1,18,30)$, $(1,17,17)$, $(1,29,24)$, $(1,24,12)$, 
$(1,9,30)$, $(1,27,15)$, $(1,13,21)$, $(1,8,10)$, $(1,5,3)$, $(1,10,27)$, $(1,26,5)$, 
$(1,5,13)$, $(1,23,3)$, $(0,1,29)$, $(1,5,16)$, $(1,5,10)$, $(1,11,18)$, $(1,15,29)$, 
$(1,21,20)$, $(1,1,4)$, $(1,29,19)$, $(1,29,28)$, $(1,23,15)$, $(1,0,2)$, $(1,22,8)$, 
$(1,27,18)$, $(1,7,21)$, $(1,25,28)$, $(1,0,3)$, $(1,0,15)$, $(1,7,17)$, $(1,30,29)$, 
$(1,6,22)$, $(1,1,9)$, $(1,9,24)$, $(1,12,30)$, $(1,13,24)$, $(1,14,27)$, $(1,26,20)$, 
$(1,19,20)$, $(1,2,29)$, $(1,7,26)$, $(1,1,17)$, $(1,0,6)$, $(1,23,26)$, $(1,6,3)$, 
$(1,20,8)$, $(1,10,20)$, $(1,23,27)$, $(1,26,27)$, $(1,26,8)$, $(1,4,28)$, $(1,0,30)$, 
$(1,2,20)$, $(1,3,23)$, $(1,10,18)$, $(1,28,11)$, $(1,16,10)$, $(1,0,7)$, $(1,3,7)$, 
$(1,22,10)$, $(1,24,4)$, $(1,16,17)$, $(1,10,22)$, $(1,17,9)$, $(1,3,22)$, $(1,8,9)$, 
$(1,7,13)$, $(1,16,20)$, $(1,25,15)$, $(1,19,25)$, $(1,27,25)$, $(1,19,21)$, $(1,0,9)$, 
$(1,29,10)$, $(1,17,2)$, $(1,8,20)$, $(1,17,1)$, $(1,4,3)$, $(1,19,13)$, $(1,22,1)$, 
$(1,0,16)$, $(1,30,28)$, $(1,5,20)$, $(1,15,13)$, $(1,21,30)$, $(1,12,18)$, $(1,8,22)$, 
$(1,0,11)$, $(1,13,1)$, $(1,20,13)$, $(1,3,25)$, $(1,29,26)$, $(1,11,20)$, $(1,16,7)$, 
$(1,11,5)$, $(1,10,15)$, $(1,10,19)$, $(1,0,20)$, $(1,14,19)$, $(1,28,1)$, $(1,5,23)$, 
$(1,9,5)$, $(1,1,5)$, $(1,2,18)$, $(1,27,10)$, $(1,6,21)$, $(1,11,15)$, $(1,8,28)$, 
$(1,3,24)$, $(1,28,6)$, $(1,6,19)$, $(1,16,18)$, $(1,3,17)$, $(1,2,13)$, $(1,20,28)$, 
$(1,16,29)$, $(1,24,1)$, $(1,1,12)$, $(1,10,16)$, $(1,6,7)$, $(1,12,28)$, $(1,15,0)$, 
$(1,21,5)$, $(1,14,29)$, $(1,5,18)$, $(1,1,15)$, $(1,18,14)$, $(1,2,2)$, $(1,10,5)$, 
$(1,18,12)$, $(1,29,29)$, $(1,30,25)$, $(1,1,16)$, $(1,8,1)$, $(1,14,17)$, $(1,15,16)$, 
$(1,21,24)$, $(1,27,23)$, $(1,15,21)$, $(1,21,11)$, $(1,4,5)$, $(1,24,2)$, $(1,11,13)$, 
$(1,9,6)$, $(1,9,2)$, $(1,13,29)$, $(1,13,19)$, $(1,1,22)$, $(1,7,9)$, $(1,15,26)$, 
$(1,21,23)$, $(1,23,1)$, $(1,23,11)$, $(1,5,8)$, $(1,19,8)$, $(1,8,2)$, $(1,28,25)$, 
$(1,7,15)$, $(1,14,1)$, $(1,29,7)$, $(1,8,29)$, $(1,20,16)$, $(1,1,24)$, $(1,5,25)$, 
$(1,16,4)$, $(1,28,26)$, $(1,9,7)$, $(1,25,17)$, $(1,13,6)$, $(1,29,23)$, $(1,20,24)$, 
$(1,4,11)$, $(1,25,6)$, $(1,17,29)$, $(1,17,4)$, $(1,1,29)$, $(1,3,10)$, $(1,2,5)$, 
$(1,12,22)$, $(1,22,9)$, $(1,25,19)$, $(1,25,13)$, $(1,11,10)$, $(1,22,30)$, $(1,6,13)$, 
$(1,19,26)$, $(1,18,19)$, $(1,4,4)$, $(1,9,4)$, $(1,30,12)$, $(1,13,5)$, $(1,11,3)$, 
$(1,2,14)$, $(1,17,18)$, $(1,11,25)$, $(1,12,2)$, $(1,5,27)$, $(1,26,7)$, $(1,17,19)$, 
$(1,6,27)$, $(1,26,24)$, $(1,11,16)$, $(1,29,0)$, $(1,4,10)$, $(1,15,28)$, $(1,21,19)$, 
$(1,30,24)$, $(1,2,23)$, $(1,14,8)$, $(1,9,11)$, $(1,20,22)$, $(1,13,9)$, $(1,30,0)$, 
$(1,18,15)$, $(1,6,6)$, $(1,18,13)$, $(1,12,12)$, $(1,19,14)$, $(1,23,12)$, $(1,24,6)$, 
$(1,28,8)$, $(1,23,13)$, $(1,3,21)$, $(1,12,27)$, $(1,26,23)$, $(1,16,3)$, $(1,7,11)$, 
$(1,24,3)$, $(1,4,26)$, $(1,14,24)$, $(1,8,5)$, $(1,22,26)$, $(1,3,30)$, $(1,18,23)$, 
$(1,30,30)$, $(1,16,11)$, $(1,27,12)$, $(1,4,15)$, $(1,5,4)$, $(1,25,16)$, $(1,24,15)$, 
$(1,30,22)$, $(1,15,9)$, $(1,21,28)$, $(1,29,3)$, $(1,12,17)$, $(1,25,4)$, $(1,7,18)$, 
$(1,19,2)$, $(1,14,6)$, $(1,10,17)$, $(1,24,30)$, $(1,10,6)$, $(1,14,10)$, $(1,27,21)$.

\subsection*{$m_{20}(2,31)\ge 567$}
$(0,0,1)$, $(1,5,2)$, $(0,1,18)$, $(1,7,1)$, $(1,11,0)$, $(1,29,20)$, $(1,28,8)$, 
$(1,6,4)$, $(1,15,23)$, $(1,14,18)$, $(0,1,1)$, $(1,12,14)$, $(1,27,9)$, $(1,16,13)$, 
$(1,26,18)$, $(1,10,27)$, $(1,4,19)$, $(1,9,17)$, $(1,15,20)$, $(1,27,17)$, $(0,1,2)$, 
$(1,27,22)$, $(1,6,3)$, $(1,11,27)$, $(1,22,7)$, $(1,7,2)$, $(1,3,13)$, $(1,20,13)$, 
$(1,15,14)$, $(1,28,5)$, $(0,1,4)$, $(1,1,4)$, $(1,26,22)$, $(1,4,28)$, $(1,2,14)$, 
$(1,4,8)$, $(1,12,5)$, $(1,25,14)$, $(1,15,26)$, $(1,8,28)$, $(0,1,5)$, $(1,14,13)$, 
$(1,5,16)$, $(1,24,3)$, $(1,24,28)$, $(1,12,23)$, $(1,26,27)$, $(1,17,0)$, $(1,15,0)$, 
$(1,26,29)$, $(0,1,6)$, $(1,29,21)$, $(1,15,10)$, $(1,20,8)$, $(0,1,26)$, $(1,18,11)$, 
$(1,29,14)$, $(1,8,23)$, $(1,15,11)$, $(1,17,13)$, $(0,1,7)$, $(1,0,20)$, $(1,25,4)$, 
$(1,19,17)$, $(1,29,3)$, $(1,2,12)$, $(1,16,29)$, $(1,4,16)$, $(1,15,24)$, $(1,29,24)$, 
$(0,1,8)$, $(1,8,16)$, $(1,4,29)$, $(1,27,7)$, $(1,20,17)$, $(1,14,19)$, $(1,5,25)$, 
$(1,26,8)$, $(1,15,21)$, $(1,9,16)$, $(0,1,9)$, $(1,30,5)$, $(1,14,23)$, $(1,5,19)$, 
$(1,17,1)$, $(1,21,5)$, $(1,0,26)$, $(1,18,25)$, $(1,15,7)$, $(1,3,26)$, $(0,1,11)$, 
$(1,13,29)$, $(1,3,11)$, $(1,18,26)$, $(1,27,13)$, $(1,27,24)$, $(1,19,16)$, $(1,23,26)$, 
$(1,15,29)$, $(1,13,30)$, $(0,1,14)$, $(1,24,8)$, $(1,2,24)$, $(1,21,30)$, $(1,21,12)$, 
$(1,8,0)$, $(1,9,18)$, $(1,10,11)$, $(1,15,18)$, $(1,10,4)$, $(0,1,16)$, $(1,16,12)$, 
$(1,22,12)$, $(1,0,2)$, $(1,19,22)$, $(1,0,16)$, $(1,30,20)$, $(1,1,3)$, $(1,15,13)$, 
$(1,18,1)$, $(0,1,20)$, $(1,3,3)$, $(1,0,19)$, $(1,13,9)$, $(1,30,29)$, $(1,25,28)$, 
$(1,21,28)$, $(1,22,1)$, $(1,15,2)$, $(1,21,27)$, $(0,1,28)$, $(1,22,9)$, $(1,18,2)$, 
$(1,28,29)$, $(1,9,10)$, $(1,24,30)$, $(1,18,10)$, $(1,11,5)$, $(1,15,25)$, $(1,0,0)$, 
$(0,1,30)$, $(1,2,19)$, $(1,7,21)$, $(1,1,24)$, $(1,4,4)$, $(1,3,10)$, $(1,23,9)$, 
$(1,2,28)$, $(1,15,4)$, $(1,6,21)$, $(1,0,3)$, $(1,4,5)$, $(1,20,4)$, $(1,8,18)$, 
$(1,10,22)$, $(1,2,22)$, $(1,19,11)$, $(1,13,23)$, $(1,18,21)$, $(1,17,8)$, $(1,0,4)$, 
$(1,12,12)$, $(1,9,23)$, $(1,2,9)$, $(1,17,23)$, $(1,23,17)$, $(1,17,22)$, $(1,25,29)$, 
$(1,12,25)$, $(1,23,19)$, $(1,0,5)$, $(1,23,10)$, $(1,29,11)$, $(1,24,11)$, $(1,4,30)$, 
$(1,14,28)$, $(1,5,26)$, $(1,29,0)$, $(1,14,3)$, $(1,5,17)$, $(1,0,7)$, $(1,9,21)$, 
$(1,7,18)$, $(1,10,21)$, $(1,0,25)$, $(1,21,16)$, $(1,26,19)$, $(1,26,14)$, $(1,19,10)$, 
$(1,6,24)$, $(1,0,12)$, $(1,28,26)$, $(1,14,20)$, $(1,5,29)$, $(1,16,14)$, $(1,11,11)$, 
$(1,20,21)$, $(1,16,9)$, $(1,13,14)$, $(1,1,20)$, $(1,0,13)$, $(1,0,17)$, $(1,3,8)$, 
$(1,21,22)$, $(1,14,27)$, $(1,5,8)$, $(1,8,25)$, $(1,1,17)$, $(1,23,28)$, $(1,8,7)$, 
$(1,0,18)$, $(1,10,18)$, $(1,10,10)$, $(1,27,0)$, $(1,12,9)$, $(1,3,7)$, $(1,12,3)$, 
$(1,20,11)$, $(1,1,22)$, $(1,24,26)$, $(1,0,22)$, $(1,26,1)$, $(1,28,24)$, $(1,28,17)$, 
$(1,13,18)$, $(1,10,26)$, $(1,11,24)$, $(1,11,22)$, $(1,30,13)$, $(1,4,10)$, $(1,0,29)$, 
$(1,22,13)$, $(1,13,2)$, $(1,16,30)$, $(1,29,7)$, $(1,25,18)$, $(1,10,14)$, $(1,12,7)$, 
$(1,4,20)$, $(1,29,30)$, $(1,1,8)$, $(1,28,18)$, $(1,10,12)$, $(1,6,16)$, $(1,8,30)$, 
$(1,11,26)$, $(1,7,19)$, $(1,25,16)$, $(1,22,18)$, $(1,10,2)$, $(1,1,9)$, $(1,21,1)$, 
$(1,27,8)$, $(1,1,29)$, $(1,20,3)$, $(1,26,21)$, $(1,8,17)$, $(1,18,23)$, $(1,30,30)$, 
$(1,17,30)$, $(1,1,10)$, $(1,9,25)$, $(1,13,4)$, $(1,24,25)$, $(1,19,13)$, $(1,2,29)$, 
$(1,23,18)$, $(1,10,0)$, $(1,3,5)$, $(1,1,28)$, $(1,1,11)$, $(1,25,24)$, $(1,30,0)$, 
$(1,18,22)$, $(1,12,21)$, $(1,24,1)$, $(1,3,27)$, $(1,27,14)$, $(1,24,21)$, $(1,25,0)$, 
$(1,1,13)$, $(1,11,21)$, $(1,2,23)$, $(1,22,24)$, $(1,17,2)$, $(1,14,25)$, $(1,5,23)$, 
$(1,28,13)$, $(1,17,26)$, $(1,26,4)$, $(1,1,16)$, $(1,17,9)$, $(1,22,11)$, $(1,9,2)$, 
$(1,22,14)$, $(1,21,2)$, $(1,13,7)$, $(1,13,28)$, $(1,4,22)$, $(1,23,23)$, $(1,1,25)$, 
$(1,30,14)$, $(1,20,6)$, $(1,11,3)$, $(1,14,1)$, $(1,5,28)$, $(1,20,24)$, $(1,12,29)$, 
$(1,19,29)$, $(1,28,12)$, $(1,1,26)$, $(1,22,30)$, $(1,6,2)$, $(1,27,11)$, $(1,30,27)$, 
$(1,7,17)$, $(1,30,4)$, $(1,9,1)$, $(1,25,7)$, $(1,2,1)$, $(1,1,27)$, $(1,16,11)$, 
$(1,23,29)$, $(1,21,8)$, $(1,9,20)$, $(1,23,22)$, $(1,22,20)$, $(1,8,2)$, $(1,26,24)$, 
$(1,16,26)$, $(1,2,0)$, $(1,6,27)$, $(1,18,20)$, $(1,16,23)$, $(1,7,5)$, $(1,17,29)$, 
$(1,24,23)$, $(1,16,4)$, $(1,17,20)$, $(1,13,19)$, $(1,2,3)$, $(1,30,7)$, $(1,30,19)$, 
$(1,27,12)$, $(1,4,1)$, $(1,13,1)$, $(1,6,13)$, $(1,3,4)$, $(1,4,24)$, $(1,7,25)$, 
$(1,2,10)$, $(1,13,16)$, $(1,27,27)$, $(1,21,18)$, $(1,10,9)$, $(1,19,12)$, $(1,14,14)$, 
$(1,5,4)$, $(1,14,9)$, $(1,5,27)$, $(1,2,13)$, $(1,9,9)$, $(1,8,26)$, $(1,30,9)$, 
$(1,26,20)$, $(1,18,5)$, $(1,3,1)$, $(1,19,4)$, $(1,25,8)$, $(1,16,16)$, $(1,2,16)$, 
$(1,18,17)$, $(1,20,25)$, $(1,28,11)$, $(1,21,3)$, $(1,3,24)$, $(1,25,27)$, $(1,11,4)$, 
$(1,21,14)$, $(1,18,14)$, $(1,2,20)$, $(1,14,10)$, $(1,5,3)$, $(1,23,16)$, $(1,6,14)$, 
$(1,16,22)$, $(1,27,4)$, $(1,7,4)$, $(1,27,5)$, $(1,2,30)$, $(1,2,25)$, $(1,27,25)$, 
$(1,25,22)$, $(1,29,10)$, $(1,16,17)$, $(1,21,26)$, $(1,19,3)$, $(1,22,4)$, $(1,29,2)$, 
$(1,12,20)$, $(1,3,2)$, $(1,17,17)$, $(1,19,24)$, $(1,18,9)$, $(1,4,17)$, $(1,22,27)$, 
$(1,14,2)$, $(1,5,18)$, $(1,10,30)$, $(1,22,23)$, $(1,3,12)$, $(1,30,10)$, $(1,21,29)$, 
$(1,29,17)$, $(1,27,1)$, $(1,26,23)$, $(1,9,22)$, $(1,17,5)$, $(1,29,22)$, $(1,20,14)$, 
$(1,3,15)$, $(1,3,16)$, $(1,9,7)$, $(1,28,0)$, $(1,23,7)$, $(1,16,10)$, $(1,28,21)$, 
$(1,11,14)$, $(1,25,17)$, $(1,24,29)$, $(1,11,20)$, $(1,3,17)$, $(1,6,11)$, $(1,22,16)$, 
$(1,20,2)$, $(1,7,23)$, $(1,6,12)$, $(1,17,21)$, $(1,12,13)$, $(1,7,28)$, $(1,9,11)$, 
$(1,3,20)$, $(1,22,0)$, $(1,4,2)$, $(1,24,19)$, $(1,22,22)$, $(1,16,2)$, $(1,20,9)$, 
$(1,24,0)$, $(1,20,16)$, $(1,19,25)$, $(1,3,23)$, $(1,8,29)$, $(1,17,19)$, $(1,4,27)$, 
$(1,28,3)$, $(1,19,30)$, $(1,4,11)$, $(1,19,8)$, $(1,11,10)$, $(1,26,10)$, $(1,4,7)$, 
$(1,11,9)$, $(1,9,5)$, $(1,21,6)$, $(1,16,5)$, $(1,30,24)$, $(1,23,13)$, $(1,28,30)$, 
$(1,20,1)$, $(1,20,20)$, $(1,4,21)$, $(1,30,26)$, $(1,16,7)$, $(1,19,23)$, $(1,12,1)$, 
$(1,17,3)$, $(1,18,18)$, $(1,10,20)$, $(1,21,9)$, $(1,23,21)$, $(1,5,11)$, $(1,7,7)$, 
$(1,29,25)$, $(1,6,9)$, $(1,14,24)$, $(1,5,12)$, $(1,13,22)$, $(1,28,7)$, $(1,22,5)$, 
$(1,19,2)$, $(1,9,28)$, $(1,11,19)$, $(1,19,28)$, $(1,27,21)$, $(1,14,11)$, $(1,5,22)$, 
$(1,9,12)$, $(1,18,13)$, $(1,25,12)$, $(1,25,19)$, $(1,6,23)$, $(1,13,25)$, $(1,17,10)$, 
$(1,7,14)$, $(1,14,30)$, $(1,6,8)$, $(1,30,8)$, $(1,13,0)$, $(1,7,24)$, $(1,8,1)$, 
$(1,8,14)$, $(1,13,3)$, $(1,21,23)$, $(1,24,10)$, $(1,17,7)$, $(1,6,18)$, $(1,10,28)$, 
$(1,26,17)$, $(1,7,12)$, $(1,20,28)$, $(1,23,5)$, $(1,27,18)$, $(1,10,7)$, $(1,18,29)$, 
$(1,18,24)$, $(1,6,28)$, $(1,25,13)$, $(1,8,3)$, $(1,7,10)$, $(1,19,18)$, $(1,10,19)$, 
$(1,20,26)$, $(1,25,26)$, $(1,23,8)$, $(1,20,27)$, $(1,8,11)$, $(1,28,1)$, $(1,12,16)$, 
$(1,12,0)$, $(1,30,17)$, $(1,11,16)$, $(1,13,8)$, $(1,19,1)$, $(1,30,25)$, $(1,24,5)$, 
$(1,9,3)$, $(1,10,25)$, $(1,29,28)$, $(1,13,27)$, $(1,25,5)$, $(1,28,20)$, $(1,30,21)$, 
$(1,29,27)$, $(1,26,25)$, $(1,11,1)$, $(1,16,18)$, $(1,11,12)$, $(1,27,3)$, $(1,29,26)$, 
$(1,21,21)$, $(1,28,28)$, $(1,24,9)$, $(1,28,9)$, $(1,29,19)$, $(1,11,13)$, $(1,23,0)$.

\subsection*{$m_{21}(2,31)\ge 597$}
$(0,1,0)$, $(1,23,16)$, $(1,2,18)$, $(1,22,29)$, $(1,26,4)$, $(1,7,1)$, $(1,18,10)$, 
$(1,1,20)$, $(1,20,12)$, $(1,26,0)$, $(0,1,1)$, $(1,2,13)$, $(1,12,12)$, $(1,21,6)$, 
$(1,16,18)$, $(1,11,17)$, $(1,0,14)$, $(1,3,8)$, $(1,14,13)$, $(1,4,26)$, $(0,1,2)$, 
$(1,19,11)$, $(1,8,2)$, $(1,28,12)$, $(1,15,7)$, $(1,19,18)$, $(1,9,12)$, $(1,17,17)$, 
$(1,15,18)$, $(1,20,24)$, $(0,1,4)$, $(1,17,24)$, $(1,7,15)$, $(1,13,8)$, $(1,3,30)$, 
$(0,1,7)$, $(1,21,29)$, $(1,17,9)$, $(1,23,21)$, $(1,7,12)$, $(1,10,13)$, $(1,13,18)$, 
$(1,28,13)$, $(1,3,20)$, $(1,30,15)$, $(0,1,9)$, $(1,29,8)$, $(1,5,10)$, $(1,24,13)$, 
$(1,1,8)$, $(1,4,20)$, $(1,14,4)$, $(1,16,23)$, $(1,11,29)$, $(0,1,27)$, $(0,1,10)$, 
$(1,10,23)$, $(1,15,4)$, $(1,6,2)$, $(1,10,14)$, $(1,18,14)$, $(1,16,7)$, $(1,11,22)$, 
$(1,1,10)$, $(1,19,28)$, $(0,1,12)$, $(1,14,28)$, $(1,22,6)$, $(1,7,25)$, $(1,14,27)$, 
$(1,20,22)$, $(1,8,26)$, $(1,18,11)$, $(1,2,15)$, $(1,17,5)$, $(0,1,14)$, $(1,24,25)$, 
$(1,19,14)$, $(1,14,0)$, $(1,13,16)$, $(1,15,2)$, $(1,21,30)$, $(1,4,2)$, $(1,4,25)$, 
$(1,29,19)$, $(0,1,19)$, $(1,5,9)$, $(1,23,24)$, $(1,8,17)$, $(1,27,15)$, $(1,29,27)$, 
$(1,27,8)$, $(1,24,6)$, $(1,5,30)$, $(1,7,14)$, $(0,1,21)$, $(1,3,22)$, $(1,27,3)$, 
$(1,20,14)$, $(1,29,6)$, $(1,25,11)$, $(1,24,19)$, $(1,29,7)$, $(1,21,17)$, $(1,24,8)$, 
$(0,1,22)$, $(1,20,20)$, $(1,25,29)$, $(1,1,11)$, $(1,17,29)$, $(1,13,25)$, $(1,6,23)$, 
$(1,8,9)$, $(1,23,27)$, $(1,21,20)$, $(0,1,23)$, $(1,30,17)$, $(1,6,28)$, $(1,25,5)$, 
$(1,9,3)$, $(1,30,0)$, $(1,22,16)$, $(0,1,25)$, $(1,18,2)$, $(1,14,17)$, $(0,1,24)$, 
$(1,26,12)$, $(1,20,1)$, $(1,4,18)$, $(1,21,11)$, $(1,5,24)$, $(1,26,22)$, $(1,30,1)$, 
$(1,25,6)$, $(1,10,2)$, $(0,1,29)$, $(1,4,0)$, $(1,26,16)$, $(1,10,1)$, $(1,20,0)$, 
$(1,3,16)$, $(1,29,11)$, $(1,26,25)$, $(1,24,1)$, $(1,27,27)$, $(0,1,30)$, $(1,27,21)$, 
$(1,16,22)$, $(1,11,24)$, $(1,28,26)$, $(1,17,10)$, $(1,2,17)$, $(1,25,0)$, $(1,27,16)$, 
$(1,25,4)$, $(1,0,6)$, $(1,22,3)$, $(1,27,29)$, $(1,27,0)$, $(1,19,16)$, $(1,30,5)$, 
$(1,12,8)$, $(1,15,3)$, $(1,22,5)$, $(1,13,20)$, $(1,0,9)$, $(1,23,6)$, $(1,17,0)$, 
$(1,20,16)$, $(1,18,15)$, $(1,9,11)$, $(1,12,14)$, $(1,20,21)$, $(1,9,28)$, $(1,7,11)$, 
$(1,0,29)$, $(1,20,28)$, $(1,28,4)$, $(1,28,11)$, $(1,25,22)$, $(1,4,8)$, $(1,12,2)$, 
$(1,18,20)$, $(1,3,10)$, $(1,5,8)$, $(1,1,1)$, $(1,7,3)$, $(1,8,6)$, $(1,14,23)$, 
$(1,26,5)$, $(1,14,1)$, $(1,17,15)$, $(1,22,9)$, $(1,23,22)$, $(1,26,14)$, $(1,1,3)$, 
$(1,5,4)$, $(1,30,9)$, $(1,23,23)$, $(1,19,10)$, $(1,3,14)$, $(1,19,2)$, $(1,15,10)$, 
$(1,27,4)$, $(1,15,24)$, $(1,1,5)$, $(1,30,7)$, $(1,20,2)$, $(1,22,23)$, $(1,30,11)$, 
$(1,29,3)$, $(1,3,13)$, $(1,24,22)$, $(1,3,19)$, $(1,17,25)$, $(1,1,6)$, $(1,18,13)$, 
$(1,19,23)$, $(1,7,23)$, $(1,25,19)$, $(1,15,28)$, $(1,14,19)$, $(1,25,13)$, $(1,28,15)$, 
$(1,1,17)$, $(1,1,9)$, $(1,23,26)$, $(1,4,28)$, $(1,2,23)$, $(1,27,22)$, $(1,17,20)$, 
$(1,8,27)$, $(1,3,25)$, $(1,18,29)$, $(1,10,6)$, $(1,1,13)$, $(1,15,30)$, $(1,10,26)$, 
$(1,30,23)$, $(1,2,0)$, $(1,18,16)$, $(1,5,0)$, $(1,4,16)$, $(1,19,9)$, $(1,23,28)$, 
$(1,1,14)$, $(1,10,17)$, $(1,29,30)$, $(1,17,23)$, $(1,24,2)$, $(1,9,21)$, $(1,1,26)$, 
$(1,9,2)$, $(1,2,8)$, $(1,13,23)$, $(1,1,15)$, $(1,14,15)$, $(1,28,20)$, $(1,5,23)$, 
$(1,10,12)$, $(1,7,29)$, $(1,7,18)$, $(1,27,26)$, $(1,10,3)$, $(1,2,2)$, $(1,1,16)$, 
$(1,24,10)$, $(1,17,3)$, $(1,18,23)$, $(1,21,13)$, $(1,21,4)$, $(1,27,12)$, $(1,13,28)$, 
$(1,17,18)$, $(1,29,0)$, $(1,1,21)$, $(1,26,9)$, $(1,23,1)$, $(1,3,23)$, $(1,16,21)$, 
$(1,11,13)$, $(1,26,3)$, $(1,16,1)$, $(1,11,14)$, $(1,27,30)$, $(1,1,22)$, $(1,25,25)$, 
$(1,5,7)$, $(1,9,23)$, $(1,4,3)$, $(1,1,30)$, $(1,8,18)$, $(1,25,21)$, $(1,20,23)$, 
$(1,5,20)$, $(1,27,11)$, $(1,9,5)$, $(1,6,29)$, $(1,4,30)$, $(1,22,12)$, $(1,2,1)$, 
$(1,6,8)$, $(1,28,28)$, $(1,27,25)$, $(1,26,2)$, $(1,3,1)$, $(1,9,24)$, $(1,13,27)$, 
$(1,18,26)$, $(1,26,19)$, $(1,2,6)$, $(1,4,17)$, $(1,13,3)$, $(1,19,3)$, $(1,4,24)$, 
$(1,20,10)$, $(1,15,9)$, $(1,23,11)$, $(1,16,26)$, $(1,11,4)$, $(1,2,9)$, $(1,23,9)$, 
$(1,23,30)$, $(1,9,22)$, $(1,21,7)$, $(1,14,5)$, $(1,19,30)$, $(1,21,8)$, $(1,22,26)$, 
$(1,15,8)$, $(1,2,10)$, $(1,8,30)$, $(1,29,9)$, $(1,23,14)$, $(1,24,4)$, $(1,8,0)$, 
$(1,6,16)$, $(1,3,12)$, $(1,25,26)$, $(1,25,18)$, $(1,2,11)$, $(1,20,7)$, $(1,7,24)$, 
$(1,5,11)$, $(1,28,0)$, $(1,21,16)$, $(1,13,14)$, $(1,5,15)$, $(1,20,26)$, $(1,14,7)$, 
$(1,2,16)$, $(1,7,19)$, $(1,30,21)$, $(1,30,10)$, $(1,13,15)$, $(1,4,7)$, $(1,25,15)$, 
$(1,8,4)$, $(1,13,26)$, $(1,7,0)$, $(1,2,22)$, $(1,5,28)$, $(1,16,8)$, $(1,11,12)$, 
$(1,16,12)$, $(1,11,18)$, $(1,30,18)$, $(1,18,19)$, $(1,6,26)$, $(1,13,6)$, $(1,2,25)$, 
$(1,3,6)$, $(1,6,12)$, $(1,29,15)$, $(1,10,18)$, $(1,24,3)$, $(1,21,25)$, $(1,28,3)$, 
$(1,7,26)$, $(1,6,30)$, $(1,2,26)$, $(1,2,27)$, $(1,14,3)$, $(1,9,17)$, $(1,3,21)$, 
$(1,8,20)$, $(1,16,17)$, $(1,11,19)$, $(1,4,29)$, $(1,5,26)$, $(1,3,27)$, $(1,2,28)$, 
$(1,30,24)$, $(1,3,7)$, $(1,26,30)$, $(1,27,1)$, $(1,5,13)$, $(1,18,17)$, $(1,17,2)$, 
$(1,26,26)$, $(1,8,1)$, $(1,3,2)$, $(1,20,6)$, $(1,21,21)$, $(1,3,11)$, $(1,21,10)$, 
$(1,18,0)$, $(1,8,16)$, $(1,9,7)$, $(1,16,4)$, $(1,11,28)$, $(1,3,24)$, $(1,22,14)$, 
$(1,22,28)$, $(1,4,12)$, $(1,30,20)$, $(1,9,18)$, $(1,14,9)$, $(1,23,29)$, $(1,18,5)$, 
$(1,4,5)$, $(1,3,29)$, $(1,15,17)$, $(1,30,22)$, $(1,14,22)$, $(1,29,12)$, $(1,6,24)$, 
$(1,24,18)$, $(1,6,20)$, $(1,19,21)$, $(1,17,30)$, $(1,4,6)$, $(1,15,5)$, $(1,17,6)$, 
$(1,28,19)$, $(1,22,2)$, $(1,17,27)$, $(1,16,30)$, $(1,11,30)$, $(1,19,0)$, $(1,28,16)$, 
$(1,4,9)$, $(1,23,13)$, $(1,6,7)$, $(1,18,3)$, $(1,28,9)$, $(1,23,19)$, $(1,28,24)$, 
$(1,29,28)$, $(1,8,21)$, $(1,24,24)$, $(1,4,27)$, $(1,24,14)$, $(1,13,12)$, $(1,27,5)$, 
$(1,18,18)$, $(1,13,22)$, $(1,22,27)$, $(1,9,13)$, $(1,25,28)$, $(1,29,14)$, $(1,5,1)$, 
$(1,14,30)$, $(1,6,10)$, $(1,9,27)$, $(1,26,10)$, $(1,28,1)$, $(1,13,4)$, $(1,14,25)$, 
$(1,30,4)$, $(1,26,13)$, $(1,5,5)$, $(1,22,4)$, $(1,18,21)$, $(1,22,8)$, $(1,6,0)$, 
$(1,17,16)$, $(1,21,28)$, $(1,10,10)$, $(1,6,5)$, $(1,5,17)$, $(1,5,18)$, $(1,19,6)$, 
$(1,27,6)$, $(1,8,7)$, $(1,15,20)$, $(1,13,13)$, $(1,29,21)$, $(1,13,29)$, $(1,24,12)$, 
$(1,14,2)$, $(1,5,25)$, $(1,21,15)$, $(1,25,30)$, $(1,30,13)$, $(1,16,5)$, $(1,11,27)$, 
$(1,7,17)$, $(1,22,24)$, $(1,16,2)$, $(1,11,7)$, $(1,6,13)$, $(1,14,24)$, $(1,15,0)$, 
$(1,29,16)$, $(1,20,3)$, $(1,24,27)$, $(1,10,11)$, $(1,10,8)$, $(1,25,27)$, $(1,19,24)$, 
$(1,6,15)$, $(1,16,29)$, $(1,11,2)$, $(1,13,24)$, $(1,25,20)$, $(1,7,5)$, $(1,28,14)$, 
$(1,9,29)$, $(1,22,20)$, $(1,21,3)$, $(1,7,7)$, $(1,24,0)$, $(1,25,16)$, $(1,16,27)$, 
$(1,11,9)$, $(1,23,15)$, $(1,30,3)$, $(1,17,12)$, $(1,9,8)$, $(1,29,18)$, $(1,7,8)$, 
$(1,21,5)$, $(1,26,21)$, $(1,25,1)$, $(1,22,21)$, $(1,19,15)$, $(1,19,19)$, $(1,26,27)$, 
$(1,15,1)$, $(1,15,25)$, $(1,7,10)$, $(1,26,7)$, $(1,22,1)$, $(1,29,17)$, $(1,9,4)$, 
$(1,24,15)$, $(1,26,6)$, $(1,29,1)$, $(1,28,22)$, $(1,24,5)$, $(1,7,13)$, $(1,9,25)$, 
$(1,8,24)$, $(1,19,8)$, $(1,10,22)$, $(1,9,15)$, $(1,15,22)$, $(1,16,0)$, $(1,11,16)$, 
$(1,8,13)$, $(1,8,22)$, $(1,10,4)$, $(1,22,17)$, $(1,20,9)$, $(1,23,25)$, $(1,13,11)$, 
$(1,13,17)$, $(1,16,10)$, $(1,11,21)$, $(1,10,5)$, $(1,10,9)$, $(1,23,10)$, $(1,29,29)$, 
$(1,30,19)$, $(1,24,30)$, $(1,20,4)$, $(1,20,19)$, $(1,16,11)$, $(1,11,11)$, $(1,14,20)$, 
$(1,10,25)$, $(1,16,28)$, $(1,11,20)$, $(1,17,8)$, $(1,27,7)$, $(1,28,29)$, $(1,14,21)$, 
$(1,18,7)$, $(1,27,14)$, $(1,30,25)$, $(1,10,27)$, $(1,28,6)$, $(1,30,14)$, $(1,21,9)$, 
$(1,23,17)$, $(1,18,28)$.

\subsection*{$m_{22}(2,31)\ge 631$}
$(0,1,0)$, $(1,25,29)$, $(1,16,15)$, $(1,7,16)$, $(1,30,0)$, $(1,23,7)$, $(0,1,17)$, 
$(1,10,12)$, $(1,4,21)$, $(1,4,4)$, $(0,1,2)$, $(1,6,24)$, $(1,2,22)$, $(1,16,21)$, 
$(1,18,0)$, $(1,4,15)$, $(1,19,4)$, $(1,5,29)$, $(1,9,11)$, $(1,2,20)$, $(0,1,5)$, 
$(1,29,17)$, $(1,27,25)$, $(1,27,3)$, $(1,8,0)$, $(1,6,6)$, $(1,12,9)$, $(1,4,20)$, 
$(1,28,4)$, $(1,28,29)$, $(0,1,7)$, $(1,22,7)$, $(1,8,19)$, $(1,1,23)$, $(1,29,0)$, 
$(1,9,8)$, $(1,21,7)$, $(1,1,24)$, $(1,15,30)$, $(1,16,1)$, $(0,1,8)$, $(1,24,1)$, 
$(1,6,20)$, $(1,11,1)$, $(1,12,0)$, $(1,12,10)$, $(1,18,18)$, $(1,20,9)$, $(1,19,22)$, 
$(1,6,19)$, $(0,1,13)$, $(1,9,15)$, $(1,0,23)$, $(1,17,25)$, $(1,1,0)$, $(1,13,21)$, 
$(1,29,19)$, $(1,0,15)$, $(1,6,17)$, $(1,29,21)$, $(0,1,14)$, $(1,4,30)$, $(1,7,4)$, 
$(1,18,29)$, $(1,13,0)$, $(1,19,25)$, $(1,9,20)$, $(1,17,13)$, $(1,25,10)$, $(1,9,26)$, 
$(0,1,15)$, $(1,21,10)$, $(1,15,0)$, $(1,29,11)$, $(1,6,0)$, $(1,22,27)$, $(1,2,25)$, 
$(1,14,17)$, $(0,1,29)$, $(1,5,27)$, $(0,1,19)$, $(1,19,16)$, $(1,1,7)$, $(1,20,6)$, 
$(1,17,0)$, $(1,11,30)$, $(1,14,12)$, $(1,6,7)$, $(1,24,12)$, $(1,18,16)$, $(0,1,20)$, 
$(1,28,20)$, $(1,30,8)$, $(1,23,18)$, $(1,10,0)$, $(1,27,20)$, $(1,24,27)$, $(1,24,14)$, 
$(1,29,2)$, $(1,1,28)$, $(0,1,21)$, $(1,15,28)$, $(1,12,17)$, $(1,10,28)$, $(1,22,0)$, 
$(1,15,12)$, $(1,5,14)$, $(1,18,22)$, $(1,22,16)$, $(1,19,8)$, $(0,1,22)$, $(1,7,21)$, 
$(1,11,2)$, $(1,30,15)$, $(1,23,0)$, $(1,8,28)$, $(1,30,5)$, $(1,23,5)$, $(1,17,26)$, 
$(1,12,2)$, $(0,1,23)$, $(0,1,28)$, $(1,20,13)$, $(1,25,26)$, $(1,3,0)$, $(1,28,0)$, 
$(1,7,17)$, $(1,3,11)$, $(1,21,18)$, $(1,27,6)$, $(0,1,24)$, $(1,8,18)$, $(1,14,16)$, 
$(1,0,19)$, $(1,2,0)$, $(1,2,24)$, $(1,4,28)$, $(1,16,4)$, $(1,0,29)$, $(1,3,12)$, 
$(0,1,25)$, $(1,0,11)$, $(1,5,5)$, $(1,13,9)$, $(1,9,0)$, $(1,24,18)$, $(1,0,22)$, 
$(1,13,8)$, $(1,8,13)$, $(1,25,22)$, $(0,1,26)$, $(1,18,19)$, $(1,21,28)$, $(1,19,2)$, 
$(1,11,0)$, $(1,10,19)$, $(1,13,26)$, $(1,30,6)$, $(1,23,14)$, $(1,7,11)$, $(0,1,27)$, 
$(1,27,23)$, $(1,22,12)$, $(1,14,13)$, $(1,25,0)$, $(1,21,16)$, $(1,17,1)$, $(1,28,19)$, 
$(1,20,20)$, $(1,15,9)$, $(0,1,30)$, $(1,12,6)$, $(1,13,1)$, $(1,24,22)$, $(1,4,0)$, 
$(1,3,4)$, $(1,15,29)$, $(1,11,21)$, $(1,1,27)$, $(1,22,15)$, $(1,0,1)$, $(1,19,18)$, 
$(1,24,25)$, $(1,5,10)$, $(1,21,5)$, $(1,0,8)$, $(1,18,6)$, $(1,5,3)$, $(1,19,6)$, 
$(1,18,23)$, $(1,0,3)$, $(1,15,1)$, $(1,22,8)$, $(1,28,30)$, $(1,22,30)$, $(1,0,27)$, 
$(1,4,24)$, $(1,27,4)$, $(1,16,29)$, $(1,15,10)$, $(1,0,4)$, $(1,7,29)$, $(1,6,27)$, 
$(1,6,23)$, $(1,20,11)$, $(1,0,30)$, $(1,10,3)$, $(1,3,17)$, $(1,11,26)$, $(1,4,14)$, 
$(1,0,6)$, $(1,8,10)$, $(1,17,12)$, $(1,30,25)$, $(1,23,24)$, $(1,0,16)$, $(1,16,13)$, 
$(1,25,18)$, $(1,12,8)$, $(1,6,2)$, $(1,0,13)$, $(1,25,28)$, $(1,20,22)$, $(1,27,17)$, 
$(1,13,22)$, $(1,0,24)$, $(1,1,19)$, $(1,16,19)$, $(1,14,3)$, $(1,28,25)$, $(1,0,14)$, 
$(1,11,15)$, $(1,18,5)$, $(1,1,20)$, $(1,19,17)$, $(1,0,17)$, $(1,14,20)$, $(1,10,30)$, 
$(1,17,11)$, $(1,8,21)$, $(1,0,18)$, $(1,29,14)$, $(1,12,16)$, $(1,29,12)$, $(1,11,3)$, 
$(1,0,21)$, $(1,24,16)$, $(1,21,15)$, $(1,10,13)$, $(1,25,12)$, $(1,0,25)$, $(1,30,26)$, 
$(1,23,1)$, $(1,21,1)$, $(1,7,27)$, $(1,0,26)$, $(1,9,22)$, $(1,14,2)$, $(1,7,5)$, 
$(1,12,28)$, $(1,1,3)$, $(1,18,28)$, $(1,7,10)$, $(1,19,24)$, $(1,12,21)$, $(1,16,5)$, 
$(1,29,15)$, $(1,5,18)$, $(1,2,30)$, $(1,3,26)$, $(1,1,4)$, $(1,1,29)$, $(1,24,24)$, 
$(1,8,9)$, $(1,16,3)$, $(1,4,26)$, $(1,21,4)$, $(1,11,29)$, $(1,7,1)$, $(1,4,6)$, 
$(1,1,5)$, $(1,27,22)$, $(1,20,28)$, $(1,4,12)$, $(1,2,4)$, $(1,20,29)$, $(1,10,16)$, 
$(1,28,24)$, $(1,17,5)$, $(1,16,14)$, $(1,1,6)$, $(1,30,20)$, $(1,23,25)$, $(1,3,5)$, 
$(1,5,6)$, $(1,14,24)$, $(1,13,24)$, $(1,16,2)$, $(1,16,17)$, $(1,9,30)$, $(1,1,9)$, 
$(1,5,16)$, $(1,3,14)$, $(1,24,28)$, $(1,9,19)$, $(1,9,25)$, $(1,18,27)$, $(1,14,19)$, 
$(1,19,12)$, $(1,10,10)$, $(1,1,10)$, $(1,4,27)$, $(1,13,4)$, $(1,2,29)$, $(1,18,25)$, 
$(1,12,12)$, $(1,19,9)$, $(1,2,28)$, $(1,29,16)$, $(1,20,27)$, $(1,1,11)$, $(1,12,1)$, 
$(1,1,16)$, $(1,10,23)$, $(1,21,27)$, $(1,18,17)$, $(1,17,14)$, $(1,19,23)$, $(1,3,18)$, 
$(1,13,12)$, $(1,1,12)$, $(1,21,26)$, $(1,6,11)$, $(1,9,16)$, $(1,7,28)$, $(1,5,1)$, 
$(1,5,13)$, $(1,25,3)$, $(1,6,13)$, $(1,25,20)$, $(1,1,21)$, $(1,14,10)$, $(1,2,15)$, 
$(1,1,22)$, $(1,28,11)$, $(1,11,6)$, $(1,6,26)$, $(1,8,8)$, $(1,27,9)$, $(1,7,8)$, 
$(1,1,25)$, $(1,2,18)$, $(1,30,18)$, $(1,23,21)$, $(1,15,23)$, $(1,28,15)$, $(1,12,11)$, 
$(1,22,13)$, $(1,25,2)$, $(1,27,11)$, $(1,1,30)$, $(1,8,14)$, $(1,22,26)$, $(1,17,10)$, 
$(1,29,22)$, $(1,7,13)$, $(1,25,25)$, $(1,17,9)$, $(1,10,27)$, $(1,15,3)$, $(1,2,3)$, 
$(1,22,2)$, $(1,9,18)$, $(1,3,3)$, $(1,9,9)$, $(1,30,9)$, $(1,23,6)$, $(1,11,17)$, 
$(1,15,18)$, $(1,17,28)$, $(1,2,6)$, $(1,24,6)$, $(1,27,13)$, $(1,25,19)$, $(1,29,30)$, 
$(1,4,5)$, $(1,10,4)$, $(1,6,29)$, $(1,28,13)$, $(1,25,14)$, $(1,2,7)$, $(1,5,30)$, 
$(1,24,19)$, $(1,12,18)$, $(1,5,11)$, $(1,17,7)$, $(1,17,17)$, $(1,4,9)$, $(1,8,4)$, 
$(1,12,29)$, $(1,2,8)$, $(1,19,27)$, $(1,17,2)$, $(1,24,7)$, $(1,11,8)$, $(1,20,17)$, 
$(1,30,19)$, $(1,23,13)$, $(1,25,7)$, $(1,18,3)$, $(1,2,9)$, $(1,13,15)$, $(1,13,10)$, 
$(1,6,8)$, $(1,10,24)$, $(1,11,18)$, $(1,16,24)$, $(1,20,14)$, $(1,16,20)$, $(1,21,21)$, 
$(1,2,10)$, $(1,12,13)$, $(1,25,17)$, $(1,8,1)$, $(1,13,7)$, $(1,14,28)$, $(1,11,28)$, 
$(1,13,6)$, $(1,21,30)$, $(1,29,7)$, $(1,2,11)$, $(1,28,14)$, $(1,21,25)$, $(1,13,30)$, 
$(1,12,23)$, $(1,27,30)$, $(1,21,20)$, $(1,9,28)$, $(1,10,8)$, $(1,24,8)$, $(1,2,12)$, 
$(1,17,23)$, $(1,14,8)$, $(1,17,16)$, $(1,18,20)$, $(1,6,22)$, $(1,19,3)$, $(1,21,24)$, 
$(1,3,25)$, $(1,16,22)$, $(1,2,13)$, $(1,25,8)$, $(1,11,14)$, $(1,9,13)$, $(1,25,1)$, 
$(1,15,21)$, $(1,3,22)$, $(1,17,15)$, $(1,22,1)$, $(1,10,17)$, $(1,2,17)$, $(1,7,3)$, 
$(1,7,22)$, $(1,29,5)$, $(1,24,17)$, $(1,28,23)$, $(1,18,10)$, $(1,24,23)$, $(1,19,26)$, 
$(1,27,26)$, $(1,3,8)$, $(1,7,23)$, $(1,13,23)$, $(1,4,3)$, $(1,5,4)$, $(1,4,29)$, 
$(1,20,4)$, $(1,14,29)$, $(1,12,26)$, $(1,11,19)$, $(1,3,10)$, $(1,22,11)$, $(1,24,3)$, 
$(1,30,12)$, $(1,23,12)$, $(1,13,20)$, $(1,3,16)$, $(1,8,16)$, $(1,15,25)$, $(1,8,2)$, 
$(1,3,13)$, $(1,25,21)$, $(1,10,20)$, $(1,18,15)$, $(1,29,25)$, $(1,19,14)$, $(1,4,8)$, 
$(1,12,4)$, $(1,3,29)$, $(1,14,5)$, $(1,3,21)$, $(1,6,30)$, $(1,5,15)$, $(1,15,8)$, 
$(1,15,5)$, $(1,21,12)$, $(1,22,19)$, $(1,22,5)$, $(1,11,16)$, $(1,24,10)$, $(1,3,24)$, 
$(1,18,8)$, $(1,16,26)$, $(1,19,7)$, $(1,13,11)$, $(1,7,26)$, $(1,15,13)$, $(1,25,27)$, 
$(1,7,7)$, $(1,22,9)$, $(1,3,28)$, $(1,27,7)$, $(1,6,16)$, $(1,27,5)$, $(1,9,23)$, 
$(1,11,22)$, $(1,8,7)$, $(1,28,18)$, $(1,9,27)$, $(1,28,12)$, $(1,4,1)$, $(1,17,4)$, 
$(1,21,29)$, $(1,30,23)$, $(1,23,28)$, $(1,28,27)$, $(1,8,22)$, $(1,11,20)$, $(1,5,20)$, 
$(1,27,2)$, $(1,4,17)$, $(1,22,4)$, $(1,29,29)$, $(1,21,9)$, $(1,20,18)$, $(1,21,13)$, 
$(1,25,5)$, $(1,4,25)$, $(1,15,4)$, $(1,9,29)$, $(1,4,18)$, $(1,11,4)$, $(1,22,29)$, 
$(1,5,22)$, $(1,10,26)$, $(1,14,30)$, $(1,20,10)$, $(1,7,14)$, $(1,14,18)$, $(1,15,20)$, 
$(1,5,7)$, $(1,12,25)$, $(1,28,21)$, $(1,27,24)$, $(1,10,15)$, $(1,16,7)$, $(1,27,1)$, 
$(1,14,14)$, $(1,5,9)$, $(1,28,17)$, $(1,5,21)$, $(1,7,12)$, $(1,27,15)$, $(1,28,26)$, 
$(1,13,16)$, $(1,5,24)$, $(1,22,20)$, $(1,16,11)$, $(1,13,27)$, $(1,11,5)$, $(1,28,3)$, 
$(1,29,3)$, $(1,17,19)$, $(1,15,11)$, $(1,27,19)$, $(1,6,9)$, $(1,18,26)$, $(1,10,5)$, 
$(1,19,19)$, $(1,7,9)$, $(1,21,14)$, $(1,8,12)$, $(1,8,24)$, $(1,21,17)$, $(1,16,16)$, 
$(1,6,10)$, $(1,14,25)$, $(1,22,6)$, $(1,20,7)$, $(1,9,5)$, $(1,8,11)$, $(1,10,11)$, 
$(1,20,23)$, $(1,12,7)$, $(1,15,14)$, $(1,7,18)$, $(1,28,10)$, $(1,10,1)$, $(1,9,10)$, 
$(1,13,19)$, $(1,27,14)$, $(1,17,21)$, $(1,20,26)$, $(1,24,15)$, $(1,17,6)$, $(1,8,30)$, 
$(1,18,21)$, $(1,28,6)$, $(1,19,15)$, $(1,20,3)$, $(1,11,27)$, $(1,27,28)$, $(1,18,13)$, 
$(1,25,11)$, $(1,19,28)$, $(1,10,14)$, $(1,15,7)$, $(1,19,20)$, $(1,29,6)$, $(1,22,22)$, 
$(1,18,9)$.

\subsection*{$m_{23}(2,31)\ge 663$}
$(0,0,1)$, $(1,1,13)$, $(1,20,19)$, $(1,24,19)$, $(1,4,25)$, $(1,18,13)$, $(1,5,28)$, 
$(1,26,24)$, $(1,7,16)$, $(1,15,2)$, $(0,1,0)$, $(1,16,9)$, $(1,22,14)$, $(1,25,24)$, 
$(1,9,11)$, $(1,5,7)$, $(1,25,0)$, $(1,15,25)$, $(1,1,10)$, $(1,4,27)$, $(0,1,1)$, 
$(1,14,24)$, $(0,1,13)$, $(1,6,22)$, $(1,25,22)$, $(1,7,27)$, $(1,23,9)$, $(1,20,2)$, 
$(1,16,25)$, $(1,29,4)$, $(0,1,2)$, $(1,7,30)$, $(1,18,24)$, $(1,20,30)$, $(1,17,1)$, 
$(1,25,21)$, $(1,27,22)$, $(1,18,5)$, $(1,14,23)$, $(1,21,25)$, $(0,1,3)$, $(1,8,7)$, 
$(1,19,6)$, $(1,14,0)$, $(1,7,29)$, $(1,26,0)$, $(1,20,7)$, $(1,14,11)$, $(1,24,2)$, 
$(1,16,11)$, $(0,1,5)$, $(1,27,4)$, $(1,4,28)$, $(1,29,13)$, $(1,8,20)$, $(1,23,1)$, 
$(1,9,10)$, $(0,1,14)$, $(1,5,14)$, $(1,27,17)$, $(0,1,7)$, $(1,22,26)$, $(1,30,25)$, 
$(1,26,29)$, $(1,10,2)$, $(1,16,24)$, $(1,24,20)$, $(1,7,6)$, $(1,25,3)$, $(1,20,16)$, 
$(0,1,8)$, $(1,12,8)$, $(1,2,2)$, $(1,16,10)$, $(1,12,15)$, $(1,15,14)$, $(1,1,15)$, 
$(1,3,12)$, $(1,10,19)$, $(0,1,9)$, $(0,1,10)$, $(1,21,18)$, $(1,6,23)$, $(1,12,21)$, 
$(1,28,26)$, $(0,1,12)$, $(1,18,25)$, $(1,17,11)$, $(1,0,23)$, $(1,11,24)$, $(1,21,12)$, 
$(1,18,16)$, $(1,25,10)$, $(1,22,0)$, $(1,6,14)$, $(0,1,15)$, $(1,25,19)$, $(1,11,26)$, 
$(1,10,11)$, $(1,29,17)$, $(1,13,25)$, $(1,28,2)$, $(1,16,8)$, $(1,15,24)$, $(1,1,0)$, 
$(0,1,16)$, $(1,26,27)$, $(1,24,9)$, $(1,17,15)$, $(1,24,0)$, $(1,19,23)$, $(1,6,8)$, 
$(1,0,1)$, $(1,2,11)$, $(1,26,8)$, $(0,1,18)$, $(1,17,17)$, $(1,29,12)$, $(1,5,17)$, 
$(1,5,16)$, $(1,11,5)$, $(1,22,29)$, $(1,4,26)$, $(1,4,13)$, $(1,30,13)$, $(0,1,19)$, 
$(1,3,29)$, $(1,13,21)$, $(1,15,5)$, $(1,1,21)$, $(1,29,30)$, $(1,7,19)$, $(1,10,17)$, 
$(1,20,29)$, $(1,23,12)$, $(0,1,20)$, $(1,0,5)$, $(1,15,16)$, $(1,1,28)$, $(1,22,18)$, 
$(1,0,19)$, $(1,8,30)$, $(1,9,3)$, $(1,27,5)$, $(1,12,6)$, $(0,1,21)$, $(1,13,16)$, 
$(1,23,27)$, $(1,22,9)$, $(1,14,28)$, $(1,27,10)$, $(1,11,1)$, $(1,12,14)$, $(1,9,18)$, 
$(1,17,20)$, $(0,1,23)$, $(1,29,20)$, $(1,3,15)$, $(1,19,25)$, $(1,2,12)$, $(1,28,20)$, 
$(1,30,24)$, $(1,17,22)$, $(1,21,30)$, $(1,25,30)$, $(0,1,26)$, $(1,4,6)$, $(1,10,13)$, 
$(1,28,8)$, $(1,21,27)$, $(1,30,9)$, $(1,19,27)$, $(1,5,9)$, $(1,13,22)$, $(1,13,15)$, 
$(0,1,29)$, $(1,2,21)$, $(1,0,7)$, $(1,21,4)$, $(1,13,6)$, $(1,2,8)$, $(1,13,23)$, 
$(1,27,7)$, $(1,8,17)$, $(1,0,22)$, $(0,1,30)$, $(1,30,28)$, $(1,21,1)$, $(1,27,3)$, 
$(1,6,7)$, $(1,9,16)$, $(1,12,12)$, $(1,22,30)$, $(1,28,6)$, $(1,19,7)$, $(1,0,10)$, 
$(1,26,25)$, $(1,3,14)$, $(1,30,22)$, $(1,24,24)$, $(1,30,5)$, $(1,10,20)$, $(1,25,12)$, 
$(1,14,2)$, $(1,16,22)$, $(1,0,11)$, $(1,17,12)$, $(1,27,1)$, $(1,17,0)$, $(1,30,2)$, 
$(1,16,4)$, $(1,7,0)$, $(1,12,11)$, $(1,20,21)$, $(1,21,29)$, $(1,0,13)$, $(1,24,29)$, 
$(1,8,10)$, $(1,13,29)$, $(1,7,14)$, $(1,5,1)$, $(1,20,4)$, $(1,3,27)$, $(1,26,9)$, 
$(1,4,30)$, $(1,0,16)$, $(1,20,6)$, $(1,22,5)$, $(1,9,27)$, $(1,14,9)$, $(1,23,20)$, 
$(1,21,21)$, $(1,2,15)$, $(1,30,1)$, $(1,7,28)$, $(1,0,17)$, $(1,30,17)$, $(1,14,30)$, 
$(1,8,11)$, $(1,18,15)$, $(1,10,8)$, $(1,27,30)$, $(1,29,29)$, $(1,6,18)$, $(1,18,0)$, 
$(1,0,18)$, $(1,14,18)$, $(1,25,15)$, $(1,29,6)$, $(1,17,29)$, $(1,2,3)$, $(1,3,25)$, 
$(1,11,30)$, $(1,23,15)$, $(1,13,24)$, $(1,0,20)$, $(1,27,23)$, $(1,10,27)$, $(1,4,9)$, 
$(1,27,13)$, $(1,3,23)$, $(1,15,12)$, $(1,1,3)$, $(1,11,8)$, $(1,5,19)$, $(1,0,21)$, 
$(1,22,2)$, $(1,16,16)$, $(1,5,25)$, $(1,22,21)$, $(1,26,18)$, $(1,23,24)$, $(1,4,8)$, 
$(1,24,13)$, $(1,25,16)$, $(1,0,24)$, $(1,6,3)$, $(1,23,29)$, $(1,25,4)$, $(1,10,3)$, 
$(1,15,15)$, $(1,1,22)$, $(1,8,25)$, $(1,21,19)$, $(1,12,4)$, $(1,0,25)$, $(1,25,27)$, 
$(1,17,9)$, $(1,3,24)$, $(1,23,7)$, $(1,0,26)$, $(1,19,8)$, $(1,19,26)$, $(1,15,30)$, 
$(1,1,5)$, $(1,7,10)$, $(1,30,19)$, $(1,20,14)$, $(1,13,4)$, $(1,30,23)$, $(1,0,28)$, 
$(1,13,20)$, $(1,24,22)$, $(1,26,20)$, $(1,9,17)$, $(1,6,21)$, $(1,24,10)$, $(1,10,18)$, 
$(1,8,14)$, $(1,26,5)$, $(1,0,29)$, $(1,29,19)$, $(1,30,11)$, $(1,6,10)$, $(1,15,26)$, 
$(1,1,14)$, $(1,12,23)$, $(1,23,19)$, $(1,15,0)$, $(1,1,1)$, $(1,1,2)$, $(1,16,21)$, 
$(1,3,19)$, $(1,9,8)$, $(1,30,16)$, $(1,8,15)$, $(1,22,10)$, $(1,21,17)$, $(1,17,14)$, 
$(1,15,28)$, $(1,1,6)$, $(1,6,26)$, $(1,17,19)$, $(1,23,10)$, $(1,14,12)$, $(1,13,14)$, 
$(1,10,30)$, $(1,6,27)$, $(1,11,9)$, $(1,15,3)$, $(1,1,7)$, $(1,23,2)$, $(1,16,19)$, 
$(1,19,5)$, $(1,23,22)$, $(1,17,7)$, $(1,17,8)$, $(1,8,5)$, $(1,28,18)$, $(1,15,29)$, 
$(1,1,8)$, $(1,28,15)$, $(1,14,19)$, $(1,5,3)$, $(1,8,26)$, $(1,14,20)$, $(1,6,16)$, 
$(1,4,18)$, $(1,22,13)$, $(1,15,22)$, $(1,1,9)$, $(1,25,1)$, $(1,8,19)$, $(1,7,21)$, 
$(1,19,21)$, $(1,20,25)$, $(1,19,15)$, $(1,11,3)$, $(1,29,24)$, $(1,15,27)$, $(1,1,11)$, 
$(1,15,6)$, $(1,1,19)$, $(1,14,22)$, $(1,11,19)$, $(1,3,16)$, $(1,28,0)$, $(1,22,6)$, 
$(1,27,12)$, $(1,15,21)$, $(1,1,16)$, $(1,30,14)$, $(1,22,19)$, $(1,2,7)$, $(1,7,18)$, 
$(1,2,10)$, $(1,20,3)$, $(1,10,14)$, $(1,21,7)$, $(1,15,10)$, $(1,1,17)$, $(1,21,3)$, 
$(1,4,19)$, $(1,13,13)$, $(1,21,6)$, $(1,11,2)$, $(1,16,20)$, $(1,17,30)$, $(1,5,4)$, 
$(1,15,8)$, $(1,1,27)$, $(1,9,9)$, $(1,28,19)$, $(1,17,18)$, $(1,12,27)$, $(1,7,9)$, 
$(1,8,23)$, $(1,20,28)$, $(1,30,30)$, $(1,15,17)$, $(1,1,30)$, $(1,2,28)$, $(1,2,19)$, 
$(1,6,12)$, $(1,29,8)$, $(1,9,21)$, $(1,7,4)$, $(1,29,22)$, $(1,20,1)$, $(1,15,18)$, 
$(1,2,1)$, $(1,18,22)$, $(1,19,3)$, $(1,12,2)$, $(1,16,3)$, $(1,14,21)$, $(1,18,17)$, 
$(1,27,2)$, $(1,16,26)$, $(1,23,17)$, $(1,2,4)$, $(1,22,15)$, $(1,23,5)$, $(1,18,26)$, 
$(1,12,0)$, $(1,11,27)$, $(1,19,9)$, $(1,10,4)$, $(1,18,11)$, $(1,8,16)$, $(1,2,5)$, 
$(1,4,0)$, $(1,8,13)$, $(1,10,25)$, $(1,5,18)$, $(1,19,11)$, $(1,21,24)$, $(1,5,21)$, 
$(1,4,23)$, $(1,25,13)$, $(1,2,6)$, $(1,5,6)$, $(1,18,18)$, $(1,4,1)$, $(1,19,13)$, 
$(1,4,10)$, $(1,3,13)$, $(1,18,14)$, $(1,28,29)$, $(1,5,22)$, $(1,2,13)$, $(1,29,26)$, 
$(1,25,6)$, $(1,3,28)$, $(1,8,28)$, $(1,5,8)$, $(1,17,25)$, $(1,9,26)$, $(1,9,1)$, 
$(1,4,24)$, $(1,2,14)$, $(1,11,11)$, $(1,13,0)$, $(1,14,10)$, $(1,29,5)$, $(1,2,18)$, 
$(1,10,5)$, $(1,17,2)$, $(1,16,18)$, $(1,26,26)$, $(1,27,26)$, $(1,7,12)$, $(1,24,6)$, 
$(1,20,27)$, $(1,27,9)$, $(1,2,20)$, $(1,14,29)$, $(1,28,23)$, $(1,28,4)$, $(1,23,16)$, 
$(1,10,29)$, $(1,22,16)$, $(1,7,8)$, $(1,22,12)$, $(1,20,23)$, $(1,2,27)$, $(1,21,9)$, 
$(1,9,29)$, $(1,11,29)$, $(1,14,17)$, $(1,26,28)$, $(1,12,3)$, $(1,21,10)$, $(1,24,28)$, 
$(1,9,14)$, $(1,2,30)$, $(1,27,14)$, $(1,14,16)$, $(1,19,30)$, $(1,22,23)$, $(1,8,2)$, 
$(1,16,2)$, $(1,16,27)$, $(1,10,9)$, $(1,11,10)$, $(1,3,2)$, $(1,16,23)$, $(1,13,18)$, 
$(1,29,28)$, $(1,23,4)$, $(1,12,5)$, $(1,30,4)$, $(1,9,2)$, $(1,16,6)$, $(1,28,1)$, 
$(1,3,4)$, $(1,27,20)$, $(1,19,10)$, $(1,7,7)$, $(1,11,14)$, $(1,20,26)$, $(1,11,21)$, 
$(1,8,8)$, $(1,12,7)$, $(1,3,22)$, $(1,3,5)$, $(1,29,11)$, $(1,10,22)$, $(1,9,23)$, 
$(1,18,3)$, $(1,13,27)$, $(1,8,9)$, $(1,30,0)$, $(1,28,3)$, $(1,18,28)$, $(1,4,2)$, 
$(1,16,13)$, $(1,11,12)$, $(1,6,5)$, $(1,21,5)$, $(1,26,1)$, $(1,5,15)$, $(1,12,29)$, 
$(1,24,8)$, $(1,11,22)$, $(1,4,4)$, $(1,20,13)$, $(1,9,7)$, $(1,10,15)$, $(1,26,21)$, 
$(1,6,0)$, $(1,5,30)$, $(1,24,7)$, $(1,28,28)$, $(1,7,11)$, $(1,4,16)$, $(1,21,13)$, 
$(1,22,24)$, $(1,18,4)$, $(1,11,4)$, $(1,25,18)$, $(1,5,12)$, $(1,30,27)$, $(1,18,9)$, 
$(1,29,25)$, $(1,4,21)$, $(1,17,13)$, $(1,13,17)$, $(1,12,20)$, $(1,10,7)$, $(1,30,26)$, 
$(1,5,23)$, $(1,14,15)$, $(1,17,4)$, $(1,14,7)$, $(1,6,2)$, $(1,16,7)$, $(1,27,29)$, 
$(1,21,20)$, $(1,8,27)$, $(1,29,9)$, $(1,26,7)$, $(1,29,27)$, $(1,28,9)$, $(1,24,26)$, 
$(1,6,4)$, $(1,19,12)$, $(1,17,28)$, $(1,6,25)$, $(1,10,28)$, $(1,14,6)$, $(1,9,24)$, 
$(1,19,17)$, $(1,9,5)$, $(1,25,8)$, $(1,6,29)$, $(1,20,24)$, $(1,29,23)$, $(1,19,0)$, 
$(1,25,20)$, $(1,12,18)$, $(1,27,6)$, $(1,13,11)$, $(1,11,25)$, $(1,7,22)$, $(1,7,3)$, 
$(1,9,20)$, $(1,20,0)$, $(1,18,8)$, $(1,7,20)$, $(1,13,8)$, $(1,20,11)$, $(1,12,28)$, 
$(1,18,12)$, $(1,23,18)$, $(1,7,15)$, $(1,21,15)$, $(1,20,20)$, $(1,26,2)$, $(1,16,5)$, 
$(1,8,3)$, $(1,24,23)$, $(1,29,1)$, $(1,11,23)$, $(1,17,6)$, $(1,8,12)$, $(1,9,22)$, 
$(1,22,28)$, $(1,19,29)$, $(1,30,15)$, $(1,27,16)$, $(1,10,1)$, $(1,12,26)$, $(1,29,15)$, 
$(1,18,30)$, $(1,30,21)$, $(1,23,6)$, $(1,26,23)$, $(1,26,10)$, $(1,18,6)$, $(1,24,14)$, 
$(1,24,12)$, $(1,13,1)$, $(1,26,12)$, $(1,25,25)$, $(1,23,23)$, $(1,24,16)$, $(1,27,24)$, 
$(1,22,25)$, $(1,13,3)$, $(1,25,17)$, $(1,19,24)$, $(1,18,7)$.

\subsection*{$m_{24}(2,31)\ge 698$}
$(0,0,1)$, $(1,29,8)$, $(1,7,4)$, $(1,11,11)$, $(1,5,30)$, $(0,1,1)$, $(1,19,8)$, 
$(1,25,14)$, $(1,9,29)$, $(1,15,4)$, $(0,1,2)$, $(1,5,8)$, $(1,29,30)$, $(1,8,7)$, 
$(1,10,17)$, $(0,1,3)$, $(1,27,8)$, $(1,22,2)$, $(1,5,3)$, $(1,29,11)$, $(0,1,4)$, 
$(1,1,8)$, $(1,28,26)$, $(0,1,22)$, $(1,23,8)$, $(0,1,5)$, $(1,21,8)$, $(1,11,20)$, 
$(1,17,19)$, $(1,7,0)$, $(0,1,7)$, $(1,2,8)$, $(1,26,18)$, $(1,13,24)$, $(1,2,13)$, 
$(0,1,8)$, $(1,18,8)$, $(1,3,19)$, $(1,28,13)$, $(1,14,19)$, $(0,1,9)$, $(1,12,8)$, 
$(1,17,13)$, $(1,6,25)$, $(1,13,3)$, $(0,1,10)$, $(1,4,8)$, $(1,23,6)$, $(1,12,2)$, 
$(1,6,15)$, $(0,1,11)$, $(1,30,8)$, $(1,16,9)$, $(1,3,21)$, $(1,20,22)$, $(0,1,13)$, 
$(1,10,8)$, $(1,18,17)$, $(1,22,5)$, $(1,1,28)$, $(0,1,14)$, $(1,16,8)$, $(1,0,7)$, 
$(1,24,18)$, $(1,19,6)$, $(0,1,15)$, $(1,3,8)$, $(1,30,3)$, $(1,20,23)$, $(1,16,20)$, 
$(0,1,17)$, $(1,24,8)$, $(1,2,15)$, $(1,21,14)$, $(1,22,23)$, $(0,1,18)$, $(1,11,8)$, 
$(1,24,10)$, $(1,7,16)$, $(1,9,1)$, $(0,1,21)$, $(1,28,8)$, $(1,13,28)$, $(1,15,6)$, 
$(1,4,14)$, $(0,1,26)$, $(1,13,8)$, $(1,21,29)$, $(1,0,17)$, $(1,3,29)$, $(0,1,27)$, 
$(1,6,8)$, $(1,12,24)$, $(1,18,10)$, $(1,18,21)$, $(0,1,28)$, $(1,26,8)$, $(1,6,0)$, 
$(1,19,1)$, $(1,12,18)$, $(0,1,29)$, $(1,0,8)$, $(1,15,5)$, $(1,10,20)$, $(1,0,12)$, 
$(0,1,30)$, $(1,22,8)$, $(1,14,1)$, $(1,16,28)$, $(1,26,25)$, $(1,0,0)$, $(1,14,7)$, 
$(1,25,7)$, $(1,17,2)$, $(1,18,4)$, $(1,0,1)$, $(1,7,21)$, $(1,30,13)$, $(1,18,16)$, 
$(1,11,18)$, $(1,0,2)$, $(1,11,13)$, $(1,15,26)$, $(1,7,17)$, $(1,12,16)$, $(1,0,3)$, 
$(1,9,17)$, $(1,20,1)$, $(1,6,3)$, $(1,1,7)$, $(1,0,5)$, $(1,21,24)$, $(1,17,16)$, 
$(1,1,26)$, $(1,19,2)$, $(1,0,6)$, $(1,30,6)$, $(1,1,3)$, $(1,22,10)$, $(1,26,19)$, 
$(1,0,10)$, $(1,1,2)$, $(1,2,29)$, $(1,28,1)$, $(1,15,10)$, $(1,0,14)$, $(1,26,14)$, 
$(1,16,21)$, $(1,5,20)$, $(1,29,13)$, $(1,0,15)$, $(1,27,12)$, $(1,5,14)$, $(1,29,15)$, 
$(1,30,11)$, $(1,0,16)$, $(1,12,11)$, $(1,3,24)$, $(1,27,18)$, $(1,23,25)$, $(1,0,19)$, 
$(1,6,23)$, $(1,9,25)$, $(1,23,24)$, $(1,4,1)$, $(1,0,21)$, $(1,18,30)$, $(1,8,30)$, 
$(1,8,0)$, $(1,28,15)$, $(1,0,22)$, $(1,16,3)$, $(1,14,0)$, $(1,2,9)$, $(1,2,5)$, 
$(1,0,24)$, $(1,13,9)$, $(1,12,10)$, $(1,14,22)$, $(1,14,12)$, $(1,0,25)$, $(1,19,28)$, 
$(1,23,17)$, $(1,19,30)$, $(1,8,24)$, $(1,0,26)$, $(1,2,0)$, $(1,4,19)$, $(1,30,29)$, 
$(1,13,14)$, $(1,0,28)$, $(1,17,1)$, $(1,13,5)$, $(1,9,14)$, $(1,6,28)$, $(1,0,29)$, 
$(1,25,16)$, $(1,28,23)$, $(1,26,4)$, $(1,3,3)$, $(1,1,0)$, $(1,15,12)$, $(1,15,25)$, 
$(1,3,13)$, $(1,9,15)$, $(1,1,1)$, $(1,3,11)$, $(1,25,28)$, $(1,27,7)$, $(1,14,4)$, 
$(1,1,4)$, $(1,23,23)$, $(1,23,15)$, $(1,25,23)$, $(1,13,0)$, $(1,1,9)$, $(1,10,9)$, 
$(1,5,22)$, $(1,29,22)$, $(1,4,26)$, $(1,1,12)$, $(1,20,15)$, $(1,22,24)$, $(1,6,20)$, 
$(1,28,29)$, $(1,1,14)$, $(1,1,16)$, $(1,17,7)$, $(1,30,14)$, $(1,10,19)$, $(1,1,15)$, 
$(1,6,19)$, $(1,20,11)$, $(1,2,21)$, $(1,27,25)$, $(1,1,17)$, $(1,30,21)$, $(1,24,6)$, 
$(1,5,28)$, $(1,29,2)$, $(1,1,18)$, $(1,27,13)$, $(1,13,12)$, $(1,28,30)$, $(1,8,11)$, 
$(1,1,19)$, $(1,13,17)$, $(1,4,0)$, $(1,22,16)$, $(1,19,24)$, $(1,1,20)$, $(1,11,22)$, 
$(1,12,21)$, $(1,7,12)$, $(1,23,9)$, $(1,1,21)$, $(1,4,24)$, $(1,14,3)$, $(1,26,15)$, 
$(1,24,13)$, $(1,1,22)$, $(1,5,6)$, $(1,29,23)$, $(1,19,9)$, $(1,16,12)$, $(1,1,23)$, 
$(1,25,18)$, $(1,11,30)$, $(1,8,4)$, $(1,7,7)$, $(1,1,24)$, $(1,24,5)$, $(1,3,9)$, 
$(1,13,26)$, $(1,3,22)$, $(1,1,25)$, $(1,12,4)$, $(1,16,16)$, $(1,4,5)$, $(1,25,17)$, 
$(1,2,1)$, $(1,20,7)$, $(1,3,1)$, $(1,22,12)$, $(1,17,21)$, $(1,2,3)$, $(1,12,13)$, 
$(1,5,18)$, $(1,29,7)$, $(1,28,5)$, $(1,2,4)$, $(1,17,17)$, $(1,24,9)$, $(1,25,1)$, 
$(1,26,22)$, $(1,2,7)$, $(1,11,6)$, $(1,6,11)$, $(1,7,5)$, $(1,27,29)$, $(1,2,10)$, 
$(1,24,4)$, $(1,26,26)$, $(1,17,20)$, $(1,25,15)$, $(1,2,14)$, $(1,13,20)$, $(1,23,16)$, 
$(1,10,25)$, $(1,22,25)$, $(1,2,16)$, $(1,3,12)$, $(1,12,0)$, $(1,12,28)$, $(1,10,3)$, 
$(1,2,17)$, $(1,25,11)$, $(1,10,14)$, $(1,27,4)$, $(1,6,6)$, $(1,2,19)$, $(1,16,10)$, 
$(1,19,13)$, $(1,19,23)$, $(1,3,16)$, $(1,2,20)$, $(1,9,23)$, $(1,14,17)$, $(1,5,2)$, 
$(1,29,12)$, $(1,2,24)$, $(1,5,26)$, $(1,29,5)$, $(1,18,6)$, $(1,23,1)$, $(1,2,25)$, 
$(1,7,9)$, $(1,21,30)$, $(1,8,22)$, $(1,11,10)$, $(1,2,26)$, $(1,23,28)$, $(1,9,21)$, 
$(1,28,21)$, $(1,21,18)$, $(1,2,28)$, $(1,30,15)$, $(1,11,7)$, $(1,23,29)$, $(1,7,13)$, 
$(1,3,0)$, $(1,3,14)$, $(1,18,1)$, $(1,9,26)$, $(1,22,6)$, $(1,3,2)$, $(1,16,11)$, 
$(1,11,17)$, $(1,27,15)$, $(1,7,29)$, $(1,3,4)$, $(1,10,10)$, $(1,6,24)$, $(1,25,30)$, 
$(1,8,13)$, $(1,3,5)$, $(1,26,23)$, $(1,4,2)$, $(1,26,7)$, $(1,9,28)$, $(1,3,10)$, 
$(1,4,9)$, $(1,28,18)$, $(1,24,22)$, $(1,30,2)$, $(1,3,15)$, $(1,15,16)$, $(1,30,9)$, 
$(1,17,28)$, $(1,21,22)$, $(1,3,20)$, $(1,21,17)$, $(1,21,3)$, $(1,18,5)$, $(1,4,15)$, 
$(1,3,26)$, $(1,5,4)$, $(1,29,29)$, $(1,12,19)$, $(1,25,20)$, $(1,4,3)$, $(1,6,21)$, 
$(1,23,10)$, $(1,15,14)$, $(1,24,3)$, $(1,4,4)$, $(1,13,13)$, $(1,25,3)$, $(1,25,6)$, 
$(1,21,4)$, $(1,4,7)$, $(1,15,24)$, $(1,11,21)$, $(1,16,7)$, $(1,7,19)$, $(1,4,10)$, 
$(1,20,5)$, $(1,17,0)$, $(1,10,18)$, $(1,4,20)$, $(1,4,18)$, $(1,7,11)$, $(1,14,26)$, 
$(1,27,23)$, $(1,11,28)$, $(1,4,21)$, $(1,12,23)$, $(1,22,29)$, $(1,20,10)$, $(1,12,7)$, 
$(1,4,22)$, $(1,18,25)$, $(1,5,11)$, $(1,29,9)$, $(1,27,2)$, $(1,4,28)$, $(1,5,0)$, 
$(1,29,20)$, $(1,26,30)$, $(1,8,29)$, $(1,4,29)$, $(1,19,15)$, $(1,4,30)$, $(1,8,1)$, 
$(1,30,1)$, $(1,5,1)$, $(1,29,14)$, $(1,20,12)$, $(1,21,13)$, $(1,26,10)$, $(1,5,7)$, 
$(1,29,25)$, $(1,16,0)$, $(1,24,26)$, $(1,28,14)$, $(1,5,10)$, $(1,29,26)$, $(1,21,15)$, 
$(1,19,25)$, $(1,9,7)$, $(1,5,12)$, $(1,29,28)$, $(1,11,16)$, $(1,18,0)$, $(1,7,3)$, 
$(1,5,15)$, $(1,29,1)$, $(1,10,13)$, $(1,12,5)$, $(1,15,19)$, $(1,5,17)$, $(1,29,19)$, 
$(1,14,25)$, $(1,6,10)$, $(1,30,18)$, $(1,5,19)$, $(1,29,21)$, $(1,17,3)$, $(1,28,2)$, 
$(1,17,23)$, $(1,5,23)$, $(1,29,24)$, $(1,15,28)$, $(1,14,24)$, $(1,20,29)$, $(1,5,24)$, 
$(1,29,6)$, $(1,22,18)$, $(1,15,18)$, $(1,21,0)$, $(1,6,7)$, $(1,18,22)$, $(1,26,20)$, 
$(1,19,18)$, $(1,18,9)$, $(1,6,12)$, $(1,10,4)$, $(1,19,19)$, $(1,26,28)$, $(1,25,12)$, 
$(1,6,16)$, $(1,24,20)$, $(1,14,5)$, $(1,30,16)$, $(1,15,21)$, $(1,6,17)$, $(1,11,14)$, 
$(1,11,9)$, $(1,7,23)$, $(1,7,22)$, $(1,6,18)$, $(1,17,12)$, $(1,16,23)$, $(1,21,12)$, 
$(1,13,29)$, $(1,6,22)$, $(1,15,23)$, $(1,28,7)$, $(1,22,9)$, $(1,24,16)$, $(1,6,26)$, 
$(1,6,29)$, $(1,22,0)$, $(1,30,25)$, $(1,14,2)$, $(1,28,0)$, $(1,7,6)$, $(1,17,10)$, 
$(1,17,29)$, $(1,11,24)$, $(1,22,17)$, $(1,7,14)$, $(1,15,0)$, $(1,26,5)$, $(1,11,26)$, 
$(1,20,24)$, $(1,7,18)$, $(1,25,19)$, $(1,19,3)$, $(1,11,4)$, $(1,14,14)$, $(1,7,28)$, 
$(1,22,4)$, $(1,9,9)$, $(1,11,25)$, $(1,25,22)$, $(1,8,2)$, $(1,19,16)$, $(1,22,3)$, 
$(1,10,26)$, $(1,9,30)$, $(1,8,3)$, $(1,23,19)$, $(1,18,26)$, $(1,16,29)$, $(1,16,30)$, 
$(1,8,6)$, $(1,15,13)$, $(1,20,30)$, $(1,8,25)$, $(1,17,30)$, $(1,8,9)$, $(1,9,24)$, 
$(1,16,22)$, $(1,17,14)$, $(1,12,30)$, $(1,8,10)$, $(1,16,6)$, $(1,13,16)$, $(1,14,28)$, 
$(1,24,30)$, $(1,8,12)$, $(1,25,5)$, $(1,23,5)$, $(1,28,4)$, $(1,30,30)$, $(1,8,17)$, 
$(1,14,20)$, $(1,15,20)$, $(1,22,1)$, $(1,23,30)$, $(1,8,21)$, $(1,26,29)$, $(1,26,11)$, 
$(1,21,16)$, $(1,27,30)$, $(1,8,23)$, $(1,18,23)$, $(1,12,14)$, $(1,30,5)$, $(1,22,30)$, 
$(1,9,2)$, $(1,9,20)$, $(1,10,7)$, $(1,16,13)$, $(1,28,24)$, $(1,9,4)$, $(1,12,12)$, 
$(1,9,6)$, $(1,25,29)$, $(1,24,21)$, $(1,9,5)$, $(1,24,11)$, $(1,22,19)$, $(1,27,5)$, 
$(1,20,18)$, $(1,9,11)$, $(1,27,3)$, $(1,27,24)$, $(1,30,0)$, $(1,21,11)$, $(1,9,13)$, 
$(1,17,9)$, $(1,14,11)$, $(1,23,22)$, $(1,27,0)$, $(1,9,18)$, $(1,26,16)$, $(1,20,17)$, 
$(1,13,18)$, $(1,16,15)$, $(1,10,1)$, $(1,17,15)$, $(1,23,3)$, $(1,16,18)$, $(1,10,2)$, 
$(1,10,11)$, $(1,28,6)$, $(1,27,26)$, $(1,26,1)$, $(1,14,23)$, $(1,10,21)$, $(1,15,11)$, 
$(1,13,23)$, $(1,27,21)$, $(1,19,26)$, $(1,10,22)$, $(1,24,29)$, $(1,10,29)$, $(1,14,9)$, 
$(1,18,13)$, $(1,10,24)$, $(1,21,23)$, $(1,24,1)$, $(1,21,25)$, $(1,28,19)$, $(1,12,6)$, 
$(1,19,7)$, $(1,21,20)$, $(1,30,28)$, $(1,28,20)$, $(1,12,20)$, $(1,27,14)$, $(1,17,11)$, 
$(1,19,10)$, $(1,27,11)$, $(1,12,26)$, $(1,13,25)$, $(1,30,17)$, $(1,26,13)$, $(1,30,7)$, 
$(1,13,1)$, $(1,18,2)$, $(1,21,9)$, $(1,20,14)$, $(1,19,22)$, $(1,13,4)$, $(1,20,20)$, 
$(1,20,2)$, $(1,25,26)$, $(1,14,16)$, $(1,13,6)$, $(1,26,12)$, $(1,18,19)$, $(1,15,2)$, 
$(1,23,2)$, $(1,13,11)$, $(1,24,25)$, $(1,27,20)$, $(1,18,3)$, $(1,24,28)$, $(1,15,9)$, 
$(1,19,29)$, $(1,30,22)$, $(1,23,18)$, $(1,20,3)$, $(1,15,15)$, $(1,16,26)$, $(1,25,10)$, 
$(1,22,20)$, $(1,16,4)$, $(1,16,24)$, $(1,23,11)$, $(1,30,12)$, $(1,27,1)$, $(1,24,17)$, 
$(1,17,22)$, $(1,20,16)$, $(1,26,6)$, $(1,28,9)$, $(1,23,0)$, $(1,18,7)$, $(1,27,16)$, 
$(1,21,10)$, $(1,23,12)$, $(1,22,11)$, $(1,22,13)$, $(1,28,17)$.

\subsection*{$m_{25}(2,31)\ge 733$}
$(0,0,1)$, $(1,10,28)$, $(1,4,29)$, $(1,22,24)$, $(1,19,19)$, $(1,30,5)$, $(1,13,7)$, 
$(1,8,14)$, $(0,1,0)$, $(1,0,17)$, $(1,18,29)$, $(1,9,10)$, $(1,1,13)$, $(1,7,9)$, 
$(0,1,15)$, $(1,15,18)$, $(0,1,1)$, $(1,13,22)$, $(1,19,29)$, $(1,4,7)$, $(1,28,22)$, 
$(1,18,30)$, $(1,19,4)$, $(1,23,27)$, $(0,1,3)$, $(1,16,16)$, $(1,6,29)$, $(1,6,2)$, 
$(1,29,12)$, $(1,28,4)$, $(1,15,6)$, $(1,17,28)$, $(0,1,4)$, $(1,22,4)$, $(1,9,29)$, 
$(1,15,26)$, $(1,15,28)$, $(1,14,28)$, $(1,30,14)$, $(1,0,5)$, $(0,1,5)$, $(0,1,29)$, 
$(1,25,29)$, $(1,20,29)$, $(1,7,15)$, $(1,12,27)$, $(1,8,25)$, $(1,10,24)$, $(0,1,6)$, 
$(1,29,21)$, $(1,17,29)$, $(1,26,14)$, $(1,0,23)$, $(1,5,8)$, $(1,14,22)$, $(1,2,15)$, 
$(0,1,7)$, $(1,4,9)$, $(1,7,29)$, $(1,23,6)$, $(1,23,10)$, $(1,6,24)$, $(1,5,11)$, 
$(1,13,8)$, $(0,1,9)$, $(1,7,3)$, $(1,8,29)$, $(1,10,23)$, $(1,3,24)$, $(1,25,18)$, 
$(1,24,17)$, $(1,29,26)$, $(0,1,10)$, $(1,20,8)$, $(1,15,29)$, $(1,8,28)$, $(1,10,16)$, 
$(1,13,12)$, $(1,26,16)$, $(1,4,25)$, $(0,1,11)$, $(1,8,1)$, $(1,10,29)$, $(1,30,4)$, 
$(1,6,25)$, $(1,20,0)$, $(1,22,18)$, $(1,21,17)$, $(0,1,12)$, $(1,26,27)$, $(1,12,29)$, 
$(1,1,30)$, $(1,8,5)$, $(1,10,26)$, $(1,29,30)$, $(1,9,19)$, $(0,1,18)$, $(1,21,6)$, 
$(1,27,29)$, $(1,2,12)$, $(1,17,8)$, $(1,3,7)$, $(1,3,12)$, $(1,18,2)$, $(0,1,19)$, 
$(1,18,12)$, $(1,29,29)$, $(1,19,16)$, $(1,16,18)$, $(1,23,17)$, $(1,28,15)$, $(1,19,7)$, 
$(0,1,21)$, $(1,17,14)$, $(1,0,29)$, $(1,14,13)$, $(0,1,24)$, $(1,11,26)$, $(1,22,29)$, 
$(1,17,21)$, $(1,21,30)$, $(1,29,20)$, $(1,21,3)$, $(1,30,0)$, $(0,1,25)$, $(1,28,23)$, 
$(1,14,29)$, $(1,21,11)$, $(1,5,4)$, $(1,11,11)$, $(1,12,23)$, $(1,20,12)$, $(0,1,26)$, 
$(1,5,7)$, $(1,30,29)$, $(1,18,3)$, $(1,14,7)$, $(1,22,1)$, $(1,18,20)$, $(1,24,1)$, 
$(0,1,27)$, $(1,1,15)$, $(1,2,29)$, $(1,24,19)$, $(1,20,9)$, $(1,27,19)$, $(1,27,0)$, 
$(1,1,10)$, $(0,1,28)$, $(1,27,25)$, $(1,26,29)$, $(1,29,22)$, $(1,11,6)$, $(1,21,16)$, 
$(1,11,8)$, $(1,25,6)$, $(0,1,30)$, $(1,12,24)$, $(1,21,29)$, $(1,28,9)$, $(1,30,2)$, 
$(1,2,22)$, $(1,7,10)$, $(1,22,22)$, $(1,0,0)$, $(1,28,7)$, $(1,5,19)$, $(1,18,16)$, 
$(1,28,11)$, $(1,3,3)$, $(1,1,3)$, $(1,29,14)$, $(1,0,1)$, $(1,17,25)$, $(1,15,12)$, 
$(1,24,28)$, $(1,27,5)$, $(1,11,30)$, $(1,5,25)$, $(1,9,14)$, $(1,0,2)$, $(1,9,24)$, 
$(1,20,24)$, $(1,29,7)$, $(1,20,25)$, $(1,18,11)$, $(1,7,5)$, $(1,7,14)$, $(1,0,3)$, 
$(1,24,22)$, $(1,23,25)$, $(1,7,25)$, $(1,3,16)$, $(1,26,7)$, $(1,2,24)$, $(1,18,14)$, 
$(1,0,4)$, $(1,7,16)$, $(1,25,5)$, $(1,28,5)$, $(1,25,24)$, $(1,9,0)$, $(1,9,16)$, 
$(1,14,14)$, $(1,0,7)$, $(1,27,3)$, $(1,18,13)$, $(1,16,12)$, $(1,30,23)$, $(1,12,14)$, 
$(1,0,13)$, $(1,22,14)$, $(1,0,8)$, $(1,11,1)$, $(1,29,27)$, $(1,6,23)$, $(1,8,15)$, 
$(1,10,15)$, $(1,23,0)$, $(1,27,14)$, $(1,0,9)$, $(1,3,0)$, $(1,7,30)$, $(1,15,10)$, 
$(1,5,28)$, $(1,23,24)$, $(1,8,26)$, $(1,10,14)$, $(1,0,10)$, $(1,23,18)$, $(1,30,17)$, 
$(1,20,20)$, $(1,23,12)$, $(1,22,9)$, $(1,11,27)$, $(1,28,14)$, $(1,0,11)$, $(1,26,30)$, 
$(1,28,6)$, $(1,26,1)$, $(1,13,14)$, $(1,0,20)$, $(1,4,4)$, $(1,21,14)$, $(1,0,15)$, 
$(1,1,23)$, $(1,26,26)$, $(1,8,27)$, $(1,10,27)$, $(1,2,19)$, $(1,28,12)$, $(1,16,14)$, 
$(1,0,16)$, $(1,30,15)$, $(1,11,21)$, $(1,30,9)$, $(1,22,6)$, $(1,5,2)$, $(1,22,10)$, 
$(1,24,14)$, $(1,0,22)$, $(1,12,5)$, $(1,17,23)$, $(1,13,6)$, $(1,9,21)$, $(1,4,18)$, 
$(1,12,17)$, $(1,11,14)$, $(1,0,24)$, $(1,2,27)$, $(1,16,2)$, $(1,19,18)$, $(1,29,17)$, 
$(1,7,1)$, $(1,24,21)$, $(1,19,14)$, $(1,0,25)$, $(1,21,10)$, $(1,19,3)$, $(1,12,4)$, 
$(1,24,18)$, $(1,6,17)$, $(1,19,9)$, $(1,6,14)$, $(1,0,26)$, $(1,19,2)$, $(1,24,15)$, 
$(1,14,8)$, $(1,12,8)$, $(1,19,26)$, $(1,21,20)$, $(1,20,14)$, $(1,0,27)$, $(1,20,6)$, 
$(1,3,8)$, $(1,5,21)$, $(1,6,3)$, $(1,17,27)$, $(1,25,11)$, $(1,25,14)$, $(1,0,30)$, 
$(1,25,26)$, $(1,6,9)$, $(1,21,22)$, $(1,16,1)$, $(1,20,10)$, $(1,20,30)$, $(1,2,14)$, 
$(1,1,0)$, $(1,14,12)$, $(1,20,21)$, $(1,25,0)$, $(1,13,19)$, $(1,13,18)$, $(1,4,17)$, 
$(1,4,30)$, $(1,1,1)$, $(1,11,24)$, $(1,11,10)$, $(1,7,20)$, $(1,28,19)$, $(1,6,27)$, 
$(1,7,24)$, $(1,17,9)$, $(1,1,2)$, $(1,15,8)$, $(1,6,28)$, $(1,11,19)$, $(1,26,19)$, 
$(1,14,30)$, $(1,24,12)$, $(1,1,11)$, $(1,1,5)$, $(1,18,27)$, $(1,30,16)$, $(1,9,4)$, 
$(1,17,19)$, $(1,22,2)$, $(1,6,1)$, $(1,2,7)$, $(1,1,7)$, $(1,13,16)$, $(1,18,22)$, 
$(1,4,13)$, $(1,6,19)$, $(1,24,26)$, $(1,11,23)$, $(1,16,13)$, $(1,1,9)$, $(1,19,23)$, 
$(1,25,3)$, $(1,16,10)$, $(1,15,19)$, $(1,19,28)$, $(1,17,6)$, $(1,29,23)$, $(1,1,12)$, 
$(1,21,15)$, $(1,26,18)$, $(1,19,17)$, $(1,3,19)$, $(1,12,6)$, $(1,25,4)$, $(1,19,1)$, 
$(1,1,17)$, $(1,22,11)$, $(1,14,24)$, $(1,26,23)$, $(1,23,19)$, $(1,17,4)$, $(1,30,26)$, 
$(1,7,18)$, $(1,1,18)$, $(1,5,17)$, $(1,12,25)$, $(1,5,5)$, $(1,16,19)$, $(1,3,22)$, 
$(1,9,8)$, $(1,28,27)$, $(1,1,20)$, $(1,6,13)$, $(1,23,4)$, $(1,3,21)$, $(1,8,19)$, 
$(1,10,13)$, $(1,15,22)$, $(1,20,28)$, $(1,1,21)$, $(1,12,20)$, $(1,9,11)$, $(1,27,15)$, 
$(1,21,19)$, $(1,25,7)$, $(1,23,20)$, $(1,18,5)$, $(1,1,22)$, $(1,17,0)$, $(1,5,13)$, 
$(1,17,2)$, $(1,7,19)$, $(1,5,15)$, $(1,20,13)$, $(1,21,24)$, $(1,1,25)$, $(1,28,18)$, 
$(1,28,17)$, $(1,8,12)$, $(1,10,19)$, $(1,21,21)$, $(1,2,2)$, $(1,14,21)$, $(1,1,27)$, 
$(1,4,21)$, $(1,29,1)$, $(1,30,22)$, $(1,29,19)$, $(1,15,11)$, $(1,22,28)$, $(1,13,25)$, 
$(1,2,4)$, $(1,26,0)$, $(1,24,4)$, $(1,8,7)$, $(1,10,21)$, $(1,18,6)$, $(1,11,15)$, 
$(1,17,26)$, $(1,2,5)$, $(1,9,22)$, $(1,6,11)$, $(1,9,6)$, $(1,30,20)$, $(1,3,15)$, 
$(1,29,5)$, $(1,22,12)$, $(1,2,6)$, $(1,2,11)$, $(1,23,13)$, $(1,14,1)$, $(1,16,30)$, 
$(1,14,27)$, $(1,27,13)$, $(1,24,25)$, $(1,2,8)$, $(1,13,15)$, $(1,28,30)$, $(1,4,11)$, 
$(1,11,7)$, $(1,6,7)$, $(1,18,18)$, $(1,18,17)$, $(1,2,9)$, $(1,29,18)$, $(1,26,17)$, 
$(1,5,10)$, $(1,2,16)$, $(1,20,26)$, $(1,13,10)$, $(1,11,4)$, $(1,4,12)$, $(1,23,3)$, 
$(1,2,20)$, $(1,12,9)$, $(1,2,17)$, $(1,15,27)$, $(1,11,28)$, $(1,15,0)$, $(1,25,28)$, 
$(1,12,22)$, $(1,14,3)$, $(1,11,18)$, $(1,2,21)$, $(1,22,7)$, $(1,15,23)$, $(1,18,28)$, 
$(1,3,26)$, $(1,20,11)$, $(1,30,1)$, $(1,25,16)$, $(1,2,23)$, $(1,19,20)$, $(1,7,2)$, 
$(1,26,20)$, $(1,17,16)$, $(1,19,24)$, $(1,15,30)$, $(1,3,28)$, $(1,2,25)$, $(1,11,3)$, 
$(1,9,15)$, $(1,13,2)$, $(1,20,5)$, $(1,4,2)$, $(1,16,26)$, $(1,9,5)$, $(1,2,26)$, 
$(1,25,25)$, $(1,27,8)$, $(1,23,23)$, $(1,24,11)$, $(1,28,0)$, $(1,17,22)$, $(1,28,20)$, 
$(1,3,2)$, $(1,17,13)$, $(1,11,0)$, $(1,23,11)$, $(1,29,16)$, $(1,29,13)$, $(1,27,26)$, 
$(1,12,1)$, $(1,3,4)$, $(1,14,15)$, $(1,5,16)$, $(1,3,23)$, $(1,30,25)$, $(1,23,30)$, 
$(1,12,30)$, $(1,22,15)$, $(1,3,5)$, $(1,6,10)$, $(1,14,23)$, $(1,9,7)$, $(1,11,9)$, 
$(1,18,8)$, $(1,21,9)$, $(1,14,10)$, $(1,3,6)$, $(1,19,22)$, $(1,26,22)$, $(1,27,21)$, 
$(1,14,5)$, $(1,19,0)$, $(1,20,1)$, $(1,9,3)$, $(1,3,9)$, $(1,13,26)$, $(1,17,15)$, 
$(1,4,10)$, $(1,13,27)$, $(1,15,1)$, $(1,4,28)$, $(1,26,2)$, $(1,3,10)$, $(1,4,1)$, 
$(1,8,8)$, $(1,10,25)$, $(1,17,1)$, $(1,14,9)$, $(1,28,3)$, $(1,13,21)$, $(1,3,13)$, 
$(1,28,16)$, $(1,27,9)$, $(1,16,9)$, $(1,8,13)$, $(1,10,10)$, $(1,25,10)$, $(1,21,26)$, 
$(1,3,18)$, $(1,11,17)$, $(1,13,5)$, $(1,30,13)$, $(1,4,8)$, $(1,9,18)$, $(1,22,17)$, 
$(1,25,13)$, $(1,3,27)$, $(1,23,9)$, $(1,20,7)$, $(1,14,4)$, $(1,18,10)$, $(1,26,6)$, 
$(1,16,0)$, $(1,29,0)$, $(1,3,30)$, $(1,16,24)$, $(1,7,21)$, $(1,11,12)$, $(1,4,0)$, 
$(1,18,15)$, $(1,27,24)$, $(1,28,2)$, $(1,21,23)$, $(1,15,21)$, $(1,13,28)$, $(1,7,7)$, 
$(1,4,5)$, $(1,20,4)$, $(1,28,13)$, $(1,9,20)$, $(1,29,8)$, $(1,24,30)$, $(1,13,13)$, 
$(1,29,6)$, $(1,4,16)$, $(1,15,16)$, $(1,21,28)$, $(1,24,27)$, $(1,5,22)$, $(1,21,27)$, 
$(1,13,30)$, $(1,27,23)$, $(1,4,20)$, $(1,19,25)$, $(1,13,23)$, $(1,6,0)$, $(1,4,24)$, 
$(1,24,13)$, $(1,25,15)$, $(1,8,3)$, $(1,10,1)$, $(1,7,13)$, $(1,13,11)$, $(1,21,12)$, 
$(1,4,26)$, $(1,27,12)$, $(1,5,18)$, $(1,7,17)$, $(1,27,4)$, $(1,29,4)$, $(1,13,24)$, 
$(1,14,25)$, $(1,5,0)$, $(1,15,5)$, $(1,26,28)$, $(1,21,18)$, $(1,9,17)$, $(1,16,7)$, 
$(1,16,11)$, $(1,18,26)$, $(1,5,1)$, $(1,28,1)$, $(1,22,0)$, $(1,12,26)$, $(1,16,6)$, 
$(1,12,12)$, $(1,15,20)$, $(1,30,28)$, $(1,5,3)$, $(1,19,30)$, $(1,11,16)$, $(1,8,2)$, 
$(1,10,11)$, $(1,19,11)$, $(1,8,21)$, $(1,10,4)$, $(1,5,6)$, $(1,24,7)$, $(1,9,2)$, 
$(1,18,0)$, $(1,8,23)$, $(1,10,30)$, $(1,6,8)$, $(1,16,5)$, $(1,5,12)$, $(1,8,0)$, 
$(1,10,9)$, $(1,24,5)$, $(1,24,20)$, $(1,26,10)$, $(1,9,12)$, $(1,25,22)$, $(1,5,23)$, 
$(1,27,18)$, $(1,20,17)$, $(1,6,21)$, $(1,26,8)$, $(1,20,2)$, $(1,12,16)$, $(1,23,1)$, 
$(1,6,16)$, $(1,12,10)$, $(1,12,0)$, $(1,26,21)$, $(1,7,11)$, $(1,16,25)$, $(1,24,9)$, 
$(1,17,24)$, $(1,6,18)$, $(1,25,17)$, $(1,14,16)$, $(1,17,20)$, $(1,8,22)$, $(1,10,7)$, 
$(1,26,15)$, $(1,18,1)$, $(1,6,20)$, $(1,27,30)$, $(1,7,22)$, $(1,30,18)$, $(1,16,17)$, 
$(1,26,24)$, $(1,22,3)$, $(1,27,11)$, $(1,6,30)$, $(1,30,3)$, $(1,15,24)$, $(1,27,28)$, 
$(1,22,21)$, $(1,23,15)$, $(1,25,12)$, $(1,7,6)$, $(1,8,4)$, $(1,10,20)$, $(1,16,27)$, 
$(1,17,3)$, $(1,23,5)$, $(1,21,8)$, $(1,17,12)$, $(1,9,25)$, $(1,8,10)$, $(1,10,8)$, 
$(1,27,20)$, $(1,15,7)$, $(1,12,28)$, $(1,18,9)$, $(1,15,25)$, $(1,29,28)$, $(1,8,11)$, 
$(1,10,22)$, $(1,22,26)$, $(1,30,8)$, $(1,9,27)$, $(1,24,10)$, $(1,30,21)$, $(1,24,23)$, 
$(1,21,4)$, $(1,22,5)$, $(1,29,9)$, $(1,23,21)$, $(1,19,13)$.

\subsection*{$m_{26}(2,31)\ge 768$}
$(0,0,1)$, $(1,23,27)$, $(1,9,18)$, $(1,6,4)$, $(1,30,15)$, $(1,7,8)$, $(0,1,0)$, 
$(1,15,30)$, $(1,13,27)$, $(1,17,3)$, $(1,18,8)$, $(1,15,25)$, $(0,1,1)$, $(1,1,12)$, 
$(1,22,24)$, $(1,25,22)$, $(1,27,21)$, $(1,12,7)$, $(0,1,3)$, $(1,7,2)$, $(1,23,3)$, 
$(1,20,14)$, $(1,4,5)$, $(1,1,3)$, $(0,1,6)$, $(0,1,19)$, $(1,10,28)$, $(1,2,10)$, 
$(1,0,13)$, $(1,2,9)$, $(0,1,7)$, $(1,9,9)$, $(1,3,20)$, $(1,28,2)$, $(1,21,2)$, 
$(1,30,22)$, $(0,1,8)$, $(1,16,18)$, $(1,16,26)$, $(1,11,12)$, $(1,1,11)$, $(1,18,12)$, 
$(0,1,10)$, $(1,4,7)$, $(1,20,4)$, $(1,29,16)$, $(0,1,29)$, $(1,17,6)$, $(0,1,11)$, 
$(1,14,11)$, $(1,12,17)$, $(1,16,20)$, $(1,26,23)$, $(1,26,29)$, $(0,1,12)$, $(1,0,24)$, 
$(1,14,6)$, $(1,21,28)$, $(1,28,19)$, $(1,11,1)$, $(0,1,13)$, $(1,21,20)$, $(1,24,13)$, 
$(1,5,21)$, $(1,8,28)$, $(1,22,5)$, $(0,1,14)$, $(1,29,17)$, $(1,25,23)$, $(0,1,15)$, 
$(1,18,25)$, $(1,4,30)$, $(1,7,18)$, $(1,23,29)$, $(1,14,19)$, $(0,1,16)$, $(1,3,19)$, 
$(1,6,19)$, $(1,22,11)$, $(1,19,6)$, $(1,6,2)$, $(0,1,17)$, $(1,2,0)$, $(1,29,1)$, 
$(1,27,19)$, $(1,25,25)$, $(1,21,30)$, $(0,1,20)$, $(1,22,8)$, $(1,15,16)$, $(1,1,27)$, 
$(1,11,22)$, $(1,10,26)$, $(0,1,21)$, $(1,20,1)$, $(1,8,8)$, $(1,15,6)$, $(1,24,27)$, 
$(1,24,17)$, $(0,1,22)$, $(1,26,22)$, $(1,30,11)$, $(1,10,29)$, $(1,9,26)$, $(1,8,14)$, 
$(0,1,23)$, $(1,24,15)$, $(1,18,15)$, $(1,13,9)$, $(1,13,18)$, $(1,0,28)$, $(0,1,25)$, 
$(1,19,13)$, $(1,7,29)$, $(1,23,25)$, $(1,7,30)$, $(1,23,11)$, $(0,1,27)$, $(1,12,4)$, 
$(1,5,9)$, $(1,18,17)$, $(1,17,10)$, $(1,19,18)$, $(0,1,28)$, $(1,28,29)$, $(1,1,0)$, 
$(1,26,5)$, $(1,16,12)$, $(1,28,10)$, $(0,1,30)$, $(1,25,3)$, $(1,17,5)$, $(1,8,1)$, 
$(1,5,3)$, $(1,27,4)$, $(1,0,2)$, $(1,24,19)$, $(1,28,8)$, $(1,15,20)$, $(1,2,27)$, 
$(1,20,9)$, $(1,0,4)$, $(1,6,8)$, $(1,15,23)$, $(1,17,27)$, $(1,6,12)$, $(1,12,6)$, 
$(1,0,5)$, $(1,12,22)$, $(1,8,12)$, $(1,30,26)$, $(1,20,6)$, $(1,9,1)$, $(1,0,6)$, 
$(1,27,26)$, $(1,30,20)$, $(1,29,7)$, $(1,9,24)$, $(1,16,23)$, $(1,0,7)$, $(1,8,23)$, 
$(1,18,10)$, $(1,7,23)$, $(1,23,18)$, $(1,22,2)$, $(1,0,8)$, $(1,15,29)$, $(1,26,27)$, 
$(1,1,2)$, $(1,1,23)$, $(1,21,21)$, $(1,0,9)$, $(1,26,3)$, $(1,7,6)$, $(1,23,17)$, 
$(1,27,3)$, $(1,24,26)$, $(1,0,10)$, $(1,21,12)$, $(1,27,2)$, $(1,8,11)$, $(1,7,16)$, 
$(1,23,14)$, $(1,0,12)$, $(1,13,14)$, $(1,1,1)$, $(1,12,25)$, $(1,26,30)$, $(1,5,15)$, 
$(1,0,16)$, $(1,20,20)$, $(1,22,3)$, $(1,4,28)$, $(1,29,11)$, $(1,4,3)$, $(1,0,17)$, 
$(1,0,25)$, $(1,11,30)$, $(1,20,22)$, $(1,21,10)$, $(1,8,20)$, $(1,0,19)$, $(1,19,28)$, 
$(1,0,26)$, $(1,5,16)$, $(1,12,5)$, $(1,27,0)$, $(1,0,20)$, $(1,3,1)$, $(1,19,16)$, 
$(1,14,1)$, $(1,25,26)$, $(1,1,29)$, $(1,0,22)$, $(1,17,13)$, $(1,3,13)$, $(1,26,12)$, 
$(1,4,4)$, $(1,26,19)$, $(1,0,23)$, $(1,22,4)$, $(1,16,29)$, $(1,24,5)$, $(1,17,25)$, 
$(1,14,30)$, $(1,0,27)$, $(1,18,5)$, $(1,14,17)$, $(1,10,18)$, $(1,5,8)$, $(1,15,11)$, 
$(1,0,29)$, $(1,28,18)$, $(1,25,21)$, $(1,25,24)$, $(1,19,2)$, $(1,18,16)$, $(1,1,4)$, 
$(1,14,0)$, $(1,20,13)$, $(1,17,7)$, $(1,7,22)$, $(1,23,24)$, $(1,1,5)$, $(1,3,29)$, 
$(1,27,7)$, $(1,2,17)$, $(1,24,2)$, $(1,12,13)$, $(1,1,8)$, $(1,15,3)$, $(1,14,27)$, 
$(1,3,6)$, $(1,3,3)$, $(1,8,9)$, $(1,1,9)$, $(1,6,7)$, $(1,1,16)$, $(1,8,13)$, 
$(1,6,25)$, $(1,29,30)$, $(1,1,10)$, $(1,4,1)$, $(1,24,14)$, $(1,30,19)$, $(1,27,24)$, 
$(1,9,10)$, $(1,1,13)$, $(1,19,15)$, $(1,17,20)$, $(1,13,20)$, $(1,11,10)$, $(1,3,4)$, 
$(1,1,15)$, $(1,30,17)$, $(1,13,19)$, $(1,9,2)$, $(1,28,21)$, $(1,16,17)$, $(1,1,18)$, 
$(1,10,19)$, $(1,29,23)$, $(1,10,22)$, $(1,18,20)$, $(1,27,28)$, $(1,1,19)$, $(1,5,4)$, 
$(1,11,3)$, $(1,26,1)$, $(1,17,23)$, $(1,2,3)$, $(1,1,21)$, $(1,7,10)$, $(1,23,6)$, 
$(1,12,0)$, $(1,2,6)$, $(1,10,11)$, $(1,1,22)$, $(1,16,6)$, $(1,2,24)$, $(1,7,24)$, 
$(1,23,5)$, $(1,19,20)$, $(1,1,24)$, $(1,27,8)$, $(1,15,4)$, $(1,18,27)$, $(1,21,11)$, 
$(1,17,18)$, $(1,1,25)$, $(1,24,30)$, $(1,19,5)$, $(1,24,23)$, $(1,4,0)$, $(1,21,22)$, 
$(1,1,28)$, $(1,9,16)$, $(1,30,0)$, $(1,1,30)$, $(1,2,26)$, $(1,18,28)$, $(1,25,12)$, 
$(1,29,18)$, $(1,24,25)$, $(1,2,1)$, $(1,22,21)$, $(1,6,23)$, $(1,3,17)$, $(1,20,11)$, 
$(1,6,22)$, $(1,2,2)$, $(1,13,29)$, $(1,10,6)$, $(1,22,7)$, $(1,11,26)$, $(1,22,15)$, 
$(1,2,5)$, $(1,21,15)$, $(1,3,28)$, $(1,13,15)$, $(1,12,14)$, $(1,13,17)$, $(1,2,8)$, 
$(1,15,10)$, $(1,16,27)$, $(1,8,16)$, $(1,16,28)$, $(1,3,2)$, $(1,2,12)$, $(1,5,12)$, 
$(1,17,15)$, $(1,4,23)$, $(1,13,2)$, $(1,25,4)$, $(1,2,13)$, $(1,25,8)$, $(1,15,8)$, 
$(1,15,27)$, $(1,29,27)$, $(1,26,21)$, $(1,2,14)$, $(1,6,18)$, $(1,20,10)$, $(1,5,29)$, 
$(1,6,24)$, $(1,18,9)$, $(1,2,15)$, $(1,14,4)$, $(1,13,1)$, $(1,30,24)$, $(1,17,16)$, 
$(1,9,11)$, $(1,2,16)$, $(1,2,25)$, $(1,8,30)$, $(1,26,0)$, $(1,14,21)$, $(1,5,5)$, 
$(1,2,18)$, $(1,17,22)$, $(1,22,17)$, $(1,21,1)$, $(1,18,4)$, $(1,19,26)$, $(1,2,19)$, 
$(1,30,7)$, $(1,30,14)$, $(1,9,22)$, $(1,9,19)$, $(1,20,12)$, $(1,2,20)$, $(1,4,6)$, 
$(1,14,20)$, $(1,18,14)$, $(1,8,0)$, $(1,24,18)$, $(1,2,23)$, $(1,28,26)$, $(1,24,24)$, 
$(1,20,26)$, $(1,25,13)$, $(1,11,14)$, $(1,2,28)$, $(1,12,23)$, $(1,7,11)$, $(1,23,13)$, 
$(1,4,17)$, $(1,28,24)$, $(1,2,30)$, $(1,16,16)$, $(1,19,22)$, $(1,12,9)$, $(1,19,23)$, 
$(1,8,25)$, $(1,3,5)$, $(1,5,24)$, $(1,24,29)$, $(1,29,29)$, $(1,16,10)$, $(1,29,19)$, 
$(1,3,8)$, $(1,15,5)$, $(1,10,27)$, $(1,12,24)$, $(1,28,1)$, $(1,4,22)$, $(1,3,10)$, 
$(1,18,21)$, $(1,30,21)$, $(1,4,18)$, $(1,9,23)$, $(1,24,1)$, $(1,3,11)$, $(1,9,4)$, 
$(1,17,28)$, $(1,11,0)$, $(1,17,17)$, $(1,30,4)$, $(1,3,12)$, $(1,25,17)$, $(1,3,26)$, 
$(1,14,10)$, $(1,27,25)$, $(1,20,30)$, $(1,3,14)$, $(1,29,28)$, $(1,18,6)$, $(1,28,5)$, 
$(1,26,18)$, $(1,19,14)$, $(1,3,18)$, $(1,21,6)$, $(1,26,16)$, $(1,22,16)$, $(1,18,24)$, 
$(1,5,7)$, $(1,3,21)$, $(1,22,1)$, $(1,29,12)$, $(1,7,28)$, $(1,23,28)$, $(1,8,24)$, 
$(1,3,23)$, $(1,30,23)$, $(1,11,5)$, $(1,9,14)$, $(1,7,9)$, $(1,23,16)$, $(1,3,25)$, 
$(1,10,30)$, $(1,28,3)$, $(1,3,27)$, $(1,7,14)$, $(1,23,20)$, $(1,5,11)$, $(1,29,8)$, 
$(1,15,12)$, $(1,3,30)$, $(1,4,29)$, $(1,5,13)$, $(1,21,23)$, $(1,20,7)$, $(1,10,25)$, 
$(1,4,2)$, $(1,11,28)$, $(1,6,9)$, $(1,25,28)$, $(1,16,22)$, $(1,6,29)$, $(1,4,8)$, 
$(1,15,7)$, $(1,12,27)$, $(1,4,10)$, $(1,22,9)$, $(1,17,11)$, $(1,26,26)$, $(1,29,0)$, 
$(1,19,4)$, $(1,4,12)$, $(1,8,5)$, $(1,28,13)$, $(1,18,11)$, $(1,22,19)$, $(1,8,18)$, 
$(1,4,13)$, $(1,10,10)$, $(1,14,2)$, $(1,16,15)$, $(1,27,1)$, $(1,13,6)$, $(1,4,14)$, 
$(1,17,12)$, $(1,10,21)$, $(1,19,9)$, $(1,21,4)$, $(1,21,24)$, $(1,4,15)$, $(1,16,25)$, 
$(1,13,30)$, $(1,29,20)$, $(1,8,26)$, $(1,27,22)$, $(1,4,20)$, $(1,12,15)$, $(1,22,26)$, 
$(1,13,21)$, $(1,28,16)$, $(1,5,19)$, $(1,4,25)$, $(1,18,30)$, $(1,11,24)$, $(1,30,18)$, 
$(1,26,17)$, $(1,9,28)$, $(1,4,27)$, $(1,30,29)$, $(1,8,15)$, $(1,9,29)$, $(1,13,8)$, 
$(1,15,26)$, $(1,5,0)$, $(1,30,10)$, $(1,24,3)$, $(1,9,12)$, $(1,11,29)$, $(1,21,3)$, 
$(1,5,1)$, $(1,10,1)$, $(1,16,11)$, $(1,20,0)$, $(1,28,23)$, $(1,29,9)$, $(1,5,2)$, 
$(1,7,26)$, $(1,23,4)$, $(1,24,21)$, $(1,10,2)$, $(1,17,0)$, $(1,5,10)$, $(1,25,0)$, 
$(1,7,20)$, $(1,23,8)$, $(1,15,13)$, $(1,22,27)$, $(1,5,14)$, $(1,28,6)$, $(1,19,8)$, 
$(1,15,28)$, $(1,27,27)$, $(1,27,23)$, $(1,5,18)$, $(1,13,7)$, $(1,6,21)$, $(1,27,29)$, 
$(1,18,1)$, $(1,6,15)$, $(1,5,30)$, $(1,22,25)$, $(1,28,30)$, $(1,30,6)$, $(1,16,9)$, 
$(1,9,25)$, $(1,6,5)$, $(1,7,19)$, $(1,23,26)$, $(1,6,6)$, $(1,11,21)$, $(1,14,12)$, 
$(1,6,20)$, $(1,10,5)$, $(1,25,5)$, $(1,6,11)$, $(1,28,14)$, $(1,27,15)$, $(1,6,13)$, 
$(1,13,22)$, $(1,24,0)$, $(1,6,16)$, $(1,24,12)$, $(1,21,16)$, $(1,6,17)$, $(1,26,13)$, 
$(1,28,20)$, $(1,6,26)$, $(1,21,26)$, $(1,10,23)$, $(1,6,28)$, $(1,14,7)$, $(1,26,10)$, 
$(1,7,4)$, $(1,23,22)$, $(1,28,28)$, $(1,17,14)$, $(1,21,7)$, $(1,14,29)$, $(1,7,7)$, 
$(1,23,7)$, $(1,17,1)$, $(1,20,28)$, $(1,24,20)$, $(1,16,4)$, $(1,7,13)$, $(1,23,0)$, 
$(1,18,26)$, $(1,12,1)$, $(1,14,18)$, $(1,28,9)$, $(1,7,15)$, $(1,23,21)$, $(1,13,25)$, 
$(1,16,30)$, $(1,25,14)$, $(1,20,16)$, $(1,7,21)$, $(1,23,19)$, $(1,26,9)$, $(1,27,9)$, 
$(1,29,21)$, $(1,11,20)$, $(1,7,25)$, $(1,23,30)$, $(1,21,8)$, $(1,15,15)$, $(1,28,27)$, 
$(1,13,26)$, $(1,8,2)$, $(1,19,1)$, $(1,11,13)$, $(1,8,6)$, $(1,17,19)$, $(1,17,29)$, 
$(1,25,16)$, $(1,28,4)$, $(1,10,4)$, $(1,8,7)$, $(1,25,9)$, $(1,10,0)$, $(1,27,14)$, 
$(1,26,24)$, $(1,8,17)$, $(1,8,29)$, $(1,20,23)$, $(1,19,24)$, $(1,26,15)$, $(1,11,19)$, 
$(1,13,0)$, $(1,9,0)$, $(1,12,3)$, $(1,25,29)$, $(1,12,28)$, $(1,30,3)$, $(1,22,20)$, 
$(1,9,3)$, $(1,21,17)$, $(1,19,21)$, $(1,29,3)$, $(1,30,13)$, $(1,14,26)$, $(1,9,7)$, 
$(1,18,2)$, $(1,26,20)$, $(1,25,18)$, $(1,30,1)$, $(1,28,0)$, $(1,9,13)$, $(1,27,16)$, 
$(1,10,9)$, $(1,20,29)$, $(1,30,9)$, $(1,12,12)$, $(1,9,21)$, $(1,26,11)$, $(1,24,7)$, 
$(1,22,6)$, $(1,30,25)$, $(1,19,30)$, $(1,9,27)$, $(1,19,7)$, $(1,27,11)$, $(1,14,5)$, 
$(1,30,8)$, $(1,15,2)$, $(1,10,3)$, $(1,10,15)$, $(1,24,16)$, $(1,13,3)$, $(1,19,11)$, 
$(1,16,24)$, $(1,10,14)$, $(1,11,4)$, $(1,25,20)$, $(1,14,25)$, $(1,17,30)$, $(1,12,10)$, 
$(1,10,16)$, $(1,21,18)$, $(1,29,5)$, $(1,20,2)$, $(1,20,17)$, $(1,18,0)$, $(1,11,2)$, 
$(1,29,6)$, $(1,29,15)$, $(1,21,0)$, $(1,13,12)$, $(1,14,22)$, $(1,11,15)$, $(1,19,10)$, 
$(1,13,5)$, $(1,18,7)$, $(1,29,13)$, $(1,12,18)$, $(1,11,18)$, $(1,27,13)$, $(1,16,3)$, 
$(1,12,21)$, $(1,18,22)$, $(1,20,3)$, $(1,11,23)$, $(1,12,19)$, $(1,12,16)$, $(1,29,2)$, 
$(1,14,14)$, $(1,24,11)$, $(1,12,8)$, $(1,15,19)$, $(1,21,27)$, $(1,16,1)$, $(1,21,9)$, 
$(1,22,10)$, $(1,13,10)$, $(1,20,15)$, $(1,26,6)$, $(1,14,8)$, $(1,15,0)$, $(1,25,27)$, 
$(1,22,13)$, $(1,20,19)$, $(1,16,14)$, $(1,14,24)$, $(1,22,12)$, $(1,24,9)$, $(1,17,2)$, 
$(1,27,5)$, $(1,22,28)$, $(1,26,4)$, $(1,20,25)$, $(1,27,30)$.

\subsection*{$m_{27}(2,31)\ge 805$}
$(0,0,1)$, $(1,0,19)$, $(1,18,6)$, $(1,26,30)$, $(1,25,9)$, $(0,1,1)$, $(1,13,22)$, 
$(1,5,23)$, $(1,25,19)$, $(1,0,17)$, $(0,1,2)$, $(1,9,2)$, $(1,7,18)$, $(1,8,18)$, 
$(1,2,30)$, $(0,1,3)$, $(1,16,6)$, $(1,25,4)$, $(1,0,23)$, $(1,6,25)$, $(0,1,4)$, 
$(1,10,7)$, $(1,26,17)$, $(1,4,5)$, $(1,18,10)$, $(0,1,5)$, $(1,27,30)$, $(1,17,24)$, 
$(1,9,29)$, $(0,1,22)$, $(0,1,6)$, $(1,19,21)$, $(1,23,9)$, $(1,11,20)$, $(1,1,8)$, 
$(0,1,7)$, $(1,12,17)$, $(1,1,2)$, $(1,24,8)$, $(1,13,24)$, $(0,1,9)$, $(1,21,0)$, 
$(1,15,29)$, $(1,23,28)$, $(1,19,1)$, $(0,1,10)$, $(1,15,1)$, $(1,12,21)$, $(1,22,17)$, 
$(1,14,15)$, $(0,1,11)$, $(1,22,5)$, $(0,1,13)$, $(1,20,26)$, $(1,21,14)$, $(0,1,12)$, 
$(1,18,16)$, $(1,24,22)$, $(1,14,22)$, $(1,26,0)$, $(0,1,14)$, $(1,7,23)$, $(1,20,1)$, 
$(1,7,7)$, $(1,12,2)$, $(0,1,15)$, $(1,8,28)$, $(1,4,10)$, $(1,1,3)$, $(1,4,12)$, 
$(0,1,16)$, $(1,30,14)$, $(1,13,3)$, $(1,30,12)$, $(1,20,23)$, $(0,1,17)$, $(1,5,13)$, 
$(1,19,19)$, $(1,29,1)$, $(1,23,27)$, $(0,1,18)$, $(1,14,27)$, $(1,10,26)$, $(1,16,13)$, 
$(1,16,28)$, $(0,1,21)$, $(1,25,20)$, $(1,0,20)$, $(1,12,0)$, $(1,8,7)$, $(0,1,23)$, 
$(1,6,18)$, $(1,3,28)$, $(1,21,6)$, $(1,3,21)$, $(0,1,25)$, $(1,2,29)$, $(1,14,16)$, 
$(1,27,10)$, $(1,30,26)$, $(0,1,27)$, $(1,26,25)$, $(1,30,7)$, $(1,18,4)$, $(1,7,16)$, 
$(0,1,28)$, $(1,1,24)$, $(1,22,27)$, $(1,2,14)$, $(1,28,13)$, $(0,1,29)$, $(1,23,10)$, 
$(1,8,0)$, $(1,19,15)$, $(1,24,18)$, $(0,1,30)$, $(1,24,15)$, $(1,2,15)$, $(1,3,25)$, 
$(1,29,4)$, $(1,0,0)$, $(1,7,26)$, $(1,14,25)$, $(1,12,15)$, $(1,25,28)$, $(1,0,1)$, 
$(1,21,9)$, $(1,11,16)$, $(1,8,24)$, $(1,25,17)$, $(1,0,2)$, $(1,10,29)$, $(1,21,15)$, 
$(1,7,3)$, $(1,25,12)$, $(1,0,3)$, $(1,24,12)$, $(1,1,17)$, $(1,15,16)$, $(1,25,18)$, 
$(1,0,4)$, $(1,2,21)$, $(1,5,29)$, $(1,24,19)$, $(1,25,0)$, $(1,0,5)$, $(1,28,16)$, 
$(1,6,1)$, $(1,28,10)$, $(1,25,6)$, $(1,0,6)$, $(1,8,27)$, $(1,29,8)$, $(1,6,13)$, 
$(1,25,1)$, $(1,0,7)$, $(1,15,3)$, $(1,20,12)$, $(1,14,26)$, $(1,25,21)$, $(1,0,8)$, 
$(1,27,15)$, $(1,16,0)$, $(1,13,5)$, $(1,25,8)$, $(1,0,12)$, $(1,30,18)$, $(1,4,26)$, 
$(1,22,8)$, $(1,25,3)$, $(1,0,13)$, $(1,22,10)$, $(1,15,28)$, $(1,17,27)$, $(1,25,2)$, 
$(1,0,15)$, $(1,13,1)$, $(1,22,18)$, $(1,1,1)$, $(1,25,14)$, $(1,0,18)$, $(1,5,24)$, 
$(1,12,19)$, $(1,10,4)$, $(1,25,29)$, $(1,0,21)$, $(1,9,28)$, $(1,3,23)$, $(1,27,20)$, 
$(1,25,30)$, $(1,0,22)$, $(1,1,20)$, $(1,19,9)$, $(1,11,25)$, $(1,25,23)$, $(1,0,24)$, 
$(1,14,2)$, $(1,17,3)$, $(1,3,12)$, $(1,25,5)$, $(1,0,25)$, $(1,20,8)$, $(1,2,20)$, 
$(1,20,28)$, $(1,25,16)$, $(1,0,26)$, $(1,4,23)$, $(1,28,5)$, $(1,21,18)$, $(1,25,15)$, 
$(1,0,27)$, $(1,16,4)$, $(1,10,13)$, $(1,4,2)$, $(1,25,22)$, $(1,0,29)$, $(1,3,22)$, 
$(1,23,21)$, $(1,2,22)$, $(1,25,7)$, $(1,0,30)$, $(1,29,17)$, $(1,9,10)$, $(1,23,29)$, 
$(1,25,10)$, $(1,1,0)$, $(1,18,3)$, $(1,6,12)$, $(1,4,0)$, $(1,26,6)$, $(1,1,5)$, 
$(1,17,28)$, $(1,5,28)$, $(1,20,19)$, $(1,23,15)$, $(1,1,6)$, $(1,10,17)$, $(1,30,0)$, 
$(1,6,14)$, $(1,16,5)$, $(1,1,9)$, $(1,11,23)$, $(1,24,3)$, $(1,19,12)$, $(1,8,29)$, 
$(1,1,10)$, $(1,21,21)$, $(1,1,30)$, $(1,16,22)$, $(1,9,26)$, $(1,1,12)$, $(1,26,20)$, 
$(1,22,4)$, $(1,18,5)$, $(1,27,3)$, $(1,1,13)$, $(1,27,26)$, $(1,3,29)$, $(1,30,27)$, 
$(1,30,25)$, $(1,1,14)$, $(1,7,30)$, $(1,10,10)$, $(1,9,4)$, $(1,15,8)$, $(1,1,15)$, 
$(1,6,24)$, $(1,21,20)$, $(1,27,6)$, $(1,7,1)$, $(1,1,16)$, $(1,13,4)$, $(1,16,7)$, 
$(1,7,21)$, $(1,28,0)$, $(1,1,19)$, $(1,1,25)$, $(1,28,1)$, $(1,14,8)$, $(1,5,7)$, 
$(1,1,21)$, $(1,8,5)$, $(1,20,5)$, $(1,12,25)$, $(1,2,16)$, $(1,1,22)$, $(1,3,6)$, 
$(1,15,23)$, $(1,22,2)$, $(1,6,4)$, $(1,1,26)$, $(1,23,2)$, $(1,29,16)$, $(1,2,17)$, 
$(1,13,14)$, $(1,1,27)$, $(1,12,29)$, $(1,4,13)$, $(1,29,20)$, $(1,17,2)$, $(1,1,28)$, 
$(1,30,13)$, $(1,12,9)$, $(1,11,18)$, $(1,19,27)$, $(1,2,1)$, $(1,24,5)$, $(1,10,30)$, 
$(1,22,3)$, $(1,9,12)$, $(1,2,2)$, $(1,16,20)$, $(1,30,1)$, $(1,5,20)$, $(1,13,19)$, 
$(1,2,3)$, $(1,12,12)$, $(1,17,9)$, $(1,11,14)$, $(1,3,17)$, $(1,2,4)$, $(1,22,1)$, 
$(1,20,0)$, $(1,4,21)$, $(1,17,26)$, $(1,2,5)$, $(1,8,4)$, $(1,3,20)$, $(1,8,17)$, 
$(1,27,28)$, $(1,2,6)$, $(1,29,15)$, $(1,5,14)$, $(1,6,19)$, $(1,19,14)$, $(1,2,7)$, 
$(1,6,0)$, $(1,12,24)$, $(1,19,6)$, $(1,30,10)$, $(1,2,8)$, $(1,26,9)$, $(1,11,27)$, 
$(1,18,7)$, $(1,14,13)$, $(1,2,9)$, $(1,11,10)$, $(1,22,25)$, $(1,27,29)$, $(1,15,7)$, 
$(1,2,10)$, $(1,10,8)$, $(1,23,22)$, $(1,21,4)$, $(1,21,2)$, $(1,2,12)$, $(1,18,24)$, 
$(1,28,7)$, $(1,28,28)$, $(1,26,3)$, $(1,2,13)$, $(1,30,17)$, $(1,7,8)$, $(1,29,27)$, 
$(1,28,22)$, $(1,2,19)$, $(1,21,30)$, $(1,14,18)$, $(1,9,16)$, $(1,5,5)$, $(1,2,24)$, 
$(1,23,3)$, $(1,16,12)$, $(1,15,10)$, $(1,16,1)$, $(1,2,25)$, $(1,9,6)$, $(1,19,3)$, 
$(1,13,12)$, $(1,6,30)$, $(1,2,27)$, $(1,15,18)$, $(1,26,13)$, $(1,10,15)$, $(1,18,20)$, 
$(1,2,28)$, $(1,7,2)$, $(1,21,28)$, $(1,16,9)$, $(1,11,0)$, $(1,3,0)$, $(1,24,13)$, 
$(1,26,29)$, $(1,28,14)$, $(1,18,27)$, $(1,3,1)$, $(1,19,22)$, $(1,16,8)$, $(1,12,18)$, 
$(1,30,5)$, $(1,3,2)$, $(1,20,14)$, $(1,24,0)$, $(1,21,8)$, $(1,4,1)$, $(1,3,4)$, 
$(1,26,28)$, $(1,29,26)$, $(1,18,1)$, $(1,15,17)$, $(1,3,5)$, $(1,16,15)$, $(1,28,27)$, 
$(1,19,24)$, $(1,20,13)$, $(1,3,8)$, $(1,3,26)$, $(1,17,7)$, $(1,4,20)$, $(1,23,23)$, 
$(1,3,9)$, $(1,11,24)$, $(1,15,9)$, $(1,11,26)$, $(1,13,0)$, $(1,3,10)$, $(1,12,16)$, 
$(1,21,3)$, $(1,5,12)$, $(1,28,19)$, $(1,3,14)$, $(1,10,1)$, $(1,8,16)$, $(1,30,29)$, 
$(1,10,21)$, $(1,3,15)$, $(1,22,29)$, $(1,7,17)$, $(1,16,17)$, $(1,19,20)$, $(1,3,16)$, 
$(1,14,0)$, $(1,23,1)$, $(1,10,3)$, $(1,29,12)$, $(1,3,18)$, $(1,13,8)$, $(1,20,4)$, 
$(1,13,10)$, $(1,14,24)$, $(1,3,19)$, $(1,4,18)$, $(1,10,14)$, $(1,15,25)$, $(1,5,25)$, 
$(1,3,30)$, $(1,6,2)$, $(1,9,15)$, $(1,29,6)$, $(1,12,7)$, $(1,4,3)$, $(1,14,12)$, 
$(1,22,22)$, $(1,19,18)$, $(1,7,20)$, $(1,4,4)$, $(1,8,15)$, $(1,26,7)$, $(1,23,16)$, 
$(1,18,15)$, $(1,4,6)$, $(1,21,24)$, $(1,10,5)$, $(1,6,9)$, $(1,11,21)$, $(1,4,7)$, 
$(1,19,25)$, $(1,14,21)$, $(1,20,2)$, $(1,28,2)$, $(1,4,9)$, $(1,11,29)$, $(1,16,29)$, 
$(1,8,8)$, $(1,22,16)$, $(1,4,14)$, $(1,30,4)$, $(1,24,30)$, $(1,15,20)$, $(1,9,5)$, 
$(1,4,15)$, $(1,20,9)$, $(1,11,9)$, $(1,11,22)$, $(1,29,10)$, $(1,4,16)$, $(1,15,27)$, 
$(1,6,20)$, $(1,29,13)$, $(1,13,6)$, $(1,4,19)$, $(1,13,28)$, $(1,13,17)$, $(1,28,29)$, 
$(1,5,4)$, $(1,4,22)$, $(1,6,16)$, $(1,23,26)$, $(1,12,6)$, $(1,27,25)$, $(1,4,24)$, 
$(1,5,1)$, $(1,29,19)$, $(1,14,5)$, $(1,19,23)$, $(1,4,25)$, $(1,7,0)$, $(1,30,23)$, 
$(1,21,17)$, $(1,23,24)$, $(1,4,28)$, $(1,27,21)$, $(1,19,10)$, $(1,27,14)$, $(1,14,14)$, 
$(1,4,29)$, $(1,9,30)$, $(1,28,15)$, $(1,9,23)$, $(1,16,30)$, $(1,5,0)$, $(1,9,19)$, 
$(1,22,7)$, $(1,22,26)$, $(1,8,30)$, $(1,5,2)$, $(1,27,19)$, $(1,27,16)$, $(1,29,2)$, 
$(1,19,0)$, $(1,5,8)$, $(1,7,19)$, $(1,26,8)$, $(1,24,28)$, $(1,18,14)$, $(1,5,9)$, 
$(1,11,19)$, $(1,5,26)$, $(1,15,19)$, $(1,24,23)$, $(1,5,15)$, $(1,21,19)$, $(1,20,22)$, 
$(1,30,3)$, $(1,27,12)$, $(1,5,17)$, $(1,18,19)$, $(1,16,21)$, $(1,21,25)$, $(1,26,26)$, 
$(1,5,21)$, $(1,14,19)$, $(1,30,9)$, $(1,11,15)$, $(1,10,2)$, $(1,5,22)$, $(1,10,19)$, 
$(1,28,24)$, $(1,27,0)$, $(1,22,20)$, $(1,5,27)$, $(1,26,19)$, $(1,15,13)$, $(1,18,22)$, 
$(1,12,5)$, $(1,5,30)$, $(1,30,19)$, $(1,9,27)$, $(1,17,21)$, $(1,7,13)$, $(1,6,7)$, 
$(1,27,7)$, $(1,24,6)$, $(1,16,18)$, $(1,24,17)$, $(1,6,8)$, $(1,14,30)$, $(1,8,6)$, 
$(1,9,22)$, $(1,27,24)$, $(1,6,15)$, $(1,30,16)$, $(1,28,6)$, $(1,24,9)$, $(1,11,28)$, 
$(1,6,17)$, $(1,22,23)$, $(1,30,6)$, $(1,14,28)$, $(1,14,4)$, $(1,6,21)$, $(1,13,27)$, 
$(1,7,6)$, $(1,17,13)$, $(1,8,21)$, $(1,6,22)$, $(1,7,9)$, $(1,11,6)$, $(1,23,14)$, 
$(1,20,18)$, $(1,6,23)$, $(1,26,4)$, $(1,17,6)$, $(1,18,8)$, $(1,21,10)$, $(1,6,26)$, 
$(1,10,18)$, $(1,14,6)$, $(1,28,20)$, $(1,7,29)$, $(1,6,28)$, $(1,24,29)$, $(1,22,6)$, 
$(1,13,2)$, $(1,23,25)$, $(1,7,5)$, $(1,7,22)$, $(1,13,29)$, $(1,29,5)$, $(1,24,26)$, 
$(1,7,10)$, $(1,18,9)$, $(1,11,1)$, $(1,27,1)$, $(1,26,23)$, $(1,7,12)$, $(1,27,4)$, 
$(1,20,3)$, $(1,17,12)$, $(1,29,3)$, $(1,7,14)$, $(1,23,20)$, $(1,15,26)$, $(1,28,3)$, 
$(1,23,12)$, $(1,7,15)$, $(1,15,21)$, $(1,27,8)$, $(1,8,25)$, $(1,15,24)$, $(1,7,28)$, 
$(1,22,24)$, $(1,17,23)$, $(1,19,16)$, $(1,10,16)$, $(1,8,1)$, $(1,17,4)$, $(1,28,26)$, 
$(1,30,22)$, $(1,8,2)$, $(1,8,3)$, $(1,10,12)$, $(1,19,8)$, $(1,17,17)$, $(1,24,25)$, 
$(1,8,9)$, $(1,11,2)$, $(1,30,30)$, $(1,23,5)$, $(1,9,17)$, $(1,8,10)$, $(1,28,18)$, 
$(1,16,2)$, $(1,14,23)$, $(1,15,14)$, $(1,8,22)$, $(1,18,25)$, $(1,22,14)$, $(1,29,24)$, 
$(1,26,24)$, $(1,8,23)$, $(1,13,13)$, $(1,15,0)$, $(1,10,0)$, $(1,20,27)$, $(1,8,26)$, 
$(1,24,27)$, $(1,20,10)$, $(1,19,13)$, $(1,28,23)$, $(1,9,0)$, $(1,14,17)$, $(1,29,30)$, 
$(1,20,30)$, $(1,9,8)$, $(1,9,7)$, $(1,10,22)$, $(1,17,1)$, $(1,13,30)$, $(1,19,28)$, 
$(1,9,18)$, $(1,20,25)$, $(1,16,27)$, $(1,23,30)$, $(1,27,13)$, $(1,9,21)$, $(1,12,4)$, 
$(1,23,0)$, $(1,28,30)$, $(1,12,14)$, $(1,9,24)$, $(1,13,26)$, $(1,26,15)$, $(1,19,30)$, 
$(1,18,26)$, $(1,9,25)$, $(1,24,20)$, $(1,28,25)$, $(1,17,30)$, $(1,21,1)$, $(1,10,24)$, 
$(1,24,10)$, $(1,24,2)$, $(1,12,26)$, $(1,20,20)$, $(1,10,25)$, $(1,13,18)$, $(1,12,10)$, 
$(1,29,22)$, $(1,22,28)$, $(1,10,27)$, $(1,27,5)$, $(1,23,13)$, $(1,23,7)$, $(1,13,23)$, 
$(1,11,3)$, $(1,26,12)$, $(1,16,25)$, $(1,18,2)$, $(1,13,9)$, $(1,11,7)$, $(1,21,5)$, 
$(1,13,7)$, $(1,17,14)$, $(1,22,9)$, $(1,11,17)$, $(1,17,18)$, $(1,17,0)$, $(1,27,18)$, 
$(1,28,9)$, $(1,12,3)$, $(1,15,12)$, $(1,21,7)$, $(1,30,8)$, $(1,19,29)$, $(1,12,8)$, 
$(1,18,23)$, $(1,23,18)$, $(1,22,21)$, $(1,26,16)$, $(1,12,22)$, $(1,24,14)$, $(1,17,16)$, 
$(1,16,23)$, $(1,17,15)$, $(1,12,28)$, $(1,28,8)$, $(1,15,5)$, $(1,29,29)$, $(1,27,23)$, 
$(1,13,20)$, $(1,18,21)$, $(1,30,15)$, $(1,17,10)$, $(1,26,10)$, $(1,14,1)$, $(1,18,13)$, 
$(1,17,25)$, $(1,23,17)$, $(1,26,27)$, $(1,14,7)$, $(1,16,24)$, $(1,29,14)$, $(1,21,29)$, 
$(1,17,20)$, $(1,15,4)$, $(1,30,24)$, $(1,16,16)$, $(1,20,7)$, $(1,29,18)$, $(1,16,3)$, 
$(1,22,12)$, $(1,24,21)$, $(1,29,23)$, $(1,29,21)$, $(1,18,18)$, $(1,18,29)$, $(1,19,2)$, 
$(1,26,5)$, $(1,26,21)$, $(1,21,13)$, $(1,24,4)$, $(1,27,27)$, $(1,21,27)$, $(1,22,13)$.

\end{document}